\documentclass{article}
\makeatletter
\def\input@path{{template-jems/}}
\makeatother
\usepackage[journal=JEMS,lang=british]{ems-journal-jems}

\usepackage{stmaryrd}
\usepackage{bbm}

\numberwithin{equation}{section}

\newcommand{\assign}{:=}
\newcommand{\backassign}{=:}
\newcommand{\mathD}{\mathrm{D}}
\newcommand{\mathd}{\mathrm{d}}

\makeatletter
\newcommand*\botimes{{\mathpalette\botimes@{1.5}}}
\newcommand*\botimes@[2]{\mathbin{\vcenter{\hbox{\scalebox{#2}{\hspace{0.0mm}$\m@th#1\otimes$\hspace{0.0mm}}}}}}
\newcommand*\Cdot{{\mathpalette\Cdot@{.7}}}
\newcommand*\Cdot@[2]{\mathbin{\vcenter{\hbox{\scalebox{#2}{\hspace{0.1mm}$\m@th#1\bullet$\hspace{0.1mm}}}}}}
\makeatother
\renewcommand{\interleave}{%
	\mathopen{}\left\vert\kern-0.28ex\middle\vert\kern-0.28ex\right\vert
}
\newcommand{\linterleave}{%
	\mathopen{}\left\vert\kern-0.28ex\middle\vert\kern-0.28ex\middle\vert
}
\newcommand{\rinterleave}{%
	\middle\vert\kern-0.28ex\middle\vert\kern-0.28ex\right\vert\mathclose{}
}
\usepackage{mathtools}

\DeclarePairedDelimiterX{\trinorm}[1]
{\lvert\!\lvert\!\lvert}
{\rvert\!\rvert\!\rvert}
{#1}

\theoremstyle{plain}
\newtheorem{theorem}{Theorem}[section]
\newtheorem{lemma}[theorem]{Lemma}
\theoremstyle{definition}
\newtheorem{definition}[theorem]{Definition}
\newtheorem{remark}[theorem]{Remark}

\begin{document}

\title{Parabolic stochastic quantisation of the fractional $\Phi^4_3$ model in the full subcritical regime}
\titlemark{Fractional $\Phi^4_3$ stochastic quantisation}

\emsauthor{1}{
    \givenname{Pawe\l{}}
    \surname{Duch}
    \mrid{}
    \zblid{}
    \orcid{0000-0002-2445-2189}}{P.~Duch}
\emsauthor{2}{
    \givenname{Massimiliano}
    \surname{Gubinelli}
    \mrid{}
    \zblid{}
    \orcid{0000-0002-4014-2949}}{M.~Gubinelli}
\emsauthor{3}{
    \givenname{Paolo}
    \surname{Rinaldi}
    \mrid{}
    \zblid{}
    \orcid{0000-0003-1578-2677}}{P.~Rinaldi}

\Emsaffil{1}{
    \pretext{}
    \department{Faculty of Mathematics and Computer Science}
    \organisation{Adam Mickiewicz University in Pozna\'{n}}
    \rorid{04g6bbq64}
    \address{ul. Uniwersytetu Pozna\'{n}skiego 4}
    \zip{61-614}
    \city{Pozna\'{n}}
    \country{Poland}
    \posttext{}
    \affemail{pawel.duch@amu.edu.pl}
    \furtheremail{}}
\Emsaffil{2}{
    \pretext{}
    \department{Mathematical Institute}
    \organisation{University of Oxford}
    \rorid{052gg0110}
    \address{Andrew Wiles Building, Radcliffe Observatory Quarter, Woodstock Road}
    \zip{OX2 6GG}
    \city{Oxford}
    \country{United Kingdom}
    \posttext{}
    \affemail{gubinelli@maths.ox.ac.uk}
    \furtheremail{}}
\Emsaffil{3}{
    \pretext{}
    \department{Institute for Applied Mathematics and Hausdorff Center for Mathematics}
    \organisation{University of Bonn}
    \rorid{041nas322}
    \address{Endenicher Allee 60}
    \zip{53115}
    \city{Bonn}
    \country{Germany}
    \posttext{}
    \affemail{rinaldi@iam.uni-bonn.de}
    \furtheremail{}}

\classification[60H17, 60H30, 81T08, 81T17]{35B45}
\keywords{stochastic quantisation, renormalisation group, fractional Laplacian, singular SPDEs}

\begin{abstract}
We present a construction of the measure of the fractional
  $\Phi^4$ Euclidean quantum field theory on $\mathbb{R}^3$ in the full
  subcritical regime via parabolic stochastic quantisation. Our approach is
  based on a truncated flow equation for the effective description
  of the model at sufficiently small scales and on coercive estimates for the
  nonlinear stochastic partial differential equation describing the
  interacting field. The constructed measure is translation invariant,
  reflection positive, and has quartic exponential tails.
\end{abstract}

\maketitle

{\small\tableofcontents}

\section{Introduction}\label{sec:intro}

For $\varepsilon \in 2^{-\mathbb{N}_0}$ and $M \in
\mathbb{N}_+$ let
\[ \mathbb{R}_{\varepsilon}^d \assign (\varepsilon \mathbb{Z})^d, \qquad
   \mathbb{T}_{\varepsilon, M}^d \assign (\varepsilon
   \mathbb{Z})^d / (M\mathbb{Z})^d, \qquad \Omega_{\varepsilon, M} \assign \{
   \varphi : \mathbb{T}_{\varepsilon, M}^d \rightarrow \mathbb{R} \} . \]
Define a probability measure $\nu_{\varepsilon, M}$ on $\Omega_{\varepsilon,
M}$ by
\begin{equation}
  \nu_{\varepsilon, M} (\mathd \varphi) \assign \frac{\exp (- S_{\varepsilon,
  M} (\varphi))}{Z_{\varepsilon, M}} \prod_{x \in \mathbb{T}_{\varepsilon,
  M}^d} \mathd \varphi (x) \label{eq:measure},
\end{equation}
where
\begin{equation}
  S_{\varepsilon, M} (\varphi) \assign 2 \varepsilon^d  \sum_{x \in
  \mathbb{T}_{\varepsilon, M}^d} \left[ \frac{1}{2} \varphi (x)  (-
  \Delta_{\varepsilon})^s \varphi (x) + \frac{m^2}{2} \varphi (x)^2 +
  \frac{\lambda}{4} \varphi (x)^4 - \frac{r_{\varepsilon, M}}{2} \varphi (x)^2
  \right]. \label{eq:action}
\end{equation}
The normalisation constant $Z_{\varepsilon, M} > 0$ is chosen so that the
total mass of $\nu_{\varepsilon, M}$ equals one. The parameters $m > 0,
  \lambda > 0, r_{\varepsilon, M} \in \mathbb{R}$ are referred to as
\textit{mass}, \textit{coupling constant} and \textit{mass
renormalisation}, respectively. The operator $(- \Delta_{\varepsilon})^s$ is
the discrete fractional Laplacian of order $s \in (0, 1]$ defined via
functional calculus as the $s$-th power of the discrete nearest-neighbour
Laplacian $- \Delta_{\varepsilon}$. In order to pass to the limit, we define a
probability measure $\hat{\nu}_{\varepsilon, M}$ on $\mathcal{S}'
(\mathbb{R}^3)$ by embedding sample paths of $\nu_{\varepsilon, M}$ in
$\mathcal{S}' (\mathbb{R}^3)$ via a suitable Fourier multiplier,
see~\eqref{eq:embedding} below.

For the sake of clarity, we restrict our considerations to $d = 3$. Then $s =
1$ corresponds to the standard $\Phi_3^4$ model, while for fractional
exponents $s > s_c \assign 3 / 4$ the model is subcritical, meaning that the
nonlinear part can be treated as a perturbation of the Gaussian measure at
small scales. Since reflection positivity is expected to hold only for $s
\leqslant 1$, we restrict our analysis to $s \in (s_c, 1]$. For further
discussion of the measure constructed in Theorem~\ref{thm:main}, we refer
to~{\cite{GH21}}. In what follows, we present
a detailed analysis of the fractional regime $s \in (3 / 4, 1)$. The case $s =
1$, corresponding to the classical Laplacian, can be handled by the same
strategy with substantial simplifications. The main result of this paper is the
following theorem.

\begin{theorem}
  \label{thm:main}Let $d = 3$ and fix $s \in (3 / 4, 1]$, $m > 0$, $\lambda >
  0$. There exists a choice of mass renormalisation
  $(r_{\varepsilon, M})_{\varepsilon \in 2^{-\mathbb{N}_0}, M \in
  \mathbb{N}_+}$ such that $(\hat{\nu}_{\varepsilon, M})_{\varepsilon \in
  2^{-\mathbb{N}_0}, M \in \mathbb{N}_+}$ is a tight family of probability
  measures on $\mathcal{S}' (\mathbb{R}^3)$. Any accumulation point $\nu$ of
  this family is non-Gaussian, invariant under translations, reflection
  positive, and satisfies
  \begin{equation}
    \int e^{\theta \| (1 + | \Cdot |)^{- b}  (1 - \Delta )^{- a} \phi
    \|_{L^2}^4} \nu (\mathd \phi) < \infty, \label{eq:int-bound}
  \end{equation}
  for sufficiently large $a, b > 0$ and sufficiently small $\theta > 0$.
\end{theorem}

\begin{proof}
  Tightness is established in Section~\ref{sec:tightness} and
  the bound~\eqref{eq:int-bound} is proved in Section~\ref{sec:integrability}.
  The reflection positivity and translation invariance of any accumulation
  point follow as in~{\cite{GH21}}, since
  these properties are inherited from the corresponding properties of the
  approximate measures $(\nu_{\varepsilon, M})_{\varepsilon, M}$.
  Non-Gaussianity follows from~\eqref{eq:int-bound}, as Gaussian measures
  cannot integrate functions exhibiting super-exponential growth.
\end{proof}

We expect that any accumulation point is invariant under all
Euclidean transformations. Unfortunately, a direct characterisation of the
limiting measures remains unavailable. As a result, we must deduce the
properties of the accumulation points indirectly, relying on certain features
of the approximating measures. Notably, since the measures on lattices are not
rotation invariant, establishing rotation invariance of the continuum
limit poses a particularly challenging problem. It is conceivable that this
issue could be addressed using the technique developed
in~{\cite{duchStochasticQuantizationTwoDimensional2024}} in the context of
the $P_2 (\Phi)$ model. However, pursuing this direction lies beyond the scope of
the present paper.

Let us also mention the recent work~{\cite{duchErgodicityInfiniteVolume2025}},
where uniqueness of the limiting measure is established in the small-coupling
regime for the standard $\Phi^4_3$ model. A similar strategy might extend to
the fractional variant, providing a potential route toward a full verification
of the OS axioms, at least in the regime of small coupling.

Our proof strategy introduces a novel combination of renormalisation group
ideas and PDE techniques which we believe may be useful more broadly in the
context of the theory of subcritical singular SPDEs. The proof also applies to
the vector version of the model where the field takes values in $\mathbb{R}^n$
for $n > 1$ and the functional $S_{\varepsilon, M} (\varphi)$ depends on
$\varphi$ in an $O (n)$-symmetric way. We briefly discuss how
to adapt the proof to the vector case in Section~\ref{sec:vector}.

Theorem~\ref{thm:main} gives a construction of a model of Euclidean quantum
field theory (EQFT) commonly referred to in the literature as the
\textit{fractional $\Phi^4_3$} model. This terminology stems from the form
of the density appearing in the formula for the measure~\eqref{eq:measure} and
the fact that the underlying space is three-dimensional. In the case $s = 1$,
the model reduces to the classical $\Phi_3^4$ theory, which has long been
regarded as a crucial benchmark in constructive quantum field theory.
Foundational results by Glimm and
Jaffe~{\cite{glimmPositivityPhiHamiltonian1973}}, Feldman and
Osterwalder~{\cite{feldmanLambdaVarphiField1974,feldmanWightmanAxiomsMass1976}}
and other pioneers of EQFT laid the groundwork for demonstrating the existence
of models satisfying the Wightman axioms for local relativistic QFT through
probabilistic methods based on the Euclidean framework~{\cite{glimmQuantumPhysicsFunctional1987}}. In the
fractional regime, that is, for $s \in (s_c, 1)$, the model we consider was
introduced by Brydges, Mitter, and Scoppola~{\cite{brydgesCriticalF42003}} as
a rigorous foundation for proving the existence of a nontrivial infrared
fixed point via an $\epsilon$-expansion for sufficiently small $\epsilon = s -
s_c > 0$. See also~{\cite{brydges_non_gaussian_1998}} for an analogous result
concerning a related four-dimensional model.

Our proof also extends to the cases $d = 1$ or $d = 2$ and $s
\in (d / 4, 1]$ without requiring any substantial modifications. We stress
that we crucially need the condition of subcriticality (\textit{i.e.}
super-renormalisability), which translates into $s > d / 4$, where $s$ is the
power of the fractional Laplacian. In particular, if $d \geqslant 4$ this
would require $s > 1$. However, it is known that the fractional Laplacian $(-
\Delta)^s$ with $s > 1$ does not satisfy the maximum principle, and as a
result, our a priori estimate does not apply in this regime. Moreover, for $s
> 1$ the measure is not expected to be reflection positive---since the
Gaussian part itself lacks this property---and is therefore of limited
interest from the perspective of constructive quantum field theory. Finally,
we note that for $d \leqslant 4$ it was shown
in~{\cite{panisTrivialityScalingLimits2025}} that at criticality $s = d /
4$ as well as in the supercritical regime $s < d / 4$, any possible continuum
limit is trivial, in the sense that it is Gaussian.

In recent years there has been a renewed interest in EQFTs due to the
development of an alternative approach to proving theorems such as
Theorem~\ref{thm:main}. This new approach is grounded in the basic ideas of
stochastic calculus and it is usually called \textit{stochastic
quantisation}. This term was introduced by Parisi and
Wu~{\cite{parisi_perturbation_1981}} to describe the quantisation of gauge
theories via the construction of a stochastic process that evolves in fictitious
time and whose stationary distribution is the target Euclidean QFT. This
stochastic evolution is a nonlinear stochastic partial differential equation
of a singular kind, for which a particular procedure of renormalisation is
needed to give it a precise meaning. The analysis of such equations requires a
mix of probabilistic and analytic arguments that go beyond the usual framework of
It\^o's stochastic differential equations. For this reason it took some time
before the SPDE community learned how to handle such singular equations and
discovered methods such as regularity structures~{\cite{hairer2014theory}},
paracontrolled
calculus~{\cite{gubinelli2015paracontrolled,catellierParacontrolledDistributions3dimensional2018}},
and the renormalisation group~{\cite{kupiainen2016renormalization}}, which finally
made it possible to tackle the problem of the stochastic quantisation of the $\Phi^4_3$
model. Gubinelli and
Hofmanov{\'a}~{\cite{GH21}} obtained the
analogue of Theorem~\ref{thm:main} with $s = 1$ and a small range of values
below that using a mix of paracontrolled calculus for the small scale
singularities of the equation and coercive estimates to tame the large-scale
fluctuations. We refer the reader to the introduction
to~{\cite{GH21}} for a deeper review of the
literature and the history of constructive QFT, and also to contextualise the
meaning and consequences of Theorem~\ref{thm:main}.

The probability measure $\nu_{\varepsilon, M}$ in~\eqref{eq:measure} is the
equilibrium distribution of the Langevin dynamics governed by the finite
system of SDEs
\begin{equation}
  \mathcal{L}_{\varepsilon} \phi^{(\varepsilon, M)} + \lambda
  {(\phi^{(\varepsilon, M)})^3}  - r_{\varepsilon, M} \phi^{(\varepsilon, M)}
  = \xi^{(\varepsilon, M)}, \label{eq:dynamics-eps-M}
\end{equation}
on $\Lambda_{\varepsilon, M} \assign \mathbb{R} \times
\mathbb{T}_{\varepsilon, M}^d$, where
\[ \mathcal{L}_{\varepsilon} \assign \partial_t + (- \Delta_{\varepsilon})^s +
   m^2 \]
and $\xi^{(\varepsilon, M)}$ is a spacetime white noise satisfying
\[ 
\mathbb{E} [\xi^{(\varepsilon, M)} (t, x) \xi^{(\varepsilon, M)} (s, y)] =
   \delta (t - s) \,\varepsilon^{-d}\mathbbm{1}_{x = y}, \qquad (t, x), (s, y) \in
   \Lambda_{\varepsilon, M} . \]
The constants $m > 0, \lambda > 0, r_{\varepsilon, M} \in \mathbb{R}$ coincide
with the parameters appearing in~\eqref{eq:action}. By standard stochastic
analysis arguments, there exists a unique stationary solution
$\phi^{(\varepsilon, M)}$ of~\eqref{eq:dynamics-eps-M} and we have $\operatorname{Law}
(\phi^{(\varepsilon, M)} (t)) = \nu_{\varepsilon, M}$ for all $t \in
\mathbb{R}$. In what follows, we identify $\phi^{(\varepsilon, M)}$ and
$\xi^{(\varepsilon, M)}$ with periodic functions on
\(\Lambda_{\varepsilon}:=\mathbb{R}\times \mathbb{R}_{\varepsilon}^d\).

The nontrivial step is now to control the solutions to the stochastic
quantisation equation~\eqref{eq:dynamics-eps-M} uniformly as $\varepsilon
\rightarrow 0$ and $M \rightarrow \infty$. At small scales, it is expected
that the nonlinear term in the dynamics is a perturbation of the linearised
equation driven by spacetime white noise. Consequently, in the continuum
limit, the solutions $\phi^{(\varepsilon, M)}$ are expected to converge to
random distributions belonging to Besov spaces of negative regularity slightly
worse than $s - d / 2$, which is the regularity of the Gaussian free field.
This low regularity presents a major analytical challenge, as it complicates
the control of the nonlinear term.

Inspired by the works of
Wilson~{\cite{wilson1971renormalization,wilson1974renormalization}} and
Polchinski~{\cite{polchinskiRenormalizationEffectiveLagrangians1984}} on the
continuous renormalisation group and by the more recent approach introduced by
one of the authors
in~{\cite{duchFlowEquationApproach2021,duchRenormalizationSingularElliptic2022}},
we use a flow equation to effectively describe the solution
$\phi^{(\varepsilon, M)}$ of the SPDE at some spatial scale much larger than
$\varepsilon$ (see Kupiainen~{\cite{kupiainen2016renormalization}} for a
discrete counterpart). Let $\phi_{\sigma}^{(\varepsilon, M)}$ denote a
description of the solution at a scale of order
\[ \llbracket \sigma \rrbracket  \assign (1 - \sigma) \gg \varepsilon > 0 \]
for some $\sigma \in (0, 1)$. The flow equation approach consists in deriving
a parabolic equation for $\phi^{(\varepsilon, M)}_{\sigma}$ of the form
\begin{equation}
  \mathcal{L}_{\varepsilon} \phi^{(\varepsilon, M)}_{\sigma} =
  F_{\sigma}^{(\varepsilon, M)} (\phi^{(\varepsilon, M)}_{\sigma}) .
  \label{eq:effective-spde}
\end{equation}
Here, $\psi \mapsto F_{\sigma}^{(\varepsilon, M)} (\psi)$ is an analytic
functional depending on the noise $\xi^{(\varepsilon, M)}$, called the
\textit{effective force}, such that
\begin{equation}
  F_1^{(\varepsilon, M)} (\psi) = - \lambda \psi^3 + r_{\varepsilon, M} \psi +
  \xi^{(\varepsilon, M)} . \label{eq:final-cond}
\end{equation}
In particular, \eqref{eq:dynamics-eps-M} can be recovered
from~\eqref{eq:effective-spde} for $\sigma = 1$ and $\phi^{(\varepsilon, M)}_1
= \phi^{(\varepsilon, M)}$. The functional $F_{\sigma}^{(\varepsilon, M)}$ can
be obtained by solving a \textit{flow equation}
\begin{equation}
  \partial_{\sigma} F_{\sigma}^{(\varepsilon, M)} =\mathbb{B}_{\sigma}
  (F_{\sigma}^{(\varepsilon, M)}, F_{\sigma}^{(\varepsilon, M)}),
  \label{eq:flow-intro}
\end{equation}
\textit{backwards} for $\sigma \in (\mu, 1]$ with the final
condition~\eqref{eq:final-cond} and where $\mathbb{B}_{\sigma}$ is an
appropriate bilinear operator. The parameter $\sigma \in [0, 1]$ does not have
any specific physical meaning and the spatial scale of the decomposition is
fixed conventionally to be of order $\llbracket \sigma \rrbracket$, that is
$\phi_{\sigma}^{(\varepsilon, M)}$ is expected to fluctuate at spatial scales
of order $\llbracket \sigma \rrbracket$ or larger, and in particular to be a
locally bounded function on $\Lambda_0 \assign \mathbb{R} \times
\mathbb{R}^d$, when extended in some reasonable way from the lattice
$\Lambda_{\varepsilon}$ to the continuum. A key ingredient is the control of
the stochastic process $(F^{(\varepsilon, M)}_{\sigma})_{\sigma}$ solving the
flow equation~\eqref{eq:flow-intro}. Following a simple but crucial
observation
of~{\cite{duchFlowEquationApproach2021,duchRenormalizationSingularElliptic2022}},
this control can be obtained by studying the evolution in the scale parameter
of the cumulants $(\mathcal{F}_{\sigma}^{(\varepsilon, M)})_{\sigma}$ of the
process $(F^{(\varepsilon, M)}_{\sigma})_{\sigma}$, which themselves satisfy a
kind of higher-order deterministic flow equation
\begin{equation}
  \partial_{\sigma} \mathcal{F}_{\sigma}^{(\varepsilon, M)}
  =\mathcal{A}_{\sigma} (\mathcal{F}_{\sigma}^{(\varepsilon, M)})
  +\mathcal{B}_{\sigma} (\mathcal{F}_{\sigma}^{(\varepsilon, M)},
  \mathcal{F}_{\sigma}^{(\varepsilon, M)}), \label{eq:flow-cumul-intro}
\end{equation}
with prescribed final condition $\mathcal{F}_1^{(\varepsilon, M)}$. Upon
choosing this initial condition appropriately by tuning the parameter
$r_{\varepsilon, M}$ in \eqref{eq:final-cond}, one can prove estimates uniform in
$\varepsilon$ and $M$ for the cumulants
$(\mathcal{F}_{\sigma}^{(\varepsilon, M)})_{\sigma}$ and therefore, by a
Kolmogorov-type argument, suitable bounds on the effective force
$(F^{(\varepsilon, M)}_{\sigma})_{\sigma}$ uniform as $\varepsilon \rightarrow
0$ and $M \rightarrow \infty$.

The flow equation~\eqref{eq:flow-intro} is bilinear and therefore solvable in
general only in a perturbative regime, e.g. in a small interval $I =
[\bar{\sigma}, 1]$ around the initial condition at $\sigma = 1$ or for small
data (or small time). The size of this perturbative region depends crucially
on the noise $\xi^{(\varepsilon, M)}$. Although this dependence can be
controlled uniformly in $\varepsilon$ and $M$, large noise fluctuations can
make the region arbitrarily small and restrict the resulting existence theory
to \textit{local} results. A similar
limitation is present in the work of
Kupiainen~{\cite{kupiainen2016renormalization}} who, instead, uses a discrete
renormalisation group (RG) iteration, and more generally in all the other
approaches which use an expansion of solutions in order to resolve the
singular terms and control the limit as $\varepsilon \rightarrow 0$, e.g. in
regularity structures and also in paracontrolled calculus. This difficulty is
the signal of the ``large field problem'', well known in constructive EQFT.

From the point of view of the stochastic quantisation equation, the
large-field problem can be addressed by exploiting the coercivity of the
nonlinear term which pulls the solution back from infinity. While this
observation is standard in PDE theory, it still requires some nontrivial
adaptation to be effective for singular SPDEs. Mourrat and Weber~{\cite{mourrat2017global}} were the first to solve this problem in their proof of
global existence for the $\Phi^4_2$ dynamics on the full space with the usual
Laplacian diffusion term and subsequently Gubinelli and Hofmanov{\'a} in the
context of paracontrolled analysis of $\Phi^4$
models~{\cite{gubinelliGlobalSolutionsElliptic2019,GH21}} including the
three-dimensional parabolic setting. Moinat and Weber~{\cite{MoinatWeber2020}}
proved the so-called spacetime localisation property for the dynamical
$\Phi^4_3$ model in the framework of regularity structures. This result was
further extended by Chandra, Moinat, and
Weber~{\cite{chandraPrioriBoundsPhi2023}} to cover the full subcritical
regime. In the latter work, the authors modify the covariance of the noise,
rather than the diffusion term, in order to explore regularities arbitrarily
close to the critical threshold. Although their estimates suffice to establish
tightness of the invariant measures associated with their SPDE, these measures
are not explicit, and it is not clear whether they are reflection positive.
Consequently, their role in the stochastic quantisation of Euclidean QFTs is,
at present, not fully understood.

In the broader context of global solutions for singular SPDEs, we also
mention recent work of Chandra, Feltes, and
Weber~{\cite{chandraPrioriBounds2d2024}}, which establishes results for the
stochastic quantisation of the two-dimensional sine-Gordon model on a periodic
domain for parameter values slightly above the first threshold. This has since
been extended up to the second threshold
in~{\cite{bringmannGlobalWellposednessDynamical2025}}.
Moreover,~{\cite{bringmannGlobalWellposednessStochastic2024}} proves long-time
well-posedness of the two-dimensional Abelian Higgs model.

During the revision of the present paper, a new
preprint~{\cite{esquivelPrioriBoundsDynamic2024}} appeared in which a priori
bounds for the fractional $\Phi_3^4$ model on the three-dimensional torus were
established in the full subcritical regime. The main distinctions between our
results and those of~{\cite{esquivelPrioriBoundsDynamic2024}} are that we
employ the flow-equation framework and work on an infinite lattice,
whereas the authors of~{\cite{esquivelPrioriBoundsDynamic2024}} use the regularity structures
approach and study the continuum model on a torus with mollified noise. As far
as we know, mollification of the noise is not a feasible strategy for
establishing stochastic quantisation of an EQFT, in particular with respect to
reflection positivity.

The main contribution of our work is the identification of a framework in
which the flow equation method is combined with PDE estimates for the
dynamics. This hybrid approach yields a powerful variant of the
renormalisation group (RG). Instead of requiring an exact solution to the flow
equation~\eqref{eq:flow-intro}, it suffices to construct a suitable
approximate solution $(F_{\sigma} )_{\sigma} = (F_{\sigma}^{(\varepsilon,
M)})_{\sigma}$ satisfying~\eqref{eq:final-cond} for which the quantity
\begin{equation}
  H_{\sigma} \assign \partial_{\sigma} F_{\sigma}  +\mathbb{B}_{\sigma}
  (F_{\sigma} , F_{\sigma} ), \label{eq:approx-flow}
\end{equation}
is small enough in an appropriate sense. The price to pay for this
approximation is a \textit{remainder} term $R^{(\varepsilon, M)}_{\sigma}$
in the SPDE which now reads as a \textit{system} of two equations:
\begin{equation}
  \left\{\begin{array}{l}
    \mathcal{L}_{\varepsilon} \phi _{\sigma} = J_{\sigma} [F_{\sigma}  (\phi
    _{\sigma}) + R_{\sigma} ],\\
    \partial_{\sigma} R_{\sigma}  = -H_{\sigma}  (\phi _{\sigma}) - \mathD
    F_{\sigma}  (\phi _{\sigma}) ((\partial_{\sigma} G _{\sigma}) R_{\sigma}
    ), \qquad R_1  = 0,
  \end{array}\right.  \label{eq:effective-spde-remainder}
\end{equation}
for the pair of scale-dependent functions
\[ (\phi _{\sigma}, R _{\sigma})_{\sigma} = (\phi^{(\varepsilon, M)}_{\sigma},
   R^{(\varepsilon, M)}_{\sigma})_{\sigma} . \]
Here $G  \assign \mathcal{L}_{\varepsilon}^{- 1}$ is the fractional heat
kernel, with $\mathcal{L}_{\varepsilon}^{- 1}$ denoting the inverse of
$\mathcal{L}_{\varepsilon}$ as defined in~\eqref{eq:ell-minus-one} below,
$(J_{\sigma})_{\sigma}$ is a family of smoothing operators (see
Definition~\ref{def:regularizing-kernels}) and
\[ (G _{\sigma})_{\sigma \in (0, 1)} \assign (J_{\sigma} G)_{\sigma \in (0,
   1)} , \]
is a scale decomposition of $G$. Moreover, $\mathD F  (\psi) \tilde{\psi}$
denotes the functional derivative of a functional $F$ at $\psi$ in the
direction of $\tilde{\psi}$.

One can prove that the term $F_{\sigma}  (\phi _{\sigma})$ retains the
coercive structure of~\eqref{eq:final-cond}, that is,
\[ J_{\sigma} F_{\sigma}  (\phi _{\sigma}) = - \lambda \phi_{\sigma}^3  + Q
   _{\sigma} (\phi _{\sigma}), \]
where $Q _{\sigma} (\phi _{\sigma})$ is ``smaller'' than the cubic
contribution provided $\llbracket \sigma \rrbracket \ll 1$. This, together with
the linearity in $R^{(\varepsilon, M)}_{\sigma}$ of the second equation
of~\eqref{eq:effective-spde-remainder}, makes this system amenable to standard
PDE techniques: by choosing $\llbracket \sigma \rrbracket \ll 1$ one can
control the noncoercive part $Q _{\sigma} (\phi _{\sigma})$ of the effective
force using the coercive part ${- \lambda \phi_{\sigma}^3} $, thereby
resolving the large-field problem. At the same time, for $\llbracket \sigma
\rrbracket > 0$, we have uniform estimates for $Q _{\sigma}, \mathD F
_{\sigma}, H _{\sigma}$ as $\varepsilon \rightarrow 0$ and $M \rightarrow
\infty$, provided the renormalisation constant $r_{\varepsilon, M}$ appearing
in the boundary condition~\eqref{eq:final-cond} for the effective force is
chosen appropriately. This allows full control of~\eqref{eq:dynamics-eps-M}
and the proof of tightness of the laws of the processes $(\phi^{(\varepsilon,
M)})_{\varepsilon, M}$, and therefore of the family of measures
$(\nu_{\varepsilon, M})_{\varepsilon, M}$.

The implementation of this plan has to deal with two main technical
difficulties:
\begin{enumerate}[(a)]
  \item We decompose the fractional heat kernel $G$ scale by scale
  (cf.~\eqref{eq:est-kernel-G}), using
  \[ (\partial_{\sigma} G_{\sigma})_{\sigma}. \]
  This produces kernels with limited spacetime decay (see
  Lemma~\ref{lem:Gdot-decay-old}), reflecting the restricted smoothness of $G$ away
  from the origin. This algebraic decay of $\partial_{\sigma} G
  _{\sigma}$ propagates to the kernels of the effective force, necessitating a
  careful choice of weighted spaces for both the solutions and the kernels. To
  handle this limited decay, we localise the relevant operators and kernels at
  several points in the proof. The main instance is the localisation procedure
  for the relevant cumulants, which governs the flow of the renormalisation
  constants; see Appendix~\ref{app:localization}. It would be
  interesting to devise an alternative strategy to bypass this problem with
  some other scale decomposition (or an additional localisation procedure).
  
  \item The natural setting of the analysis provides only negative spacetime
  regularity for the solutions of the SPDE. Such regularity is insufficient
  for stochastic quantisation, since one must be able to compute the marginal
  of the solution at a given time. To address this, it is necessary to work
  with distributional norms defined via smoothing operators $(K_{\mu})_{\mu}$
  that provide only limited time regularisation. Suitable Schauder estimates
  for these operators can be established (see Lemma~\ref{lem:schauder-full}).
  These smoothing operators also influence the definition of the norms for the
  kernels of the effective force (cf. Definition~\ref{def:norm-kernels}) and require
  a careful adaptation of the Kolmogorov-type argument used in
  Lemma~\ref{lem:v-and-K}.
\end{enumerate}

{\noindent}\textbf{Comparison with other approaches.} The possibility of
working with an approximate flow equation makes it easier to compare the RG
approach advocated in this paper (and originally proposed
in~{\cite{kupiainen2016renormalization}}
and~{\cite{duchFlowEquationApproach2021,duchRenormalizationSingularElliptic2022}})
with regularity structures~{\cite{hairer2014theory}} and paracontrolled
distributions~{\cite{gubinelli2015paracontrolled}}. There is a clear parallel
among the various approaches. The flow equation constructs a random object
$F_{\sigma}^{(\varepsilon, M)}$ --- the scale-dependent {\textit{effective
force}} --- which encapsulates the influence of the noise and takes the form
of a finite polynomial built from the noise and the linear part of the
equation. This object corresponds respectively to the {\textit{model}} in the
theory of regularity structures, the {\textit{enhanced noise}} in paracontrolled
calculus, or the {\textit{rough path}} in rough path theory.

While $F_{\sigma}^{(\varepsilon, M)}$ is obtained through a probabilistic
construction, the remainder term $R_{\sigma}^{(\varepsilon, M)}$ is defined
analytically in terms of $F_{\sigma}^{(\varepsilon, M)}$. This deterministic
component mirrors the analytic machinery of regularity structures, the
paracontrolled operators in the paracontrolled calculus, and the sewing lemma
in rough path analysis.

When the parameter $s$ is near its critical value of $s_c = 3 / 4$ the number
of terms which have to be accounted for in the approximation
$F_{\sigma}^{(\varepsilon, M)}$ of the solution of the flow equation grows in
an unbounded manner. A notable advantage of the flow equation approach,
however, lies in its relative insensitivity to this increasing complexity: the
analysis remains compact and largely independent of the distance to
criticality. This analytical efficiency --- the ability to capture the
nonlinear propagation of randomness with minimal combinatorial overhead ---
was first observed by
Polchinski~{\cite{polchinskiRenormalizationEffectiveLagrangians1984}} in his
proof of perturbative renormalisability of the Euclidean $\phi^4_4$ QFT. For a
modern account of this approach to perturbation theory of QFTs, the reader can
consult the book of
Salmhofer~{\cite{salmhoferRenormalizationIntroduction2007}} or
Kopper~{\cite{kopperRenormalizationTheoryBased2007}}.

As we already noted, the application of RG ideas to SPDEs is made efficient by
the observation of one of the
authors~{\cite{duchFlowEquationApproach2021,duchRenormalizationSingularElliptic2022}}
that flow equations can also be used to estimate cumulants directly, thereby
avoiding explicit and cumbersome inductive arguments on trees --- much as
Polchinski's method circumvents the inductive structure of BPHZ
renormalisation. There are further conceptual similarities with recent
developments by Otto, Weber, and
collaborators~{\cite{ottoPrioriBoundsQuasilinear2021,linaresDiagramfreeApproachStochastic2021}},
who use PDE-based arguments to derive probabilistic estimates for the modes in
regularity structures. The flow equation framework, however, offers an
additional advantage: renormalisation conditions naturally appear as boundary
conditions for the corresponding flow equation.

The combination of the flow equation approach with stochastic quantisation in
the context of the construction of EQFT has recently been exploited by Meyer
and one of the authors in~{\cite{gubinelliFBSDEApproachSineGordon2024}} to
study the sine--Gordon model in the full space up to the second threshold and,
by De~Vecchi, Fresta and one of the authors,
in~{\cite{devecchiStochasticAnalysisSubcritical2022}}, to develop a new
approach to Euclidean Fermionic theories. In both papers the stochastic
quantisation is obtained by using a forward-backward stochastic differential
equation together with an approximate analysis of Polchinski's flow equation.
Finally, we mention the work {\cite{duch2024construction}},
where the Polchinski equation was employed to construct the Gross--Neveu
model, a critical fermionic model of Euclidean quantum field theory.

{\noindent}{\textbf{Conclusions. }}Despite the technical difficulties due to
the analysis of the fractional heat equation, we would like to stress that the
present paper is self-contained and presents complete arguments for all
proofs. Compared with the few existing works on singular SPDEs in the full
subcritical regime, we believe our work is the first to address, in a unified
framework, several intricate aspects of the problem:
\begin{enumerate}[(a)]
  \item we present the entire argument --- both analytic and probabilistic ---
  in full detail;
  
  \item we work with an extremely nonlocal SPDE posed on a (semi-)discrete
  space;
  
  \item we establish an a priori bound in the full subcritical regime, valid
  globally in spacetime;
  
  \item we obtain optimal tail estimates for the solutions.
\end{enumerate}
We hope that this work illustrates the strength and flexibility of the
stochastic quantisation approach in tackling, in a genuinely nonperturbative
manner, the construction of Euclidean quantum field models in the subcritical
regime.{\medskip}

{\noindent}{\textbf{Plan of the paper. }}In Section~\ref{sec:sq} we introduce
the main objects of our analysis: the scale decomposition, the spacetime
weighted norms which will be used to control the large fields and all the
intermediate results which are needed in the proof of Theorem~\ref{thm:main}.
The coercive estimates will be proven in Section~\ref{sec:apriori}, while the
approximate flow equation for the effective force will be analysed in
Section~\ref{sec:flow} via the flow equation for the cumulants.
Appendix~\ref{app:technical} contains some technical lemmas and Schauder
estimates tailored to our setting while Appendix~\ref{app:flow-eq} contains
the detailed definition of the various contributions to the flow equations for
the cumulants and their analytic estimates.{\medskip}

\subsection{Preliminaries and notation}\label{sec:notations}

In this section we introduce the main notation that we are going to use
throughout the paper.

We let $\mathbb{R}_0 \assign \mathbb{R}$, $\mathbb{R}_{\varepsilon} \assign
\varepsilon \mathbb{Z}$, $\mathbb{T}_{\varepsilon, M} \assign
(\varepsilon \mathbb{Z}) / (M\mathbb{Z})$ and define the corresponding
spacetime domains
\[ \Lambda_{\varepsilon} \assign \mathbb{R} \times
   \mathbb{R}_{\varepsilon}^d, \qquad \Lambda_{\varepsilon, M} \assign
   \mathbb{R} \times \mathbb{T}_{\varepsilon, M}^d, \qquad d = 3, \quad
   \varepsilon \in 2^{-\mathbb{N}_0}, \quad M \in \mathbb{N}_+ . \]
The assumption that $\varepsilon \in 2^{-\mathbb{N}_0}$ ensures that
$M\mathbb{Z} \subset \varepsilon \mathbb{Z}$ for every $M \in \mathbb{N}_+$.
These domains are continuous in the time direction and discrete in the $d$
spatial dimensions. For a measurable function $f : \Lambda_{\varepsilon}
\rightarrow \mathbb{R}^n$, we write
\[ \| f \| \assign \| f \|_{L^{\infty} (\Lambda_{\varepsilon})} . \]
Given a nonnegative weight $w \in C (\Lambda_0)$ we denote by $C
(\Lambda_{\varepsilon}, w)$ the space of continuous functions $f :
\Lambda_{\varepsilon} \rightarrow \mathbb{R}$ such that
\[ \| f \|_{L^{\infty} (w)} \assign \sup_{z \in \Lambda_{\varepsilon}} | w (z)
   f (z) | < \infty . \]
We always identify functions $f : \Lambda_{\varepsilon, M} \rightarrow
\mathbb{R}$ with their spatially periodic extensions $f :
\Lambda_{\varepsilon} \rightarrow \mathbb{R}$. In particular, we write $C
(\Lambda_{\varepsilon, M}, w)$ for the subspace of $C (\Lambda_{\varepsilon},
w)$ consisting of functions that are periodic in space with period $M$. We
denote by $\mathcal{S} (\Lambda_{\varepsilon})$ the space of smooth rapidly
decreasing functions over $\Lambda_{\varepsilon}$. The Fourier transform of $f
\in \mathcal{S} (\Lambda_{\varepsilon})$ is defined as
\[ \hat{f} (\omega, k) \assign \int_{\Lambda_{\varepsilon}} e^{- i (\omega t +
   k \cdot x)} f (t, x)\,\mathd t \mathd x, \qquad (\omega, k) \in
   \Lambda_{\varepsilon}^{\ast} \assign \mathbb{R} \times
   (\mathbb{R}_{\varepsilon}^d)^{\ast}, \]
where
\[ (\mathbb{R}_{\varepsilon}^d)^{\ast} \assign (- \pi / \varepsilon, \pi /
   \varepsilon]^d \]
is the dual of the group $\mathbb{R}_{\varepsilon}^d$. The integral over
$\mathbb{R}_{\varepsilon}^d$ is understood with respect to the counting measure
(still written as $\mathd x$) with normalisation
\[ \int_{\mathbb{R}_{\varepsilon}^d} \, f (x)\,\mathd x \assign \varepsilon^d
   \sum_{x \in \mathbb{R}^d_{\varepsilon}} f (x), \]
which ensures weak convergence to the Lebesgue measure as $\varepsilon
\rightarrow 0$. As usual, the Fourier transform extends to $\mathcal{S}'
(\Lambda_{\varepsilon})$ by duality. The inverse Fourier transform is given by
\[ f (t, x) = \int_{\Lambda^{\ast}_{\varepsilon}}  \hat{f} (\omega, k) e^{i
   (\omega t + k \cdot x)} \frac{\mathd \omega \mathd k}{(2 \pi)^{d + 1}},
   \qquad (t, x) \in \Lambda_{\varepsilon} . \]
For $\varepsilon \geqslant 0$, we define the Laplacian $\Delta_{\varepsilon}$
as the Fourier multiplier on $\mathcal{S}' (\Lambda_{\varepsilon})$ with
symbol
\[ (\mathbb{R}_{\varepsilon}^d)^{\ast} \ni k \mapsto - q_{\varepsilon}^2 (k)
   \in \mathbb{R}, \]
where
\begin{equation}
  q_0 (k) = | k |, \qquad q_{\varepsilon} (k) \assign \left[ \sum_{i = 1}^d
  \left( \frac{2}{\varepsilon} \sin \left(\frac{\varepsilon k_i}{2}\right)
  \right)^2 \right]^{1 / 2}, \qquad \varepsilon > 0. \label{eq:discrete-symbol}
\end{equation}
Note that for $\varepsilon > 0$, the operator $\Delta_{\varepsilon}$ coincides
with the standard nearest-neighbour discrete Laplacian.

\paragraph*{The fractional Laplacian.} For $s \in (0, 1)$ the fractional
Laplacian $(- \Delta_{\varepsilon})^s$ is defined as the Fourier multiplier
with symbol $(\mathbb{R}_{\varepsilon}^d)^{\ast} \ni k \mapsto
q_{\varepsilon}^{2 s} (k)$. In particular it is self-adjoint and positive in
$L^2 (\mathbb{R}^d_{\varepsilon})$ and for $s \in (0, 1)$ it has the
(discrete, when $\varepsilon > 0$) heat--kernel
representation~{\cite{kwasnickiTenEquivalentDefinitions2017}}
\begin{equation}
  (- \Delta_{\varepsilon})^s f = C_s  \int_{\mathbb{R}_+} (f - e^{\theta
  \Delta_{\varepsilon}} f) \theta^{- 1 - s} \mathd \theta,
  \label{eq:fractional-laplacian-heat-kernel}
\end{equation}
with the constant $C_s = | \Gamma (- s) |^{- 1}$. In the continuum, the
fractional Laplacian has, for $s \in (0, 1)$, the integral
representation~{\cite{kwasnickiTenEquivalentDefinitions2017}}:
\begin{equation}
  (- \Delta_0)^s f (x) = C_{d, s}  \text{PV} \int_{\mathbb{R}^d} \frac{f (x) -
  f (y)}{| x - y |^{d + 2 s}} \mathd y, \qquad x \in \mathbb{R}^d,
  \label{eq:integral-repr-frac-lapl-cont}
\end{equation}
where $C_{d, s}$ is a universal constant. In the discrete setting a similar
formula holds~{\cite{ciaurriFractionalDiscreteLaplacian2015}}:
\begin{equation}
  (- \Delta_{\varepsilon})^s f (x) = \varepsilon^d  \sum_{y \in
  \mathbb{R}_{\varepsilon}^d | y \neq x} H^{(\varepsilon)}_s (x - y)  (f (x) -
  f (y)), \qquad x \in \mathbb{R}_{\varepsilon}^d,
  \label{eq:integral-rep-frac-lap-discrete}
\end{equation}
where the kernel $H^{(\varepsilon)}_s : \mathbb{R}_{\varepsilon}^d \rightarrow
\mathbb{R}$ is positive and such that $H^{(\varepsilon)}_s (0) = 0$,
$H^{(\varepsilon)}_s (x) = H^{(\varepsilon)}_s (- x)$ and
\[ | H^{(\varepsilon)}_s (x) | \leqslant C_{d, s}'  | x |^{- d - 2 s}, \qquad
   x \in \mathbb{R}_{\varepsilon}^d, \]
uniformly in $\varepsilon > 0$ for some constant $C_{d, s}' > 0$.
Note that in our
notation~\eqref{eq:integral-rep-frac-lap-discrete} can be equivalently written
as
\[ (- \Delta_{\varepsilon})^s f (x) =
   \int_{\mathbb{R}_{\varepsilon}^d} H^{(\varepsilon)}_s (x - y)  (f (x) - f
   (y))\,\mathd y. \]
For $\varepsilon > 0$, we can encode the
representation~\eqref{eq:integral-rep-frac-lap-discrete} of the fractional
Laplacian, via a positive measure $\mu^{(\varepsilon)}_s$ on
$\Lambda_{\varepsilon} \times \Lambda_{\varepsilon}$ for which
\[ \langle f, (- \Delta_{\varepsilon})^s g \rangle =
   \int_{\Lambda_{\varepsilon} \times \Lambda_{\varepsilon}} f (z)  (g (z) - g
   (z')) \mu^{(\varepsilon)}_s (\mathd z \mathd z') . \]
We also define the kernel $\mu^{(\varepsilon)}_s (z, \mathd z')$ arising from
the disintegration of $\mu^{(\varepsilon)}_s$ with respect to the measure
$\mathd z$ on $\Lambda_{\varepsilon}$. Specifically, we have
\begin{equation}\begin{aligned}
  \mu^{(\varepsilon)}_s (\mathd z \mathd z') &= \delta (t - u)
  H^{(\varepsilon)}_s (x - y)\,\mathd z \mathd z', \\ \mu^{(\varepsilon)}_s
  (z, \mathd z') &= \delta (t - u) H^{(\varepsilon)}_s (x - y)\,\mathd z',
  \end{aligned}\label{eq:int-kernel}\end{equation}
where $z = (t, x)$ and $z' = (u, y) .$ The kernel $\mu^{(\varepsilon)}_s$ is
symmetric, that is,
\[ \mu^{(\varepsilon)}_s (z, \mathd z')\,\mathd z = \mu^{(\varepsilon)}_s (z',
   \mathd z)\,\mathd z' . \]
With this notation in place, the following Leibniz-type
formula with remainder holds:
\begin{equation}
  \begin{aligned}
    I^{(\varepsilon)}_s (f, g)
    &\assign f (- \Delta_{\varepsilon})^s g + g (- \Delta_{\varepsilon})^s f
    - (- \Delta_{\varepsilon})^s (f g)\\
    &= \int_{\Lambda_{\varepsilon}} (f (\Cdot) - f (z))
    (g (\Cdot) - g (z))
    \mu^{(\varepsilon)}_{ s} (\Cdot, \mathd z) .
  \end{aligned}
  \label{def-I}
\end{equation}
Let us also introduce the fractional difference operator
\begin{equation}
  \mathfrak{D}^{(\varepsilon)}_s (f) (z) \assign
  I^{(\varepsilon)}_s (f, f)^{1 / 2} (z) = \left[ \int_{\Lambda_{\varepsilon}}
  [f (z') - f (z)]^2 \mu^{(\varepsilon)}_s (z, \mathd z') \right]^{1 / 2},
  \label{eq:D}
\end{equation}
which behaves like a derivative of order $s$. In particular, using that, for
all $\delta > 0$, one has the estimates
\begin{align*}
  \int_{B (z, \delta)} [f (z') - f (z)]^2
  \mu^{(\varepsilon)}_s (z, \mathd z')
  &\lesssim \| \nabla_{\varepsilon} f \|^2 \delta^{2 - 2 s},\\
  \int_{B (z, \delta)^c} [f (z') - f (z)]^2
  \mu^{(\varepsilon)}_s (z, \mathd z')
  &\lesssim \| f \|^2 \delta^{- 2 s},
\end{align*}
where
\[ B (z, \delta) \assign \{ z' \in \Lambda_{\varepsilon} |  | z - z' |_s
   \leqslant \delta \} \]
denotes the ball of radius $\delta > 0$. The discrete gradient is defined as
\[ \nabla_{\varepsilon} f = (\partial^{1+}_\varepsilon f,
   \partial^{1-}_\varepsilon f, \ldots,
   \partial^{d+}_\varepsilon f, \partial^{d-}_\varepsilon f), \]
where, for \(i \in \{1, \ldots, d\}\),
\[
  \partial^{i\pm}_\varepsilon f (x)
  \assign
  \pm \frac{f (x \pm \varepsilon e_i) - f (x)}{\varepsilon},
  \qquad x \in \mathbb{R}^d_\varepsilon .
\]
For a spacetime multi-index
\[
  A \in \mathbb{N}_0^{\{ 0, 1\pm, \ldots, d\pm\}},
\]
we write its spatial part as
\[
  \bar{A} \assign (A_{i,+}, A_{i,-})_{i = 1}^d
  \in \mathbb{N}_0^{\{1\pm, \ldots, d\pm\}},
\]
and set
\[
  |\bar{A}| \assign \sum_{i = 1}^d (A_{i,+} + A_{i,-}),
  \qquad
  |A|_s \assign 2 s A_0 + |\bar{A}|.
\]
We then define
\[
  \partial_\varepsilon^{\bar{A}}
  \assign
  \prod_{i = 1}^d
  (\partial^{i+}_\varepsilon)^{A_{i,+}}
  (\partial^{i-}_\varepsilon)^{A_{i,-}},
  \qquad
  \partial_\varepsilon^A
  \assign
  \partial_t^{A_0}\partial_\varepsilon^{\bar{A}}.
\]
Whenever the lattice spacing \(\varepsilon\) is fixed and no confusion can
arise, we suppress the subscript \(\varepsilon\) and write
\(\partial^{i\pm}\), \(\partial^{\bar{A}}\), and \(\partial^A\). Finally,
\(\| \nabla_{\varepsilon} f \|\) denotes the supremum norm. Choosing
\(\delta = \| f \| / \| \nabla_{\varepsilon} f
\|\) yields
\begin{equation}
  \mathfrak{D}^{(\varepsilon)}_s (f) (z) \lesssim \|
  \nabla_{\varepsilon} f \| \delta^{1 - s} + \| f \| \delta^{- s} \lesssim \|
  \nabla_{\varepsilon} f \|^s  \| f \|^{1 - s} . \label{eq:D-est}
\end{equation}
\begin{remark}
  Although the kernel representation of the fractional Laplacian in the
  continuum case $\varepsilon = 0$ requires the presence of the principal
  value, the above considerations and the results presented below extend to
  the continuum with only minor modifications in the proofs, or simply by
  taking the $\varepsilon \rightarrow 0$ limit in the relevant inequalities.
  We emphasise, however, that the main results of the paper do not rely on
  these continuum extensions. In particular, the limit $\varepsilon
  \rightarrow 0$ in Theorem~\ref{thm:main} is obtained via tightness, using
  the a priori estimates established uniformly for $\varepsilon \in
  2^{-\mathbb{N}_0}$.
\end{remark}

\begin{remark}
  A basic observation is that the fractional Laplacian (whether in the
  continuum or on the lattice) satisfies an inequality under the action of
  convex functions. Let $s \in (0, 1],$ $\Phi : \mathbb{R} \rightarrow
  \mathbb{R}$ be a convex function and $\Phi'$ one of its sub-differentials,
  then, for any $\varepsilon \geqslant 0$ and $u \in C
  (\Lambda_{\varepsilon})$, we have
  \[ (- \Delta_{\varepsilon})^s \Phi (u) \leqslant \Phi' (u)  (-
     \Delta_{\varepsilon})^s u. \]
  Indeed, let $\Phi : \mathbb{R} \rightarrow \mathbb{R}$ be a convex function,
  then
  \[ \Phi (a) - \Phi (b) \leqslant \Phi' (a)  (a - b), \qquad a, b \in
     \mathbb{R}, \]
  so if $u : \Lambda_{\varepsilon} \rightarrow \mathbb{R}$ is a continuous and
  bounded function, we have
  \[ \Phi (u) - e^{\theta \Delta_{\varepsilon}} \Phi (u) \leqslant \Phi' (u) 
     (u - e^{\theta \Delta_{\varepsilon}} u), \]
  since $e^{\theta \Delta_{\varepsilon}}$ has a positive definite probability
  kernel. The claimed inequality now follows
  from~\eqref{eq:fractional-laplacian-heat-kernel} in the case
  $s \in (0, 1)$. The case $s = 1$ is elementary. For $s > 1$ the result is
  not true. It is clear that the same proof works for $\varepsilon = 0$
  with some additional regularity assumption. We will
  incorporate this idea in the proof of the key
  Lemma~\ref{lem:est-apriori-laplacian}, below.
\end{remark}

The operator $G = G^{(\varepsilon)} =\mathcal{L}^{- 1}_{\varepsilon}$ is
defined as
\begin{equation}
  (\mathcal{L}^{- 1}_{\varepsilon} f) (t, \Cdot) \assign \int_{- \infty}^t
  e^{- (m^2 + (- \Delta_{\varepsilon})^s)  (t - u)} f (u, \Cdot)\,\mathd u,
  \qquad t \in \mathbb{R}, \label{eq:ell-minus-one}
\end{equation}
and will be applied to continuous functions on $\Lambda_{\varepsilon}$ with at
most a limited polynomial growth in spacetime. Indeed, on account of Lemma~5.4
of~{\cite{Grigoryan-heat-kernel-metric-meas-sp}} together with the argument of
Section~1 of~{\cite{Grigoryan-Telcs-Sub-Gaussian-estimates}}, the kernel $G  (t,
x)$ of $G$ satisfies
\begin{equation}\label{eq:est-kernel-G}
	G(t,x)
	\lesssim
	\mathbbm{1}_{{t\geq0}}\,\varepsilon^{-d}\mathbbm{1}_{{x=0}}\,
	e^{-m^2t-ct/\varepsilon^{2s}}
	+
	\frac{\mathbbm{1}_{{t\geq0}}\,t\,e^{-m^2t}}
	{\bigl(t^{1/(2s)}+\varepsilon+|x|\bigr)^{d+2s}},
\end{equation}
uniformly in $\varepsilon \geqslant 0$. Here
$\mathbbm{1}_{t \geqslant 0} = 0$ if $t < 0$ and $1$
if $t \geqslant 0$. If $s = 1$, then the above estimate is not optimal and the
following bound
\begin{equation}\label{eq:est-kernel-G-s-1}
	G(t,x)
	\lesssim
	\mathbbm{1}_{{t\geq0}}\,\varepsilon^{-d}\mathbbm{1}_{{x=0}}
	e^{-m^2t-ct/\varepsilon^2}
	+
	\frac{\mathbbm{1}_{{t\geq0}}\,t\, e^{-m^2t}}{( t^{1/2}+\varepsilon)^{d+2}}
	\exp\left(
	-c\frac{|x|^2}{t+\varepsilon|x|}
	\right).
\end{equation}
holds true uniformly in $\varepsilon \geqslant 0$.

\paragraph*{Parameters.} During the subsequent analysis, we
shall introduce several parameters
\begin{equation}
  \alpha, \beta, \gamma, \delta, \theta, \vartheta, \varrho, \nu, \kappa,
  \bar{\kappa}, \kappa_{\mathfrak{o}}, \bar{\ell}, \bar{k}, a, \flat .
  \label{eq:parameters}
\end{equation}
Although their precise values are not specified at this stage, all these
parameters are to be regarded as fixed once and for all. Their choice depends
solely on the power $s \in (3 / 4, 1)$ of the fractional Laplacian. Since the
specific constraints determining these values will emerge later in the
analysis, we postpone their detailed specification to
Sections~\ref{sec:bounds-on-parameters} and~\ref{sec:post-processing}
below.

\paragraph*{Space-time weights.}

\begin{definition}
  \label{def:distance}We define the fractional parabolic norm by
  \begin{equation}
    | z |_s \assign | z_0 |^{1 / 2 s} + {| \bar{z} | } , \qquad z = (z_0,
    \bar{z}) \in \Lambda_0 =\mathbb{R} \times \mathbb{R}^3,
    \label{eq:dist-parabolic}
  \end{equation}
  where $| \Cdot |$ denotes the usual Euclidean norm on $\mathbb{R}$ and
  $\mathbb{R}^3$. For $\varkappa > 0$ and $z = (z_0, \bar{z}) \in \Lambda_0$
  we denote by
  \[ \varkappa .z \assign (\varkappa^{2 s} z_0, \varkappa \bar{z}) \in
     \Lambda_0 \]
  the fractional parabolic rescaling, satisfying
  \[ | \llbracket \mu \rrbracket^a .z |_s = \llbracket \mu \rrbracket^a | z
     |_s . \]
\end{definition}

\begin{definition}
  \label{def:weights}We introduce the following Japanese brackets:
  \begin{equation}
    \langle z \rangle_s \assign (1 + | z_0 |^{1 / s} + | \bar{z} |^2)^{1 / 2}
    \qquad \text{and} \qquad \langle \bar{z} \rangle \assign (1 + | \bar{z}
    |^2)^{1 / 2}
    \label{eq:Japanese}
  \end{equation}
  for $z = (z_0, \bar{z}) \in \Lambda_0$, where, as in the previous definition, $| \Cdot |$ denotes the usual
  Euclidean norm on $\mathbb{R}$ and $\mathbb{R}^3$.
  \begin{enumerate}[(a)]
    \item Let $a > 1$ and $\nu \in (0, 1 / 3)$.
    
    \item Let $(\chi_i : \Lambda_0 \rightarrow \mathbb{R}_+)_{i \geqslant -
    1}$ be a dyadic partition of unity on $\Lambda_0$ with $\chi_i$ supported
    on an annulus of radius $\sim 2^{ai}$ for $i \geqslant 0$, $\chi_{- 1}$
    supported in a ball of radius $\approx 1$ and $\sum_{i \geqslant - 1}
    \chi_i = 1$.
    
    \item Define $(\mu_i)_{i \geqslant - 1} \subseteq [1 / 2, 1)$ by
    $\llbracket \mu_i \rrbracket = 2^{- i - 2}$, $i \geqslant - 1$.
    
    \item Let $\zeta : \Lambda_0 \rightarrow \mathbb{R}$ be a weight defined
    by
    \[ \zeta (z) \assign \langle z \rangle_s^{- 1}, \qquad z \in \Lambda_0 .
    \]
    For $\mu \in [0, 1]$, we also introduce the associated rescaled weights
    \[ \zeta_{\mu} (z) \assign \langle \llbracket \mu \rrbracket^a .z
       \rangle_s^{- 1}, \qquad \rho_{\mu} \assign \zeta_{\mu}^{\nu} = \langle
       \llbracket \mu \rrbracket^a .z \rangle_s^{- \nu}, \qquad z \in
       \Lambda_0 . \]
  \end{enumerate}
\end{definition}

\begin{remark}
  Concerning the parameters $\nu$ and $a$ introduced above, we impose the
  relation
  \begin{equation}
   a \nu = \gamma, \label{eq:a-nu-gamma}
  \end{equation}
  where the parameter $\gamma \in (0,s)$ will be fixed in
  Section~\ref{sec:post-processing}. In particular, throughout our analysis we
  shall be concerned with very small values of $\nu > 0$ and with $a > 1$.
\end{remark}

\begin{remark}
  \label{rem:reg-weight}The form of the Japanese bracket~\eqref{eq:Japanese}
  is motivated by the requirement that the weight function $\zeta$ (and its
  powers) be $C^1$ in the time variable and $C^2$ in the spatial variables.
  This regularity is essential in the a priori estimates established below
  (see Section~\ref{sec:apriori}). Indeed, one readily verifies that
  \[ \partial_t \zeta (t, x) = - \frac{| t |^{(1 - s) / s} \partial_t | t |}{2
     s (1 + | t |^{1 / s} + | x |^2)^{3 / 2}}, \qquad (t, x) \in \Lambda_0, \]
  which is continuous for $s < 1$. An analogous argument applies to all
  spatial derivatives.
\end{remark}

\begin{remark}
  \label{rem:weight-properties}Our weights satisfy the following properties.
  \begin{enumerate}[(a)]
    \item We have
    \[ \zeta  (z) \zeta ^{- 1} (z_1) \lesssim \zeta ^{- 1} (z - z_1), \]
    uniformly over $z, z_1 \in \Lambda_0$. Consequently,
    \begin{align*}
      \zeta_{\mu} (z) \zeta_{\mu}^{- 1} (z_1) &\lesssim  \zeta_{\mu}^{- 1}
      (z - z_1) \leqslant \zeta^{- 1} (z - z_1),\\
      \rho_{\mu} (z) \rho_{\mu}^{- 1} (z_1) &\lesssim  \rho_{\mu}^{- 1} (z -
      z_1) \lesssim \rho^{- 1}_0 (z - z_1),
    \end{align*}
    uniformly over $z, z_1 \in \Lambda_0$ and $\mu \in [0, 1]$.
    
    \item We have
    \[ \chi_i \zeta_{\mu}^{- 1} \lesssim (1 + (\llbracket \mu \rrbracket^a
       \llbracket \mu_i \rrbracket^{- a})^2)^{1 / 2} \lesssim 1, \]
    uniformly in $i \geqslant - 1$ and $\mu \geqslant \mu_i$.
  \end{enumerate}
\end{remark}

\paragraph*{Scale decomposition.} Let us introduce a scale decomposition of
spacetime functions parametrised by $\sigma \in [0, 1]$, where we set
\[ \llbracket \sigma \rrbracket \assign (1 - \sigma) \]
for convenience. The value $\sigma = 1$ corresponds to allowing fluctuations
at all scales, whereas $\sigma \in (1/2,1)$ corresponds to allowing fluctuations only at spatial scales $\gtrsim \llbracket
\sigma \rrbracket$ or equivalently at Fourier scales $\lesssim \llbracket
\sigma \rrbracket^{- 1}$.

\begin{definition}
  \label{def:function-j}
  Consider a smooth and compactly supported function $j
  : \mathbb{R} \rightarrow \mathbb{R}_+$ such that
  \[ j (\eta) = \left\{ \begin{aligned}
       1 \qquad \operatorname{if} \qquad | \eta | \leqslant 1,\\
       0 \qquad \operatorname{if} \qquad | \eta | \geqslant 2.
     \end{aligned} \right. \]
  For $\ell = 0, 1, 2, \ldots$ denote
  \[ j_{\sigma, \ell} (\eta) \assign j (2^{- \ell} \sigma^{- 1}  \llbracket
     \sigma \rrbracket \eta), \qquad \eta \in \mathbb{R}, \]
  and let $j_{\sigma} \assign j_{\sigma, 0}$ and $\tilde{j}_{\sigma} \assign
  j_{\sigma, 1}$.
\end{definition}

\begin{remark}
  Note that $j_{\sigma, \ell} (\eta) j_{\sigma, \ell'} (\eta) = j_{\sigma,
  \ell} (\eta)$ for $0 \leqslant \ell < \ell'$.
\end{remark}

\begin{definition}
  The family $(J_{\sigma} )_{\sigma \in (0, 1)} =
  (J_{\sigma}^{(\varepsilon)})_{\sigma \in (0, 1)}$ of Fourier multipliers
  acting on distributions is defined as
  \begin{equation}
    \label{eq:cutoff}
    \begin{aligned}
      J_{\sigma}  f (t, x) &\assign
      \int_{\Lambda^{\ast}_{\varepsilon}} j_{\sigma} (| \omega |^{1 / 2 s})
      j_{\sigma} (q_{\varepsilon} (k))  \hat{f} (\omega, k) e^{i (\omega t + k
      \cdot x)}  \frac{\mathd \omega \mathd k}{(2 \pi)^{d + 1}},
    \end{aligned}
  \end{equation} for $(t, x) \in \Lambda_{\varepsilon}$ and 
  $f \in \mathcal{S}'(\Lambda_{\varepsilon})$,
  where $q_{\varepsilon}^2 (k)$ is the symbol of the Laplacian introduced
  in~\eqref{eq:discrete-symbol}. In addition, we define $(\tilde{J}_{\sigma,
  \ell})_{\sigma \in (0, 1), \ell = 0,1,2, \ldots}$ by~\eqref{eq:cutoff} with
  the function $j_{\sigma}$ replaced by $j_{\sigma, \ell}$. In particular, $J_\sigma=\tilde J_{\sigma,0}$. Moreover, we let
  $\tilde{J}_{\sigma}  \assign \tilde{J}_{\sigma, 1} $ and $\dot{J}_{\sigma}
  \assign \partial_{\sigma} J_{\sigma}$.
\end{definition}

Note that $J_{\sigma}  f \rightarrow f$ as $\sigma \nearrow 1$ in
$\mathcal{S}' (\Lambda_{\varepsilon})$. We let $q_0 (k) \assign | k |$, so
that $q_{\varepsilon} (k) \rightarrow q_0 (k)$ pointwise for $k \in
\mathbb{R}^d$ as $\varepsilon \rightarrow 0$. We observe that, on account of
the above definitions, for any $\ell \in \mathbb{N}_+$ and $\sigma \in (0, 1)$,
we have
\begin{equation}
  \tilde{J}_{\sigma, \ell} J_{\sigma} = J_{\sigma}, \qquad \qquad
  \tilde{J}_{\sigma, \ell + 1} \tilde{J}_{\sigma, \ell} = \tilde{J}_{\sigma,
  \ell} . \label{eq:cutoff-symbol-properties}
\end{equation}
Furthermore, since $J_{\sigma}$, $\tilde{J}_{\sigma, \ell}$ as well as
$\mathcal{L}_{\varepsilon}$ are Fourier multipliers, they all commute. The
operators $(J_{\sigma} )_{\sigma}$ are used to define the scale decomposition
\[ \sigma \mapsto G_{\sigma} \assign J_{\sigma}  G \]
of the Green function $G$. The operators $(\tilde{J}_{\sigma} )_{\sigma}$ play
an auxiliary role and will be used frequently in the estimates. We extend the notation continuously to the endpoint by setting
\[
J_1=\tilde J_1=1,
\qquad
G_1=G.
\]

\begin{remark}
  \label{rmk:js}We shall use the fact that if $\sigma < \mu_i$, where $\mu_i$
  is as in Definition~\ref{def:weights}, then $J_{\sigma} J_{\mu_{i + 1}} =
  J_{\sigma}$. This is a consequence of the definitions of $J_{\sigma}$ and
  $\mu_i$, and of the fact that, for $\sigma < \mu_i$,
  \[ 2 \sigma \llbracket \sigma \rrbracket^{- 1} \leqslant 2 \mu_i 
     \llbracket \mu_i \rrbracket^{- 1} = 2 (2^{i + 2} - 1) < (2^{i + 3} - 1) =
     \mu_{i + 1}  \llbracket \mu_{i + 1} \rrbracket^{- 1}, \]
  which implies $j_{\sigma} (\eta) j_{\mu_{i + 1}} (\eta) = j_{\sigma} (\eta)$
  for all $\eta \in \mathbb{R}$.
\end{remark}

\paragraph*{Smoothing operators.} To establish a suitable Schauder estimate, we
will employ smoothing operators $(K_{\sigma})_{\sigma}$ that possess only
limited regularising effects in time. It is convenient to choose these
operators as inverses of differential operators.

\begin{definition}
  \label{def:regularizing-kernels}For $\sigma, \eta \in (0, 1)$ let
  \begin{align*}
    L_{\sigma} &\assign (1 + \llbracket \sigma \rrbracket^{2 s}
    \partial_t) (1 - \llbracket \sigma \rrbracket^2 \Delta_\varepsilon)^2,\\
    K_{\sigma} &\assign L_{\sigma}^{- 1}
    = {(1 + \llbracket \sigma \rrbracket^{2 s} \partial_t)^{- 1}}
    (1 - \llbracket \sigma \rrbracket^2 \Delta_\varepsilon)^{- 2},\\
    K_{\eta, \sigma} &\assign L_{\sigma} K_{\eta} .
  \end{align*}
\end{definition}

In the tightness argument, we will also employ a Littlewood--Paley
decomposition acting solely on the spatial variables.

\begin{definition}
  \label{def:space-LP}(Spatial LP blocks) Let $(\hat{\Delta}_i :
  [0, \infty) \rightarrow \mathbb{R}_+)_{i \geqslant - 1}$ be
  a dyadic partition of unity on $[0, \infty)$, where
  $\hat{\Delta}_{- 1}$ is supported in $[0, 1]$, 
  $\hat{\Delta}_0$ is supported in $[1 / 2, 3 / 2]$ and
  $\hat{\Delta}_i (\Cdot) \assign \hat{\Delta}_0 (2^{- i} \Cdot)$ for $i
  \geqslant 1$. We define spatial Littlewood-Paley blocks $(\bar{\Delta}_i)_{i
  \geqslant - 1}$ as the Fourier multipliers on $\mathcal{S}'
  (\mathbb{R}_{\varepsilon}^d)$ associated with the symbols $k \mapsto
  \hat{\Delta}_i (q_{\varepsilon} (k))$.
\end{definition}

\begin{remark}
  By extension, we identify the fractional Laplacian $(-
  \Delta_{\varepsilon})^s$ and the spatial Littlewood-Paley blocks
  $(\bar{\Delta}_i)_{i \geqslant - 1}$ with the corresponding operators on $C
  (\Lambda_{\varepsilon}, \zeta)$ acting trivially on the time variable.
\end{remark}

\paragraph*{Convolution operators.} We call $T$ a convolution
operator on $C (\Lambda_{\varepsilon}, w)$ if there exists a signed measure $m
(\mathd z)$ such that
\[ (T f) (z) = \int_{\Lambda_{\varepsilon}} m (\mathd z') f (z - z') . \]
By slight abuse of notation, we usually denote the measure $m (\mathd z)$
associated to $T$ by $T (\mathd z)$. We call $T (\mathd z)$ the kernel of the
operator $T$. We denote by $| T (\mathd z) |$ the variation of the measure $T
(\mathd z)$. Given a nonnegative weight $w \in C (\Lambda_0)$ we write
\[ \| T \|_{\operatorname{TV} (w)} \assign \int_{\Lambda_{\varepsilon}} w (z)  | T
   (\mathd z) | \]
for the weighted total variation norm of the kernel of $T$. If $T (\mathd z)$
is absolutely continuous with respect to the measure $\mathd z$ on
$\Lambda_{\varepsilon}$, we write $T (z)$ for its density, \textit{i.e.} $T
(\mathd z) = T (z) \text{ } \mathd z$. For $p \in [1, \infty] $, we denote by
$\| T \|_{L^p (w)}$ the weighted $L^p$ norm of the density of the kernel of
$T$, \textit{i.e.}
\[ \| T \|_{L^p (w)} \assign \| z \mapsto w (z) T (z) \|_{L^p}, \]
where $\| \Cdot \|_{L^p}$ is the usual $L^p$ norm on $\Lambda_{\varepsilon}$
with respect to the measure $\mathd z$. Given a weight $w$, a convolution
operator $T$ and a function $f$ we write
\begin{equation}
  (w T f) (z) \assign w (z)  (T f) (z) . \label{eq:notation-w-T}
\end{equation}

Note that the operators
\[ G, G_{\sigma}, K_{\sigma}, K_{\eta, \sigma}, J_{\sigma},
   \tilde{J}_{\sigma}, (- \Delta_{\varepsilon})^s, \bar{\Delta}_i \]
introduced above are all convolution operators on $C (\Lambda_{\varepsilon},
\zeta^{\alpha})$ for every $\alpha \in [0, 2 s)$. The kernels of the operators
$G_{\sigma}, J_{\sigma}, \tilde{J}_{\sigma}$ have smooth densities. The
kernels of $G, K_{\sigma}$ have densities of limited regularity. The kernels
of $K_{\eta, \sigma}, (- \Delta_{\varepsilon})^s, \bar{\Delta}_i$ do not
possess densities. In particular, the kernels of $(- \Delta_{\varepsilon})^s,
\bar{\Delta}_i$ are proportional to the Dirac delta in time.

The following lemmas collect the fundamental properties of these convolution
operators, which will be used repeatedly throughout our analysis.

\begin{definition}
  \label{def:weight-w}For $\omega \in \mathbb{R}$ and $\mu \in [0, 1)$, we
  write
  \[ w_{\mu}^{\omega} (z) \assign (1 + \llbracket \mu \rrbracket^{- 1} | z
     |_s)^{\omega}, \qquad z \in \Lambda_0 . \]
\end{definition}

\begin{lemma}
  \label{lem:bounds-K-J}For all $\flat \in (0, 2 s)$, we have
  \begin{enumerate}[(a)]
    \item $\| {K_{\sigma}}  \|_{\operatorname{TV} ( {w^2_{\sigma}} 
    )} \vee \| {K_{\sigma, \eta}}  \|_{\operatorname{TV} (
    {w^2_{\sigma}}  )} \lesssim 1,$
    
    \item $\| J_{\sigma} \|_{\operatorname{TV} ( {w^2_{\sigma}}  )} \vee \|
    \tilde{J}_{\sigma} \|_{\operatorname{TV} ( {w^2_{\sigma}}  )} \vee \|
    L^2_{\sigma} J_{\sigma} \|_{\operatorname{TV} ( {w^2_{\sigma}}  )} \vee
    \| L^2_{\sigma} \tilde{J}_{\sigma} \|_{\operatorname{TV} ( {w^2_{\sigma}} 
    )} \lesssim 1,$
    
    \item $\| \partial_t \tilde{J}_{\sigma} \|_{\operatorname{TV} (
    {\zeta^{-1}}  )} \lesssim \llbracket \sigma \rrbracket^{- 2 s}$
    and $\| \nabla_{\varepsilon} \tilde{J}_{\sigma} \|_{\operatorname{TV} (
    {\zeta^{-1}}  )} \lesssim \llbracket \sigma \rrbracket^{- 1}$,
    
    \item $\| G \|_{\operatorname{TV} ( {w^{\flat}_0}  )} \lesssim 1$ and
    $\| \dot{G}_{\sigma} \|_{\operatorname{TV} ( {w^{\flat}_{\sigma}}  )}
    \vee \| L_{\sigma}^3 \dot{G}_{\sigma} \|_{\operatorname{TV} (
    {w^{\flat}_{\sigma}}  )} \lesssim \llbracket \sigma \rrbracket^{2 s
    - 1},$
  \end{enumerate}
  uniformly over $\varepsilon \in 2^{-\mathbb{N}_0}$ and $1 / 2 \leqslant
  \sigma \leqslant \eta < 1$.
\end{lemma}

\begin{remark}
  \label{rem:decay-rate}The kernels of the operators $K_{\sigma}$,
  $K_{\sigma, \eta}$, $J_{\sigma}$ and $\tilde{J}_{\sigma}$ exhibit rapid decay
  in spacetime. In particular, these operators have finite total variation
  norms $\| \Cdot \|_{\operatorname{TV} (w_{\sigma}^{\omega})}$ for all $\omega
  \geqslant 0$. In contrast, for fractional Laplacians of order $s \in (0,
  1)$, the kernels of the operators $G$ and $\dot{G}_{\sigma}$ decay much more
  slowly. Their total variation norms $\| \Cdot \|_{\operatorname{TV}
  (w_{\sigma}^{\omega})}$ are finite only for $\omega \in [0, 2 s)$. This slow
  decay of $G$ and $\dot{G}_{\sigma}$ constitutes a major analytical
  difficulty and is one of the main reasons why PDEs involving fractional
  Laplacians are particularly challenging to study.
\end{remark}

\begin{remark}
  \label{rem:weights-zeta-w}Note that all of the above bounds remain valid if
  one replaces the weights by $\zeta^{- \alpha}$ or $\zeta^{- \alpha}_{\mu}$
  with any $\mu \in (0, 1)$ and $\alpha \in [0, 1]$, since $1 < 2 s < 2$ and
  \[ \zeta^{- \alpha}_{\mu} \leqslant \zeta^{- \alpha} {\lesssim w^1_{\sigma}}
     {{\leqslant w^2_{\sigma}} } , \]
  for all $\mu, \sigma \in (0, 1)$.
\end{remark}

\begin{proof}
  Item a) is proved in Lemma~\ref{lem:K-kernel}. Items b) and c) follow from
  Lemma~\ref{lem:J-decay}. Item~d) follows from the estimates for
  $\dot{G}_{\sigma}$ established in Lemma~\ref{lem:Gdot-decay}.
\end{proof}

\begin{lemma}
  \label{lem:bound-L-J}We have
  \[ \| \mathcal{L}_{\varepsilon} \tilde{J}_{\sigma} \|_{\operatorname{TV} \left(
     {\zeta^{- 1}}  \right)} \lesssim \llbracket \sigma \rrbracket^{- 2 s} \]
  uniformly over $\varepsilon \in 2^{-\mathbb{N}_0}$ and $1 / 2 \leqslant
  \sigma < 1$.
\end{lemma}

\begin{proof}
  First note that
  \[ \| (- \Delta_{\varepsilon})^s \tilde J_{\sigma} \|_{\operatorname{TV} (\zeta^{- 1})}
     \leqslant \sum_{i \geqslant - 1} \| (- \Delta_{\varepsilon})^s
     \bar{\Delta}_i  \tilde J_{\sigma} \|_{\operatorname{TV} (\zeta^{- 1})}. \]
  Next, observe that $\bar{\Delta}_i 
  \tilde J_{\sigma} = 0$, unless $\llbracket \mu_i \rrbracket \geqslant 2^{- 10}
  \llbracket \sigma \rrbracket .$ Consequently, we have
  \begin{align*}
    \| (- \Delta_{\varepsilon})^s \tilde J_{\sigma} \|_{\operatorname{TV}(\zeta^{- 1})}
    &\leqslant
    \sum_{i\mid \llbracket \mu_i \rrbracket \geqslant 2^{-10}
    \llbracket \sigma \rrbracket}
    \| (- \Delta_{\varepsilon})^s \bar{\Delta}_i
    \|_{\operatorname{TV}(\zeta^{- 1})}
    \| \tilde J_{\sigma} \|_{\operatorname{TV}(\zeta^{- 1})}\\
    &\lesssim
    \sum_{i\mid \llbracket \mu_i \rrbracket \geqslant 2^{-10}
    \llbracket \sigma \rrbracket}
    \llbracket \mu_i \rrbracket^{- 2 s}
    \lesssim \llbracket \sigma \rrbracket^{- 2 s},
  \end{align*}
  where the second estimate follows from Lemmas~\ref{lem:bounds-K-J}
  and~\ref{lem:LP-block}. Since $\| \tilde{J}_{\sigma} \|_{\operatorname{TV} (
  	{\zeta^{-1}}  )} \lesssim 1$ and $\| \partial_t \tilde{J}_{\sigma} \|_{\operatorname{TV} (
  	{\zeta^{-1}}  )} \lesssim \llbracket \sigma \rrbracket^{- 2 s}$ by Lemma~\ref{lem:bounds-K-J}, the proof is complete.
\end{proof}

The following lemma combines the estimates from the preceding lemmas with
Young's convolution inequality. The resulting bounds will be used frequently
in the subsequent analysis, often without explicit reference.

\begin{lemma}
  For all $\alpha \in [0,1]$, the following estimates hold:
  \begin{enumerate}[(a)]
    \setlength{\itemsep}{0.25em}
    \setlength{\parsep}{0pt}

    \item
    \[
      \bigl\| \zeta_{\mu}^{\alpha} K_{\sigma} f \bigr\|
      \lesssim
      \bigl\| \zeta_{\mu}^{\alpha} K_{\eta} f \bigr\|
      \lesssim
      \bigl\| \zeta_{\mu}^{\alpha} f \bigr\|.
    \]
    \item
    \[
      \bigl\| \zeta_{\mu}^{\alpha} J_{\sigma} f \bigr\|
      \lesssim
      \bigl\| \zeta_{\mu}^{\alpha} \tilde J_{\sigma} f \bigr\|
      \lesssim
      \bigl\| \zeta_{\mu}^{\alpha} K_{\sigma} f \bigr\|
      \lesssim
      \bigl\| \zeta_{\mu}^{\alpha} f \bigr\|.
    \]
    \item
    \[
      \begin{aligned}
        \bigl\| \zeta_{\mu}^{\alpha}
        \mathcal L_{\varepsilon} \tilde J_{\sigma} f \bigr\|
        &\lesssim
        \llbracket \sigma \rrbracket^{-2s}
        \bigl\| \zeta_{\mu}^{\alpha} f \bigr\|, \\
        \bigl\| \zeta_{\mu}^{\alpha}
        \partial_t \tilde J_{\sigma} f \bigr\|
        &\lesssim
        \llbracket \sigma \rrbracket^{-2s}
        \bigl\| \zeta_{\mu}^{\alpha} f \bigr\|, \\
        \bigl\| \zeta_{\mu}^{\alpha}
        \nabla_{\varepsilon} \tilde J_{\sigma} f \bigr\|
        &\lesssim
        \llbracket \sigma \rrbracket^{-1}
        \bigl\| \zeta_{\mu}^{\alpha} f \bigr\|.
      \end{aligned}
    \]
    \item
    \[
      \begin{gathered}
        \bigl\| \zeta_{\mu}^{\alpha} G f \bigr\|
        \lesssim
        \bigl\| \zeta_{\mu}^{\alpha} f \bigr\|,
        \\
        \bigl\| \zeta_{\mu}^{\alpha} \dot G_{\sigma} f \bigr\|
        \vee
        \bigl\| \zeta_{\mu}^{\alpha}
        L_{\sigma}^{3} \dot G_{\sigma} f \bigr\|
        \lesssim
        \llbracket \sigma \rrbracket^{2s-1}
        \bigl\| \zeta_{\mu}^{\alpha} f \bigr\|.
      \end{gathered}
    \]
  \end{enumerate}
  The estimates are uniform over
  \[
    \varepsilon \in 2^{-\mathbb N_0},
    \qquad
    0 \leqslant \mu \leqslant 1,
    \qquad
    \frac12 \leqslant \sigma \leqslant \eta < 1,
    \qquad
    f \in \mathcal S'(\Lambda_{\varepsilon}).
  \]
  Here $\| \Cdot \|$ denotes the $L^{\infty}$ norm over
  $\Lambda_{\varepsilon}$, and we use the notation introduced
  in~\eqref{eq:notation-w-T}.
\end{lemma}

\begin{proof}
  The stated bounds are consequences of the previous lemmas, together with
  Young's inequality for convolution, the identities
  \[ K_{\sigma} = K_{\sigma, \eta} {K_{\eta}} , \qquad J_{\sigma} = J_{\sigma}
     \tilde{J}_{\sigma}, \qquad \tilde{J}_{\sigma} = L^2_{\sigma}
     \tilde{J}_{\sigma} K_{\sigma} K_{\sigma} \]
  and Remarks~\ref{rem:weight-properties} a) and~\ref{rem:weights-zeta-w}. For
  example, to prove the first bound, we note that
  \[ \begin{aligned}
       \| \zeta^{\alpha}_{\mu} K_{\sigma} f \| & \leqslant  \sup_z
       \int_{\Lambda_{\varepsilon}} \zeta^{\alpha}_{\mu} (z)  | (K_{\eta} f)
       (z - z_1) |  | K_{\sigma, \eta} (\mathd z_1) |\\
       & \lesssim  \sup_z \int_{\Lambda_{\varepsilon}} \zeta^{\alpha}_{\mu}
       (z - z_1)  | (K_{\eta} f) (z - z_1) | \zeta^{- 1} (z_1)  | K_{\sigma,
       \eta} (\mathd z_1) |\\
       & \leqslant  \| \zeta^{\alpha}_{\mu} (K_{\eta} f) \|  \left\|
       {K_{\sigma, \eta}}  \right\|_{\operatorname{TV} \left( {\zeta^{- 1}}
       \right)}\\
       & \lesssim  \| \zeta^{\alpha}_{\mu} f \| .
     \end{aligned} \]
  The proof of the remaining bounds proceeds analogously.
\end{proof}

\section{Stochastic quantisation}\label{sec:sq}

In this section we lay out the main steps in the proof of
Theorem~\ref{thm:main}, starting from the effective equation at (fractional,
parabolic) spacetime scale $\llbracket \sigma \rrbracket \assign 1 - \sigma$,
$\sigma \in [0, 1]$, obtained through a scale decomposition and the
introduction of the approximate effective force. Since our analysis primarily
concerns scales close to $1$, we shall restrict all scale parameters, such as
$\sigma$, to the interval $[1 / 2, 1]$.

\subsection{Scale decomposition}\label{sec:scale-decomposition}

Let $\phi^{(\varepsilon, M)}$ be a stationary solution to the finite system of
SDEs~\eqref{eq:dynamics-eps-M} and define the scale-dependent field
\[ \phi^{(\varepsilon, M)}_{\sigma} \assign J_{\sigma} \phi^{(\varepsilon,
   M)} \]
localised at (fractional, parabolic) spacetime scales $\gtrsim \llbracket
\sigma \rrbracket$. Then we have
\[ \mathcal{L}_{\varepsilon} {\phi_{\sigma}^{(\varepsilon, M)}}  = J_{\sigma}
   F^{(\varepsilon, M)} (\phi^{(\varepsilon, M)}), \qquad \sigma \in [1 / 2,
   1], \]
where
\begin{equation}
  F^{(\varepsilon, M)} (\phi) {\assign - \lambda \phi^3}  + r_{\varepsilon, M}
  \phi + \xi^{(\varepsilon, M)} . \label{eq:force}
\end{equation}
We call $F^{(\varepsilon, M)}$ the force. Let $\mathcal{E} 
\assign \bigcap_{\alpha > 0} C (\Lambda_{\varepsilon, M}, \zeta^{\alpha})$ be
the space of continuous functions on $\Lambda_{\varepsilon}$ that are periodic in space
and exhibit subpolynomial growth in time at infinity. Denote by
$\hat{\mathcal{E}} \subset \mathcal{S}' (\Lambda_{\varepsilon, M})$ the sum of $\mathcal{E}$ and the image
of $\mathcal{E}$ under $\partial_t$, interpreted as the time derivative
operator on $\mathcal{S}' (\Lambda_{\varepsilon, M})$. Note that the white
noise $\xi^{(\varepsilon, M)}$ on $\Lambda_{\varepsilon, M}$ belongs a.s. to
$\hat{\mathcal{E}}$ and the stochastic convolution $G \xi^{(\varepsilon, M)}$
belongs a.s. to $\mathcal{E}$. Consider a family of functionals, referred to
as the {\textit{effective force}},
\[ \left( {F_{\sigma}^{(\varepsilon, M)}}  : \mathcal{E} \rightarrow
   \hat{\mathcal{E}} \right)_{\sigma \in [1 / 2, 1]}, \]
which is differentiable in $\sigma \in (1 / 2, 1)$ and satisfies the final
condition
\begin{equation}
  F _1 (\psi) = F^{(\varepsilon, M)} (\psi) . \label{eq:final-cond-f}
\end{equation}
Using the identity $\phi_1^{(\varepsilon, M)} = \phi^{(\varepsilon, M)}$, we
obtain
\begin{equation}
  F^{(\varepsilon, M)} (\phi^{(\varepsilon, M)}) = F _1 (\phi^{(\varepsilon,
  M)}_1) = F _{\mu} (\phi^{(\varepsilon, M)}_{\mu}) + R_{\mu}^{(\varepsilon,
  M)}, \label{eq:force-F-R}
\end{equation}
for all $\mu \in [1 / 2, 1]$, where
\[ R_{\mu}^{(\varepsilon, M)} \assign \int_{\mu}^1 [\partial_{\sigma}
   F^{(\varepsilon, M)}_{\sigma} (\phi^{(\varepsilon, M)}_{\sigma}) + \mathD
   F^{(\varepsilon, M)}_{\sigma} (\phi^{(\varepsilon, M)}_{\sigma})
   (\partial_{\sigma} \phi^{(\varepsilon, M)}_{\sigma})]\,\mathd \sigma . \]
Here $\mathD F (\psi) \psi'$ denotes the Fr{\'e}chet derivative of a
functional $F $ in the direction of $\psi' \in \mathcal{E} $ at the point
$\psi \in \mathcal{E}$. Moreover, we have
\[ \partial_{\sigma} \phi^{(\varepsilon, M)}_{\sigma} = \dot{G}_{\sigma} 
   (F^{(\varepsilon, M)} (\phi^{(\varepsilon, M)})) = \dot{G}_{\sigma}
   (F^{(\varepsilon, M)}_{\sigma} (\phi^{(\varepsilon, M)}_{\sigma}) +
   R_{\sigma}^{(\varepsilon, M)}), \]
where $\dot{G}_{\sigma}  \assign \mathcal{L}^{- 1}_{\varepsilon}
\dot{J}_{\sigma}$ and $\dot{J}_{\sigma} \assign \partial_{\sigma} J_{\sigma}$.
Consequently, we deduce that for any choice of the effective force
\[ (F_{\sigma})_{\sigma \in [1 / 2, 1]} = (F_{\sigma}^{(\varepsilon,
   M)})_{\sigma \in [1 / 2, 1]}, \]
satisfying the conditions specified above, the pair
\[ (\phi_{\mu}, R_{\mu})_{\mu \in [1 / 2, 1]} = (\phi^{(\varepsilon,
   M)}_{\mu}, R^{(\varepsilon, M)}_{\mu})_{\mu \in [1 / 2, 1]} \]
satisfies the system of equations
\begin{equation}
  \left\{\begin{array}{l}
    \mathcal{L}_{\varepsilon} \phi _{\mu} = J_{\mu} (F _{\mu} (\phi _{\mu}) +
    R _{\mu})\\
    R_{\mu}  = \int_{\mu}^1 H_{\sigma}  (\phi _{\sigma})\,\mathd \sigma +
    \int_{\mu}^1 [\mathD F _{\sigma} (\phi _{\sigma})  (\dot{G}_{\sigma} 
    R_{\sigma})]\,\mathd \sigma,
  \end{array}\right. \text{} \label{eq:rg-spde}
\end{equation}
where
\begin{equation}
  H_{\sigma}  (\psi) \assign \partial_{\sigma} F _{\sigma} (\psi) + \mathD F
  _{\sigma} (\psi)  (\dot{G}_{\sigma}  F _{\sigma} (\psi)) . \label{eq:def-H}
\end{equation}
The first equation in~\eqref{eq:rg-spde} is the effective stochastic PDE at
scale $\mu$, while the second propagates the error made by replacing the
original force by the scale-dependent effective force. The term $H_{\sigma}$
measures the error of $F_{\sigma}$ in satisfying the exact flow equation. If
$H_{\sigma}=0$, the second equation generates no remainder and one recovers
$R_{\mu}=0$.

Our main goal will be to show that this system allows for good a priori
estimates for a suitably chosen effective force $(F_{\sigma} )_{\sigma \in [1
/ 2, 1]}$.

\begin{remark}
  Except for Sections~\ref{sec:tightness}
  and~\ref{sec:integrability}, we shall almost always suppress the explicit
  dependence on $\varepsilon \in 2^{-\mathbb{N}_0}$ and $M \in \mathbb{N}_+ .$
  To avoid repetition, we will not restate this each time, but \textbf{all
  estimates should be understood as uniform in $\varepsilon \in
  2^{-\mathbb{N}_0}$ and $M \in \mathbb{N}_+$}, unless explicitly stated
  otherwise. In particular, we shall write $\Lambda$ for
  $\Lambda_{\varepsilon}$ and $\Lambda_M$ for $\Lambda_{\varepsilon, M}$.
  Likewise, we write $(- \Delta)^s$ for the discrete fractional Laplacian.
\end{remark}

\subsection{Overview of the strategy}

The goal of this section is to provide a blueprint that guides the reader
through the technical aspects of the proofs and highlights the underlying
heuristics. We will not present rigorous arguments here. Our intention is
instead to convey intuition and the global structure of the analysis. Readers
who find the informal reasoning confusing can safely skip this section without
loss of logical continuity.

Our approach to obtaining a priori global spacetime estimates for solutions to
the system~\eqref{eq:rg-spde} relies on several new conceptual and technical
ideas:
\begin{enumerate}[(a)]
  \item the introduction of a remainder term $R_{\sigma}$, which circumvents
  the need to explicitly solve the flow equation for $F_{\sigma}$ (unlike in
  the original approach
  of~{\cite{duchFlowEquationApproach2021,duchRenormalizationSingularElliptic2022}});
  
  \item the use of new weighted norms that control the solution across the
  entire space, inspired by the spatial decomposition introduced
  in~{\cite{gubinelliGlobalSolutionsElliptic2019}};
  
  \item a stopping argument for the effective force, allowing us to close
  nonlinear estimates;
  
  \item a collection of technical innovations addressing the difficulties
  arising from the limited spacetime decay of the fractional heat kernel.
\end{enumerate}
In what follows, we focus primarily on the new ideas required to handle the
large-field problem. For a pedagogical overview of the flow equation approach,
we refer the reader to the lecture notes~{\cite{duch-notes}}.

To keep the exposition clear, let us ignore the remainder $R_{\sigma}$ and
model the original equation by
\begin{equation}
  \mathcal{L} \phi_{\sigma} - \lambda \phi_{\sigma}^3 \approx F_{\sigma}
  (\phi_{\sigma}) - \lambda \phi_{\sigma}^3, \label{eq:model-coerc}
\end{equation}
where $F_{\sigma} (\phi_{\sigma})$ is a polynomial in the field
$\phi_{\sigma}$. To measure the spacetime growth of fields, we use the weight
\[ \qquad \zeta_{\sigma} (z) \approx (1 + \llbracket \sigma \rrbracket^a | z
   |_s )^{- 1}, \qquad z \in \Lambda, \]
where $| z |_s$ is the fractional parabolic norm on $\Lambda$ and $a > 1$
is an exponent chosen to balance the scale behaviour at large distances. Its
precise value will later be crucial to closing our nonlinear estimates.

Pathwise bounds for the random effective force $(F_{\sigma})_{\sigma}$ are
obtained by analysing the flow equation for its probabilistic cumulants. The
outcome is that $F_{\sigma} (\phi_{\sigma})$ behaves as a random nonlocal
polynomial in $\phi_{\sigma}$ whose coefficients are localised in regions of
size $\approx \llbracket \sigma \rrbracket$ and scale roughly as

\[ \llbracket \sigma \rrbracket^{(k - 3) \beta + \delta \ell}, \]

where $k \in \mathbb{N}_0$ is the monomial degree, $\beta > \gamma$ encodes
the field scaling, $\delta > 0$ measures the distance to criticality, and
$\ell$ is the perturbative order. Due to a Kolmogorov-type argument needed
to extract the almost sure behaviour of the force $F_{\sigma}$ from its
moments, we also lose a bit in the spacetime growth, which will be modelled
by a weight $\zeta_{\sigma}^{- \kappa_{\mathfrak{o}} (\ell + 1)}$, where
$\kappa_{\mathfrak{o}} > 0$ is an arbitrarily small exponent. Overall we have,
schematically,
\begin{equation}
  F_{\sigma} (\phi_{\sigma}) \approx \xi + \bar{r} \phi_{\sigma}  - \lambda
  \phi_{\sigma}^3 + \sum_{\ell = 1}^{\bar{\ell}} \sum_{k = 0}^{\bar{k}}
  \zeta_{\sigma}^{- \kappa_{\mathfrak{o}} (\ell + 1)}  \llbracket \sigma
  \rrbracket^{(k - 3) \beta + \delta \ell} \phi_{\sigma}^k,
  \label{eq:model-post}
\end{equation}
where the sums over $k$ and $\ell$ are finite. The cutoffs $\bar{k}$ and
$\bar{\ell}$ are chosen so that the equation for the remainder $R_{\sigma}$
can be solved, a technical aspect we omit in this heuristic discussion.

To estimate the size of the solution to~\eqref{eq:model-coerc}, we introduce
a constant $C_{\Phi} = \trinorm{ \phi }$ (cf.
Definition~\ref{def:triple-norms}) such that
\begin{equation}
  | \phi_{\sigma} (z) | = | (J_{\sigma} \phi) (z) | \leqslant C_{\Phi}
  \zeta_{\sigma}^{- 1 / 3} (z)  \llbracket \sigma \rrbracket^{- \gamma},
  \qquad z \in \Lambda, \label{eq:trivial-est}
\end{equation}
valid for all $\sigma \geqslant \bar{\mu}$, where $\bar{\mu}$, chosen so that
$\llbracket \bar{\mu} \rrbracket \ll 1$, is a random scale to be chosen later.

For moderate distances $| z | \lesssim \llbracket \sigma \rrbracket^{- a}$,
the spatial weight $\zeta_{\sigma}^{- 1 / 3} (z)$ is of order one and we are
describing the distributional nature of the solution, growing like $\llbracket
\sigma \rrbracket^{- \gamma}$ as $\llbracket \sigma \rrbracket \searrow 0$ for
some $\gamma > 0$. For large distances $| z | \gg \llbracket \sigma
\rrbracket^{- a}$, the spatial growth can be improved as follows. Let
$\llbracket\hat{\mu} \rrbracket \ll \llbracket
\sigma \rrbracket$ be such that $| z | \approx \llbracket \hat{\mu}
\rrbracket^{- a}$. Observing that $\phi_{\sigma} = J_{\sigma}
\phi_{\hat{\mu}}$ due to the properties of the smoothing operators, we obtain
\begin{equation}
  \begin{aligned}
    | \phi_{\sigma} (z) | & \approx | (J_{\sigma} \phi_{\hat{\mu}})  (z) |
    \approx | \phi_{\hat{\mu}} (z) |\\
    & \lesssim C_{\Phi}  (1 + \llbracket \hat{\mu} \rrbracket^a | z |)^{1 /
    3}  \llbracket \hat{\mu} \rrbracket^{- \gamma}\\
    & \approx C_{\Phi}  | z |^{\gamma / a} \approx C_{\Phi}  (1 + \llbracket
    \sigma \rrbracket^a | z |)^{\gamma / a}  \llbracket \sigma \rrbracket^{-
    \gamma} .
  \end{aligned} \label{eq:model-better}
\end{equation}
This improved spatial growth replaces $\zeta_{\sigma}^{- 1 / 3}$ with
$\zeta_{\sigma}^{- \gamma / a}$, and by taking $a$ large we can make the
effective growth arbitrarily mild, an essential feature for closing nonlinear
bounds.

Since higher-order monomials ($k > 3$) in~\eqref{eq:model-post} are
accompanied by small coefficients $\llbracket \sigma \rrbracket^{(k - 3) \beta
+ \delta \ell} \ll 1$, a coercive bound for~\eqref{eq:model-coerc} yields,
schematically,
\begin{equation}
  | \phi_{\sigma} |^3 \approx | F_{\sigma} (\phi_{\sigma}) + \lambda
  \phi_{\sigma}^3 | \lesssim \zeta_{\sigma}^{- (\bar{\ell} + 1)
  \kappa_{\mathfrak{o}}} \sum_{k, \ell} \llbracket \sigma \rrbracket^{(k - 3)
  \beta + \delta \ell} | \phi_{\sigma} |^k . \label{eq:model-size-force}
\end{equation}
We estimated the spacetime growth by taking the worst possible weight
$\zeta_{\sigma}^{- (\bar{\ell} + 1) \kappa_{\mathfrak{o}}}$, where
$\bar{\ell}$ denotes the largest perturbative order that needs to be
considered. Substituting~\eqref{eq:model-better}
in~\eqref{eq:model-size-force}, and ignoring the mild nonlocality of the
effective force, we obtain
\[ | \phi_{\sigma} |^3 \approx | F_{\sigma} (\phi_{\sigma}) + \lambda
   \phi_{\sigma}^3 | \lesssim \zeta_{\sigma}^{- (\bar{\ell} + 1)
   \kappa_{\mathfrak{o}}} \left[ \sum_{k, \ell} \llbracket \sigma
   \rrbracket^{(k - 3) \beta -  k \gamma + \delta \ell} \zeta_{\sigma}^{- k
   \gamma / a} C_{\Phi}^k \right] . \]
The terms with $\ell = 0$, namely the first two terms on the right-hand side
of~\eqref{eq:model-post}, are explicit and yield improved estimates, so that
in total we arrive at
\begin{equation}
  | \phi_{\sigma} |^3 \lesssim \llbracket \sigma \rrbracket^{- 3 \beta +
  \delta} \zeta_{\sigma}^{- (\bar{\ell} + 1) \kappa_{\mathfrak{o}} - \bar{k}
  \gamma / a}  (1 + C_{\Phi})^{\bar{k}}, \label{eq:rough-estimate}
\end{equation}
where $\bar{k}$ is the maximal degree of the monomials. Here we used the fact
that $\beta > \gamma$, which allows the field amplitude to be compensated by
the kernel size. The constant $C_{\Phi}$ can then be estimated as
\[ C_{\Phi} \approx \sup_{\sigma \geqslant \bar{\mu}}  \llbracket \sigma
   \rrbracket^{\gamma} \zeta_{\sigma}^{1 / 3}  | \phi_{\sigma} | \lesssim \sup_{\sigma \geqslant \bar{\mu}}
   \llbracket \sigma \rrbracket^{\gamma - \beta + \delta / 3}
   \zeta_{\sigma}^{(1 - (\bar{\ell} + 1) \kappa_{\mathfrak{o}} - \bar{k}
   \gamma / a) / 3}  (1 + C_{\Phi})^{\bar{k} / 3} . \]
Choosing $\gamma \leqslant \beta$ such that $\gamma - \beta + \delta / 3
\geqslant \vartheta > 0$ and then taking $a$ large enough such that
\[ 1 - (\bar{\ell} + 1) \kappa_{\mathfrak{o}} - \bar{k} \gamma / a \geqslant
   0, \]
we obtain the bound
\[ C_{\Phi} \lesssim \llbracket \bar{\mu} \rrbracket^{\vartheta}  (1 +
   C_{\Phi})^{\bar{k} / 3} . \]
It follows that, for sufficiently small $\llbracket \bar{\mu} \rrbracket \ll
1$, the nonlinear estimate closes (after selecting a small branch with the use of a stopping-time argument) and yields $C_{\Phi} \approx 1$.

Let us now address the treatment of the nonlocality of the effective force
kernels. The main technical challenges in this work stem from the limited
decay of the slice propagator $\dot{G}_{\sigma}$ associated with the
fractional parabolic operator $\mathcal{L}$. Roughly speaking, we only have
algebraic decay of the form (see Lemma~\ref{lem:Gdot-decay-old}):
\begin{equation}
  | \dot{G}_{\sigma} (z) | \lesssim   \llbracket \sigma \rrbracket^{- d - 1} 
  (1 + | z |_s / \llbracket \sigma \rrbracket)^{- d - 2 s + \epsilon}, \qquad
  z \in \Lambda, \label{eq:limited-decay}
\end{equation}
where $\epsilon > 0$ is arbitrarily small. This behaviour contrasts sharply
with the standard heat kernel, which has stretched-exponential decay, and also
with the fractional Laplacian case appearing in the usual (``static'')
renormalisation group analysis of the fractional $\Phi^4$
model {\cite{brydgesCriticalF42003}}. The difference arises from the limited
smoothness of the symbol of the fractional heat operator $\mathcal{L}$.

The main consequence of~\eqref{eq:limited-decay} is that the effective force
kernels, obtained by solving a flow equation driven by $\dot{G}_{\sigma}$,
inherit a similar algebraic decay. The monomials which appear in $F_{\sigma}
(\phi_{\sigma})$ take the schematic form
\[ F_{\sigma} (\phi_{\sigma}) (z) \approx \sum_{k, \ell} \int F^{[\ell],
   (k)}_{\sigma} (z ; z_1, \ldots, z_k) \prod_{j = 1}^k \phi_{\sigma} (z_j)
   \mathd z_j, \]
where $F^{[\ell], (k)}_{\sigma}$ are random kernels. Ignoring their
distributional nature and thinking of them as bona fide functions, their
spatial nonlocality and spatial growth can be modelled as (cf.
Definitions~\ref{def:weights-kernels} and~\ref{def:norm-kernels})
\[ F^{[\ell], (k)}_{\sigma} (z, z_1, \ldots, z_k) \approx \zeta_{\sigma}^{-
   (\ell + 1) \kappa_{\mathfrak{o}}} (z)  (1 + \llbracket \sigma \rrbracket^{-
   1} \operatorname{St} (z, z_1, \ldots, z_k) )^{- (\flat - \ell
   \kappa_{\mathfrak{o}})}, \]
where $\operatorname{St} (z, z_1, \ldots, z_k)$ measures the diameter of the set $\{
z, z_1, \ldots, z_k \}$. The initial decay exponent $\flat \approx 2 s$
follows from~\eqref{eq:limited-decay}, while the additional loss $\ell
\kappa_{\mathfrak{o}}$ reflects the growth at spacetime infinity of the kernel
in its output variable, induced by the similar growth of the noise. It is now
clear that, to proceed as in~\eqref{eq:rough-estimate}, we must be able to
compensate for the spatial growth of the fields $\phi_{\sigma}$ by exploiting
the limited decay of the kernels away from the diagonal.

The remainder of this section makes this argument rigorous and establishes
further properties of the solution. Section~\ref{sec:apriori} develops the
coercive estimate required in ~\eqref{eq:model-size-force}, while
Section~\ref{sec:flow} contains the detailed analysis of the random force
coefficients and the derivation of the precise form of
estimate~\eqref{eq:model-post}.

\subsection{Main estimate}

In this section we introduce a family of weighted norms that will be used to
measure the size of fields over spacetime uniformly across scales. These norms
depend on an exponent $\gamma > 0$, to be fixed in
Section~\ref{sec:post-processing}, and a terminal scale $\bar{\mu} \in [1 / 2,
1)$, which will later be tuned according to the noise amplitude.

\begin{definition}
  \label{def:triple-norms}For $\psi \in \mathcal{S}' (\Lambda)$ and $f \in C
  ([1 / 2, 1), C (\Lambda, \zeta))$, we set
  \begin{equation}
    \trinorm{ \psi } = \trinorm{ \psi }_{\bar{\mu}}
    \assign \sup_{\bar{\mu} \leqslant \sigma < 1} \llbracket \sigma
    \rrbracket^{\gamma}  \| \zeta^{1 / 3}_{\sigma} J_{\sigma} \psi \|,
    \label{eq:sol-norm}
  \end{equation}
  \begin{equation}
    \trinorm{ f }_{\#} = \trinorm{ f }_{\#, \bar{\mu}}
    \assign \sup_{\bar{\mu} \leqslant \sigma < 1}  \llbracket \sigma \rrbracket^{3
    \gamma}  \| \zeta_{\sigma} f_{\sigma} \| . \label{eq:sharp-norm}
  \end{equation}
\end{definition}

The powers in these norms are adapted to the cubic coercivity. Roughly, a
bound of size $\llbracket\sigma\rrbracket^{-\gamma}$ for $\phi_{\sigma}$
produces a bound of size $\llbracket\sigma\rrbracket^{-3\gamma}$ for its
cubic force. Likewise, the weight $\zeta_{\sigma}^{1/3}$ on the field becomes
the full weight $\zeta_{\sigma}$ after cubing.

The above ``triple norms'' exhibit a specific behaviour with respect to the
weight. In particular, the following lemma shows that changing the weight
yields equivalent norms.

\begin{lemma}
  \label{lem:triple-phi}For $\alpha \in [3 \nu, 1]$ the following bounds
  \begin{equation}
    \sup_{\sigma, \mu | \sigma \geqslant \mu \geqslant \bar{\mu}} \llbracket
    \sigma \rrbracket^{\gamma}  \| \zeta_{\mu}^{\alpha / 3} \phi_{\sigma} \|
    \lesssim \trinorm{ \phi }_{\bar{\mu}} \lesssim \sup_{\sigma,
    \mu | \sigma \geqslant \mu \geqslant \bar{\mu}} \llbracket \sigma
    \rrbracket^{\gamma}  \| \zeta_{\mu}^{\alpha / 3} \phi_{\sigma} \|,
    \label{eq:equiv-norms}
  \end{equation}
  \begin{equation}
    \sup_{\sigma, \mu | \sigma \geqslant \mu \geqslant \bar{\mu}} \llbracket
    \sigma \rrbracket^{3 \gamma}  \| \zeta_{\mu}^{\alpha} \mathcal{L}
    \phi_{\sigma} \| \lesssim \trinorm{ \sigma \mapsto \mathcal{L}
    \phi_{\sigma} }_{\#, \bar{\mu}} \lesssim \sup_{\sigma, \mu
    | \sigma \geqslant \mu \geqslant \bar{\mu}} \llbracket \sigma
    \rrbracket^{3 \gamma}  \| \zeta^{\alpha}_{\mu} \mathcal{L} \phi_{\sigma}
    \|, \label{eq:equiv-norms-L}
  \end{equation}
  hold uniformly in $\bar{\mu} \in [1 / 2, 1)$ and $\phi \in \mathcal{S}'
  (\Lambda)$, where $\phi_{\sigma} \assign J_{\sigma} \phi$ for $\sigma \in [1
  / 2, 1)$.
\end{lemma}

\begin{proof}
  Let us prove the first inequality in~\eqref{eq:equiv-norms}. We start by
  observing that, on account of the support properties of $\{ \chi_i \}_i$, we
  have the following decomposition
  \[ \| \zeta_{\mu}^{\alpha / 3} \phi_{\sigma} \| \leqslant \sup_i \|
     \zeta_{\mu}^{\alpha / 3} (\chi_i + \chi_{i - 1} + \chi_{i + 1})
     \phi_{\sigma} \| \lesssim \sup_i \| \zeta_{\mu}^{\alpha / 3} \chi_i
     \phi_{\sigma} \| . \]
  Let us consider separately the cases $\mu_i \leqslant \sigma$ and $\mu_i
  \geqslant \sigma$. In the first case, we use the bound
  \[ \| \zeta_{\mu}^{\alpha / 3} \chi_i
     \zeta_{\mu_i}^{- 1 / 3} \| \lesssim 1 \]
  to obtain
  \begin{align*}
    \sup_{i | \mu_i \leqslant \sigma} \llbracket \sigma \rrbracket^{\gamma} 
     \| \zeta_{\mu}^{\alpha / 3} \chi_i \phi_{\sigma} \|
     &\lesssim \sup_{i | \mu_i \leqslant \sigma}
     \llbracket \sigma \rrbracket^{\gamma}
     \| \zeta_{\mu}^{\alpha / 3} \chi_i \zeta_{\mu_i}^{- 1 / 3} \|
     \| \zeta^{1 / 3}_{\mu_i} \phi_{\sigma} \|\\
     &\lesssim \llbracket \sigma \rrbracket^{\gamma}
     \sup_{i | \mu_i \leqslant \sigma}  \| \zeta^{1 / 3}_{\mu_i}
     \phi_{\sigma} \|\\
     &\leqslant \llbracket \sigma \rrbracket^{\gamma}
     \| \zeta^{1 / 3}_{\sigma} \phi_{\sigma} \|
     \lesssim \trinorm{ \phi } .
  \end{align*}
  Let us turn to the case $\mu_i > \sigma$. Note that
  \[ \| \zeta_{\mu}^{\alpha / 3} \chi_i \zeta_{\mu_{i +
     1}}^{- 1 / 3} \| \lesssim (1 + \llbracket \mu \rrbracket^a \llbracket
     \mu_{i + 1} \rrbracket^{- a})^{- \alpha / 3} \lesssim
     \llbracket \mu \rrbracket^{- a \alpha / 3}  \llbracket
     \mu_{i + 1} \rrbracket^{\alpha a / 3} . \]
  Moreover, on account of Remarks~\ref{rmk:js},~\ref{rem:weight-properties}
  and Lemma~\ref{lem:bounds-K-J}, we have
  \[ \| \zeta_{\mu_{i + 1}}^{1 / 3} \phi_{\sigma} \| = \| \zeta_{\mu_{i +
     1}}^{1 / 3} J_{\sigma} \phi_{\mu_{i + 1}} \| \lesssim \| J_{\sigma}
     \|_{\operatorname{TV} (\zeta ^{- 1})}  \| \zeta_{\mu_{i + 1}}^{1 / 3}
     \phi_{\mu_{i + 1}} \| \lesssim \| \zeta_{\mu_{i + 1}}^{1 / 3}
     \phi_{\mu_{i + 1}} \| . \]
  Combining the above estimates, we deduce that
  \begin{align*}
    &\sup_{i | \mu_i > \sigma} \llbracket \sigma \rrbracket^{\gamma}
    \| \zeta_{\mu}^{\alpha / 3} \chi_i \phi_{\sigma} \|\\
    &\quad \lesssim  \sup_{i | \mu_i > \sigma}
    \llbracket \sigma \rrbracket^{\gamma}
    \| \zeta_{\mu}^{\alpha / 3} \chi_i
    \zeta_{\mu_{i + 1}}^{- 1 / 3} \|
    \| \zeta^{1 / 3}_{\mu_{i + 1}} \phi_{\mu_{i + 1}} \|\\
 &\quad \lesssim  \llbracket \sigma \rrbracket^{\gamma} \llbracket \mu
    \rrbracket^{- a \alpha / 3} \left(\sup_{i | \mu_i > \sigma} \llbracket
    \mu_{i + 1} \rrbracket^{\alpha a / 3} \llbracket \mu_{i +
    1} \rrbracket^{- \gamma}\right)  \left(\sup_{i | \mu_i > \sigma} \llbracket \mu_{i +
    1} \rrbracket^{\gamma} \| \zeta^{1 / 3}_{\mu_{i + 1}}
    \phi_{\mu_{i + 1}} \|\right)\\[3pt]
 &\quad \lesssim  \llbracket \sigma \rrbracket^{\alpha a / 3} \llbracket \mu
    \rrbracket^{- a \alpha / 3}  \trinorm{ \phi }\\[3pt]
 &\quad \lesssim  \trinorm{ \phi },
  \end{align*}
  where we used that $\alpha a / 3 \geqslant \nu a = \gamma$ and $\sigma
  \geqslant \mu .$ This proves the first inequality in~\eqref{eq:equiv-norms}.
  The second inequality in~\eqref{eq:equiv-norms} is a direct consequence of
  the fact that $\zeta^{1 / 3}_{\sigma} \leqslant \zeta_{\sigma}^{\alpha /
  3}$. The inequalities~\eqref{eq:equiv-norms-L} are proved similarly.
\end{proof}

The two norms just defined fix the analytical setting for the global analysis
of the SPDE~\eqref{eq:rg-spde}. In Section~\ref{sec:apriori}, we will prove suitable
coercive estimates for the fractional parabolic equation with cubic nonlinearity,
which allow us to control the large values of the fields. Below we state a
direct consequence of this estimate.

\begin{theorem}
  \label{thm:apriori}(Coercive estimate) For $\phi \in \mathcal{S}' (\Lambda)$
  and $\sigma \in [1 / 2, 1)$ define
  \[ f_{\sigma} \assign \mathcal{L} \phi_{\sigma} + \lambda \phi_{\sigma}^3,
     \qquad \phi_{\sigma} \assign J_{\sigma} \phi, \qquad \mathcal{L} \assign
     \partial_t + (- \Delta)^s + m^2 . \]
  Then the following bound
  \[ \trinorm{ \phi }_{\bar{\mu}} \lesssim \lambda^{- 1 / 2} 
     \llbracket \bar{\mu} \rrbracket^{\gamma} + \lambda^{- 1 / 3}  \trinorm{
     f }_{\#, \bar{\mu}}^{1 / 3} \]
  holds uniformly in $\phi \in \mathcal{S}' (\Lambda)$, $\bar{\mu} \in [1 / 2,
  1)$ and $\lambda \in (0, \infty)$.
\end{theorem}

\begin{proof}
  We apply the a priori estimates from
  Theorem~\ref{thm:partial-coercive-estimates} to $u = \phi_{\sigma}$ and
  $\rho = \zeta_{\sigma}^{1 / 3}$. The constants $A$ and $B$ appearing in the
  theorem can be bounded as follows:
  \[ A = \| (- \Delta)^s \zeta_{\sigma}^{2 / 3} \| + \| \zeta_{\sigma}^{1 / 3}
     (\partial_t \zeta_{\sigma}^{1 / 3}) \| \lesssim \llbracket \sigma
     \rrbracket^{2 s a}, \]
  and
  \begin{align*}
    B &=  \| \zeta_{\sigma}^{1 / 3} u \|
    \left( \| \zeta_{\sigma}^{1 / 3}(\partial_t \zeta_{\sigma}^{1 / 3}) \|
    + \| \zeta_{\sigma}^{1 / 3} (-\Delta)^s \zeta_{\sigma}^{1 / 3} \|
    + \big\| \zeta_{\sigma}^{2 / 3}
    \mathfrak{D}_s \left( {\zeta_{\sigma}^{- 1 / 3}}  \right)
    \mathfrak{D}_s(\zeta_{\sigma}^{1 / 3}) \big\| \right)\\
    &\quad + \big\| {\zeta_{\sigma}^{1 / 3}}
    \mathfrak{D}_s(\zeta_{\sigma}^{1 / 3})
    \mathfrak{D}_s(\zeta_{\sigma}^{1 / 3} u) \big\|\\[3pt]
    &\lesssim Q^2_{\sigma}  \| \zeta_{\sigma}^{1 / 3} \phi_{\sigma} \|
    + Q_{\sigma}  \| \mathfrak{D}_s(\zeta_{\sigma}^{1 / 3} \phi_{\sigma}) \| .
  \end{align*}
  Here we define
  \[ Q_{\sigma}^2 \assign \| \partial_t \zeta_{\sigma}^{1 / 3} \| + \| (-
     \Delta)^s \zeta_{\sigma}^{1 / 3} \| + \| \mathfrak{D}_s
     (\zeta_{\sigma}^{1 / 3}) \|^2 + \| \zeta_{\sigma}^{1 / 3} \mathfrak{D}_s
     (\zeta_{\sigma}^{- 1 / 3}) \|^2 \lesssim \llbracket \sigma \rrbracket^{2
     s a} . \]
  The bound $\| (- \Delta)^s \zeta_{\sigma}^{\alpha} \|
  \lesssim \llbracket \sigma \rrbracket^{2 s a}$ is proved in
  Lemma~\ref{lem:fractional-bound}. To estimate $\| \mathfrak{D}_s
  (\zeta_{\sigma}^{1 / 3}) \|^2$ we used~\eqref{eq:D-est}, while the bound on
  $\| \zeta_{\sigma}^{1 / 3} \mathfrak{D}_s (\zeta_{\sigma}^{- 1 / 3}) \|^2$
  follows from Lemma~\ref{lem:D-zeta}. Moreover, using the
  bound~\eqref{eq:D-est}, we obtain
  \begin{equation}
    \| \mathfrak{D}_s (\zeta_{\sigma}^{1 / 3} \phi_{\sigma}) \| \lesssim \|
    \nabla_{\varepsilon} (\zeta_{\sigma}^{1 / 3} \phi_{\sigma}) \|^s  \|
    \zeta_{\sigma}^{1 / 3} \phi_{\sigma} \|^{1 - s} . \label{eq:D-proof-est}
  \end{equation}
  We note that
  \[ \| \nabla_{\varepsilon} (\zeta_{\sigma}^{1 / 3} \phi_{\sigma}) \|
     \lesssim \| (\nabla_{\varepsilon} \zeta_{\sigma}^{1 / 3}) \phi_{\sigma}
     \| + \| \zeta_{\sigma}^{1 / 3}  (\nabla_{\varepsilon} \phi_{\sigma}) \| .
  \]
  Hence, using the fact that $| \nabla_{\varepsilon} \zeta_{\sigma}^{1 / 3} |
  \lesssim \zeta_{\sigma}^{1 / 3}$ and $\phi_{\sigma} = \tilde{J}_{\sigma}
  \phi_{\sigma}$ as well as $\| \nabla_{\varepsilon} \tilde{J}_{\sigma}
  \|_{\operatorname{TV} (\zeta^{- 1})} \lesssim \llbracket \sigma \rrbracket^{- 1}$,
  we arrive at
  \begin{align*}
    \| \nabla_{\varepsilon}(\zeta_{\sigma}^{1 / 3} \phi_{\sigma}) \|
    &\lesssim \| \zeta_{\sigma}^{1 / 3} \phi_{\sigma} \|
    + \| \zeta_{\sigma}^{1 / 3}
    (\nabla_{\varepsilon} \tilde{J}_{\sigma, 1})\phi_{\sigma} \|\\
    &\lesssim \| \zeta_{\sigma}^{1 / 3} \phi_{\sigma} \|
    + \llbracket \sigma \rrbracket^{- 1}
    \| \zeta_{\sigma}^{1 / 3}\phi_{\sigma} \|\\
    &\lesssim \llbracket \sigma \rrbracket^{- 1}
    \| \zeta_{\sigma}^{1 / 3}\phi_{\sigma} \| .
  \end{align*}
  Together with~\eqref{eq:D-proof-est}, this proves that
  \[ \| \mathfrak{D}_s (\zeta_{\sigma}^{1 / 3} \phi_{\sigma}) \| \lesssim
     \llbracket \sigma \rrbracket^{- s}  \| \zeta_{\sigma}^{1 / 3}
     \phi_{\sigma} \|, \]
  uniformly in $\sigma \geqslant \bar{\mu}$. Therefore,
  \[ B \lesssim (\llbracket \sigma \rrbracket^{2 s a} + \llbracket \sigma
     \rrbracket^{s a} \llbracket \sigma \rrbracket^{- s}) \| \zeta_{\sigma}^{1
     / 3} \phi_{\sigma} \| \lesssim \llbracket \sigma \rrbracket^{s a - s}  \|
     \zeta_{\sigma}^{1 / 3} \phi_{\sigma} \| . \]
  By Theorem~\ref{thm:partial-coercive-estimates} and Young's inequality, we
  obtain
  \[ \begin{aligned}
       \| \zeta_{\sigma}^{1 / 3} \phi_{\sigma} \| & \leqslant  \lambda^{- 1 /
       2} C \llbracket \sigma \rrbracket^{a s} + \lambda^{- 1 / 3} (\|
       \zeta_{\sigma} f_{\sigma} \| + C \llbracket \sigma \rrbracket^{a s - s}
       \| \zeta_{\sigma}^{1 / 3} \phi_{\sigma} \|)^{1 / 3} \\
       & \leqslant  \lambda^{- 1 / 2} C \llbracket \sigma \rrbracket^{a s} +
       \lambda^{- 1 / 3} \| \zeta_{\sigma} f_{\sigma} \|^{1 / 3}  + C
       \lambda^{- 1 / 3} \llbracket \sigma \rrbracket^{(a s - s) / 3} \|
       \zeta_{\sigma}^{1 / 3} \phi_{\sigma} \|^{1 / 3}\\
       & \leqslant  \lambda^{- 1 / 2} C \llbracket \sigma \rrbracket^{a s} +
       \lambda^{- 1 / 3} \| \zeta_{\sigma} f_{\sigma} \|^{1 / 3}  + C'
       \lambda^{- 1 / 2} \llbracket \sigma \rrbracket^{(a s - s) / 2} +
       \frac{1}{2} \| \zeta_{\sigma}^{1 / 3} \phi_{\sigma} \| .
     \end{aligned} \]
  As a result, we arrive at
  \[ \| \zeta_{\sigma}^{1 / 3} \phi_{\sigma} \| \lesssim \lambda^{- 1 / 2}
     \llbracket \sigma \rrbracket^{(a s - s) / 2} + \lambda^{- 1 / 3} \|
     \zeta_{\sigma} f_{\sigma} \|^{1 / 3}  . \]
  This allows us to deduce that
  \[ \trinorm{ \phi }_{\bar{\mu}} = \sup_{\sigma \geqslant
     \bar{\mu}}  \llbracket \sigma \rrbracket^{\gamma} \| \zeta_{\sigma}^{1 /
     3} \phi_{\sigma} \| \lesssim \sup_{\sigma \geqslant \bar{\mu}} 
     \llbracket \sigma \rrbracket^{\gamma + (a s - s) / 2} \lambda^{- 1 / 2} +
     \lambda^{- 1 / 3} \sup_{\sigma \geqslant \bar{\mu}}  (\llbracket \sigma
     \rrbracket^{3 \gamma}  \| \zeta_{\sigma} f_{\sigma} \|)^{1 / 3}  . \]
  Hence, provided $\gamma + (a s - s) / 2 \geqslant \gamma \geqslant 0$, which
  follows from $a > 1$, we conclude our claim.
\end{proof}

The coercivity estimate derived above provides the key ingredient for deriving
a priori bounds for solutions to system~\eqref{eq:rg-spde}, as stated in the
following theorem.

\begin{theorem}
  \label{thm:apriori-system}Let $\phi = \phi^{(\varepsilon,
  M)}$ be a solution of the stochastic quantisation
  equation~\eqref{eq:dynamics-eps-M} and
  \[ \left( F_{\sigma} {= F_{\sigma}^{(\varepsilon, M)}}  : \mathcal{E}
     \rightarrow \hat{\mathcal{E}} \right)_{\sigma \in [1 / 2, 1]} \]
  be a family of functionals differentiable in $\sigma \in (1 / 2, 1)$ and
  satisfying the final condition~\eqref{eq:final-cond-f}. For $\sigma \in [1 /
  2, 1]$ define
  \[ \phi_{\sigma} \assign J_{\sigma} \phi, \qquad R_{\sigma} \assign F_1 (\phi)
     - F_{\sigma} (\phi_{\sigma}) . \]
  
  Suppose that there exist constants $S\geqslant 3, C_F \geqslant 1$ and
  $\vartheta, \bar{\kappa} > 0$ such that
  \begin{equation}
    \bar{\kappa} \in [\nu, 1/3), \qquad \vartheta \leqslant \gamma \wedge (2 s -
    \gamma), \qquad (s + \gamma)  \bar{\kappa} / (1 - \bar{\kappa}) \leqslant
    \vartheta / 4, \label{eq:thm-constraints}
  \end{equation}
  and $(F_{\sigma})_{\sigma}$ satisfies the following estimates
  \begin{equation}
    \begin{aligned}
      &\| \zeta_{\mu} [J_{\sigma} F_{\sigma} (\psi_{\sigma})
      - (- \lambda \psi_{\sigma}^3)] \|\\
      &\quad \leqslant \llbracket \sigma \rrbracket^{- 3 \gamma + \vartheta}
      [C_F  (1 + \trinorm{ \psi }_{\bar{\mu}})^S
      + (1 + \trinorm{ \psi }_{\bar{\mu}})^2
      \trinorm{ \mathcal{L} \psi_{\Cdot} }_{\#, \bar{\mu}}],\\[3pt]
      &\| \zeta_{\mu} K_{\sigma} H_{\sigma} (\psi_{\sigma}) \|
      \leqslant
      C_F  \llbracket \sigma \rrbracket^{\vartheta - 1}
      (1 + \trinorm{ \psi }_{\bar{\mu}})^S,\\[3pt]
      &\| \zeta^{\bar{\kappa}}_{\mu} K_{\sigma}
      F_{\sigma} (\psi_{\sigma}) \|
      \leqslant C_F  \llbracket \sigma \rrbracket^{- 3 \gamma}
      (1 + \trinorm{ \psi }_{\bar{\mu}})^S,\\[3pt]
      &\| \zeta_{\mu} K_{\sigma}
      (\mathD F_{\sigma} (\psi_{\sigma})
      (\dot{G}_{\sigma}  \hat{\psi})) \| \leqslant C_F  \llbracket \sigma \rrbracket^{\vartheta - 1}
      (1 + \trinorm{ \psi }_{\bar{\mu}})^S
      \| \zeta^{1 - \bar{\kappa}}_{\mu}
      \tilde{J}_{\sigma}^2 \hat{\psi} \|,
    \end{aligned} \label{eq:est-F}
  \end{equation}
  for all $\bar{\mu} \in [1 / 2, 1), \mu \in [\bar{\mu}, 1), \sigma \in [\mu,
  1)$ and $\psi, \hat{\psi} \in \mathcal{S}' (\Lambda)$, where we denote
  $\psi_{\sigma} \assign J_{\sigma} \psi$ and $(\mathcal{L} \psi)_{\sigma}
  \assign \mathcal{L} \psi_{\sigma}$, the functional $H_\sigma$ is defined
  by~\eqref{eq:def-H}.
  
  Then there exists a universal constant $\hat{C} > 0$ such that for all
  $\bar{\mu} \in [1 / 2, 1)$ satisfying
  \begin{equation}
    \llbracket \bar{\mu} \rrbracket^{\vartheta} \leqslant
    \hat{C}  ((\lambda^{- 1} + 1) C_F)^{- 2}, \label{eq:bound-mubar}
  \end{equation}
  we have
  \[ \trinorm{ \phi }_{\bar{\mu}} \leqslant 1, \qquad \trinorm{
     \mathcal{L} \phi_{\Cdot} }_{\#, \bar{\mu}} \leqslant \lambda
     \llbracket \bar{\mu} \rrbracket^{- \vartheta / 2}, \qquad \trinorm{
     K_{\Cdot} R_{\Cdot} }_{\#, \bar{\mu}} \leqslant \lambda
     \llbracket \bar{\mu} \rrbracket^{\vartheta / 2} . \]
\end{theorem}

\begin{remark}
  By the elementary argument presented in Section~\ref{sec:scale-decomposition},
  our assumptions on $\phi$ and $(F_{\mu})_{\mu}$ imply that $(\phi_{\mu},
  R_{\mu})_{\mu}$ is a solution of~\eqref{eq:rg-spde}.
  The four estimates in~\eqref{eq:est-F} have distinct roles. The first extracts
  the coercive cubic part of the effective force, the second controls the error
  $H_{\sigma}$ in the flow equation, the third bounds the effective
  force in the weaker weight required in the field equation, and the fourth is
  the Lipschitz estimate used to propagate the remainder $R_{\sigma}$.
\end{remark}

\begin{proof}
  Define
  \begin{equation}
    \Phi = \Phi_{\bar{\mu}} \assign 1 + \trinorm{ \phi
    }_{\bar{\mu}} + \theta_{\mathcal{L}}^{- 1}  \trinorm{
    \mathcal{L} \phi_{\Cdot} }_{\#, \bar{\mu}} + \theta_R^{- 1} 
    \trinorm{ K_{\Cdot} R_{\Cdot} }_{\#, \bar{\mu}},
    \label{eq:def-Phi}
  \end{equation}
  for constants $\theta_{\mathcal{L}}, \theta_R > 0$ to be fixed later. 
  These constants normalise the three quantities entering the bootstrap so that
  the field equation and the remainder equation can be closed in a single
  inequality. 
  Our
  goal is to bound each term on the right-hand side by an increasing function
  of $\Phi$ itself and then apply a continuity argument to establish the
  uniform bounds claimed. We omit the index $\bar{\mu} \in [1 / 2,
  1)$, as the estimates hold uniformly in this parameter.
  
  Let $C_{\mathcal{Q}} = C_F (1 + \trinorm{ \phi })^S$. We first
  prove a bound for $\trinorm{ K_{\Cdot} R_{\Cdot} }_{\#}$. To
  this end, we use Lemma~\ref{lem:weighted-gronwall} below, which shows that
  we can control $\| \zeta_{\mu} K_{\sigma} R_{\sigma} \|$ in terms of $\|
  \zeta_{\mu}^{\bar{\kappa}} \tilde{J}_{\sigma}^2 R_{\sigma} \|$. Note that
  \[ \tilde{J}_{\sigma}^2 R_{\sigma} = \tilde{J}_{\sigma}^2 (\mathcal{L} \phi
     - F_{\sigma} (\phi_{\sigma})), \]
  which follows from~\eqref{eq:dynamics-eps-M} and~\eqref{eq:force-F-R}. Since
  $\bar{\kappa} \in [\nu, 1)$ and $\vartheta \in (0, \gamma]$, by
  Lemmas~\ref{lem:bound-L-J} and~\ref{lem:triple-phi}, we have
  \begin{align*}
    \| \zeta_{\mu}^{\bar{\kappa}} \tilde{J}_{\sigma}^2
    \mathcal{L} \phi \|
    &= \| \zeta_{\mu}^{\bar{\kappa}}
    (\mathcal{L} \tilde{J}_{\sigma})
    \tilde{J}_{\sigma} \phi \|\\
    &\lesssim \| \mathcal{L} \tilde{J}_{\sigma}
    \|_{\operatorname{TV}(\zeta^{- 1})}
    \| \zeta_{\mu}^{\bar{\kappa}}
    \tilde{J}_{\sigma} \phi \|\\
    &\lesssim \llbracket \sigma \rrbracket^{- 2 s}
    \llbracket \sigma \rrbracket^{- \gamma}
    \trinorm{ \phi }\\
    &\lesssim C_{\mathcal{Q}}
    \llbracket \bar{\mu} \rrbracket^{\vartheta}
    \llbracket \sigma \rrbracket^{- 2 (s + \gamma)} .
  \end{align*}
  Since $\vartheta \in (0, 2 s - \gamma]$, by Lemma~\ref{lem:bounds-K-J} and
  the estimates~\eqref{eq:est-F}, we have
  \[ \| \zeta_{\mu}^{\bar{\kappa}} \tilde{J}_{\sigma}^2 F_{\sigma}
     (\phi_{\sigma}) \| \lesssim \| \zeta_{\mu}^{\bar{\kappa}} K_{\sigma}
     F_{\sigma} (\phi_{\sigma}) \| \lesssim C_{\mathcal{Q}} \llbracket \sigma
     \rrbracket^{- 3 \gamma} \lesssim C_{\mathcal{Q}} \llbracket \bar{\mu}
     \rrbracket^{\vartheta} \llbracket \sigma \rrbracket^{- 2 (s + \gamma)} .
  \]
  Altogether, this yields
  \begin{equation}
    \sup_{\sigma, \mu | \sigma \geqslant \mu \geqslant \bar{\mu}} \llbracket
    \sigma \rrbracket^{2 (s + \gamma)} \| \zeta_{\mu}^{\bar{\kappa}}
    \tilde{J}_{\sigma}^2 R_{\sigma} \| \lesssim C_{\mathcal{Q}} \llbracket
    \bar{\mu} \rrbracket^{\vartheta} . \label{eq:J-R-bound}
  \end{equation}
  Using~\eqref{eq:est-F}, the factorisation
  $K_{\mu} = K_{\mu, \eta} K_{\eta}$ for $\eta \geqslant \mu$, and
  Lemma~\ref{lem:bounds-K-J} we also have
  \[ \begin{aligned}
  	\sup_{\sigma, \mu | \sigma \geqslant \mu \geqslant \bar{\mu}} \left\|
     \zeta_{\mu} K_{\mu} \int_{\sigma}^1 H_{\eta} (\phi_{\eta})\,\mathd \eta
     \right\| &\lesssim \int_{\bar{\mu}}^1 \| \zeta_{\mu} K_{\mu, \eta}
     K_{\eta} H_{\eta} (\phi_{\eta}) \|\,\mathd \eta
     \\
     &\lesssim \int_{\bar{\mu}}^1 \| \zeta_{\mu} K_{\eta} H_{\eta}
     (\phi_{\eta}) \|\,\mathd \eta \lesssim C_{\mathcal{Q}} \llbracket
     \bar{\mu} \rrbracket^{\vartheta} .
     \end{aligned} \]
  By Lemma~\ref{lem:weighted-gronwall}, it then follows that
  \[ \sup_{\sigma, \mu | \sigma \geqslant \mu \geqslant \bar{\mu}} \|
     \zeta_{\mu} K_{\sigma} R_{\sigma} \| \lesssim \exp (C C_{\mathcal{Q}}
     \llbracket \bar{\mu} \rrbracket^{\vartheta / 2}), \]
  possibly after adjusting the constant $C$. Hence, by~\eqref{eq:sharp-norm},
  \begin{equation}
    \trinorm{ K_{\Cdot} R_{\Cdot} }_{\#} \lesssim \llbracket
    \bar{\mu} \rrbracket^{3 \gamma} \exp (C C_F \llbracket \bar{\mu}
    \rrbracket^{\vartheta / 2} (1 + \trinorm{ \phi })^S) \lesssim
    \llbracket \bar{\mu} \rrbracket^{3 \gamma} \exp (C C_F \llbracket
    \bar{\mu} \rrbracket^{\vartheta / 2} \Phi^S) . \label{eq:bound-R-sharp}
  \end{equation}
  The control of $\trinorm{ \mathcal{L} \phi_{\Cdot} }_{\#}$ is
  obtained from~\eqref{eq:est-F}
  \[ \begin{aligned}
       \trinorm{ \mathcal{L} \phi_{\Cdot} + \lambda \phi^3_{\Cdot}
       }_{\#} & \leqslant  \trinorm{ (J_{\Cdot} F_{\Cdot}
       (\phi_{\Cdot}) + \lambda \phi^3_{\Cdot}) }_{\#} +
       \trinorm{ J_{\Cdot} R_{\Cdot} }_{\#}\\
       & \lesssim  \llbracket \bar{\mu} \rrbracket^{\vartheta} C_F (1 +
       \trinorm{ \phi })^S + \llbracket \bar{\mu}
       \rrbracket^{\vartheta} (1 + \trinorm{ \phi })^2
       \trinorm{ \mathcal{L} \phi_{\Cdot} }_{\#} + \trinorm{
       K_{\Cdot} R_{\Cdot} }_{\#},
     \end{aligned} \]
  using in particular Lemma~\ref{lem:bounds-K-J} to bound the contribution
  from $R$. From this it follows that
  \begin{equation}
    \trinorm{ \mathcal{L} \phi_{\Cdot} + \lambda \phi^3_{\Cdot}
    }_{\#} \lesssim (\llbracket \bar{\mu} \rrbracket^{\vartheta} C_F
    + \theta_{\mathcal{L}} \llbracket \bar{\mu} \rrbracket^{\vartheta} +
    \theta_R) \Phi^S \label{eq:bound-force}
  \end{equation}
  and
  \begin{equation}
    \trinorm{ \mathcal{L} \phi_{\Cdot} }_{\#} \leqslant
    \trinorm{ \lambda \phi^3_{\Cdot} }_{\#} + \trinorm{
    \mathcal{L} \phi_{\Cdot} + \lambda \phi^3_{\Cdot} }_{\#}
    \lesssim (\lambda + \llbracket \bar{\mu} \rrbracket^{\vartheta} C_F +
    \theta_{\mathcal{L}} \llbracket \bar{\mu} \rrbracket^{\vartheta} +
    \theta_R) \Phi^S, \label{eq:bound-Lphi-sharp}
  \end{equation}
  by the triangle inequality and
  \[ \trinorm{ \lambda \phi^3_{\Cdot} }_{\#} = \lambda
     \sup_{\sigma \geqslant \bar{\mu}} \llbracket \sigma \rrbracket^{3 \gamma}
     \| \zeta_{\sigma} \phi_{\sigma}^3 \| \lesssim \lambda [\sup_{\sigma
     \geqslant \bar{\mu}} \llbracket \sigma \rrbracket^{\gamma}  \|
     \zeta_{\sigma}^{1 / 3} \phi_{\sigma} \|]^3 = \lambda \trinorm{ \phi
     }_{\bar{\mu}}^3 . \]
  Next, the a~priori estimates of Theorem~\ref{thm:apriori}
  and~\eqref{eq:bound-force} give
  \begin{equation}
    \begin{aligned}
      \trinorm{ \phi } & \lesssim \lambda^{- 1 / 2}  \llbracket
      \bar{\mu} \rrbracket^{\gamma} + \lambda^{- 1 / 3} 
      \trinorm{ \mathcal{L} \phi_{\Cdot} + \lambda \phi^3_{\Cdot}
      }_{\#}^{1 / 3}\\
      & \lesssim \lambda^{- 1 / 2}  \llbracket \bar{\mu}
      \rrbracket^{\gamma} + \lambda^{- 1 / 3}
      ( \llbracket \bar{\mu} \rrbracket^{\vartheta} C_F +
      \theta_{\mathcal{L}} \llbracket \bar{\mu} \rrbracket^{\vartheta} +
      \theta_R)^{1 / 3} \Phi^{S / 3} .
    \end{aligned} \label{eq:bound-phi-triple}
  \end{equation}
  Gathering~\eqref{eq:bound-R-sharp},~\eqref{eq:bound-Lphi-sharp}
  and~\eqref{eq:bound-phi-triple}, we obtain
  \begin{equation}
    \Phi_{\bar{\mu}} \leqslant 1 + C [\tau (\bar{\mu}) + \tau (\bar{\mu})^{1 /
    3} \Phi_{\bar{\mu}}^{S / 3} + \tau (\bar{\mu}) \Phi_{\bar{\mu}}^S + \tau
    (\bar{\mu}) \exp (C \tau (\bar{\mu}) \Phi_{\bar{\mu}}^S)],
    \label{eq:bound-Phi}
  \end{equation}
  where
  \[ \begin{aligned}
       \tau (\bar{\mu}) & \assign  \lambda^{- 1 / 2}  \llbracket \bar{\mu}
       \rrbracket^{\gamma / 2} + \lambda^{- 1} ( \llbracket
       \bar{\mu} \rrbracket^{\vartheta} C_F + \theta_{\mathcal{L}} \llbracket
       \bar{\mu} \rrbracket^{\vartheta} + \theta_R)\\
       &   + \theta_{\mathcal{L}}^{- 1} (\lambda + \llbracket \bar{\mu}
       \rrbracket^{\vartheta} C_F + \theta_{\mathcal{L}} \llbracket \bar{\mu}
       \rrbracket^{\vartheta} + \theta_R) + \theta_R^{- 1} \llbracket
       \bar{\mu} \rrbracket^{3 \gamma} + C_F \llbracket \bar{\mu}
       \rrbracket^{\vartheta / 2} .
     \end{aligned} \]
  Choose
  \begin{equation}
    \theta_R = \llbracket \bar{\mu} \rrbracket^{\vartheta / 2} \lambda, \qquad
    \theta_{\mathcal{L}} = \llbracket \bar{\mu} \rrbracket^{- \vartheta / 2}
    \lambda . \label{eq:choice-theta-LR}
  \end{equation}
  Then $\bar{\mu} \mapsto \tau (\bar{\mu})$ is a decreasing
  function and
  \begin{align*}
     \tau (\bar{\mu})
     &\leqslant  \llbracket \bar{\mu} \rrbracket^{\vartheta / 2}
     (\lambda^{- 1 / 2}  \llbracket \bar{\mu}
     \rrbracket^{(\gamma - \vartheta) / 2}
     + (2 \lambda^{- 1} + 1) C_F + 5 + \lambda^{- 1} \llbracket \bar{\mu}
     \rrbracket^{3 \gamma - \vartheta })
     \\
     &\leqslant 10 \llbracket
     \bar{\mu} \rrbracket^{\vartheta / 2} (\lambda^{- 1} + 1) C_F .
  \end{align*}
  Fix $\tau_{\ast} > 0$ small enough such that
  \begin{equation}
    C [\tau_{\ast} + \tau_{\ast}^{1 / 3} 4^{S / 3} + \tau_{\ast} 4^S +
    \tau_{\ast} \exp (C \tau_{\ast} 4^S)] \leqslant 1, \label{eq:bound-tau}
  \end{equation}
  and define $\bar{\mu}_{\ast} = \bar{\mu}_{\ast} (\lambda, C_F) \in (0, 1)$
  as the (unique) solution to
  \[ 10\, \llbracket \bar\mu_{\ast} \rrbracket^{\vartheta / 2} 
     (\lambda^{- 1} + 1) C_F = \tau_{\ast} . \]
  Then for all $\bar{\mu} \in [\bar{\mu}_{\ast}, 1)$ we have $\tau (\bar{\mu})
  \leqslant \tau_{\ast}$ and as a consequence of~\eqref{eq:bound-Phi}
  and~\eqref{eq:bound-tau},
  \[ \Phi_{\bar{\mu}} \leqslant 4 \Rightarrow \Phi_{\bar{\mu}} \leqslant 2. \]
  Define the set
  \[ A \assign \{ \bar{\mu} \in [\bar{\mu}_{\ast}, 1)  | \Phi_{\bar{\mu}}
     \leqslant 4 \} \subseteq [\bar{\mu}_{\ast}, 1) . \]
  Note that $A \neq \varnothing$ since for $\bar{\mu} \nearrow 1$ we have
  $\Phi_{\bar{\mu}} \searrow 1$. As the map $\bar{\mu} \mapsto
  \Phi_{\bar{\mu}}$ is continuous, the set $A$ is closed in
  $[\bar{\mu}_{\ast}, 1)$. Hence, to prove that $A = [\bar{\mu}_{\ast}, 1)$,
  it is enough to show that $A$ is open in $[\bar{\mu}_{\ast}, 1)$. If $\mu
  \in A$, then $\Phi_{\mu} \leqslant 2$ and by continuity, there exists a
  neighbourhood of $\mu$ within $[\bar{\mu}_{\ast}, 1)$ on which the function
  $\bar{\mu} \mapsto \Phi_{\bar{\mu}}$ takes values not exceeding $4$. We
  conclude that $A = [\bar{\mu}_{\ast}, 1)$ and therefore that
  $\Phi_{\bar{\mu}_{\ast}} \leqslant 2$. Using~\eqref{eq:choice-theta-LR} this
  implies that for $\bar{\mu} \geqslant \bar{\mu}_{\ast}$, we have
  \[ \trinorm{ \phi }_{\bar{\mu}} \leqslant 1, \qquad \trinorm{
     \mathcal{L} \phi_{\Cdot} }_{\#, \bar{\mu}} \leqslant \lambda
     \llbracket \bar{\mu} \rrbracket^{- \vartheta / 2}, \qquad \trinorm{
     K_{\Cdot} R_{\Cdot} }_{\#, \bar{\mu}} \leqslant \lambda
     \llbracket \bar{\mu} \rrbracket^{\vartheta / 2}, \]
  with
  \[ \llbracket \bar{\mu}_{\ast} \rrbracket^{- 1} = (\tau_{\ast}^{- 1} 10 (\lambda^{-
     1} + 1) C_F)^{2 / \vartheta} \lesssim ((\lambda^{- 1} + 1) C_F)^{2 /
     \vartheta}, \]
  uniformly in $C_F$ and $\lambda > 0$.
\end{proof}

To apply the previous theorem to the analysis of~\eqref{eq:rg-spde}, we need
to construct a suitable approximate solution to the flow equation
\begin{equation}
  \partial_{\sigma} F _{\sigma} + \mathD F _{\sigma}  (\dot{G}_{\sigma} F
  _{\sigma}) = 0, \label{eq:flow-eq-sec2}
\end{equation}
for the effective force
\[ (F_{\sigma})_{\sigma \in [1 / 2, 1]} = (F^{(\varepsilon,
   M)}_{\sigma})_{\sigma \in [1 / 2, 1]}, \]
subject to the final condition~\eqref{eq:final-cond-f}. Section~\ref{sec:flow}
will be devoted to the construction of such an effective force
$(F_{\sigma})_{\sigma}$, as formulated in the following theorem.

\begin{theorem}
  \label{thm:effective-force}There exists a choice of the deterministic
  parameters
  \[ (r_{\varepsilon, M})_{\varepsilon \in 2^{-\mathbb N_0}, M \in \mathbb{N}_+} \]
  and a family of random, scale-dependent functionals
  \[ \left( F_{\sigma} {= F_{\sigma}^{(\varepsilon, M)}}  : \mathcal{E}
     \rightarrow \hat{\mathcal{E}} \right)_{\sigma \in [1 / 2, 1]} \]
  such that:
  \begin{enumerate}[(a)]
    \item The map $\sigma \mapsto F_{\sigma}$ is differentiable for $\sigma
    \in (1 / 2, 1)$.
    
    \item The final condition~\eqref{eq:final-cond-f} holds.
    
    \item The estimates~\eqref{eq:est-F} are satisfied with \[C_F = 1 + \|
    F^{\mathfrak{A}} \|^2\] and $\| F^{\mathfrak{A}} \|$ introduced in
    Definition~\ref{def:norm-kernels}.
    
    \item For all $N \geqslant 1$, we have
    \begin{equation}
      \sup_{\varepsilon \in 2^{-\mathbb{N}_0}, M \in \mathbb{N}_+} \mathbb{E}
      [\| F^{\mathfrak{A}} \|^N] < \infty .
      \label{eq:stochastic-estimates-kernels}
    \end{equation}
  \end{enumerate}
\end{theorem}

\begin{proof}
  Choose $\bar{r}$ and $(r^{[\ell]}_{\varepsilon, M})_{\ell = 1, \ldots,
  \hat{\ell}}$ as in Lemma~\ref{lem:cumul-est-1}, set
  \[ r_{\varepsilon, M} \assign \bar{r} + \sum_{\ell = 1}^{\hat{\ell}}
     r^{[\ell]}_{\varepsilon, M}, \]
  and let
  \[ F_{\sigma} \assign \sum_{\ell = 0}^{\bar{\ell}} F_{\sigma}^{[\ell]} \]
  be the corresponding truncated solution of the perturbative flow
  equation~\eqref{eq:pert-flow} introduced in
  Section~\ref{sec:random-flow}. By construction, $\sigma \mapsto F_{\sigma}$
  is differentiable on $(1 / 2, 1)$, and the final data prescribed
  in~\eqref{eq:initial-cond-induction-lemma} yield the final
  condition~\eqref{eq:final-cond-f}. Next, Lemmas~\ref{lem:cumul-est-1}
  and~\ref{lem:kolmogorov} yield Item~d). Finally, Item~c) follows from the post-processing
  estimates established later in the paper: the first and third bounds
  in~\eqref{eq:est-F} are given by
  Lemmas~\ref{lem:overall-difference-bound},~\ref{lem:control-of-forcing}
  and~\ref{lem:l-0-kernels}, the fourth one by
  Lemma~\ref{lem:estimates-DF}, and the bound on $H_{\sigma}$ by
  Lemma~\ref{lem:packaging}. Thus all four items hold with
  $C_F = 1 + \| F^{\mathfrak{A}} \|^2$.
\end{proof}

We complete this section with the proof of the following auxiliary lemma, used
in the proof of Theorem~\ref{thm:apriori-system}.

\begin{lemma}
  \label{lem:weighted-gronwall}Suppose that there exist constants $S \in
  \mathbb{N}_+, C_F \geqslant 1$ and $\vartheta, \bar{\kappa} > 0$ such that
  \[ \bar{\kappa} \in [\nu, 1/2), \qquad (s + \gamma)  \bar{\kappa} / (1 -
     \bar{\kappa}) \leqslant \vartheta / 4 \]
  and $(F_{\sigma})_{\sigma}$ satisfies the following estimate
  \begin{equation}
    \begin{aligned}
      \| \zeta_{\mu} K_{\sigma} (\mathD F_{\sigma} (\psi_{\sigma})
      (\dot{G}_{\sigma}  \hat{\psi})) \|_{L^{\infty}} & \leqslant C_F 
      \llbracket \sigma \rrbracket^{\vartheta - 1}  (1 + \trinorm{ \psi
      }_{\bar{\mu}})^S  \| \zeta^{1 - \bar{\kappa}}_{\mu}
      \tilde{J}_{\sigma}^2 \hat{\psi} \|_{L^{\infty}},
    \end{aligned} \label{eq:est-DF}
  \end{equation}
  for all $\bar{\mu} \in [1 / 2, 1)$, $\mu \in [\bar{\mu}, 1)$, $\sigma \in [\mu,
  1)$, and $\psi, \hat{\psi} \in \mathcal{S}' (\Lambda)$. Here we set
  $\psi_{\sigma} \assign J_{\sigma} \psi$. Fix a function
  $[\bar{\mu}, 1) \ni \mu \mapsto \mathcal{H}_{\mu} \in \mathcal{S}'
  (\Lambda)$ and consider the linear equation
  \begin{equation}
    R_{\mu} =\mathcal{H}_{\mu} + \int_{\mu}^1 \mathD F_{\sigma}
    (\phi_{\sigma}) (\dot{G}_{\sigma} R_{\sigma})\,\mathd \sigma, \qquad \mu
    \geqslant \bar{\mu} . \label{eq:lin-R}
  \end{equation}
  Then there exists a constant $C > 0$ such that
  \begin{equation}
    \begin{aligned}
      \| \zeta_{\eta} K_{\mu} R_{\mu} \|
      \leqslant \exp (C C_{\mathcal{Q}}
      \llbracket \bar{\mu} \rrbracket^{\vartheta / 2})
      \,
      \left[\sup_{\sigma \geqslant \mu}
      \| \zeta_{\eta} K_{\sigma} \mathcal{H}_{\sigma} \|
      + \sup_{\sigma \geqslant \mu}
      \llbracket \sigma \rrbracket^{2 (s + \gamma)}
      \| \zeta_{\eta}^{\bar{\kappa}} \tilde{J}_{\sigma}^2 R_{\sigma} \|\right],
    \end{aligned}
    \label{eq:gronwall}
  \end{equation}
  for
  $\mu \geqslant \eta \geqslant \bar{\mu}$,
  where $C_{\mathcal{Q}} = C_F  (1 + \trinorm{ \phi })^S$.
\end{lemma}

\begin{proof}
  Start by observing that, by interpolation, for any $\alpha \in [0, 1]$ and
  $\beta \in \mathbb{R}$, we have by Young's inequality
  \begin{align*}
    &\llbracket \sigma \rrbracket^{\beta (1 - \alpha)}
    \| \zeta_{\eta}^{\alpha + \bar{\kappa} (1 - \alpha)}
    \tilde{J}_{\sigma}^2 R_{\sigma} \|\\
    &\quad \lesssim
    \llbracket \sigma \rrbracket^{\beta (1 - \alpha)}
    \| \zeta_{\eta} \tilde{J}_{\sigma}^2 R_{\sigma} \|^{\alpha}
    \| \zeta_{\eta}^{\bar{\kappa}}
    \tilde{J}_{\sigma}^2 R_{\sigma} \|^{1 - \alpha}\\
    &\quad \lesssim
    \| \zeta_{\eta} \tilde{J}_{\sigma}^2 R_{\sigma} \|
    + \llbracket \sigma \rrbracket^{\beta}
    \| \zeta_{\eta}^{\bar{\kappa}} \tilde{J}_{\sigma}^2 R_{\sigma} \| .
  \end{align*}
  Choose $\alpha = (1 - 2 \bar{\kappa}) / (1 - \bar{\kappa})$. Then $\alpha +
  \bar{\kappa} (1 - \alpha) = 1 - \bar{\kappa}$ and $1 - \alpha = \bar{\kappa}
  / (1 - \bar{\kappa})$. Consequently,
  \begin{equation}
    \llbracket \sigma \rrbracket^{\beta \bar{\kappa} / (1 - \bar{\kappa})} \|
    \zeta_{\eta}^{1 - \bar{\kappa}} \tilde{J}_{\sigma}^2 R_{\sigma} \| \lesssim
    \| \zeta_{\eta}  \tilde{J}_{\sigma}^2 R_{\sigma} \| + \llbracket \sigma
    \rrbracket^{\beta}  \| \zeta_{\eta}^{\bar{\kappa}} \tilde{J}_{\sigma}^2
    R_{\sigma} \| . \label{eq:interp}
  \end{equation}
  Then, by~\eqref{eq:est-DF}, we have
  \[ \| \zeta_{\eta} K_{\sigma} (\mathD F_{\sigma} (\phi_{\sigma})
     (\dot{G}_{\sigma} R_{\sigma})) \| \lesssim \llbracket \sigma
     \rrbracket^{\vartheta - 1} C_{\mathcal{Q}}  \| \zeta^{1 -
     \bar{\kappa}}_{\eta} \tilde{J}_{\sigma}^2 R_{\sigma} \|, \]
  and by~\eqref{eq:interp} with $\beta = 2 (s + \gamma)$, and
  Lemma~\ref{lem:bounds-K-J}, we have
  \begin{equation}
    \begin{aligned}
      &\| \zeta_{\eta} K_{\sigma}
      (\mathD F_{\sigma}(\phi_{\sigma})(\dot{G}_{\sigma} R_{\sigma})) \|\\
      &\quad \lesssim C_{\mathcal{Q}}
      \llbracket \sigma \rrbracket^{\vartheta - 1}
      \llbracket \sigma \rrbracket^{- 2 (s + \gamma)\bar{\kappa}
      / (1 - \bar{\kappa})}\\
      &\qquad \times
      \bigl(\| \zeta_{\eta} K_{\sigma} R_{\sigma} \|
      + \llbracket \sigma \rrbracket^{2 (s + \gamma)}
      \| \zeta_{\eta}^{\bar{\kappa}}
      \tilde{J}_{\sigma}^2 R_{\sigma} \|\bigr) .
    \end{aligned}
    \label{eq:grad-interp}
  \end{equation}
  At this point, to deduce the desired estimate, we apply Gronwall's
  inequality to
  \begin{align*}
    \| \zeta_{\eta} K_{\mu} R_{\mu} \|
    &\leqslant \| \zeta_{\eta} K_{\mu} \mathcal{H}_{\mu} \|\\
    &\quad + C C_{\mathcal{Q}} \int_{\mu}^1
    \llbracket \sigma \rrbracket^{
    \vartheta - 1 - 2 (s + \gamma)\bar{\kappa} / (1 - \bar{\kappa})}\\
    &\qquad \times
    \bigl(\| \zeta_{\eta} K_{\sigma} R_{\sigma} \|
    + \llbracket \sigma \rrbracket^{2 (s + \gamma)}
    \| \zeta_{\eta}^{\bar{\kappa}}
    \tilde{J}_{\sigma}^2 R_{\sigma} \|\bigr)\,\mathd \sigma,
  \end{align*}
  which follows from~\eqref{eq:lin-R} and~\eqref{eq:grad-interp}, and use the
  inequality $\vartheta - 2 (s + \gamma) \bar{\kappa} / (1 - \bar{\kappa})
  \geqslant \vartheta / 2$.
\end{proof}

\subsection{Tightness}\label{sec:tightness}

In this section, we apply Theorem~\ref{thm:apriori-system}, together with the
effective force constructed in Theorem~\ref{thm:effective-force}, to prove
tightness of the sequence of measures ${(\hat{\nu}_{\varepsilon,
M})_{\varepsilon, M}} $. This constitutes the first step in the proof of our
main result, Theorem~\ref{thm:main}. In order to pass to the limit
$\varepsilon \rightarrow 0$ and $M \rightarrow \infty$, we embed all the random
spatially periodic fields $\phi = \phi^{(\varepsilon, M)}$ in the same space
by extending them from $\Lambda_{\varepsilon}$ to $\Lambda_0$. Let
\begin{equation}
  \phi^{[\varepsilon, M]} (t, x) \assign \int_{\Lambda^{\ast}_{\varepsilon}} 
  \hat{\theta} (\varepsilon | k |)  \hat{\phi}^{(\varepsilon, M)} (\omega, k)
  e^{i (\omega t + k \cdot x)} \frac{\mathd \omega \mathd k}{(2 \pi)^{d + 1}},
  \qquad (t, x) \in \Lambda_0, \label{eq:embedding}
\end{equation}
where \(\hat{\theta}:\mathbb{R}_+\to\mathbb{R}_+\) is a smooth bump function such that
\(\hat{\theta}(0)=1\) and \(\hat{\theta}(\eta)=0\) for \(|\eta|>1\). The
random fields $\phi^{[\varepsilon, M]}$ now all live in the continuum domain
$\Lambda_0$ for any $\varepsilon \in 2^{-\mathbb{N}_0}$. To extract
information about the EQFT, we must evaluate the fixed-time marginal of
the SPDE solution. This is not possible when controlling only spacetime
distributional norms such as~$\trinorm{ \Cdot }_{\bar{\mu}}$.
The necessary temporal regularity can, however, be recovered by means of a
Schauder estimate adapted to our weighted norms. Set $\rho = \rho_0$.
Indeed, by Lemma~\ref{lem:schauder-full} and using $\rho \leqslant \rho_{\bar{\mu}}$, we
obtain
\[ \sup_i 2^{- i \gamma} \| \rho  \bar{\Delta}_i  \phi \| \lesssim \sup_i 2^{-
   i \gamma} \| \rho_{\bar{\mu}} \bar{\Delta}_i  \phi \|  \lesssim \llbracket
   \bar{\mu} \rrbracket^{- \gamma} [\trinorm{ \phi }_{\bar{\mu}} +
   \trinorm{ K_{\Cdot} \mathcal{L} \phi }_{\#, \bar{\mu}}] . \]
By Theorem~\ref{thm:apriori-system} there exists $\bar{\mu} \in [1 / 2, 1)$
such that $\llbracket \bar{\mu} \rrbracket^{- \vartheta / 2}$ is of order $C_F
= 1 + \| F^{\mathfrak{A}} \|^2$ and
\[ \trinorm{ \phi }_{\bar{\mu}} \leqslant 1, \qquad \trinorm{
   \mathcal{L} \phi_{\Cdot} }_{\#, \bar{\mu}} \leqslant \lambda
   \llbracket \bar{\mu} \rrbracket^{- \vartheta / 2}, \qquad \trinorm{
   K_{\Cdot} R_{\Cdot} }_{\#, \bar{\mu}} \leqslant \lambda
   \llbracket \bar{\mu} \rrbracket^{\vartheta / 2} . \]
Recalling that $\mathcal{L} \phi = F_{\mu} (\phi_{\mu}) + R_{\mu}$, we deduce
from the above bounds that
\begin{equation}
  {\trinorm{ K_{\Cdot} \mathcal{L} \phi }_{\#, \bar{\mu}}} 
  \leqslant \trinorm{ K_{\Cdot} F_{\Cdot} (\phi_{\Cdot})
  }_{\#, \bar{\mu}} + \trinorm{ K_{\Cdot} R_{\Cdot}
  }_{\#, \bar{\mu}} \lesssim C_F  (1 + \trinorm{ \phi
  }_{\bar{\mu}})^S + \lambda \llbracket \bar{\mu}
  \rrbracket^{\vartheta / 2} \lesssim C_F . \label{eq:est-Lphi}
\end{equation}
Combining this with the previous Schauder estimate yields
\begin{equation}
  \sup_i 2^{- i \gamma} \| \rho \bar{\Delta}_i \phi^{(\varepsilon, M)}
  \|_{L^{\infty} (\Lambda_{\varepsilon})} = \sup_i 2^{- i \gamma} \| \rho 
  \bar{\Delta}_i  \phi \| \lesssim \llbracket \bar{\mu} \rrbracket^{- \gamma}
  C_F \lesssim  C_F^{2 \gamma / \vartheta + 1} . \label{eq:schauder-est}
\end{equation}
To derive a meaningful estimate, we need to control a suitable weighted Besov
norm of $\phi^{[\varepsilon, M]}$ in terms of the left-hand side of the bound
above. The argument of \cite[Lemma~2.24]{MartinPerkowski2019} gives, uniformly in $\varepsilon$ and $M$,
\begin{equation}
  \sup_{t \in \mathbb{R}}
  \|\phi^{[\varepsilon,M]}(t)\|_{\mathcal C^{-\gamma}
  (\mathbb R^d,\rho(t))}
  \lesssim
  \sup_i 2^{-i\gamma}
  \|\rho\bar\Delta_i\phi^{(\varepsilon,M)}\|_{L^\infty
  (\Lambda_\varepsilon)}.
  \label{eq:weighted-extension}
\end{equation}
Although the lattice blocks in that reference are obtained by restricting the
continuum dyadic partition, the same proof applies to the present blocks: on
the support of the embedding multiplier one has
$q_\varepsilon(k)\asymp |k|$, and Lemma~\ref{lem:q-epsilon-bounds} provides
the required uniform multiplier bounds.
As a result, we have
\[ \sup_{t \in \mathbb{R}} \| \phi^{[\varepsilon, M]} (t) \|_{\mathcal{C}^{-
   \gamma} (\mathbb{R}^d, \rho (t))} \lesssim
   C_F^{2 \gamma / \vartheta + 1} \lesssim 1 + \| F^{\mathfrak{A}} \|^{4
   \gamma / \vartheta + 2}, \]
where $\mathcal{C}^{\alpha} (\mathbb{R}^d, w)$ stands for the
usual H{\"o}lder--Besov norm with the weight $w$ and the regularity index
$\alpha$ 
and here $\rho=\rho_0$, so that $\rho(t)=\rho(t,\Cdot)$. Let $\bar{\rho}(x) \assign \rho(0, x)$.
Combining the above estimate with
Theorem~\ref{thm:effective-force}, we arrive at
\begin{equation}
  \begin{aligned}
    &\sup_{\varepsilon \in 2^{-\mathbb{N}_0}, M \in \mathbb{N}_+}
    \int \| \varphi \|_{\mathcal{C}^{- \gamma}
    (\mathbb{R}^d, \bar{\rho})}^N
    \hat{\nu}_{\varepsilon, M}(\mathd \varphi)\\
    &\quad\leqslant
    \sup_{\varepsilon \in 2^{-\mathbb{N}_0}, M \in \mathbb{N}_+}
    \mathbb{E} \left[\sup_{t \in \mathbb{R}}
    \| \phi^{[\varepsilon, M]} (t) \|_{\mathcal{C}^{- \gamma}
    (\mathbb{R}^d, \rho (t))}^N\right]\\
    &\quad \lesssim
    \sup_{\varepsilon \in 2^{-\mathbb{N}_0}, M \in \mathbb{N}_+}
    \mathbb{E} \left[1 + \| F^{\mathfrak{A}} \|^{4 \gamma N / \vartheta + 2 N}\right]
    < \infty,
  \end{aligned}
  \label{eq:tightness}
\end{equation}
for all $N \in \mathbb{N}_+$. The bound~\eqref{eq:tightness} implies the
tightness of the family $(\hat{\nu}_{\varepsilon, M})_{\varepsilon \in
2^{-\mathbb{N}_0}, M \in \mathbb{N}_+}$ in $\mathcal{S}' (\mathbb R^d)$. This
proves the first part of Theorem~\ref{thm:main}.

\subsection{Integrability}\label{sec:integrability}

In order to complete the proof of Theorem~\ref{thm:main} it remains to
establish the integrability property of the measures
$\hat{\nu}_{\varepsilon, M}$ uniformly in $\varepsilon \in 2^{-\mathbb N_0}$ and $M
\in \mathbb{N}_+$ and obtain the bound~\eqref{eq:int-bound} for any limit
point. We aim to estimate quantities of the form
\[ Z_{\varepsilon, M, \theta} \assign \int \exp\Bigl(\theta \| hQ_{\varepsilon}
\varphi \|_{L^2 (\mathbb{T}^d_{\varepsilon, M})}^4 \Bigr)\,\nu_{\varepsilon, M}
(\mathd \varphi), \]
where $\theta > 0$ is a small parameter.
Here
\begin{equation}
	Q_{\varepsilon} \assign (1 - \Delta_{\varepsilon} )^{- A / 2} \label{eq:Q-A}
\end{equation}
is a regularising kernel. For \(M\geq1\), let
\begin{equation}
	\label{eq:h-B}
	|x|_M\assign\min_{k\in M\mathbb Z^d}|x+k|,
	\qquad
	h(x)\assign(1+|x|_M)^{-B}.
\end{equation}
Thus \(h\) is the \(M\)-periodic extension of
\(x\mapsto(1+|x|)^{-B}\) from the centred fundamental domain
\([-M/2,M/2)^d\).
Whenever a non-periodic weight such as \(\rho_\sigma\) occurs inside an
integral over \(\mathbb T_{\varepsilon,M}^d\), the integral is understood over
this centred fundamental domain.
The constants $A$
and $B$ are chosen big enough according to Lemma~\ref{lem:bound-o} below.

The main tool for this purpose is the Hairer--Steele argument
{\cite{hairerPhiMeasureHas2022}}, which provides optimal estimates with
respect to the growth of the function. We introduce a new \textit{tilted}
probability measure
\begin{equation}
  \nu_{\varepsilon, M, \theta} (\mathd \varphi) \assign \frac{\exp\Bigl(\theta \| hQ_{\varepsilon}
  	\varphi \|_{L^2 (\mathbb{T}^d_{\varepsilon, M})}^4 \Bigr)\,
  \nu_{\varepsilon, M} (\mathd \varphi)}{Z_{\varepsilon, M, \theta}} .
  \label{eq:exp-tilt}
\end{equation}
To prove that this measure is well defined we take advantage
of the presence of the coercive term $\| \varphi \|_{L^4
(\mathbb{T}^d_{\varepsilon, M})}^4$ in the action functional~\eqref{eq:action}
defining the original measure $\nu_{\varepsilon, M}$ and use the inequality
\[ \theta \| hQ_{\varepsilon} \varphi \|_{L^2
   (\mathbb{T}^d_{\varepsilon, M})}^4 \leqslant \frac{\lambda}{4}  \| \varphi
   \|_{L^4 (\mathbb{T}^d_{\varepsilon, M})}^4, \]
for all $\varepsilon$ and $M$ as soon as $\theta > 0$ is small
enough and $A, B$ are big enough. Observe that Jensen's inequality gives
\begin{align*}
  1 &= \int \nu_{\varepsilon, M} (\mathd \varphi)\\
    &= Z_{\varepsilon, M, \theta}
   \int \exp\Bigl(- \theta \| hQ_{\varepsilon} \varphi \|_{L^2
   (\mathbb{T}^d_{\varepsilon, M})}^4\Bigr)\, \nu_{\varepsilon, M, \theta} (\mathd
   \varphi)\\
    &\geqslant Z_{\varepsilon, M, \theta} \exp \left[ - \theta \int \|
   hQ_{\varepsilon} \varphi \|_{L^2 (\mathbb{T}^d_{\varepsilon, M})}^4
   \nu_{\varepsilon, M, \theta} (\mathd \varphi) \right] .
\end{align*}
Hence,
\begin{align}
  &\log \int \exp [\theta \| hQ_{\varepsilon} \varphi \|_{L^2
  (\mathbb{T}^d_{\varepsilon, M})}^4] \nu_{\varepsilon, M} (\mathd \varphi)
  \notag\\
  &\qquad = \log Z_{\varepsilon, M, \theta}
  \leqslant \theta \int \| hQ_{\varepsilon}
  \varphi \|_{L^2 (\mathbb{T}^d_{\varepsilon, M})}^4
  \nu_{\varepsilon, M, \theta} (\mathd \varphi) . \label{eq:hs-tilt}
\end{align}
The problem of controlling the size of $Z_{\varepsilon, M, \theta}$ is, by
virtue of~\eqref{eq:hs-tilt}, reduced to estimating certain polynomial moments
of $\varphi$ under the tilted measure $\nu_{\varepsilon, M, \theta} (\mathd
\varphi)$. Through stochastic quantisation, this measure can be identified
with the marginal law of a stationary solution to the SPDE
\begin{equation}
  \mathcal{L}_{\varepsilon} \phi  + \lambda \phi ^3  - r_{\varepsilon, M}
  \phi  = O (\phi ) + \xi^{(\varepsilon, M)}, \label{eq:tilted-spde}
\end{equation}
where the additional perturbation $O (\phi)$ is given by
\begin{align*}
	O (\phi) (t, \Cdot) &=  \frac12 \theta \frac{\delta}{\delta
		\phi} \| h Q_{\varepsilon} \phi (t, \Cdot) \|_{L^2
		(\mathbb{T}^d_{\varepsilon, M})}^4
		\\
		&= 2 \theta \| h Q_{\varepsilon} \phi
	(t, \Cdot) \|_{L^2 (\mathbb{T}^d_{\varepsilon, M})}^2 
	(Q_{\varepsilon} h^2 Q_{\varepsilon} \phi) (t, \Cdot) .
\end{align*}
We use the notation
\[ Q_{\varepsilon} h^2 Q_{\varepsilon} \phi = 
   Q_{\varepsilon} (h^2 Q_{\varepsilon} \phi), \qquad (h^i Q_{\varepsilon}
   \phi)  (t, x) = (h  (x))^i  (Q_{\varepsilon} \phi)
    (t, x), \quad i \in \{ 1, 2 \} . \]
We exploited above the fact that the operator $Q_{\varepsilon}$ is
self-adjoint and bounded. Since $Q_{\varepsilon}$ has a uniformly rapidly
decaying periodised lattice kernel, $Q_{\varepsilon} \phi$ is a well-defined
periodic lattice function with mild polynomial growth. Hence, $h
Q_{\varepsilon} \phi (t, \Cdot)$ is in $L^2 (\mathbb{T}^d_{\varepsilon, M})$.

Uniform estimates for this new equation in $\varepsilon \in 2^{-\mathbb{N}_0}$
and $M \in \mathbb{N}_+$ can be obtained by modifying our previous arguments.
For $\mu \in [1/2,1]$ set
\[ R^O_{\mu} \assign F_1 (\phi) + O (\phi) - F_{\mu} (\phi_{\mu}) . \]
As before we rewrite~\eqref{eq:tilted-spde} as a system of equations
\begin{equation}
  \left\{\begin{array}{l}
    \mathcal{L}_{\varepsilon} \phi _{\mu} = J_{\mu} (F _{\mu} (\phi _{\mu}) +
    R^O_{\mu})\\
    R^O_{\mu} = \int_{\mu}^1 H_{\sigma}  (\phi _{\sigma})\,\mathd \sigma + O
    (\phi ) + \int_{\mu}^1 \mathD F _{\sigma} (\phi _{\sigma}) 
    (\dot{G}_{\sigma}  R^O_{\sigma})\,\mathd \sigma
  \end{array}\right. \label{eq:rem-O}
\end{equation}
for
\[ (\phi _{\mu}, R^O_{\mu})_{\mu} = (\phi_{\mu}^{(\varepsilon, M)},
   R^{(\varepsilon, M), O}_{\mu})_{\mu} . \]
Applying estimate~\eqref{eq:gronwall} with $\eta = \mu$ to the equation for
$R_{\Cdot}^O$, we obtain
\begin{align*}
   &\llbracket \mu \rrbracket^{3 \gamma}
   \| \zeta_{\mu} K_{\mu} R^O_{\mu} \|\\
   &\quad \leqslant \exp (C C_{\mathcal{Q}} \llbracket \bar{\mu}
   \rrbracket^{\vartheta / 2})\\
   &\qquad \times
   \left[\sup_{\sigma \geqslant \mu}  \llbracket \mu
   \rrbracket^{3 \gamma}  \| \zeta_{\mu} K_{\sigma} \mathcal{H}_{\sigma} \|
   + \llbracket \bar{\mu} \rrbracket^{3 \gamma} \sup_{\sigma \geqslant \mu}
   \llbracket \sigma \rrbracket^{2 (s + \gamma)}
   \| \zeta_{\mu}^{\bar{\kappa}} \tilde{J}_{\sigma}^2 R_{\sigma}^O \|\right],
\end{align*}
with
\[ \mathcal{H}_{\mu} = O (\phi) + \int_{\mu}^1 H_{\sigma} (\phi_{\sigma})
   \mathd \sigma . \]
The term with $H_{\sigma}$ can be estimated as previously
from~\eqref{eq:est-F}, while to estimate $O (\phi)$ we use the uniform
boundedness of $\| K_{\sigma} \|_{\operatorname{TV} (\zeta^{- 1})}$ and
Lemma~\ref{lem:bound-o} below. We obtain that
\[ \sup_{\sigma | \sigma \geqslant \mu}  \| \zeta_{\mu} K_{\sigma}
   \mathcal{H}_{\sigma} \| \lesssim C_F  \llbracket \bar{\mu}
   \rrbracket^{\vartheta}  (1 + \trinorm{ \phi })^S + \theta
   \llbracket \mu \rrbracket^{- 3 \gamma} [\trinorm{ \phi } +
   \trinorm{ K_{\Cdot} \mathcal{L} \phi }_{\#}]^3, \]
uniformly in $\mu \geqslant \bar{\mu} .$ We also observe that
\[ \sup_{\sigma, \mu | \sigma \geqslant \mu \geqslant \bar{\mu}} \llbracket
   \sigma \rrbracket^{2 (s + \gamma)}  \| \zeta_{\mu}^{\bar{\kappa}}
   \tilde{J}_{\sigma}^2 R_{\sigma}^O \| \lesssim \llbracket \bar{\mu}
   \rrbracket^{\vartheta} C_F (1 + \trinorm{ \phi })^S, \]
uniformly in $\mu \geqslant \bar{\mu}$, by an argument analogous to the one
leading to~\eqref{eq:J-R-bound}. Therefore, combining the above estimates, we
obtain
\[ \trinorm{ K_{\Cdot} R^O_{\Cdot} }_{\#} \lesssim \exp (C
   C_{\mathcal{Q}} \llbracket \bar{\mu} \rrbracket^{\vartheta / 2})  [C_F 
   \llbracket \bar{\mu} \rrbracket^{3 \gamma + \vartheta}  (1 + \trinorm{
   \phi })^S + \theta [\trinorm{ \phi } + \trinorm{
   K_{\Cdot} \mathcal{L} \phi }_{\#}]^3] . \]
We need a good bound for $\trinorm{ K_{\Cdot} \mathcal{L} \phi
}_{\#}$ in terms of $\Phi = \Phi_{\bar{\mu}}$ defined
by~\eqref{eq:def-Phi} with $R$ replaced by $R^O$. Observe that,
by~\eqref{eq:K-force-coercive-remainder} below, the decomposition of the force, and
Lemma~\ref{lem:bounds-K-J}, we have
\begin{equation}
  \begin{aligned}
    \trinorm{ K_{\Cdot} \mathcal{L} \phi }_{\#} & \lesssim
    \trinorm{ K_{\Cdot} F (\phi_{\Cdot}) }_{\#} +
    \trinorm{ K_{\Cdot} R^O_{\Cdot} }_{\#}\\
    & \lesssim \lambda \Phi^3 + \llbracket \bar{\mu} \rrbracket^{\vartheta}
    C_F \Phi^S + \trinorm{ K_{\Cdot} R^O_{\Cdot} }_{\#}\\
    & \lesssim \lambda \Phi^3 + \llbracket \bar{\mu} \rrbracket^{\vartheta}
    C_F \Phi^S + \theta_R \Phi .
  \end{aligned} \label{eq:K-L-phi-O}
\end{equation}
As a result, for $\theta
\in (0, 1]$, we obtain
\begin{equation}
  \begin{aligned}
    \trinorm{ K_{\Cdot} R^O_{\Cdot} }_{\#}
    &\lesssim
    \exp (C C_F \Phi^S \llbracket \bar{\mu} \rrbracket^{\vartheta / 2})\\
    &\quad \times [C_F \Phi^S \llbracket \bar{\mu}
    \rrbracket^{3 \gamma + \vartheta} + \theta (1 + \lambda)^3 \Phi^{3S}
    + \theta (\llbracket \bar{\mu} \rrbracket^{\vartheta}
    C_F \Phi^S)^3
    + \theta \theta_R^3 \Phi^3]\\
    &\lesssim \exp (C C_F \Phi^S \llbracket \bar{\mu}
    \rrbracket^{\vartheta / 2})\\
    &\quad \times [\llbracket \bar{\mu} \rrbracket^{\vartheta / 2}
    + \theta (1 + \lambda)^3 + \theta \theta_R^3 ] \Phi^{3S},
  \end{aligned} \label{eq:bound-ro}
\end{equation}
with possibly different constants $C$. In the last inequality we absorbed
polynomial factors of
$C_F \llbracket \bar{\mu} \rrbracket^{\vartheta / 2} \Phi^S$ into the
exponential. Indeed, $x + x^3 \lesssim \exp (c x)$ for $x \geqslant 0$.
Gathering~\eqref{eq:bound-ro},~\eqref{eq:bound-Lphi-sharp}
and~\eqref{eq:bound-phi-triple} we arrive at
\begin{equation}
  \Phi_{\bar{\mu}} \leqslant 1 + C [\tau (\bar{\mu}) + \tau (\bar{\mu})^{1 /
  3} \Phi_{\bar{\mu}}^{S / 3} + \tau (\bar{\mu}) \Phi_{\bar{\mu}}^S + \tau
  (\bar{\mu}) \Phi_{\bar{\mu}}^{3S} \exp (C \tau (\bar{\mu})
  \Phi_{\bar{\mu}}^S)], \label{eq:bound-Phi-O}
\end{equation}
where
\[ \begin{aligned}
     \tau (\bar{\mu}) \assign & \lambda^{- 1 / 2}  \llbracket \bar{\mu}
     \rrbracket^{\gamma / 2} + \lambda^{- 1} (\llbracket \bar{\mu}
     \rrbracket^{\vartheta} C_F + \theta_{\mathcal{L}} \llbracket \bar{\mu}
     \rrbracket^{\vartheta} + \theta_R)\\
     &   + \theta_{\mathcal{L}}^{- 1} (\lambda
     + \llbracket \bar{\mu} \rrbracket^{\vartheta} C_F + \theta_{\mathcal{L}}
     \llbracket \bar{\mu} \rrbracket^{\vartheta} + \theta_R)\\
     &   + \theta_R^{- 1} [\llbracket \bar{\mu} \rrbracket^{\vartheta / 2} +
     \theta (1 + \lambda)^3 + \theta \theta_R^3]
     + C_F \llbracket \bar{\mu} \rrbracket^{\vartheta / 2} 
   \end{aligned} \]
and
\begin{equation}
  \theta_R = \lambda \theta^{1 / 2}, \qquad \theta_{\mathcal{L}} = \llbracket
  \bar{\mu} \rrbracket^{- \vartheta / 2} \lambda .
  \label{eq:choice-theta-LR-O}
\end{equation}
Then $\bar{\mu} \mapsto \tau (\bar{\mu})$ is non-increasing
and, for $\theta \in (0, 1]$, we have
\[ \tau (\bar{\mu}) \leqslant \llbracket \bar{\mu} \rrbracket^{\vartheta / 2}
   [(3 + \lambda^{- 1 / 2}) + (1 + 2 \lambda^{- 1}) C_F + \lambda^{- 1}
   \theta^{- 1 / 2}] + \theta^{1 / 2} [2 + \lambda^{- 1} (1 + \lambda)^3 +
   \lambda^2] . \]
Fix $\tau_{\ast} > 0$ so that
\[ C [\tau_{\ast} + \tau_{\ast}^{1 / 3} 4^{S / 3} + \tau_{\ast} 4^S +
   \tau_{\ast} 4^{3S} \exp (C \tau_{\ast} 4^S)] \leqslant 1. \]
Let $\theta_{\ast} = \theta_{\ast} (\lambda)$ be such that
\[ \theta_{\ast}^{1 / 2} [2 + \lambda^{- 1} (1 + \lambda)^3 + \lambda^2] =
   \tau_{\ast} / 2. \]
If needed we reduce $\tau_\ast$ further, to ensure that $\theta_\ast \leqslant
1$ and that the tilted measure with $\theta = \theta_\ast$ is well defined.
Define $\bar{\mu}_{\ast} = \bar{\mu}_{\ast} (\lambda, \theta_{\ast}, C_F)$
by
\[ \llbracket \bar{\mu}_{\ast} \rrbracket^{\vartheta / 2} [(3 + \lambda^{- 1 /
   2}) + (1 + 2 \lambda^{- 1}) C_F + \lambda^{- 1} \theta_{\ast}^{- 1 / 2}] =
   \tau_{\ast} / 2. \]
Then, for $\theta = \theta_{\ast}$ and $\bar{\mu} \geqslant
\bar{\mu}_{\ast}$, we have $\tau (\bar{\mu}) \leqslant \tau_{\ast}$,
and the continuity argument from Theorem~\ref{thm:apriori-system} can be
applied to obtain the desired estimates. We conclude that
$\Phi_{\bar{\mu}_{\ast}} \leqslant 2$. By~\eqref{eq:choice-theta-LR-O} this
implies that for $\bar{\mu} = \bar{\mu}_{\ast}$ and $\theta = \theta_{\ast}$, we have
\[ \trinorm{ \phi }_{\bar{\mu}} \leqslant 1, \qquad \trinorm{
   \mathcal{L} \phi_{\Cdot} }_{\#, \bar{\mu}} \leqslant \lambda
   \llbracket \bar{\mu} \rrbracket^{- \vartheta / 2}, \qquad \trinorm{
   K_{\Cdot} R^O_{\Cdot} }_{\#, \bar{\mu}} \leqslant \lambda
   \theta_{\ast}^{1 / 2}, \]
with
\begin{align*}
   \llbracket \bar{\mu} \rrbracket^{- 1}
   &= [2 \tau_{\ast}^{- 1} ((3 + \lambda^{- 1 / 2})
   + (1 + 2 \lambda^{- 1}) C_F
   + \lambda^{- 1} \theta_{\ast}^{- 1 / 2})]^{2 / \vartheta}\\
   &\lesssim [(1 + \lambda^{- 1}) C_F
   + \lambda^{- 1} \theta_{\ast}^{- 1 / 2}]^{2 / \vartheta},
\end{align*}
uniformly in $C_F = 1 + \| F^{\mathfrak{A}} \|^2$. By
Lemma~\ref{lem:bound-o},~\eqref{eq:K-L-phi-O}
and~\eqref{eq:choice-theta-LR-O}, we also have
\[ \begin{aligned}
     \sup_{t \in \mathbb{R}} \rho (t, 0)^4  \| hQ_{\varepsilon}
     \phi^{(\varepsilon,M)} (t) \|_{L^2 (\mathbb{T}^d_{\varepsilon, M})}^4 &
     \lesssim  \llbracket \bar{\mu} \rrbracket^{- 4 \gamma} [\trinorm{ \phi
     } + \trinorm{ K_{\Cdot} \mathcal{L} \phi
     }_{\#}]^4\\
     & \lesssim  \llbracket \bar{\mu} \rrbracket^{- 4 \gamma} [\Phi +
     \lambda \Phi^3 + \llbracket \bar{\mu} \rrbracket^{\vartheta}
     C_F \Phi^S + \theta_R \Phi]^4\\
     & \lesssim_{\lambda}  \llbracket \bar{\mu} \rrbracket^{- 4 \gamma} C_F
     \lesssim  1 + {\| F^{\mathfrak{A}} \|^{16 \gamma / \vartheta + 2}}  .
   \end{aligned} \]
As a consequence of Theorem~\ref{thm:effective-force}, we obtain
\begin{align*}
  &\sup_{\varepsilon \in 2^{-\mathbb{N}_0}, M \in \mathbb{N}_+}
  \int \| hQ_{\varepsilon} \varphi \|_{L^2 (\mathbb{T}^d_{\varepsilon, M})}^{4
  N} \nu_{\varepsilon, M, \theta_\ast} (\mathd \varphi)\\
  &\qquad =  \sup_{\varepsilon
  \in 2^{-\mathbb{N}_0}, M \in \mathbb{N}_+} \mathbb{E} [\| hQ_{\varepsilon}
  \phi^{(\varepsilon, M)} (t) \|_{L^2 (\mathbb{T}^d_{\varepsilon, M})}^{4
  N}]\\
 &\qquad \lesssim  \sup_{\varepsilon \in 2^{-\mathbb{N}_0}, M \in \mathbb{N}_+}
  \mathbb{E} [1 + \| F^{\mathfrak{A}} \|^{16 \gamma N / \vartheta + 2 N}] <
  \infty,
\end{align*}
for all $N \in \mathbb{N}_+$ and $t \in \mathbb{R}$. From this
and~\eqref{eq:hs-tilt}, we deduce easily that any accumulation point $\nu$ of
the sequence $(\hat{\nu}_{\varepsilon, M})_{\varepsilon, M}$
satisfies~\eqref{eq:int-bound} provided $\theta > 0$ is small enough. This
proves the exponential integrability, required for the Osterwalder--Schrader
reconstruction, completing the proof of Theorem~\ref{thm:main}.

\begin{remark}
  The choice of the norm used to verify exponential integrability is quite
  arbitrary. Since the corresponding tilt must be incorporated into an SPDE, we
  choose a differentiable norm. In general, we could replace the $L^2$ norm by any
  finite $L^{2 n}$ norm and similarly use a different
  spatial weight $h$ and smoothing operator $Q_{\varepsilon}$, as long as
  they remain compatible with our Schauder estimate.
\end{remark}

\begin{lemma}
  \label{lem:bound-o}Suppose that the parameters $A$ and $B$ in~\eqref{eq:Q-A}
  and~\eqref{eq:h-B} satisfy $A > \gamma$ and $B \geqslant d + 1$. Then
  \[ \| \zeta_{\sigma} O (\phi) \| \lesssim \theta \llbracket \sigma
     \rrbracket^{- 3 \gamma} [\trinorm{ \phi } + \trinorm{
     K_{\Cdot} \mathcal{L} \phi }_{\#}]^3 \]
  and
  \[ \sup_{t \in \mathbb{R}} \rho_0 (t, 0)^4  \| hQ_{\varepsilon} \phi (t)
     \|_{L^2 (\mathbb{T}^d_{\varepsilon, M})}^4 \lesssim \llbracket \bar{\mu}
     \rrbracket^{- 4 \gamma} [\trinorm{ \phi } + \trinorm{
     K_{\Cdot} \mathcal{L} \phi }_{\#}]^4 . \]
\end{lemma}

\begin{proof}
  Let us first observe that
  \begin{align*}
    \| \zeta_{\sigma} O (\phi) \| &\lesssim  \theta \| (t, x) \mapsto
    \zeta_{\sigma} (t, x)  \| h (Q_{\varepsilon} \phi) (t, \Cdot) \|_{L^2
    (\mathbb{T}^d_{\varepsilon, M})}^2  (Q_{\varepsilon} h^2 Q_{\varepsilon}
    \phi) (t, x) \|\\
 &=  \theta \left\| (t, x) \mapsto
    \left\| h \rho_{\sigma}^{- 1} (t, \Cdot)
    \rho_{\sigma} (t, \Cdot)  (Q_{\varepsilon} \phi)
    \left( t, \Cdot \right) \right\|_{L^2 (\mathbb{T}^d_{\varepsilon, M})}^2
    \right.\\
 &\qquad \left. \times
    (\zeta_{\sigma} \rho_{\sigma}^{- 1} \rho_{\sigma}
    Q_{\varepsilon} h^2 Q_{\varepsilon} \phi) (t, x) \right\|\\
 &\lesssim  \theta \| (t, x) \mapsto \| h \rho_{\sigma}^{- 1} (t,
    \Cdot) \|_{L^2 (\mathbb{T}^d_{\varepsilon, M})}^2
    (\zeta_{\sigma} \rho_{\sigma}^{- 1}) (t, x)  \|\\
 &\quad \times
    \| \rho_{\sigma} Q_{\varepsilon} \phi \|^2
    \| \rho_{\sigma} Q_{\varepsilon} h^2 Q_{\varepsilon} \phi \| .
  \end{align*}
  By Remark~\ref{rem:weight-properties} (a) we have
  \[ \rho^{- 1}_{\sigma} (t, x) \lesssim \rho^{- 1}_{\sigma} (t, 0) \rho^{-
     1}_{\sigma} (0, x) . \]
  Hence, restricting to the centred fundamental domain, where $h (x) = (1 +
  | x |)^{- B}$, and using $B \geqslant d + 1$, we arrive at
  \[ \| h \rho^{- 1}_{\sigma} (t, \Cdot) \|_{L^2 (\mathbb{T}^d_{\varepsilon,
     M})}^2 \lesssim \rho_{\sigma}^{- 2} (t, 0)  \| h \rho_{\sigma}^{- 1} (0,
     \Cdot)  \|_{L^2 (\mathbb{T}^d_{\varepsilon, M})}^2 \lesssim \rho^{-
     2}_{\sigma} (t, 0) . \]
  Since $\rho^{- 2}_{\sigma} (t, 0) \leqslant \rho^{- 2}_{\sigma} (t, x)$ for
  any $x \in \mathbb{R}^d_{\varepsilon}$ and $\rho_{\sigma} = \zeta_{\sigma}^{\nu}$ with $\nu \in
  (0, 1 / 3)$ we obtain
  \begin{align*}
     &\| (t, x) \mapsto \| h \rho_{\sigma}^{- 1} (t, \Cdot) \|_{L^2
     (\mathbb{T}^d_{\varepsilon, M})}^2
     (\zeta_{\sigma} \rho_{\sigma}^{- 1})(t, x)  \|\\
     &\quad \lesssim \| (t, x) \mapsto \rho^{- 2}_{\sigma} (t, 0)
     (\zeta_{\sigma} \rho_{\sigma}^{- 1})(t, x)  \|
     \lesssim \| \zeta_{\sigma} \rho^{- 3}_{\sigma} \| \lesssim 1.
  \end{align*}
  Since $A > \gamma$, Lemma~\ref{lem:schauder-full}, together with $\|\rho_{\sigma} Q_{\varepsilon}f\|\lesssim \|\rho_\sigma f\|$, yields
  \begin{align*}
    \| \rho_{\sigma} Q_{\varepsilon} h^2 Q_{\varepsilon} \phi \| &\lesssim 
    \|\rho_{\sigma} h^2 Q_{\varepsilon} \phi \|  \lesssim \| \rho_{\sigma} Q_{\varepsilon}
    \phi \|\\
 &\lesssim  \sum_i 2^{- i A}  \| \rho_{\sigma}  \bar{\Delta}_i  \phi \|
    \lesssim \llbracket \sigma \rrbracket^{- \gamma}  [\trinorm{ \phi
    } + \llbracket \bar{\mu} \rrbracket^{2 s - 2 \gamma} \trinorm{
    K_{\Cdot} \mathcal{L} \phi }_{\#}],
  \end{align*}
  uniformly over $\sigma \geqslant \bar{\mu}$. Combining the above estimates
  we obtain
  \[ \| \zeta_{\sigma} O (\phi) \| \lesssim \theta \llbracket \sigma
     \rrbracket^{- 3 \gamma}  [\trinorm{ \phi } + \trinorm{
     K_{\Cdot} \mathcal{L} \phi }_{\#}]^3 . \]
  Using Remark~\ref{rem:weight-properties} a) and $B \geqslant d + 1$ we get
  \[ \sup_{t \in \mathbb{R}} \rho_0 (t, 0)^4  \| h \rho^{- 1}_0 (t, \Cdot)
     \|_{L^2 (\mathbb{T}^d_{\varepsilon, M})}^4 \lesssim \| h \rho^{- 1}_0 (0,
     \Cdot) \|_{L^2 (\mathbb{T}^d_{\varepsilon, M})}^4 \lesssim 1. \]
  Since
  \[ \| \rho_0 Q_{\varepsilon} \phi \| \lesssim \| \rho_{\bar{\mu}}
     Q_{\varepsilon} \phi \| \lesssim \sum_i 2^{- i A} \| \rho_{\bar{\mu}} \bar{\Delta}_i 
     \phi \| \lesssim \llbracket \bar{\mu} \rrbracket^{- \gamma}  [\trinorm{
     \phi } + \llbracket \bar{\mu} \rrbracket^{2 s - 2 \gamma}
     \trinorm{ K_{\Cdot} \mathcal{L} \phi }_{\#}], \]
  we obtain
  \begin{align*}
       &\sup_{t \in \mathbb{R}} \rho_0 (t, 0)^4
       \| hQ_{\varepsilon} \phi (t) \|_{L^2
       (\mathbb{T}^d_{\varepsilon, M})}^4\\
       &\quad \lesssim \sup_{t \in \mathbb{R}} \rho_0 (t, 0)^4
       \| h \rho^{- 1}_0 (t, \Cdot) \|_{L^2
       (\mathbb{T}^d_{\varepsilon, M})}^4
       \| \rho_0 Q_{\varepsilon} \phi \|^4\\
       &\quad \lesssim \llbracket \bar{\mu} \rrbracket^{- 4 \gamma}
       [\trinorm{ \phi } + \trinorm{ K_{\Cdot}
       \mathcal{L} \phi }_{\#}]^4 .
  \end{align*}
  This finishes the proof.
\end{proof}

\subsection{The vector model}\label{sec:vector}

In this final subsection, we discuss the modifications required to extend our
results to the vector-valued model, where the field $\phi =
\phi^{(\varepsilon, M)}$ takes values in the Euclidean space $\mathbb{R}^n$
for some $n > 1$. We denote by $(\phi^a)_{a \in \{ 1, \ldots, n \}}$ the
components of the field $\phi$ in the canonical basis. The dynamics reads
\begin{equation}
  \mathcal{L}_{\varepsilon} \phi^a + \lambda | \phi  |^2 {\phi^a}  - r_{\varepsilon,M}  \phi^a
  = \xi^a, \qquad a = 1, \ldots, n, \label{eq:dynamics-eps-vector}
\end{equation}
on $\Lambda_{\varepsilon, M}$, where
\[ | \phi |^2  \assign \sum_a (\phi^a)^2 . \]
Here
\begin{equation}\label{eq:cov-noise-vector}
	\xi = (\xi^a)_{a \in \{ 1, \ldots, n \}} = (\xi^{(\varepsilon, M), a})_{a
		\in \{ 1, \ldots, n \}}
\end{equation}
is a vector-valued, spacetime white noise on $\Lambda_{\varepsilon, M}$ such
that
\[
  \mathbb{E} [\xi^a (t, x) \xi^b (s, y)] = \delta (t - s) \delta_{a, b}
  \,\varepsilon^{-d}\mathbbm{1}_{x = y}, \qquad (t, x), (s, y) \in \Lambda_{\varepsilon, M},
  \quad a, b = 1, \ldots, n.
\]
As before, we identify $\phi$ and $\xi$ with spatially
periodic functions on $\Lambda_{\varepsilon}$.

Our proof of the main result in the scalar-valued case extends
straightforwardly to the vector-valued case. The main subtlety lies in
obtaining the analogue of Theorem~\ref{thm:apriori} in this setting, which,
however, follows directly from Theorem~\ref{thm:partial-coercive-estimates}
below. Beyond this adjustment, the subsequent analysis carries over with
minimal changes, as it does not depend on the scalar nature of the equation
until the classification of the relevant cumulants in
Section~\ref{sec:class-rel-cumulants}. At that stage, one must exploit the
$O (n)$-symmetry of the noise~\eqref{eq:cov-noise-vector} to deduce that the
cumulants are likewise symmetric. Consequently, the only contribution to the
mass renormalisation is diagonal in the vector indices, and can therefore be
absorbed into a redefinition of the renormalisation constant $r_{\varepsilon, M} \in \mathbb{R}$.

\section{A~priori estimates}\label{sec:apriori}

This section is devoted to establishing weighted estimates for classical
solutions of a fractional parabolic equation with a cubic coercive term. To
the best of our knowledge, this result is new, although the proof closely
follows the argument developed in the case of the standard Laplacian (see, for
instance,~{\cite{gubinelliGlobalSolutionsElliptic2019}}). Furthermore, we
present the estimates directly in a vector-valued setting, as no additional
difficulties arise in this more general formulation.

\begin{theorem}
  \label{thm:partial-coercive-estimates}Let \(n\in\mathbb N_+\), $\lambda > 0$, $\bar{\nu} < 2 s$
  and $u \in C_t^1 C_x^2 (\Lambda_{\varepsilon}, \mathbb{R}^n)$ be such that
  $| \zeta^{\bar{\nu}} u |$ is a~bounded function, where the weight $\zeta$
  was introduced in Definition~\ref{def:weights}. Then $u$ belongs to the domain of
  the operator $\partial_t + (- \Delta)^s$. Define $f \in C 
  (\Lambda_{\varepsilon}, \mathbb{R}^n)$ by
  \begin{equation}
    f^a \assign \partial_t u^a + (- \Delta)^s u^a + m^2 u^a +
    \lambda | u |^2 u^a, \qquad a \in \{ 1, \ldots, n \},
    \label{eq:cubic-fractional}
  \end{equation}
  where $| u |^2  \assign \sum_a (u^a)^2$.
  
  For any positive weight $\rho \in C_t^1 C_x^2
  (\Lambda_{0})$ such that for some $\tilde{\nu} \in
  (\bar{\nu}, 2 s)$ the functions $\zeta^{- \tilde{\nu}} \rho$,
  $\zeta^{\tilde{\nu}} \rho^{- 1}$ are bounded we have
  \begin{equation}
    \| \rho u \| \leqslant 2 \lambda^{- 1 / 2} A^{1 / 2} + \lambda^{- 1 / 3} 
    (\| \rho^3 f \| + B)^{1 / 3}, \label{eq:coercive-bound}
  \end{equation}
  where
  \[ A \assign \| (- \Delta)^s \rho^2 \| + \| \rho (\partial_t \rho) \| \]
  and
  \[ B \assign \| \rho u \|  (\| \rho (\partial_t \rho) \| + \| \rho  (-
     \Delta)^s \rho \| + \| \rho^2 \mathfrak{D}_s (\rho^{- 1}) \mathfrak{D}_s
     (\rho) \|) + \| \rho  \mathfrak{D}_s (\rho) \mathfrak{D}_s (\rho u) \| .
  \]
\end{theorem}

\begin{remark}
If the weight $\rho$ were constant, both $A$ and $B$ would vanish and the
estimate would reduce to the usual maximum-principle bound with cubic damping.
Here $A$ measures the cost of applying the time derivative and the fractional
Laplacian to the weight, whereas $B$ contains the genuinely nonlocal
commutator errors produced by the interaction of $(-\Delta)^s$ with $\rho$.
\end{remark}

\begin{proof}
  We assume that $\| \rho^3 f \| + A + B < \infty$ as otherwise there is
  nothing to prove. For $N \in \mathbb{N}_+$ and $L > 0$, we define a convex
  function $\Phi = \Phi_{L, N} \in C (\mathbb{R})$ by
  \[ \Phi (\eta) \assign (\eta - L)^N_+, \qquad \Phi' (\eta)
     \assign N (\eta - L)^{N - 1}_+, \qquad \eta \in \mathbb{R}, \]
  where $\Phi'$ denotes the derivative of $\Phi$ and $(\eta)_+ \assign
  \mathbbm{1}_{\eta \geqslant 0} \eta$. We fix the parameter $N > 1$ such that
  $2 N (\tilde{\nu} - \bar{\nu}) > d + 2 s$. The parameter $L$ will be fixed
  later. Let $\hat{u} \assign \rho u$. After
  testing~\eqref{eq:cubic-fractional} with $\Phi' (| \hat{u} |^2) \hat{u}^a
  \rho^3$ and summing over $a \in \{ 1, \ldots, n \}$, we obtain
  \[ 0 = \int_{{\Lambda_{\varepsilon}} } \Phi' (| \hat{u} |^2)  [\rho^3 
     \hat{u} \cdot \partial_t u + \rho^3  \hat{u} \cdot (- \Delta)^s u + m^2
     \rho^2  | \hat{u} |^2 + \lambda (| \hat{u} |^2)^2 - \rho^3 f \cdot
     \hat{u}] . \]
  Observe that
  \begin{align*}
  	\int_{\Lambda_{\varepsilon}}
  	\Phi'(|\hat u|^2)\rho^3 \hat u\cdot \partial_t u
  	&=
  	\frac12
  	\int_{\Lambda_{\varepsilon}}
  	\Phi'(|\hat u|^2)\rho^2 \partial_t |\hat u|^2
  	-
  	\int_{\Lambda_{\varepsilon}}
  	\Phi'(|\hat u|^2)(\rho\partial_t\rho)|\hat u|^2
  	\\
  	&=
  	\frac12
  	\int_{\Lambda_{\varepsilon}}
  	\rho^2 \partial_t\Phi(|\hat u|^2)
  	-
  	\int_{\Lambda_{\varepsilon}}
  	(\rho\partial_t\rho)\Phi'(|\hat u|^2)|\hat u|^2
  	\\
  	&=
  	-\frac12
  	\int_{\Lambda_{\varepsilon}}
  	\partial_t(\rho^2)\Phi(|\hat u|^2)
  	-
  	\int_{\Lambda_{\varepsilon}}
  	(\rho\partial_t\rho)\Phi'(|\hat u|^2)|\hat u|^2
  	\\
  	&=
  	-
  	\int_{\Lambda_{\varepsilon}}
  	(\rho\partial_t\rho)\Phi(|\hat u|^2)
  	-
  	\int_{\Lambda_{\varepsilon}}
  	(\rho\partial_t\rho)\Phi'(|\hat u|^2)|\hat u|^2
  	\\
  	&\geqslant
  	-\|\rho(\partial_t\rho)\|
  	\int_{\Lambda_{\varepsilon}}\Phi(|\hat u|^2)
  	-
  	\|\hat u\|_{L^\infty}^2
  	\|\rho(\partial_t\rho)\|
  	\int_{\Lambda_{\varepsilon}}\Phi'(|\hat u|^2).
  \end{align*}
  By Lemma~\ref{lem:est-apriori-laplacian} below,
  \[ \sum_a \langle \Phi' (| \hat{u} |^2)  \hat{u}^a, \rho^3  (- \Delta)^s u^a
     \rangle \geqslant - \frac{1}{2} \| (- \Delta)^s \rho^2 \|
     \int_{\Lambda_{\varepsilon}} \Phi (| \hat{u} |^2) -
     \overline{B}  \| \hat{u} \| 
     \int_{\Lambda_{\varepsilon}} \Phi' (| \hat{u} |^2) . \]
  We also have
  \[ \begin{aligned}
       \langle \Phi' (| \hat{u} |^2), m^2 \rho^2  | \hat{u} |^2 + \lambda (|
       \hat{u} |^2)^2 \rangle & \geqslant  \langle \Phi' (| \hat{u} |^2),
       \lambda (| \hat{u} |^2)^2 \rangle\\
       & \geqslant  \lambda N \langle (| \hat{u} |^2 - L)^{N - 1}_+, ((|
       \hat{u} |^2 - L)_+ + L) L \rangle\\
       & \geqslant  \lambda L \int_{\Lambda_{\varepsilon}} \Phi  (| \hat{u}
       |^2) + \lambda L^2 \int_{\Lambda_{\varepsilon}} \Phi' (| \hat{u} |^2)
     \end{aligned} \]
  and
  \[ - \langle \Phi' (| \hat{u} |^2), \rho^3 f \cdot \hat{u} \rangle \geqslant
     - \| \rho^3 f \|  \| \hat{u} \| \int_{\Lambda_{\varepsilon}} \Phi' (|
     \hat{u} |^2) . \]
  We conclude that
  \[ 0 \geqslant [\lambda L - A] \int_{\Lambda_{\varepsilon}} \Phi (| \hat{u}
     |^2) + [\lambda L^2 - (\| \rho^3 f \| + B)  \| \hat{u} \|]
     \int_{\Lambda_{\varepsilon}} \Phi' (| \hat{u} |^2) . \]
  Taking
  \[ L > L_{\ast} \assign \max (\lambda^{- 1} A, \lambda^{- 1 / 2} [\| \rho^3
     f \| + B]^{1 / 2}  \| \hat{u} \|^{1 / 2}), \]
  we deduce that
  \[ \int_{\Lambda_{\varepsilon}} \Phi ' (| \hat{u} |^2) =
     \int_{\Lambda_{\varepsilon}} \Phi (| \hat{u} |^2) = 0, \]
  which implies that $| \hat{u} |^2 \leqslant L$ a.e. on
  $\Lambda_{\varepsilon}$. Thus,
  \[ \| \hat{u} \|^2 \leqslant \inf_{L > L_{\ast}} L = L_{\ast} \leqslant
     \lambda^{- 1} A + \lambda^{- 1 / 2} [\| \rho^3 f \| + B]^{1 / 2} \|
     \hat{u} \|^{1 / 2} . \]
  This implies that
  \[ \| \hat{u} \| \leqslant \lambda^{- 1 / 2} A^{1 / 2} + \lambda^{- 1 / 4} 
     \| \hat{u} \|^{1 / 4}  [\| \rho^3 f \| + B]^{1 / 4} . \]
  By Young's inequality, we have
  \[ \lambda^{- 1 / 4}  \| \hat{u} \|^{1 / 4}  [\| \rho^3 f \| + B]^{1 / 4}
     \leqslant \frac{\| \hat{u} \|}{4} + \frac{3}{4} \lambda^{- 1 / 3}  [\|
     \rho^3 f \| + B]^{1 / 3}, \]
  and hence
  \[ \| \hat{u} \| \leqslant 2 \lambda^{- 1 / 2} A^{1 / 2} + \lambda^{- 1 / 3}
     [\| \rho^3 f \| + B]^{1 / 3}, \]
  as claimed.
\end{proof}

The following lemma completes the proof.

\begin{lemma}
  \label{lem:est-apriori-laplacian}Suppose that $u \in C_t^1 C_x^2
  (\Lambda_{\varepsilon}, \mathbb{R}^n)$ is such that $| \zeta^{\bar{\nu}} u
  |$ is a bounded function for some $\bar{\nu} < 2 s$ and $\rho \in C_t^1
  C_x^2 (\Lambda_{\varepsilon})$ is a positive weight such that for some
  $\tilde{\nu} \in (\bar{\nu}, 2 s)$ the functions $\zeta^{- \tilde{\nu}}
  \rho$, $\zeta^{\tilde{\nu}} \rho^{- 1}$ are bounded. Let $\Phi = \Phi_{L, N}
  \in C (\mathbb{R})$ with $N > 1$ fixed as in the proof of the previous theorem
  and arbitrary $L > 0$. Then we have
  \[ \sum_a \langle \Phi' (| \hat{u} |^2)  \hat{u}^a, \rho^3  (- \Delta)^s u^a
     \rangle \geqslant - \frac{1}{2} \| (- \Delta)^s \rho^2 \|
     \int_{\Lambda_{\varepsilon}} \Phi (| \hat{u} |^2) -
     \bar{B}  \| \hat{u} \|  \int_{\Lambda_{\varepsilon}}
     \Phi' (| \hat{u} |^2) , \]
  where $\hat{u} \assign \rho u$ and
  \[ \bar{B} \assign \| \hat{u} \|  (\| \rho  (- \Delta)^s
     \rho \| + \| \rho^2 \mathfrak{D}_s (\rho^{- 1}) \mathfrak{D}_s (\rho) \|)
     + \| \rho  \mathfrak{D}_s (\rho) \mathfrak{D}_s (\hat{u}) \| . \]
\end{lemma}

\begin{proof}
  The following decomposition separates the favourable convexity contribution
  from the terms in which the fractional operator interacts with the weight.
  Leaving the sum over $a$ implicit, we have
  \[ \begin{aligned}
       &  \langle \Phi' (| \hat{u} |^2)  \hat{u}^a,
       \rho^3 (- \Delta)^s u^a \rangle \\
       & =  \int_{\Lambda_{\varepsilon}\times\Lambda_{\varepsilon}} \mu_s (\mathd z \mathd z') \Phi' (|
       \hat{u} (z) |^2)  \hat{u}^a (z) \rho^3 (z)  (u^a (z) - u^a (z')) \\
       & =  \int_{\Lambda_{\varepsilon}\times\Lambda_{\varepsilon}} \mu_s (\mathd z \mathd z') \Phi' (|
       \hat{u} (z) |^2)  \hat{u}^a (z) \rho^2 (z)  (\hat{u}^a (z) - \hat{u}^a
       (z')) & (\backassign (\mathbb{I}))\\
       & \quad  + \int_{\Lambda_{\varepsilon}\times\Lambda_{\varepsilon}} \mu_s (\mathd z \mathd z') \Phi'
       (| \hat{u} (z) |^2)  \hat{u}^a (z) \rho^2 (z)  (\rho (z') - \rho (z))
       u^a (z') & (\backassign (\mathbb{I}\mathbb{I}))
     \end{aligned} \]
  Let $V_{\tau} \assign \tau \hat{u} (z') + (1 - \tau)  \hat{u} (z)$. Since
  $\tau \in [0, 1] \mapsto \Phi (| V_{\tau} |^2)$ is a convex function (as a~composition of a convex function with an affine one),
  \[ \Phi (| V_1 |^2) - \Phi (| V_0 |^2) \geqslant \partial_{\tau}
     \Phi (| V_{\tau} |^2) |_{\tau = 0}, \]
  that is,
  \[ \Phi (| \hat{u} (z') |^2) - \Phi (| \hat{u} (z) |^2) \geqslant 2 \Phi' (|
     \hat{u} (z) |^2)  \sum_a \hat{u}^a (z)  (\hat{u}^a (z') - \hat{u}^a (z))
     . \]
  Using the above inequality we obtain
  \[ \begin{aligned}
       \mathbb{I} & \geqslant  \frac{1}{2} \int_{\Lambda_{\varepsilon} \times
       \Lambda_{\varepsilon}} \mu_s (\mathd z \mathd z') \rho^2 (z)  [\Phi (|
       \hat{u} (z) |^2) - \Phi (| \hat{u} (z') |^2)]\\
       & =  \frac{1}{2} \int_{\Lambda_{\varepsilon} \times
       \Lambda_{\varepsilon}} \mu_s (\mathd z \mathd z') \Phi (| \hat{u} (z)
       |^2)  [\rho^2 (z) - \rho^2 (z')]\\
       & \geqslant  - \frac{1}{2} \| (- \Delta)^s \rho^2 \|
       \int_{\Lambda_{\varepsilon}} \Phi (| \hat{u} |^2),
     \end{aligned} \]
  where the equality in the second line follows by an integration by parts of
  the fractional Laplacian, which is a symmetric operator in $L^2$.
  Let us now consider $(\mathbb{I}\mathbb{I})$. We first split the factor
  \(\rho^{-1}(z')\) as
  \[
  \rho^{-1}(z')
  =
  \bigl(\rho^{-1}(z')-\rho^{-1}(z)\bigr)
  +
  \rho^{-1}(z).
  \]
  This gives
  \begin{align*}
  	(\mathbb{I}\mathbb{I})
  	&=
  	\int_{\Lambda_{\varepsilon}\times\Lambda_{\varepsilon}}
  	\mu_s(\mathd z\,\mathd z')\,
  	\Phi'(|\hat u(z)|^2)
  	\hat u^a(z)\rho^2(z)
  	\bigl(\rho(z')-\rho(z)\bigr)\rho^{-1}(z')\hat u^a(z')
  	\\
  	&=
  	\int_{\Lambda_{\varepsilon}\times\Lambda_{\varepsilon}}
  	\mu_s(\mathd z\,\mathd z')\,
  	\Phi'(|\hat u(z)|^2)
  	\hat u^a(z)\rho^2(z)
  	\bigl(\rho(z')-\rho(z)\bigr)
  	\bigl(\rho^{-1}(z')-\rho^{-1}(z)\bigr)\hat u^a(z')
  	\\
  	&\quad
  	+
  	\int_{\Lambda_{\varepsilon}\times\Lambda_{\varepsilon}}
  	\mu_s(\mathd z\,\mathd z')\,
  	\Phi'(|\hat u(z)|^2)
  	\hat u^a(z)\rho^2(z)
  	\bigl(\rho(z')-\rho(z)\bigr)\rho^{-1}(z)\hat u^a(z').
  \end{align*}
  We call the first term \((\mathbb{I}\mathbb{I}_2)\). In the second term we
  split
  \[
  \hat u^a(z')
  =
  \hat u^a(z)
  +
  \bigl(\hat u^a(z')-\hat u^a(z)\bigr).
  \]
  Hence
  \begin{align*}
  	(\mathbb{I}\mathbb{I})
  	&=
  	\int_{\Lambda_{\varepsilon}\times\Lambda_{\varepsilon}}
  	\mu_s(\mathd z\,\mathd z')\,
  	\Phi'(|\hat u(z)|^2)
  	|\hat u(z)|^2\rho(z)
  	\bigl(\rho(z')-\rho(z)\bigr)
  	\\
  	&
  	+
  	\int_{\Lambda_{\varepsilon}\times\Lambda_{\varepsilon}}
  	\mu_s(\mathd z\,\mathd z')\,
  	\Phi'(|\hat u(z)|^2)
  	\hat u^a(z)\rho^2(z)
  	\bigl(\rho(z')-\rho(z)\bigr)
  	\bigl(\rho^{-1}(z')-\rho^{-1}(z)\bigr)\hat u^a(z')
  	\\
  	&
  	+
  	\int_{\Lambda_{\varepsilon}\times\Lambda_{\varepsilon}}
  	\mu_s(\mathd z\,\mathd z')\,
  	\Phi'(|\hat u(z)|^2)
  	\hat u^a(z)\rho(z)
  	\bigl(\rho(z')-\rho(z)\bigr)
  	\bigl(\hat u^a(z')-\hat u^a(z)\bigr)
  	\\
  	&=
  	(\mathbb{I}\mathbb{I}_1)
  	+
  	(\mathbb{I}\mathbb{I}_2)
  	+
  	(\mathbb{I}\mathbb{I}_3).
  \end{align*}
  For $(\mathbb{I}\mathbb{I}_1)$ we have
  \[ \begin{aligned}
       (\mathbb{I}\mathbb{I}_1) & \geqslant  - \| \hat{u} \|^2  \| \rho (-
       \Delta)^s \rho \|  \int_{\Lambda_{\varepsilon}} \Phi' (| \hat{u} (z)
       |^2) .
     \end{aligned} \]
  Next,
  \[
  \begin{aligned}
  	(\mathbb{I}\mathbb{I}_2)
  	&\geqslant
  	-\|\hat u\|^2
  	\int_{\Lambda_\varepsilon}
  	\Phi'(|\hat u(z)|^2)\rho^2(z)
  	\\
  	&\qquad\qquad \times
  	\left(
  	\int_{\Lambda_\varepsilon}
  	|\rho(z')-\rho(z)|
  	|\rho^{-1}(z')-\rho^{-1}(z)|
  	\,\mu_s(z,\mathd z')
  	\right)\mathd z
  	\\
  	&\geqslant
  	-\|\hat u\|^2
  	\|\rho^2\mathfrak D_s(\rho^{-1})\mathfrak D_s(\rho)\|
  	\int_{\Lambda_\varepsilon}
  	\Phi'(|\hat u|^2).
  \end{aligned}
  \]
  Finally, using the Cauchy--Schwarz inequality, we obtain
  \[ \begin{aligned}
       (\mathbb{I}\mathbb{I}_3) & \geqslant  - \| \hat{u} \|  \| \rho
       \mathfrak{D}_s (\rho) \mathfrak{D}_s (\hat{u}) \|
       \int_{\Lambda_{\varepsilon}} \Phi' (| \hat{u} |^2) .
     \end{aligned} \]
  Combining the lower bounds for \((\mathbb I)\) and
     \((\mathbb{I}\mathbb{I}_j)\), \(j=1,2,3\), gives the desired estimate.
\end{proof}

\section{Analysis of the flow equation}\label{sec:flow}

In this section, we establish the estimates used in the proof of
Theorem~\ref{thm:effective-force}. More precisely, we analyse the existence of
an approximate solution to the flow equation~\eqref{eq:flow-eq-sec2} for the
effective force $(F_{\sigma})_{\sigma}$, with well-controlled bounds encoded
by~\eqref{eq:est-F}. These bounds are achievable only because we can ``tune''
the boundary condition~\eqref{eq:final-cond-f} using the $\varepsilon$- and
$M$-dependent renormalisation term $r_{\varepsilon, M}$.

Conceptually, we are dealing with a random bilinear equation whose solution is
analysed via the evolution equation for its cumulants. The flow equation for
cumulants has a similar structure to the flow equation for the effective force
and propagates comparable bounds backwards from the final condition at $\sigma
= 1$, except in a low-dimensional (so-called relevant) subspace, where the
bounds must be propagated forwards from small to large $\sigma$. This
procedure requires tuning an appropriate final condition so that the solution
lies on a trajectory with controlled bounds. To simplify this tuning, we
decompose the flow equation to reduce the relevant subspace to one dimension.

Once bounds for the cumulants are established, a Kolmogorov-type argument
allows us to deduce pathwise bounds on the effective force. The section
concludes with a technical ``post-processing'' step, which extracts the
coercive term essential for the global a priori estimates and verifies the
conditions~\eqref{eq:est-F}.

\subsection{Random flow equation}\label{sec:random-flow}

To study approximate solutions $(F_{\sigma})_{\sigma \in [1 / 2, 1]}$ of the
flow equation~\eqref{eq:flow-eq-sec2} we need to set up the appropriate
spaces. Recall that
\[ \mathcal{E}  \assign \bigcap_{\alpha > 0} C (\Lambda_M, \zeta^{\alpha}),
   \qquad \hat{\mathcal{E}} \assign \mathcal{E}+\partial_t \mathcal{E}. \]
As we shall argue, the flow equation can be approximately solved in the
space $\mathcal{P} (\mathcal{E})$ of polynomial functions on $\mathcal{E}$
with values in $\hat{\mathcal{E}}$. We say that a functional $F$ belongs to
$\mathcal{P} (\mathcal{E})$ if for some $\bar{k} \in \mathbb{N}_+$ there exist
kernels $(F^{(k)})_{k \in \{ 0, \ldots, \bar{k} \}}$ of operators
$\mathcal{E}^{\otimes k} \rightarrow \hat{\mathcal{E}}$ such that
\[ F (\phi) = \sum_{k = 0}^{\bar{k}} F^{(k)} (\phi) = \sum_{k = 0}^{\bar{k}}
   \int_{\Lambda^k} F^{(k)} (\Cdot ; \mathd z_1, \ldots, \mathd z_k) \phi
   (z_1) \cdots \phi (z_k) \in \hat{\mathcal{E}}, \]
for all $\phi \in \mathcal{E}$. In order to construct a suitable approximate
solution of~\eqref{eq:flow-eq-sec2} we introduce a formal parameter $\hbar$
and make the ansatz
\[ {}^{\hbar} F_{\sigma} = \sum_{\ell \geqslant 0} \hbar^{\ell}
   F_{\sigma}^{[\ell]} . \]
Moreover, we assume that the final condition is of the form
\[ {}^{\hbar} F (\phi) = - \lambda \phi^3  + {}^{\hbar}
   r_{\varepsilon, M} \phi + \xi^{(\varepsilon, M)}, \qquad {}^{\hbar}
   r_{\varepsilon, M} = \bar{r} + \sum_{\ell \geqslant 1} \hbar^{\ell}
   {r^{[\ell]}_{\varepsilon, M}}  . \]
We are led to look for solutions of the \textit{perturbative flow equation}
\begin{equation}
  \partial_{\sigma} {}^{\hbar} F_{\sigma} + \hbar \mathD {}^{\hbar} F_{\sigma}
  (\dot{G}_{\sigma} {}^{\hbar} F_{\sigma}) = 0, \qquad {}^{\hbar} F_1 =
  {}^{\hbar} F, \label{eq:pert-flow}
\end{equation}
in the space $\mathcal{P} (\mathcal{E}) \llbracket \hbar \rrbracket$ of
formal power series in $\hbar$ with coefficients in $\mathcal{P}
(\mathcal{E})$. This setup has the advantage that now the flow equation has a
unique global solution which can be determined by induction on the order in
$\hbar$. An approximate solution to~\eqref{eq:flow-eq-sec2} is obtained by
fixing an integer $\bar{\ell} \geqslant 0$ and letting
\[ F_{\sigma} \assign \sum_{\ell = 0}^{\bar{\ell}} F_{\sigma}^{[\ell]} . \]
The choice of the value of $\bar{\ell}$ will be discussed in
Section~\ref{sec:post-processing} below. We observe that, thanks
to~\eqref{eq:pert-flow}, this truncation implies the existence of a maximal
polynomial order $\bar{k} \in \mathbb{N}_+$ in the fields for the kernels. We
decompose the force as
\[ F_{\mu} (\phi) = \sum^{\bar{\ell}}_{\ell = 0} F_{\mu}^{[\ell]} (\phi) =
   \sum^{\bar{\ell}}_{\ell = 0} \sum_{k = 0}^{\bar{k}} F_{\mu}^{[\ell],(k)}
   (\phi), \]
where $\ell$ measures the perturbative order in $\hbar$ while $k$ is the
polynomial degree in $\phi$.

Let us now introduce a condensed notation to manipulate these kernels. Let
\[
  \mathfrak{A} \assign \{ (\ell, k) | 0 \leqslant \ell \leqslant \bar{\ell},
  0 \leqslant k \leqslant \bar{k} \}, \qquad
  \mathfrak{A}_{\mathrm{ext}} \assign \mathfrak{A} \cup
  \{ (\ell, k) | \bar{\ell} + 1 \leqslant \ell \leqslant 2 \bar{\ell} + 1,
  0 \leqslant k \leqslant 2 \bar{k} \} .
\]
For $\mathfrak{a} \in \mathfrak{A}_{\mathrm{ext}}$ with $\mathfrak{a}=
(\ell, k)$ we let $k (\mathfrak{a}) \assign k$, $\ell (\mathfrak{a}) \assign
\ell$ and write
\begin{equation}
  [\mathfrak{a}] \assign - \alpha + \delta \ell (\mathfrak{a}) + \beta k
  (\mathfrak{a}), \label{eq:homogeneity-of-kernels}
\end{equation}
for suitable positive parameters $\alpha, \delta$ and $\beta$ whose values will
be fixed later. The set $\mathfrak{A}$ indexes the kernels retained in the
truncated effective force, whereas the additional indices in
$\mathfrak{A}_{\mathrm{ext}} \setminus \mathfrak{A}$ will be used only to
label the terms in its flow-equation error. For $\mathfrak{a} \in
\mathfrak{A}$, we say that a kernel

\[ F^{\mathfrak{a}} (z ; \mathd z_1, \cdots, \mathd z_k)
   \assign F^{[\ell], (k)} (z ; \mathd z_1, \cdots, \mathd z_k) \]
is \textit{relevant} if $[\mathfrak{a}] < 0$, \textit{marginal} if
$[\mathfrak{a}] = 0$ and \textit{irrelevant} if $[\mathfrak{a}] > 0$. We
refer to $z$ as the \textit{output} variable, and to $z_1, \ldots, z_k$ as
the \textit{input} variables of the kernel $F^{\mathfrak{a}} (z ; \mathd
z_1, \ldots, \mathd z_k)$. Throughout, the kernels \(F^{[\ell],(k)}(z;\mathd z_1,\ldots,\mathd z_k)\) are taken
to be symmetric in their \(k\) input variables.
To simplify the notation, we usually ignore the
fact that the kernel $F^{\mathfrak{a}} (z ; \mathd z_1, \ldots, \mathd z_k)$
is generally not well-defined pointwise in the output variable $z$.

\subsection{Norms for kernels}

To introduce suitable norms for the effective force kernels
\[ (F_{\sigma}^{\mathfrak{a}})_{\mathfrak{a} \in \mathfrak{A}, \sigma \in [1 /
   2, 1]}, \]
we first need some notation and preliminary definitions for weights and
smoothing operators.

\begin{definition}
  We denote by $\operatorname{St} (z_1, \ldots, z_n)$ the Steiner diameter of the set
  $\{ z_1, \ldots, z_n \}$ with respect to the parabolic distance introduced
  in Definition~\ref{def:distance}, \textit{i.e.} the minimum of the
  lengths of trees with nodes at the points $\{ z_1, \ldots, z_n \}$ and
  possibly other points $\{ z_1', \ldots, z_m' \}$.
\end{definition}

\begin{remark}
  \label{rem:steiner}Our fractional parabolic
  distance~\eqref{eq:dist-parabolic} satisfies the triangle inequality
  \[ | z_1 + z_2 |_s \leqslant | z_1 |_s + | z_2 |_s, \qquad z_1, z_2 \in
     \Lambda_0 . \]
  Moreover, for $n \in \mathbb{N}_+, m \in \{ 1, \ldots n - 1 \}$ and $z, z_1,
  \ldots, z_n \in \Lambda_0$, the Steiner diameter satisfies
  \[ \operatorname{St} (z_1, \ldots, z_n) \leqslant \operatorname{St} (z_1, \ldots, z_m, z) +
     \operatorname{St} (z, z_{m + 1}, \ldots, z_n), \]
  which, noting that $\operatorname{St} (z_1, z_2) = | z_1 - z_2 |_s$, can be seen as
  a generalisation of the triangle inequality for the fractional parabolic
  distance.
\end{remark}

\begin{definition}
  \label{def:weights-kernels}Let $\flat \in (1, 2 s)$ be a constant close to
  $2 s$ and $\kappa_{\mathfrak{o}} \in (0, (2-\flat) / (2 \bar{\ell} + 2))$
  be a small constant, to be fixed later.
  \begin{enumerate}[(a)]
    \item For $m \in \mathbb{N}_+$ and $\omega \in \mathbb{R}$, the weight
    $w^{(1 + m), \omega}_{\mu} \in C  (\Lambda_0^{(1 + m)})$ is defined by
    \[ w^{(1 + m), \omega}_{\mu} (z, z_1, \ldots, z_m) \assign (1 + \llbracket
       \mu \rrbracket^{- 1} \operatorname{St} (z, z_1, \ldots, z_m) )^{\omega}, \qquad
       z, z_1, \ldots, z_m \in \Lambda_0 . \]
    We write $w^{(1), \omega}_{\mu} = 1 \in C  (\Lambda_0 )$.
    
    \item For $m \in \mathbb{N}_+$, the weight $v^{(1 + m)}_{\mu} {\in
    C^{\infty}}  (\Lambda_0^{(1 + m)})$ is defined by
    \[ v^{(m + 1)}_{\mu} (z, z_1, \ldots, z_m) \assign v^{(m + 1)}
       (\llbracket \mu \rrbracket^{- 1} .z, \llbracket \mu \rrbracket^{- 1}
       .z_1, \ldots, \llbracket \mu \rrbracket^{- 1} .z_m), \]
    where $v^{(m + 1)}_{} {\in C^{\infty}}  (\Lambda_0^{(1 + m)}, [0, 1])$ is
    a fixed function such that
    \[ v^{(m + 1)} (z, z_1, \ldots, z_m) = \left\{ \begin{aligned}
         1 & ~\operatorname{if}~  \operatorname{St} (z, z_1, \ldots, z_m) \leqslant 1,\\
         0 & ~\operatorname{if}~  \operatorname{St} (z, z_1, \ldots, z_m) \geqslant 2,
       \end{aligned} \right. \]
       and
       \[
       v^{(m + 1)} (z, z_1, \ldots, z_m)=
       v^{(m + 1)} (0, z_1-z, \ldots, z_m-z).
       \]
    We write $v^{(1)}_{\mu} = 1 \in C  (\Lambda_0 )$.
    
    \item The weight $\mathfrak{o} \in C  (\Lambda_0 )$ is defined by
    \[ \mathfrak{o} (z) \assign \zeta (z)^{\kappa_{\mathfrak{o}}} = \langle z
       \rangle_s^{- \kappa_{\mathfrak{o}}}, \qquad z \in \Lambda_0 . \]
    \item For $\mathfrak{a} \in \mathfrak{A}_{\mathrm{ext}}$, we write
    \begin{align*}
       w^{\mathfrak{a}}_{\mu} &\assign
       w^{(1 + k (\mathfrak{a})), \flat}_{\mu},
       & \tilde{w}^{\mathfrak{a}}_{\mu} &\assign
       w^{(1 + k (\mathfrak{a})), \flat - \ell (\mathfrak{a})
       \kappa_{\mathfrak{o}}}_{\mu},\\
       v^{\mathfrak{a}}_{\mu} &\assign v^{(1 + k (\mathfrak{a}))}_{\mu},
       & \mathfrak{o}^{\mathfrak{a}} &{\assign
       \mathfrak{o}^{1 + \ell (\mathfrak{a})}}  .
    \end{align*}
    \item The weight $h_{\mu} \in C  (\Lambda^2_0)$ is
    defined by
    \[ h_{\mu} (z, z') \assign (1 + \llbracket \mu \rrbracket^{- 2} | z - z'
       |_s^2)^{- 1} . \]
  \end{enumerate}
\end{definition}

Recall that
\[ \langle z \rangle_s \assign (1 + | z_0 |^{1 / s} + | \bar{z} |^2)^{1 / 2} .
\]
The weights $w^{\mathfrak{a}}_{\mu}$ and ${\tilde{w}^{\mathfrak{a}}_{\mu}}$
will control the approximate localisation of the effective force kernels
$F^{\mathfrak{a}}_{\mu} (z ; \mathd z_1, \ldots, \mathd z_k)$ near the
diagonal $\{ z = z_1 = \cdots = z_k \}$. The weights
$\mathfrak{o}^{\mathfrak{a}}$ will instead be used to control the growth at
spacetime infinity of the kernels in the output variable $z$, reflecting the
corresponding growth inherited from the white noise. These three families of
weights are constructed to be compatible with the estimates for the bilinear
terms appearing in the flow equation, as detailed in the lemma below. The
additional weights ${v_{\mu}^{(m + 1)}}$ and $h_{\mu}$ will be used only in
intermediate estimates and play an auxiliary role. Some of their useful
properties are summarised in Appendix~\ref{sec:weight-properties}. Finally,
observe that
\[ w_{\mu}^{(2), \flat} (z, z') = w^{\flat}_{\mu} (z - z') = (1 + \llbracket
   \mu \rrbracket^{- 1} | z - z' |_s)^{\flat}, \]
where $w^{\omega}_{\mu}$ is the weight introduced in Definition~\ref{def:weight-w}.
The exponent $\flat - \ell (\mathfrak{a}) \kappa_{\mathfrak{o}}$ in the
definition of the weight $\tilde{w}^{\mathfrak{a}}_{\mu}$ is designed in such
a way that the bound~\eqref{eq:weight-w-o} below holds true.

\begin{lemma}
  \label{lem:weight-factorization}For all $\mathfrak{a} \in
  \mathfrak{A}_{\mathrm{ext}}$ and $\mathfrak{b}, \mathfrak{c} \in
  \mathfrak{A}$ such that $k (\mathfrak{b}) + k
  (\mathfrak{c}) = k (\mathfrak{a}) + 1$ and $\ell (\mathfrak{b}) + \ell
  (\mathfrak{c}) = \ell (\mathfrak{a}) - 1$, the bounds
  \begin{equation}
    w_{\mu}^{\mathfrak{a}} (z, z_{1,} \ldots, z_{k (\mathfrak{a})})
    {\lesssim w_{\mu}^{\mathfrak{b}}}  (z, z_{1,} \ldots, z_{k
    (\mathfrak{b}) - 1}, \hat{z}) w^{\flat}_{\mu} (\hat{z}  - \tilde{z})
    w_{\mu}^{\mathfrak{c}} (\tilde{z}, z_{k (\mathfrak{b})}, \ldots, z_{k
    (\mathfrak{a})}), \label{eq:weight-w}
  \end{equation}
  \begin{equation}
    \begin{aligned}
      \mathfrak{o}^{\mathfrak{a}} (z)
      \tilde{w}_{\mu}^{\mathfrak{a}}
      (z, z_{1,} \ldots, z_{k (\mathfrak{a})})
      & \lesssim
      \mathfrak{o}^{\mathfrak{b}} (z)
      {\tilde{w}_{\mu}^{\mathfrak{b}}}
      (z, z_{1,} \ldots, z_{k (\mathfrak{b}) - 1}, \hat{z})
      w^{\flat}_{\mu} (\hat{z}  - \tilde{z})\\
      &\qquad \times
      \mathfrak{o}^{\mathfrak{c}} (\tilde{z})
      \tilde{w}_{\mu}^{\mathfrak{c}}
      (\tilde{z}, z_{k (\mathfrak{b})}, \ldots, z_{k (\mathfrak{a})}),
    \end{aligned}
    \label{eq:weight-w-o}
  \end{equation}
  hold uniformly in $\mu \in (0, 1)$ and $z, \hat{z}, \tilde{z}, z_{1,}
  \ldots, z_{k (\mathfrak{a})} \in \Lambda_0$. For all $\mathfrak{a} \in
  \mathfrak{A}_{\mathrm{ext}}$, the bounds
  \begin{equation}
    \mathfrak{o}^{\mathfrak{a}} (z) \lesssim (\mathfrak{o}^{\mathfrak{a}} (z -
    z'))^{- 1} \mathfrak{o}^{\mathfrak{a}} (z') \lesssim w^{2 - \flat}_{\mu}
    (z - z') \mathfrak{o}^{\mathfrak{a}} (z'), \label{eq:weight-o}
  \end{equation}
  \begin{equation}
  	\begin{aligned}
    w^{\mathfrak{a}}_{\mu} (z, z_1, \ldots, z_{k (\mathfrak{a})}) &\lesssim
    w^{\flat}_{\mu} (z - z') w^{\flat}_{\mu} (z_1 - z_1') \ldots
    w^{\flat}_{\mu} (z_{k (\mathfrak{a})} - z_{k (\mathfrak{a})}')
    \\
    &\qquad\times
    w^{\mathfrak{a}}_{\mu} (z', z'_1, \ldots, z'_{k (\mathfrak{a})}),
   \end{aligned}
    \label{eq:weight-w-2}
  \end{equation}
  \begin{equation}
    \begin{aligned}
      \mathfrak{o}^{\mathfrak{a}} (z)
      \tilde{w}^{\mathfrak{a}}_{\mu} (z, z_1, \ldots, z_{k (\mathfrak{a})})
      &\lesssim w^2_{\mu} (z - z') w^{\flat}_{\mu}
      (z_1 - z_1') \ldots w^{\flat}_{\mu} (z_{k (\mathfrak{a})}
      - z_{k (\mathfrak{a})}') \\
      &\qquad \times
      \mathfrak{o}^{\mathfrak{a}} (z')
      \tilde{w}^{\mathfrak{a}}_{\mu} (z', z'_1, \ldots, z'_{k (\mathfrak{a})}),
    \end{aligned}
    \label{eq:weight-w-o-2}
  \end{equation}
  hold uniformly in $\mu \in (0, 1)$ and $z, z_{1,} \ldots, z_{k
  (\mathfrak{a})}, z', z'_{1,} \ldots, z'_{k (\mathfrak{a})} \in \Lambda_0$.
\end{lemma}

\begin{proof}
  To prove the claim we apply iteratively the inequalities from
  Remark~\ref{rem:steiner} and use the bounds
  \[ 1 + \sum_{i \in \{ 1, \ldots, n \}} a_i \leqslant \prod_{i \in \{ 1,
     \ldots, n \}} (1 + a_i), \qquad (1 + a_i)^{\omega} \leqslant (1 +
     a_i)^{\bar{\omega}}, \]
  valid for all $a_1, \ldots, a_n \geqslant 0$ and $0 \leqslant \omega
  \leqslant \bar{\omega}$. To prove~\eqref{eq:weight-o}
  and~\eqref{eq:weight-w-2}, we use furthermore that
  \[ (1 + \ell (\mathfrak{a})) \kappa_{\mathfrak{o}} \leqslant (2
     \bar{\ell} + 2) \kappa_{\mathfrak{o}} \leqslant  2 - \flat, \]
  which follows from the assumption $\kappa_{\mathfrak{o}} \in
  (0, (2-\flat) / (2 \bar{\ell} + 2))$.
\end{proof}

\begin{definition}
  \label{def:regularizing-kernels-for-kernels}For $n, m \in \mathbb{N}_0$ and
  $\sigma, \eta \in (1 / 2, 1)$ we set
  \[ K^{n, m}_{\sigma} \assign 1^{\otimes n} \otimes (K_{\sigma})^{\otimes m},
     \qquad L_{\sigma}^{n, m} \assign 1^{\otimes n} \otimes
     (L_{\sigma})^{\otimes m}, \qquad K^{n, m}_{\eta, \sigma} \assign
     L_{\sigma}^{n, m} K^{n, m}_{\eta}, \]
  \[ \tilde{K}^{n, m}_{\sigma} \assign K^{\otimes n}_{\sigma} \otimes
     (K^2_{\sigma})^{\otimes m}, \qquad \tilde{L}_{\sigma}^{n, m} \assign
     L_{\sigma}^{\otimes n} \otimes (L^2_{\sigma})^{\otimes m}, \qquad
     \tilde{K}^{n, m}_{\eta, \sigma} \assign \tilde{L}_{\sigma}^{n, m} 
     \tilde{K}^{n, m}_{\eta}, \]
  where $L_{\sigma}, K_{\sigma}, K_{\eta, \sigma}$ are introduced in
  Definition~\ref{def:regularizing-kernels}. Given $\mathfrak{a} \in
  \mathfrak{A}_{\mathrm{ext}}$ and $\mu\in (1/2,1)$,
  we write
  \[ \tilde{K}^{\mathfrak{a}}_{\mu} \assign \tilde{K}^{1, k
     (\mathfrak{a})}_{\mu} = K_{\mu}^{\otimes (1 + k (\mathfrak{a}))} K^{1, k
     (\mathfrak{a})}_{\mu} . \]
\end{definition}

\begin{definition}
  \label{def:norm-kernels}The norm of a kernel $F$ associated with an operator
  $\mathcal{E}^{\otimes k} \rightarrow \mathcal{E}$ is defined by
  \begin{equation}
    \| F \| \assign \sup_{z \in \Lambda} \int_{\Lambda^k} | F (z ; \mathd z_1,
    \cdots, \mathd z_k) | . \label{eq:basic-norm}
  \end{equation}
  Given weights $\mathfrak{o} \in C (\Lambda_0)$ and $w \in C (\Lambda_0^{1 +
  k})$, we write
  \[ (\mathfrak{o} \cdot  F \cdot w) (z ; \mathd z_1, \cdots,
     \mathd z_k) \assign \mathfrak{o} (z) F (z ; \mathd z_1, \cdots, \mathd
     z_k) w (z, z_1, \cdots, z_k) . \]
  Let $\mathfrak{a} \in \mathfrak{A}_{\mathrm{ext}}$. The norm of a kernel
  $F^{\mathfrak{a}}$
  of an operator $\mathcal{E}^{\otimes k (\mathfrak{a})} \rightarrow
  \hat{\mathcal{E}}$ at scale $\sigma \in (1 / 2, 1)$ is defined by
  \[ \| F^{\mathfrak{a}} \|_{\sigma} \assign \| \mathfrak{o}^{\mathfrak{a}}
     \cdot [\tilde{K}^{\mathfrak{a}}_{\sigma} F^{\mathfrak{a}}] \cdot
     \tilde{w}^{\mathfrak{a}}_{\sigma} \|, \]
  where $\tilde{K}^{\mathfrak{a}}_{\sigma} F^{\mathfrak{a}}$ denotes the
  convolution of $F^{\mathfrak{a}}$, viewed as an element of $\mathcal{S}'
  (\Lambda^{1 + k (\mathfrak{a})})$, with the kernel of
  $\tilde{K}^{\mathfrak{a}}_{\sigma}$. We further introduce a norm for the
  family of effective force kernels
  \[ F^{\mathfrak{A}} \assign (F^{\mathfrak{a}}_{\sigma})_{\mathfrak{a} \in
     \mathfrak{A}, \sigma \in (1 / 2, 1)}, \]
  defined by
  \begin{equation}
    \| F^{\mathfrak{A}} \| \assign \bigg(\sup_{\mathfrak{a} \in \mathfrak{A}}
    \sup_{\sigma \in (1 / 2, 1)}  \llbracket \sigma \rrbracket^{-
    [\mathfrak{a}]}  \| F^{\mathfrak{a}}_{\sigma} \|_{\sigma}\bigg) 
    \vee
    \bigg(\sup_{\sigma \in (1 / 2, 1)}  \llbracket \sigma \rrbracket^{d / 2 + s + 2
    \kappa}  \| F^{[0], (0)}_{\sigma} \|_{\sigma}\bigg), \label{eq:norm-family-F}
  \end{equation}
  where $[\mathfrak{a}]$ is defined by~\eqref{eq:homogeneity-of-kernels}.
\end{definition}

\begin{remark}
  Inspecting the parameter choices discussed in
  Sections~\ref{sec:bounds-on-parameters} and~\ref{sec:post-processing}, and using
  \eqref{eq:homogeneity-of-kernels}, shows that $[(0, 0)] < - d / 2 - s -
  2 \kappa$. Consequently, the second term on the right-hand side
  of~\eqref{eq:norm-family-F} provides a stronger bound on $F^{[0],
  (0)}_{\sigma} = \xi^{(\varepsilon, M)}$ than the first term. For technical
  reasons, this stronger control will be required in
  Section~\ref{sec:post-processing} to establish the estimates stated
  in~\eqref{eq:est-F}.
\end{remark}

\begin{remark}
  \label{rem:norm0-scale-monotonicity}By Young's inequality,
  \eqref{eq:weight-w-o-2} and Lemma~\ref{lem:bounds-K-J}, we have
  \begin{equation}
    \begin{aligned}
      \| F^{\mathfrak{a}} \|_{\mu}
      &= \| \mathfrak{o}^{\mathfrak{a}} \cdot
      [\tilde{K}^{1, k (\mathfrak{a})}_{\mu} F^{\mathfrak{a}}] \cdot
      \tilde{w}^{\mathfrak{a}}_{\mu} \|\\
      &= \| \mathfrak{o}^{\mathfrak{a}}
      \cdot [\tilde{K}^{1, k (\mathfrak{a})}_{\mu, \sigma}
      \tilde{K}^{1, k(\mathfrak{a})}_{\sigma} F^{\mathfrak{a}}] \cdot
      \tilde{w}^{\mathfrak{a}}_{\mu} \|\\
      &\leqslant \| K_{\mu, \sigma} \|_{\operatorname{TV} (w^2_{\mu})}
      \| K^2_{\mu, \sigma} \|_{\operatorname{TV}
      (w^{\flat}_{\mu})}^{k (\mathfrak{a})}\,
      \| \mathfrak{o}^{\mathfrak{a}} \cdot [\tilde{K}^{1, k
      (\mathfrak{a})}_{\sigma} F^{\mathfrak{a}}] \cdot
      \tilde{w}^{\mathfrak{a}}_{\mu} \|\\
      &\lesssim \| \mathfrak{o}^{\mathfrak{a}} \cdot [\tilde{K}^{1, k
      (\mathfrak{a})}_{\sigma} F^{\mathfrak{a}}] \cdot
      \tilde{w}^{\mathfrak{a}}_{\sigma} \|\\
      &= \| F^{\mathfrak{a}} \|_{\sigma},
    \end{aligned} \label{eq:norm-mon}
  \end{equation}
  uniformly in $\mu \in [1 / 2, 1)$ and $\sigma \in [\mu, 1)$.
\end{remark}

Given kernels $F '$ and $F''$ of operators
$\mathcal{E}^{\otimes k'} \rightarrow \mathcal{E}$ and $\mathcal{E}^{\otimes
k''} \rightarrow \mathcal{E}$ respectively, we denote by
\begin{equation}
  F =\mathfrak{C} (\dot{G}_{\sigma}) (F'  \otimes F'')
  \label{eq:frak-B-operator}
\end{equation}
the kernel of the operator $\mathcal{E}^{\otimes k} \rightarrow \mathcal{E}$,
with $k = k' + k'' - 1$, defined by
\begin{equation}
  \begin{aligned}
    F (z ; \mathd z_1, \cdots, \mathd z_k)
    &\assign \frac{1}{k!} \sum_{\pi \in \mathcal{P}_k}
    \int_{\Lambda^2}
    F'(z ; \mathd z_{\pi (1)}, \cdots,
    \mathd z_{\pi (k' - 1)}, \mathd z')\\
    &\quad \times \dot{G}_{\sigma} (z' - z'')
    F''(z'' ; \mathd z_{\pi (k')}, \cdots, \mathd z_{\pi (k)})
    \mathd z' \mathd z'',
  \end{aligned}
  , \label{eq:def-B-general}
\end{equation}
where $\mathcal{P}_k$ denotes the set of permutations of $\{ 1, \ldots, k \}$
and the integral is over the variables $z'$ and~$z''$. Note that
$\mathfrak{C}  (\dot{G}_{\sigma})$ applies $\dot{G}_{\sigma}$ to the
output variable of the kernel $F''$ and plugs the result into the last input
variable of the kernel $F'$. The above definition extends naturally to
kernels of operators $\mathcal{E}^{\otimes k'} \rightarrow \hat{\mathcal{E}}$
and $\mathcal{E}^{\otimes k''} \rightarrow \hat{\mathcal{E}}$.

The flow equation~\eqref{eq:pert-flow} for the effective force can be
rewritten as the following flow equation for the kernels
\begin{equation}
  \partial_{\sigma} F_{\sigma}^{\mathfrak{a}} = -\sum_{\mathfrak{b},
  \mathfrak{c}} B^{\mathfrak{a}}_{\mathfrak{b}, \mathfrak{c}}
  (\dot{G}_{\sigma}, F_{\sigma}^{\mathfrak{b}}, F_{\sigma}^{\mathfrak{c}}),
  \qquad F_1^{\mathfrak{a}} = F^{\mathfrak{a}}, \label{eq:flow-kernels}
\end{equation}
where the operators $B^{\mathfrak{a}}_{\mathfrak{b}, \mathfrak{c}}$ are
implicitly defined by
\begin{equation}
  \sum_{\mathfrak{b}, \mathfrak{c}} B^{\mathfrak{a}}_{\mathfrak{b},
  \mathfrak{c}} (\dot{G}_{\sigma}, F^{\mathfrak{b}}_{\sigma},
  F^{\mathfrak{c}}_{\sigma}) \assign \sum_{\ell' = 0}^{\ell (\mathfrak{a} ) -
  1} \sum^{k (\mathfrak{a})}_{k' = 0} (k' + 1) \mathfrak{C} (\dot{G}_{\sigma})
  (F^{[\ell (\mathfrak{a}) - 1 - \ell'], (k' + 1)}_{\sigma} \otimes
  F_{\sigma}^{[\ell'], (k (\mathfrak{a}) - k')}) , \label{eq:B-form}
\end{equation}
and in particular $B^{\mathfrak{a}}_{\mathfrak{b}, \mathfrak{c}} = 0$ unless

\begin{equation}
  \left\{\begin{aligned}
    \ell (\mathfrak{a} ) = \, \ell (\mathfrak{b}) + \ell (\mathfrak{c} ) +
    1,\\
    k (\mathfrak{a}) = \, k (\mathfrak{b}) + k (\mathfrak{c}) - 1,\\
    - [\mathfrak{a}] + [\mathfrak{b}] + [\mathfrak{c}] = - \alpha + \beta -
    \delta .
  \end{aligned}\right. \label{eq:link-between-homogeneities-kernels}
\end{equation}

\begin{lemma}
  \label{lem:estimate-B-kernels}For all $\mathfrak{a} \in
  \mathfrak{A}_{\mathrm{ext}}$ and $\mathfrak{b}, \mathfrak{c} \in
  \mathfrak{A}$, the following bound
  \[ \llbracket \sigma \rrbracket^{- [\mathfrak{a}]}  \|
     B^{\mathfrak{a}}_{\mathfrak{b}, \mathfrak{c}} (\dot{G}_{\sigma},
     F^{\mathfrak{b}}_{\sigma}, F^{\mathfrak{c}}_{\sigma}) \|_{\sigma}
     \lesssim \llbracket \sigma \rrbracket^{- [\mathfrak{b}] - [\mathfrak{c}]
     - 1}  \| F^{\mathfrak{b}}_{\sigma} \|_{\sigma}  \|
     F^{\mathfrak{c}}_{\sigma} \|_{\sigma} \]
  holds uniformly in $\sigma \in [1 / 2, 1)$, provided
  \begin{equation}
    2 s - \alpha + \beta - \delta \geqslant 0.
    \label{eq:constraint-flow-kernels}
  \end{equation}
\end{lemma}

\begin{proof}
  By Definitions~\ref{def:norm-kernels}
  and~\ref{def:regularizing-kernels-for-kernels}, together with
  \eqref{eq:B-form}, \eqref{eq:def-B-general} and the identity $L_{\sigma}
  K_{\sigma} = 1$, we obtain
  \begin{align*}
    \| B^{\mathfrak{a}}_{\mathfrak{b}, \mathfrak{c}} (\dot{G}_{\sigma},
    F^{\mathfrak{b}}_{\sigma}, F^{\mathfrak{c}}_{\sigma}) \|_{\sigma} &=  \|
    \mathfrak{o}^{\mathfrak{a}} \cdot [\tilde{K}_{\sigma}^{\mathfrak{a}}
    B^{\mathfrak{a}}_{\mathfrak{b}, \mathfrak{c}} (\dot{G}_{\sigma},
    F^{\mathfrak{b}}_{\sigma}, F^{\mathfrak{c}}_{\sigma})] \cdot
    \tilde{w}_{\sigma}^{\mathfrak{a}} \|\\
 &=  \| \mathfrak{o}^{\mathfrak{a}} \cdot
    [B^{\mathfrak{a}}_{\mathfrak{b}, \mathfrak{c}} (L_{\sigma}^3
    \dot{G}_{\sigma}, \tilde{K}_{\sigma}^{\mathfrak{b}}
    F^{\mathfrak{b}}_{\sigma}, \tilde{K}_{\sigma}^{\mathfrak{c}}
    F^{\mathfrak{c}}_{\sigma})] \cdot \tilde{w}_{\sigma}^{\mathfrak{a}}
    \|\\
 &\leqslant  \| \mathfrak{o}^{\mathfrak{a}} \cdot
    [B^{\mathfrak{a}}_{\mathfrak{b}, \mathfrak{c}} (| L_{\sigma}^3
    \dot{G}_{\sigma} |, | \tilde{K}_{\sigma}^{\mathfrak{b}}
    F^{\mathfrak{b}}_{\sigma} |, | \tilde{K}_{\sigma}^{\mathfrak{c}}
    F^{\mathfrak{c}}_{\sigma} |)] \cdot \tilde{w}_{\sigma}^{\mathfrak{a}}
    \| .
  \end{align*}
  Hence, by \eqref{eq:weight-w-o}, \eqref{eq:B-form},
  \eqref{eq:def-B-general}, Definition~\ref{def:norm-kernels} and
  Lemma~\ref{lem:bounds-K-J}, we arrive at
  \begin{align*}
    \| B^{\mathfrak{a}}_{\mathfrak{b}, \mathfrak{c}} (\dot{G}_{\sigma},
    F^{\mathfrak{b}}_{\sigma}, F^{\mathfrak{c}}_{\sigma}) \|_{\sigma} &\lesssim  \| [B^{\mathfrak{a}}_{\mathfrak{b}, \mathfrak{c}} (|
    L_{\sigma}^3 \dot{G}_{\sigma} | \cdot w_{\sigma}^{\flat},
    \mathfrak{o}^{\mathfrak{b}} \cdot | \tilde{K}_{\sigma}^{\mathfrak{b}}
    F^{\mathfrak{b}}_{\sigma} | \cdot \tilde{w}_{\sigma}^{\mathfrak{b}},
    \mathfrak{o}^{\mathfrak{c}} \cdot | \tilde{K}_{\sigma}^{\mathfrak{c}}
    F^{\mathfrak{c}}_{\sigma} | \cdot \tilde{w}_{\sigma}^{\mathfrak{c}})]
    \|\\
 &\lesssim  \| L_{\sigma}^3 \dot{G}_{\sigma} \|_{\operatorname{TV}
    (w_{\sigma}^{\flat})}  \| \mathfrak{o}^{\mathfrak{b}} \cdot |
    \tilde{K}_{\sigma}^{\mathfrak{b}} F^{\mathfrak{b}}_{\sigma} | \cdot
    \tilde{w}_{\sigma}^{\mathfrak{b}} \|   \| \mathfrak{o}^{\mathfrak{c}}
    \cdot | \tilde{K}_{\sigma}^{\mathfrak{c}} F^{\mathfrak{c}}_{\sigma} |
    \cdot \tilde{w}_{\sigma}^{\mathfrak{c}} \| \\
 &\lesssim  \| L_{\sigma}^3 \dot{G}_{\sigma} \|_{\operatorname{TV}
    (w_{\sigma}^{\flat})}  \| F^{\mathfrak{b}}_{\sigma} \|_{\sigma}  \|
    F^{\mathfrak{c}}_{\sigma} \|_{\sigma}\\
 &\lesssim  \llbracket \sigma \rrbracket^{2 s - 1}  \|
    F^{\mathfrak{b}}_{\sigma} \|_{\sigma}  \| F^{\mathfrak{c}}_{\sigma}
    \|_{\sigma} .
  \end{align*}
  
  Moreover, combining~\eqref{eq:link-between-homogeneities-kernels}
  and~\eqref{eq:constraint-flow-kernels} yields
  \[ - [\mathfrak{a}] + 2 s - 1 = - [\mathfrak{b}] - [\mathfrak{c}] - \alpha +
     \beta - \delta + 2 s - 1 \geqslant - [\mathfrak{b}] - [\mathfrak{c}] - 1.
  \]
  This finishes the proof.
\end{proof}

\subsection{Norms for cumulants}

In order to prove probabilistic bounds for the moments of the norm $\|
F^{\mathfrak{A}} \|$, we start by establishing bounds for the cumulants of the
family $(F_{\sigma}^{\mathfrak{a}})_{\mathfrak{a} \in \mathfrak{A}, \sigma \in
[1 / 2, 1]}$ of the effective force kernels. To this end, let us first
introduce a useful notation. For
\[ \boldsymbol{a} \in \boldsymbol{A} \assign \{ (\mathfrak{a}_1, \ldots,
   \mathfrak{a}_n) | \mathfrak{a}_k \in \mathfrak{A}, L (\boldsymbol{a})
   \leqslant 2 \bar{\ell} \}, \]
we write
\[ n (\boldsymbol{a}) \assign n, \qquad L (\boldsymbol{a}) \assign \ell
   (\mathfrak{a}_1) + \cdots + \ell (\mathfrak{a}_n), \qquad K (\boldsymbol{a})
   \assign k (\mathfrak{a}_1) + \cdots + k (\mathfrak{a}_n) . \]
For random variables $X_1, \ldots, X_n$, we denote by
$\mathfrak{K}_n (X_1, \ldots, X_n)$ their $n$-th order joint cumulant. For
$\boldsymbol{a}= (\mathfrak{a}_1, \ldots, \mathfrak{a}_n) \in \boldsymbol{A}$ and
$\sigma \in [1 / 2, 1]$, we use the following shorthand notation
\[ \mathcal{F}_{\sigma}^{\boldsymbol{a}} (z_1, \mathd Z_1, \ldots, z_n, \mathd
   Z_n) \assign \mathfrak{K}_n (F_{\sigma}^{\mathfrak{a}_1} (z_1, \mathd Z_1),
   \cdots, F_{\sigma}^{\mathfrak{a}_n} (z_n, \mathd Z_n)), \]
where
\begin{equation}
  Z_i = (z_{i, 1}, \ldots, z_{i, k (\mathfrak{a}_i)}) {\in \Lambda^{k
  (\mathfrak{a}_i)}} , \qquad i \in \{ 1, \ldots n \}, \quad n = n
  (\boldsymbol{a}), \label{eq:notation-Z}
\end{equation}
are the input variables of the $i$-th effective force kernel involved in the
joint cumulant. Since $F_{\sigma}^{\mathfrak{a}_i}$ is a kernel of an operator
$\mathcal{E}^{\otimes k (\mathfrak{a}_i)} \rightarrow \hat{\mathcal{E}}$, it
is natural to view $\mathcal{F}_{\sigma}^{\boldsymbol{a}}$ as a kernel of an
operator
\begin{equation}
  \mathcal{E}^{\otimes K (\boldsymbol{a})} \rightarrow
  \hat{\mathcal{E}}^{\otimes n (\boldsymbol{a})} . \label{eq:operator-cumulant}
\end{equation}
We define the global homogeneity of the kernel
$\mathcal{F}_{\sigma}^{\boldsymbol{a}}$ as
\[ [\boldsymbol{a}] \assign - \varrho + n (\boldsymbol{a}) (\theta + \alpha) +
   [\mathfrak{a}_1] + \cdots + [\mathfrak{a}_n], \]
for suitable parameters $\varrho$ and $\theta$ whose values will be fixed in
Section~\ref{sec:bounds-on-parameters}. By~\eqref{eq:homogeneity-of-kernels}, we
have
\begin{equation}
  [\boldsymbol{a}] = - \varrho + \theta n (\boldsymbol{a}) + \delta L
  (\boldsymbol{a}) + \beta K (\boldsymbol{a}) . \label{eq:homogeneity-cum}
\end{equation}
We say that a cumulant $\mathcal{F}^{\boldsymbol{a}}$ is \textit{relevant} if
$[\boldsymbol{a}] < 0$, \textit{marginal} if $[\boldsymbol{a}] = 0$ and
\textit{irrelevant} if $[\boldsymbol{a}] > 0$.

\begin{definition}
  For $\boldsymbol{a} \in \boldsymbol{A}$ and $\sigma, \mu \in (0, 1)$ we define
  \[ {w^{\boldsymbol{a}}_{\mu}}  \assign \bigotimes_{i = 1}^{n (\boldsymbol{a})}
     {w^{\mathfrak{a}_i}_{\mu}} , \qquad K^{\boldsymbol{a}}_{\mu} \assign
     \bigotimes_{i = 1}^{n (\boldsymbol{a})} K^{1, k (\mathfrak{a}_i)}_{\mu},
     \qquad L^{\boldsymbol{a}}_{\sigma} \assign \bigotimes_{i = 1}^{n
     (\boldsymbol{a})} L^{1, k (\mathfrak{a}_i)}_{\sigma}, \qquad
     K^{\boldsymbol{a}}_{\eta, \sigma} \assign L^{\boldsymbol{a}}_{\sigma}
     K^{\boldsymbol{a}}_{\eta}, \]
  where ${w^{\mathfrak{a}_i}_{\mu}}$ is introduced in
  Definition~\ref{def:weights-kernels}, and $K^{1, k}_{\mu}$ and $L^{1, k}_{\sigma}$
  are introduced in Definition~\ref{def:regularizing-kernels-for-kernels}
\end{definition}

The space of cumulants is endowed with the norm $\trinorm{\Cdot}$ defined by taking the $L^{\infty}$ norm on the first output variable
and the total variation norm in all other variables, with the output variables
restricted to the first period in space. The restriction of the integration of
the output variables to the first period is natural since cumulants are
kernels of operators acting between spaces of periodic functions,
cf.~\eqref{eq:operator-cumulant}.

\begin{definition}
  \label{def:norm-cumulants}Let $\boldsymbol{a} \in \boldsymbol{A}$, $n = n
  (\boldsymbol{a})$ and $k = K (\boldsymbol{a})$. The norm of a kernel
  $\mathcal{F} ^{\boldsymbol{a}}$ associated with an operator
  $\mathcal{E}^{\otimes k} \rightarrow \hat{\mathcal{E}}^{\otimes n}$ is defined by
  \begin{equation}
    \trinorm{ \mathcal{F} ^{\boldsymbol{a}} } \assign \sup_{z_1 \in
    \Lambda} \int_{\Lambda^{n - 1}_M} \left( \int_{\Lambda^k} |
    \mathcal{F}^{\boldsymbol{a}} (z_1, \mathd Z_1, \ldots, z_n, \mathd Z_n) |
    \right)\,\mathd z_2 \cdots \mathd z_n, \label{eq:basic-norm-cumulant}
  \end{equation}
  with the notation as in~\eqref{eq:notation-Z}. Given a weight $w \in C
  (\Lambda_0^{n + k})$, we write
  \[ (\mathcal{F} \cdot w) (z_1, \mathd Z_1, \ldots, z_n, \mathd Z_n)
     \assign \mathcal{F} (z_1, \mathd Z_1, \ldots, z_n, \mathd Z_n) w (z_1,
     Z_1, \cdots, z_n, Z_n) . \]
  The norm of $\mathcal{F}^{\boldsymbol{a}}$ at scale $\sigma \in (1 / 2, 1)$ is
  defined by
  \begin{equation}
    \trinorm{ \mathcal{F}^{\boldsymbol{a}} }_{\sigma} =
    \trinorm{ [K^{\boldsymbol{a}}_{\sigma} \mathcal{F}^{\boldsymbol{a}}]
    {\cdot w^{\boldsymbol{a}}_{\sigma}}  },
    \label{eq:cumulant-norms-mu}
  \end{equation}
  where $K^{\boldsymbol{a}}_{\sigma} \mathcal{F}^{\boldsymbol{a}}$ denotes
  the convolution of $\mathcal{F}^{\boldsymbol{a}}$, viewed as an element of
  $\mathcal{S}' (\Lambda^{n (\boldsymbol{a}) + K (\boldsymbol{a})})$, with the
  kernel of $K^{\boldsymbol{a}}_{\sigma}$. We further introduce a norm
  for the family of cumulants of effective force kernels
  \[ \mathcal{F}^{\boldsymbol{A}} \assign
     (\mathcal{F}^{\boldsymbol{a}}_{\sigma})_{\boldsymbol{a} \in \boldsymbol{A},
     \sigma \in (1 / 2, 1)} \]
  defined by

  \begin{equation}
    \trinorm{ \mathcal{F}^{\boldsymbol{A}} } \assign
    \sup_{\boldsymbol{a} \in \boldsymbol{A}}  \biggl(\sup_{\sigma \in (1 / 2, 1)}
    \llbracket \sigma \rrbracket^{- [\boldsymbol{a}]}  \trinorm{
    \mathcal{F}_{\sigma}^{\boldsymbol{a}} }_{\sigma}\biggr)^{1/n
    (\boldsymbol{a})} . \label{eq:norm-family-cumulants}
  \end{equation}
  Here \(\boldsymbol A\) denotes the set of all finite tuples satisfying
  \(L(\boldsymbol a)\le 2\bar\ell\). Although this set is infinite, the induction
  below shows that, at each level \(L\le 2\bar\ell\), only finitely many cumulants
  are nonzero. All suprema below may therefore be restricted to this finite set.
\end{definition}

\begin{remark}
  \label{rem:norm-scale-monotonicity}Analogously to
  Remark~\ref{rem:norm0-scale-monotonicity}, by $K^{\boldsymbol{a}}_{\mu} =
  K^{\boldsymbol{a}}_{\mu, \sigma} K_{\sigma}^{\boldsymbol{a}}$, Young's
  inequality,~\eqref{eq:weight-w-2} and Lemma~\ref{lem:bounds-K-J}, we obtain
  \begin{align*}
    \trinorm{ \mathcal{F}^{\boldsymbol{a}} }_{\mu} =
    \trinorm{ [K^{\boldsymbol{a}}_{\mu} \mathcal{F}^{\boldsymbol{a}}]
    {\cdot w^{\boldsymbol{a}}_{\mu}}  } &=  \trinorm{
    [K^{\boldsymbol{a}}_{\mu,\sigma} K^{\boldsymbol{a}}_{\sigma}
    \mathcal{F}^{\boldsymbol{a}}] {\cdot w^{\boldsymbol{a}}_{\mu}} 
    }\\
 &\lesssim  \|
    K_{\mu, \sigma} \|_{\operatorname{TV} (w^{\flat}_{\mu})}^{K (\boldsymbol{a})} 
    \trinorm{ [K^{\boldsymbol{a}}_{\sigma} \mathcal{F}^{\boldsymbol{a}}]
    {\cdot w^{\boldsymbol{a}}_{\mu}}  }\\
 &\lesssim  \trinorm{ [K^{\boldsymbol{a}}_{\sigma}
    \mathcal{F}^{\boldsymbol{a}}] {\cdot w^{\boldsymbol{a}}_{\sigma}} 
    }\\
 &=  \trinorm{ \mathcal{F}^{\boldsymbol{a}} }_{\sigma},
  \end{align*}
  uniformly in $\mu \in [1 / 2, 1)$ and $\sigma \in [\mu, 1)$.
\end{remark}

In Section~\ref{sec:concl-of-cumul-analysis}, we will pass from estimates on the
norm $\trinorm{ \mathcal{F}^{\boldsymbol{A}} }$ for the family of
cumulants to estimates on the norm $\| F^{\mathfrak{A}} \|$ for the
corresponding family of kernels, using a Kolmogorov-type argument.
Following~{\cite{duchFlowEquationApproach2021,duchRenormalizationSingularElliptic2022}},
we introduce a flow equation for cumulants, which will allow us to control the
norm $\trinorm{ \mathcal{F}^{\boldsymbol{A}} }$.

\begin{lemma}
  \label{lem:flow-cumulants}The cumulants satisfy the following flow equation:
  \begin{equation}
    \partial_{\sigma} \mathcal{F}^{\boldsymbol{a}}_{\sigma} =
    -\sum_{\boldsymbol{b}} \mathcal{A}^{\boldsymbol{a}}_{\boldsymbol{b}}
    (\dot{G}_{\sigma}, \mathcal{F}^{\boldsymbol{b}}_{\sigma}) -
    \sum_{\boldsymbol{b}, \boldsymbol{c}}
    \mathcal{B}^{\boldsymbol{a}}_{\boldsymbol{b}, \boldsymbol{c}} (\dot{G}_{\sigma},
    \mathcal{F}^{\boldsymbol{b}}_{\sigma}, \mathcal{F}^{\boldsymbol{c}}_{\sigma}),
    \label{eq:cumulant-flow}
  \end{equation}
  where the multilinear operators $\mathcal{A}^{\boldsymbol{a}}_{\boldsymbol{b}}$
  and $\mathcal{B}^{\boldsymbol{a}}_{\boldsymbol{b}, \boldsymbol{c}}$ are defined in
  Appendix~\ref{sec:A-B-cumulants}. We have
  $\mathcal{A}^{\boldsymbol{a}}_{\boldsymbol{b}} = 0$ unless
  \begin{equation}
    \left\{\begin{aligned}
      n (\boldsymbol{a}) = n (\boldsymbol{b}) - 1, \, \,\\
      L (\boldsymbol{a}) = L (\boldsymbol{b}) + 1,\\
      K (\boldsymbol{a}) = K (\boldsymbol{b}) - 1,\\
      {}[\boldsymbol{a}] = [\boldsymbol{b}] - \theta + \delta - \beta,
    \end{aligned}\right. \label{eq:characterization-a}
  \end{equation}
  and $\mathcal{B}^{\boldsymbol{a}}_{\boldsymbol{b}, \boldsymbol{c}}
  = 0$ unless
  \begin{equation}
    \left\{\begin{aligned}
      n (\boldsymbol{a}) = n (\boldsymbol{b}) + n (\boldsymbol{c}) - 1, \, \,\\
      L (\boldsymbol{a}) = L (\boldsymbol{b}) + L (\boldsymbol{c}) + 1,\\
      K (\boldsymbol{a}) = K (\boldsymbol{b}) + K (\boldsymbol{c}) - 1,\\
      {}[\boldsymbol{a}] = [\boldsymbol{b}] + [\boldsymbol{c}] + \varrho - \theta +
      \delta - \beta .
    \end{aligned}\right. \label{eq:characterization-b}
  \end{equation}
  Moreover, the operators $\mathcal{A}^{\boldsymbol{a}}_{\boldsymbol{b}}$ and
  $\mathcal{B}^{\boldsymbol{a}}_{\boldsymbol{b}, \boldsymbol{c}}$ satisfy the
  following bounds:
  \[ \begin{aligned}
       \llbracket \sigma \rrbracket^{- [\boldsymbol{a}]}  \trinorm{
       \mathcal{A}^{\boldsymbol{a}}_{\boldsymbol{b}} (\dot{G}_{\sigma},
       \mathcal{F}^{\boldsymbol{b}}_{\sigma}) }_{\sigma} & \lesssim
       \llbracket \sigma \rrbracket^{- [\boldsymbol{b}] - 1}  \trinorm{
       {\mathcal{F}^{\boldsymbol{b}}}  }_{\sigma},\\
       \llbracket \sigma \rrbracket^{- [\boldsymbol{a}]}  \trinorm{
       \mathcal{B}^{\boldsymbol{a}}_{\boldsymbol{b}, \boldsymbol{c}}
       (\dot{G}_{\sigma}, \mathcal{F}^{\boldsymbol{b}}_{\sigma},
       \mathcal{F}^{\boldsymbol{c}}_{\sigma}) }_{\sigma} & \lesssim
       \llbracket \sigma \rrbracket^{- [\boldsymbol{b}] - [\boldsymbol{c}] - 1}
       \trinorm{ \mathcal{F}^{\boldsymbol{b}}_{\sigma} }_{\sigma}
       \trinorm{ \mathcal{F}^{\boldsymbol{c}}_{\sigma} }_{\sigma},
     \end{aligned} \]
  provided the following compatibility conditions hold:
  \begin{equation}
    \theta + \beta - \delta - d \geqslant 0, \qquad - \varrho + \theta + \beta
    - \delta + 2 s \geqslant 0. \label{eq:compatibility}
  \end{equation}
\end{lemma}

\begin{proof}
  The derivation of the flow equation is a direct consequence of the
  definition of cumulants,
  see~{\cite{duchFlowEquationApproach2021,duchRenormalizationSingularElliptic2022}}.
  The explicit form of the operators is not essential for the subsequent
  discussion and is provided in Appendix~\ref{sec:A-B-cumulants}, where the
  stated bounds are also proved (see Lemma~\ref{lem:estimates-A} and
  Lemma~\ref{lem:estimates-B}).
\end{proof}

This general structure of the flow equation~\eqref{eq:cumulant-flow} allows us
to propagate estimates for the cumulants of the form
\[ \sup_{\sigma \in (1 / 2, 1)} \llbracket \sigma \rrbracket^{-
   [\boldsymbol{a}]}  \trinorm{ \mathcal{F}^{\boldsymbol{a}}_{\sigma}
   }_{\sigma} \lesssim 1 . \]
However, depending on the sign of $[\boldsymbol{a}]$, we shall handle
the cumulants differently: in particular, for $[\boldsymbol{a}] > 0$, namely for
irrelevant cumulants, the flow equation can be solved \textit{backward}
starting from the final condition at $\sigma = 1$. On the other hand, this
approach does not work for cumulants for which $[\boldsymbol{a}] < 0$ as in this
case the flow equation cannot be integrated close to $\sigma = 1$. As we shall
see in Section~\ref{sec:induction}, we will solve the flow equation for this
class of cumulants, called {\textit{relevant cumulants}}, by integrating it
{\textit{forward}}.

\begin{remark}
  \label{rem:symmetries}Before proceeding with the analysis of the flow
  equations for cumulants, we record a few remarks about symmetries. First, we
  observe that the SPDE under consideration, namely~\eqref{eq:dynamics-eps-M},
  is invariant under the transformation
  \[ \phi \mapsto - \phi, \qquad \xi \mapsto - \xi, \]
  which also preserves the law of the noise $\xi = \xi^{(\varepsilon, M)}$. By
  the discussion in Section~\ref{sec:bounds-on-parameters}, this invariance
  implies that $\mathcal{F}^{\boldsymbol{a}}_{\sigma} = 0$ for all $\sigma \in
  [1/2, 1]$ and $\boldsymbol{a} \in \boldsymbol{A}$ such that $L (\boldsymbol{a}) = 0$
  and
  \[ n (\boldsymbol{a}) + K (\boldsymbol{a}) \in 2\mathbb{N}_0 + 1. \]
  To extend this property to higher levels $L
  (\boldsymbol{a}) > 0$, we exploit the fact that it is preserved by the flow
  equation~\eqref{eq:cumulant-flow} thanks to the compatibility
  conditions~\eqref{eq:characterization-a} and~\eqref{eq:characterization-b}.
  Indeed, fix an arbitrary $\ell \in \mathbb{N}_+$ and assume inductively that
  \[ \mathcal{F}^{\boldsymbol{a}}_{\sigma} = 0, \qquad \forall \sigma \in [1/2,
     1],\, \boldsymbol{a} \in \boldsymbol{A} \quad \text{such that}
     \quad L (\boldsymbol{a}) < \ell \quad \text{and} \quad n (\boldsymbol{a}) + K
     (\boldsymbol{a}) \in 2\mathbb{N}_0 + 1. \]
  Then from the flow equation~\eqref{eq:cumulant-flow} and the
  conditions~\eqref{eq:characterization-a} and~\eqref{eq:characterization-b},
  it follows that
  \begin{align*}
     \partial_{\sigma} \mathcal{F}^{\boldsymbol{a}}_{\sigma} = 0,
     &\quad \forall \sigma \in [1/2, 1), \,\boldsymbol{a} \in \boldsymbol{A}
     \quad \text{such that} \quad
     L (\boldsymbol{a}) = \ell \quad \text{and} \quad
     n (\boldsymbol{a}) + K (\boldsymbol{a}) \in 2\mathbb{N}_0 + 1.
  \end{align*}
  Since by construction $\mathcal{F}^{\boldsymbol{a}}_{\sigma = 1} = 0$ for all
  $\boldsymbol{a} \in \boldsymbol{A}$ such that $n (\boldsymbol{a}) + K
  (\boldsymbol{a}) \in 2\mathbb{N}_0 + 1$, we conclude that
  \[ \mathcal{F}^{\boldsymbol{a}}_{\sigma} = 0, \qquad \forall \sigma \in [1/2,
     1],\, \boldsymbol{a} \in \boldsymbol{A} \quad \text{such that}
     \quad n (\boldsymbol{a}) + K (\boldsymbol{a}) \in 2\mathbb{N}_0 + 1. \]
  
  A further symmetry is given by spatial reflection, that is, the
  transformation
  \[ \phi (t, x) \mapsto \phi (t, - x), \qquad \xi (t, x) \mapsto \xi (t, -
     x), \]
  which also leaves the law of $\xi$ invariant. By an argument analogous to
  the one above, we infer that for any $\boldsymbol{a} \in
  \boldsymbol{A}$, the cumulant $\mathcal{F}^{\boldsymbol{a}}_{\sigma}$ is symmetric
  under spatial reflections.
  Finally, the law of the noise and the deterministic kernels in the final
  condition are invariant under spacetime translations, and the operators in
  the cumulant flow equation~\eqref{eq:cumulant-flow} preserve this invariance;
  see also the explicit formulas~\eqref{eq:specific-operator-A}
  and~\eqref{eq:specific-operator-B}. Consequently, all cumulants are invariant
  under simultaneous spacetime translations of their output and input
  variables. In particular, if $n(\boldsymbol a)=K(\boldsymbol a)=1$, then
  \begin{equation}
    z\longmapsto
    \int_\Lambda \mathcal F^{\boldsymbol a}_\sigma(z,\mathd z_1)
    \quad\text{is constant.}
    \label{eq:one-point-mass-translation}
  \end{equation}
\end{remark}

\subsection{Bounds on parameters}\label{sec:bounds-on-parameters}

We shall now fix the parameters
\[ \beta, \theta, \varrho, \alpha, \delta, \flat \]
introduced in the analysis of the kernels. We have to choose these parameters
so that~\eqref{eq:compatibility} is satisfied. Another constraint comes from
the requirement that $\trinorm{ \mathcal{F}^{\boldsymbol{a}}_{\sigma}
}_{\sigma} \lesssim \llbracket \sigma \rrbracket^{ [\boldsymbol{a}]}$
for $\boldsymbol{a} \in \boldsymbol{A}$ with $L (\boldsymbol{a}) = 0$. Note that for
$\boldsymbol{a} \in \boldsymbol{A}$ such that $L (\boldsymbol{a}) = 0$, we have
$\partial_{\sigma} \mathcal{F}^{\boldsymbol{a}}_{\sigma} = 0$ and
$\mathcal{F}^{\boldsymbol{a}}_{\sigma} =\mathcal{F}^{\boldsymbol{a}}_1$. Using the
equality
\[ F^{[0]}_{\sigma} (\varphi) = F^{[0]}_1 (\varphi) = - \lambda \varphi^3 +
   \bar{r}  \varphi  + \xi^{(\varepsilon, M)} \]
and the fact that the noise $\xi^{(\varepsilon, M)}$ is Gaussian, one shows
that for $\boldsymbol{a} \in \boldsymbol{A}$ such that $L (\boldsymbol{a}) = 0$ the
cumulant $\mathcal{F}^{\boldsymbol{a}}_{\sigma}$ is nonzero only if:
\begin{enumerate}[(a)]
  \item it is the covariance of the noise, that is $n (\boldsymbol{a}) = 2$, $K(\boldsymbol{a}) = 0$ and
  \[ \mathcal{F}^{\boldsymbol{a}}_{\sigma} (z, z') =\mathfrak{K}_2 (F^{[0], (0)}
     (z), F^{[0], (0)} (z')) =\mathbb{E} [\xi^{(\varepsilon, M)} (z)
     \xi^{(\varepsilon, M)} (z')] = \delta^M (z - z'), \]
  where $\delta^M$ is the periodisation in space of the Dirac
  delta with period $M$, or
  
  \item it is the expected value of the (deterministic) kernel $F^{[0],
  (3)}_1$, that is $n (\boldsymbol{a}) = 1$, $K (\boldsymbol{a}) = 3$ and
  \[ \mathcal{F}^{\boldsymbol{a}}_{\sigma} (z, z_1, z_2, z_3) = - \lambda \delta
     (z - z_1) \delta (z - z_2) \delta (z - z_3), \]
  or
  
  \item it is the expected value of the (deterministic) kernel $F^{[0],
  (1)}_1$, that is $n (\boldsymbol{a}) = 1$, $K (\boldsymbol{a}) = 1$ and
  \[ \mathcal{F}^{\boldsymbol{a}}_{\sigma} (z, z_1) = \bar{r} \delta (z - z_1) .
  \]
\end{enumerate}
In case (a), using the fact that $w_{\sigma}^{\boldsymbol{a}}
= 1$, we conclude that
\begin{equation}
  {\trinorm{ \mathcal{F}^{\boldsymbol{a}}_1 }_{\sigma}}  = \int_{\Lambda_M} \delta^M
  (\mathd z) = 1. \label{eq:est-cov-noise}
\end{equation}
As a consequence, we have to require that
\begin{equation}
  [\boldsymbol{a}] = - \varrho + 2 \theta \leqslant 0. \label{eq:two-point}
\end{equation}
In case (b), we observe that

\[ (K_{\sigma}^{\boldsymbol{a}} \mathcal{F}^{\boldsymbol{a}}_{\sigma}) (z ; z_1,
   z_2, z_3) = - \lambda \prod_{i = 1}^3 K_{\sigma} (z - z_i), \]
and thus $\trinorm{ \mathcal{F}^{\boldsymbol{a}}_{\sigma} }_{\sigma}
\lesssim 1$ uniformly in $\sigma \in (1/2, 1)$ by Lemma~\ref{lem:bounds-K-J}.
Consequently, we have to require that
\begin{equation}
  [\boldsymbol{a}] = - \varrho + \theta + 3 \beta \leqslant 0
  \label{eq:initial-c} .
\end{equation}
In case (c), we have

\[ (K_{\sigma}^{\boldsymbol{a}} \mathcal{F}^{\boldsymbol{a}}_{\sigma}) (z ; z_1) =
   \bar{r} K_{\sigma} (z - z_1), \]
and thus $\trinorm{ \mathcal{F}^{\boldsymbol{a}}_{\sigma} }_{\sigma}
\lesssim 1$ uniformly in $\sigma \in (1/2, 1)$. Consequently, we have to require that
\[ [\boldsymbol{a}] = - \varrho + \theta + \beta \leqslant 0. \]
Note that since $\beta > 0$ the last condition is implied by the
bound~\eqref{eq:initial-c}.

We collect all constraints here and fix the parameters once and for all, so
that the later arguments can refer directly to this table.
\begin{equation}
  \begin{array}{ll|l}
    {}[B] & \text{Flow kernels~$\eqref{eq:constraint-flow-kernels}$} \quad
    & 2 s - \alpha + \beta - \delta \geqslant 0\\
    {}[\mathbb{A}] & \text{Flow cumulants
    $\mathcal{A}$~\eqref{eq:compatibility}}  & \theta + \beta - \delta - d
    \geqslant 0\\
    {}[\mathbb{B}] & \text{Flow cumulants
    $\mathcal{B}$~\eqref{eq:compatibility}} & 2 s - \varrho + \theta + \beta -
    \delta \geqslant 0\\
    {}[\Xi] & \text{Initial condition $F^{[0], (0)}$~\eqref{eq:two-point}} & -
    \varrho + 2 \theta \leqslant 0\\
    {}[\Phi^3] & \text{Initial condition $F^{[0], (3)}$~\eqref{eq:initial-c}}
    & - \varrho + \theta + 3 \beta \leqslant 0\\
    {}[K] & \text{Kolmogorov~$\eqref{eq:wanted-bound}$} \quad & \alpha -
    \frac{\varrho}{2} + \theta - \frac{d + 2 s}{2} - \kappa > 0\\
    {}[\bar{F}] & \text{Kolmogorov $\bar{F}$~$\quad
    \eqref{eq:constraint-f-bar}$} & \alpha - \varrho + \theta - \kappa
    \geqslant 0
  \end{array} \label{eq:summary-ineq}
\end{equation}

The parameter $\kappa > 0$ quantifies the loss of regularity in the
Kolmogorov-type argument used to estimate the pathwise behaviour of the random
kernels and will be chosen sufficiently small.

From the constraints $[\Phi^3]$ and $[\mathbb{B}]$ we have
\[ \beta \leqslant s - \delta / 2, \]
while from $[\Xi]$, $[\mathbb{B}]$ and this last inequality we deduce that
\[ \theta \leqslant \beta - \delta + 2 s \leqslant 3 s - 3 \delta / 2. \]
Using now $[\mathbb{B}]$ again we have
\[ \varrho \leqslant \theta + \beta - \delta + 2 s \leqslant 6 s - 3 \delta .
\]
And $[\mathbb{A}]$ now gives
\begin{equation}
  0 \leqslant \theta + \beta - \delta - d \leqslant 4 s - d - 3 \delta = 3
  (\delta_{\star} - \delta), \label{eq:sub-crit}
\end{equation}
where the strict positivity of
\[ \delta_{\star} \assign \frac{4 s - d}{3} > 0, \]
defines the subcritical regime of this model. We now choose the remaining
parameters so as to saturate most of the inequalities above, namely:
\begin{equation}
  \begin{aligned}
    \beta = s - \delta / 2 - \kappa / 2, \quad \theta = 3 \beta, \quad \varrho
    = 2 \theta, \quad \alpha = 3 \beta + \kappa, \quad \delta
    = \delta_{\star} / 2 .\\
  \end{aligned} \label{eq:values-param}
\end{equation}
These choices are not intended to be optimal. Their purpose is to make the
power counting transparent while leaving positive margins, proportional to
$\delta_{\star}$ and $\kappa$, that can absorb the losses introduced by
localisation and by the lattice-uniform Kolmogorov argument.
By substituting these values into the inequalities~\eqref{eq:summary-ineq}
above, we check the following.
\begin{equation}
  \begin{array}{llll}
    {}[B] & 0 \leqslant 2 s - \alpha + \beta - \delta & = 2 s - 2 \beta -
    \kappa - \delta & = 0\\
    {}[\mathbb{A}] & 0 \leqslant \theta + \beta - \delta - d & = 4 \beta -
    \delta - d & = \frac{3}{2} \delta_{\star} - 2 \kappa\\
    {}[\mathbb{B}] & 0 \leqslant - \varrho + \theta + \beta - \delta + 2 s & =
    - 2 \beta - \delta + 2 s & = \kappa\\
    {}[\Xi] & 0 \geqslant - \varrho + 2 \theta & = 0 & = 0\\
    {}[\Phi^3] & 0 \geqslant - \varrho + \theta + 3 \beta & = 0 & = 0\\
    {}[K] & 0 < \alpha - \frac{\varrho}{2} + \theta - \frac{d + 2 s}{2} -
    \kappa & = 3 \beta - \frac{d + 2 s}{2} & = \frac{3}{4}
    \delta_{\star} - \frac{3}{2} \kappa\\
    {}[\bar{F}] & 0 \leqslant \alpha - \varrho + \theta - \kappa & = 0 & = 0
  \end{array} \label{eq:ver-param}
\end{equation}
We used that
\begin{equation}
  \frac{d + 2 s}{2} = \frac{6 s - 3 \delta_{\star}}{2} = 3 \beta - \frac{3}{2}
  (\delta_{\star} - \delta) + \frac32 \kappa . \label{eq:dplustwos}
\end{equation}
All the inequalities are satisfied provided
\begin{equation}
  \kappa \in (0, \delta_{\star} / 2) . \label{eq:cond-kappa}
\end{equation}
In addition, we choose $\kappa$ so that
\begin{equation}
  \frac{2 \beta}{\delta} \notin \mathbb{N}_+ .
  \label{eq:no-marginal-cumulants}
\end{equation}
Equivalently,
\[
  \delta \ell - 2 \beta \neq 0
  \qquad \text{for every } \ell \in \mathbb{N}_+ .
\]
Since $2 \beta = 2 s - \delta - \kappa$, the excluded values are
\[
  \kappa = 2 s - \delta (\ell + 1),
  \qquad \ell \in \mathbb{N}_+ .
\]
They form a discrete set. Hence this condition is compatible
with~\eqref{eq:cond-kappa} and with choosing $\kappa > 0$ arbitrarily small
subject to the remaining strict parameter inequalities.
Moreover, in Lemma~\ref{lem:loc-kernels} below we use the fact that the
parameter of the weight $\mathfrak{o}$ introduced in
Definition~\ref{def:weights-kernels} satisfies
\begin{equation}
  \kappa_{\mathfrak{o}} < \delta \wedge (2 - \flat) / (2 \bar{\ell} + 2) .
  \label{eq:cond-kappa-o}
\end{equation}

Finally, we explain how the value of the parameter $\flat$, introduced in
Definition~\ref{def:weights-kernels}, is fixed. The role of this parameter is to
compensate for the loss of spatial weight arising in the localisation
procedure for the relevant cumulants discussed below, where the Taylor
expansion allows us to gain in homogeneity and establish that only a local
renormalisation is needed. In Lemma~\ref{lem:cumul-est-1} and
Lemma~\ref{lem:loc-kernels}, we shall need $\flat > 2 s - \delta$. On the
other hand, as argued in Remark~\ref{rem:decay-rate}, we must also impose
$\flat < 2 s$. Thus, {\underline{we fix $\flat$}} such that
\begin{equation}
  2 s - \delta < \flat < 2 s. \label{eq:choice-of-flat} 
\end{equation}

In Section~\ref{sec:post-processing}, we shall introduce additional constraints
on the parameters $\kappa$ and $\kappa_{\mathfrak{o}}$, and we will fix all
remaining parameters
\[ \gamma, \vartheta, \bar{\kappa}, \bar{\ell}, \bar{k}, a, \nu,
   \kappa_{\mathfrak{o}}, \kappa, \]
which play a role in our analysis. We postpone this discussion because, apart
from $\kappa$ and $\kappa_{\mathfrak{o}}$, which must satisfy
conditions~\eqref{eq:cond-kappa} and~\eqref{eq:cond-kappa-o}, these parameters
are only relevant for the results established in
Section~\ref{sec:post-processing}. We emphasise once again that all these
parameters are to be regarded as fixed once and for all, and their choice
depends solely on the exponent $s \in (3 / 4, 1)$ of the fractional
Laplacian.

\subsection{Classification of
cumulants}\label{sec:class-rel-cumulants}

Given the bounds on the parameters from the previous section, we can now
examine the class of cumulants that are relevant or marginal, \textit{i.e.},
those satisfying $[\boldsymbol{a}] \leqslant 0$. Recalling that
\[ \varrho = 2 \theta, \quad \text{and} \quad \theta = \varrho - 3 \beta, \]
observe that the condition $[\boldsymbol{a}] \leqslant 0$ can be written as
\[
  [\boldsymbol{a}] = \theta (n (\boldsymbol{a}) - 2) + \beta K
  (\boldsymbol{a}) + \delta L (\boldsymbol{a}) = \beta (3 n
  (\boldsymbol{a}) - 6 + K (\boldsymbol{a})) + \delta L (\boldsymbol{a})
  \leqslant 0.
\]
Then
\begin{enumerate}[(a)]
  \item if $n (\boldsymbol{a}) > 2$, then all cumulants are irrelevant;
  
  \item if $n (\boldsymbol{a}) = 2$, then
  \[
    [\boldsymbol{a}] = \beta K (\boldsymbol{a}) + \delta L (\boldsymbol{a}) .
  \]
  Hence the only non-positive case is
  $L (\boldsymbol{a}) = K (\boldsymbol{a}) = 0$, namely the covariance of the
  noise
  \[
    \mathcal{F}^{\boldsymbol{a}}_{\sigma} =
    \mathfrak{K}_2 (F^{[0], (0)}_{\sigma}, F^{[0], (0)}_{\sigma}) .
  \]
  This cumulant has homogeneity zero and is independent of $\sigma$.
  
  \item if $n (\boldsymbol{a}) = 1$, then
  \[
    [\boldsymbol{a}] = \beta (K (\boldsymbol{a}) - 3) + \delta L
    (\boldsymbol{a}) .
  \]
  By Remark~\ref{rem:symmetries}, the cumulants with $K (\boldsymbol{a}) = 0$
  and $K (\boldsymbol{a}) = 2$ vanish. For $K (\boldsymbol{a}) = 3$, we have
  \[
    [\boldsymbol{a}] = \delta L (\boldsymbol{a}) .
  \]
  Thus the only non-positive case is $L (\boldsymbol{a}) = 0$, namely the
  deterministic cubic kernel
  \[
    \mathfrak{K}_1 (F_{\sigma}^{[0], (3)}) ,
  \]
  which is marginal and independent of $\sigma$. For $K (\boldsymbol{a}) = 1$,
  we have
  \[
    [\boldsymbol{a}] = \delta L (\boldsymbol{a}) - 2 \beta .
  \]
  At $L (\boldsymbol{a}) = 0$, this is the deterministic linear kernel
  $\mathfrak{K}_1 (F_{\sigma}^{[0], (1)}) = \bar{r}\,\delta$,
  which is relevant and independent of $\sigma$. At
  $L (\boldsymbol{a}) = \ell \geqslant 1$, these are the cumulants
  \[
    \mathbb{E} F_{\sigma}^{[\ell], (1)} .
  \]
  By~\eqref{eq:no-marginal-cumulants}, none of them is marginal.
\end{enumerate}
Summarising, the only level-zero marginal cumulants are
\[
  \mathfrak{K}_2 (F^{[0], (0)}_{\sigma}, F^{[0], (0)}_{\sigma}),
  \qquad
  \mathfrak{K}_1 (F_{\sigma}^{[0], (3)}),
\]
and the deterministic linear kernel
\[
  \mathfrak{K}_1 (F_{\sigma}^{[0], (1)})
\]
is relevant. All three have trivial flow. By~\eqref{eq:no-marginal-cumulants},
there are no marginal cumulants with $L (\boldsymbol{a}) > 0$, so the only
positive-level cumulants that must be analysed in detail are
\begin{equation}
  \mathfrak{K}_1 (F_{\sigma}^{[\ell],(1)}) =\mathbb{E}F_{\sigma}^{[\ell],(1)},
  \qquad \ell \in \{ 1, \ldots, \hat{\ell} \}, \label{eq:relevant-cum}
\end{equation}
where
\[
  \hat{\ell} \assign \max \left\{ \ell \in \mathbb{N}_+ : \delta \ell - 2
  \beta < 0 \right\}.
\]
By~\eqref{eq:no-marginal-cumulants}, we then have
\[
  \delta \hat{\ell} - 2 \beta < 0 < \delta (\hat{\ell} + 1) - 2 \beta .
\]
Equivalently,
\[
  \hat{\ell} = \left\lfloor \frac{2 \beta}{\delta} \right\rfloor .
\]
Thus the positive-level relevant cumulants are exactly
\[
  \mathbb{E} F_{\sigma}^{[\ell], (1)},
  \qquad 1 \leqslant \ell \leqslant \hat{\ell} .
\]
Since the parameter $\hat{\ell}$ appears only in this and
the following section, it is not included in the list~\eqref{eq:parameters}.
Note that $\bar\ell$ is chosen so that $\hat\ell\leq\bar\ell$.

Thus, the renormalisation problem is purely of
mass type: the only relevant components are the one-point kernels that are
linear in the field. All other positive-level cumulants are irrelevant and
can be reconstructed by integrating the flow from their terminal condition.

\subsection{Inductive procedure}\label{sec:induction}

The aim of this section is to derive bounds for the cumulants
$\mathcal{F}^{\boldsymbol{a}}_{\sigma}$ using an induction on $L
(\boldsymbol{a})$, based on the flow equation~\eqref{eq:cumulant-flow}. More
precisely, Lemma~\ref{lem:flow-cumulants} shows that the flow equation defines
a triangular system with respect to $L (\boldsymbol{a})$. For irrelevant
cumulants, we have $\mathcal{F}^{\boldsymbol{a}}_{\sigma = 1} = 0$, and the
desired bounds follow by directly integrating the flow equation from $\sigma =
1$. For relevant cumulants $\mathcal{F}^{\boldsymbol{a}}_{\sigma}$, which for $L
(\boldsymbol{a}) > 0$ coincide with \eqref{eq:relevant-cum}, the bound for
$\partial_{\sigma} \mathcal{F}^{\boldsymbol{a}}_{\sigma}$ is not integrable at
$\sigma = 1$. To overcome this, we decompose relevant cumulants into a local
relevant part and a~nonlocal irrelevant part, using the procedure detailed in
Appendix~\ref{app:localization}. The irrelevant part can again be integrated
from $\sigma = 1$ with zero boundary condition. For the local relevant part, we
instead impose the boundary condition at $\sigma = 1 / 2$, which acts as the
renormalisation condition fixing the mass counterterm $r_{\varepsilon, M}$.
This allows the integration of the flow equation from $\sigma = 1 / 2$,
thereby avoiding the non-integrable singularity at $\sigma = 1$.

\begin{lemma}
  \label{lem:cumul-est-1}For any $\bar{r} \in \mathbb{R}$ there exists
  a choice of constants $(r^{[\ell]}_{\varepsilon, M})_{\ell = 1,
  \ldots, \hat{\ell}}$ such that the solution
  \[ F^{\mathfrak{A}} \assign (F^{\mathfrak{a}}_{\sigma})_{\mathfrak{a} \in
     \mathfrak{A}, \sigma \in (1 / 2, 1)}, \]
  of the approximate flow equation with the final condition
  \begin{equation}
    \begin{aligned}
      F^{[0]}_1 (\phi) & = - \lambda \phi^3 + \bar{r} \phi +
      \xi^{(\varepsilon, M)},& \\
      F^{[\ell]}_1 (\phi) & = r^{[\ell]}_{\varepsilon, M} \phi, &\ell \in \{
      1, \ldots, \hat{\ell} \},\\
      F^{[\ell]}_1 (\phi) & = 0, &\ell \in \{ \hat{\ell} + 1, \ldots,
      \bar{\ell} \},
    \end{aligned} \label{eq:initial-cond-induction-lemma}
  \end{equation}
  satisfies
  \[ \trinorm{ \mathcal{F}^{\boldsymbol{A}} } = \sup_{\boldsymbol{a}
     \in \boldsymbol{A}} \biggl[\sup_{\sigma \in (1 / 2, 1)} \llbracket \sigma
     \rrbracket^{- [\boldsymbol{a}]}  \trinorm{
     \mathcal{F}^{\boldsymbol{a}}_{\sigma} }_{\sigma}\biggr]^{1 / n
     (\boldsymbol{a})} \lesssim 1. \]
\end{lemma}

\begin{remark}
  Note that
  \[ F  (\phi) = \sum_{\ell = 0}^{\bar{\ell}} F^{[\ell]}_1 (\phi), \qquad
     r_{\varepsilon, M} = \bar{r} + \sum_{\ell = 1}^{\hat{\ell}}
     r^{[\ell]}_{\varepsilon, M}, \]
  where $F = F^{(\varepsilon, M)}$ is the force~\eqref{eq:force} and $r
  _{\varepsilon, M}$ is the mass renormalisation, which appears in the
  expression~\eqref{eq:action} for the action. The parameter $\hat{\ell} \in
  \mathbb{N}_+$ was fixed at the end of Section~\ref{sec:class-rel-cumulants},
  whereas the parameter $\bar{\ell} \geqslant \hat{\ell}$ will be fixed in
  Section~\ref{sec:post-processing}.
\end{remark}

\begin{proof}
  For $\ell \in \{ 0, \ldots, 2 \bar{\ell} \}$, define
  \[ M_{\ell} \assign 1 + \sup_{\substack{
       \boldsymbol{a} \mid L (\boldsymbol{a}) \leqslant \ell
     }} \biggl[\sup_{\sigma \in (1 / 2, 1)} \llbracket \sigma
     \rrbracket^{- [\boldsymbol{a}]}  \trinorm{
     \mathcal{F}^{\boldsymbol{a}}_{\sigma} }_{\sigma}\biggr]^{1 / n
     (\boldsymbol{a})} . \]
  We will prove by induction on $\ell$ that
  \begin{equation}
    M_{\ell} < \infty, \quad \operatorname{and} \quad
    \mathcal{F}^{\boldsymbol{a}}_{\sigma} = 0, \qquad \text{for all
    $\boldsymbol{a}$ such that $L (\boldsymbol{a}) \leqslant \ell$ and $n
    (\boldsymbol{a}) > 2^{\ell + 1}$,} \label{eq:induction}
  \end{equation}
  for any $\ell \in \{ 0, \ldots, 2 \bar{\ell} \}$.
  
  Let us remark that the second part of the above statement implies that there
  are only finitely many nonzero cumulants
  $\mathcal{F}^{\boldsymbol{a}}_{\sigma}$ such that $L (\boldsymbol{a}) \leqslant
  \ell$. We first note that the case $\ell = 0$ was already discussed in Section
  \ref{sec:bounds-on-parameters}, where it was proved that $M_0 < \infty$.
  
  Let us now consider the induction step. Assume that the
  conditions expressed by~\eqref{eq:induction} hold true at order $\ell -
  1$, for a fixed $\ell \in \mathbb{N}_+$. We shall prove that then the same
  is true at order~$\ell$. The proof is based on the flow
  equation~\eqref{eq:cumulant-flow}. We first note that, by Lemma
  \ref{lem:flow-cumulants}, we have
  \[ \begin{aligned}
       \trinorm{ \mathcal{A}^{\boldsymbol{a}}_{\boldsymbol{b}}
       (\dot{G}_{\sigma}, \mathcal{F}^{\boldsymbol{b}}_{\sigma})
       }_{\sigma} & \lesssim  \llbracket \sigma
       \rrbracket^{[\boldsymbol{a}] - 1} M_{\ell - 1}^{n (\boldsymbol{b})},\\
       \trinorm{ \mathcal{B}^{\boldsymbol{a}}_{\boldsymbol{b}, \boldsymbol{c}}
       (\dot{G}_{\sigma}, \mathcal{F}^{\boldsymbol{b}}_{\sigma},
       \mathcal{F}^{\boldsymbol{c}}_{\sigma}) }_{\sigma} & \lesssim
       \llbracket \sigma \rrbracket^{[\boldsymbol{a}] - 1} M_{\ell - 1}^{n
       (\boldsymbol{b})} M_{\ell - 1}^{n (\boldsymbol{c})},
     \end{aligned} \]
  for all $\boldsymbol{a} \in \boldsymbol{A}$ such that $L (\boldsymbol{a}) = \ell$.
  Recall that $n (\boldsymbol{a}) = n (\boldsymbol{b}) - 1$ and $L (\boldsymbol{b})
  = L (\boldsymbol{a}) - 1 = \ell - 1$ in the first line and $n (\boldsymbol{a}) =
  n (\boldsymbol{b}) + n (\boldsymbol{c}) - 1$ and $L (\boldsymbol{b}) + L
  (\boldsymbol{c}) = L (\boldsymbol{a}) - 1 = \ell - 1$ in the second line above.
  Hence, using the induction hypothesis we conclude that
  \[ \mathcal{A}^{\boldsymbol{a}}_{\boldsymbol{b}} (\dot{G}_{\sigma},
     \mathcal{F}^{\boldsymbol{b}}_{\sigma}) = 0, \qquad
     \mathcal{B}^{\boldsymbol{a}}_{\boldsymbol{b}, \boldsymbol{c}}
     (\dot{G}_{\sigma}, \mathcal{F}^{\boldsymbol{b}}_{\sigma},
     \mathcal{F}^{\boldsymbol{c}}_{\sigma}) = 0, \]
  if $n (\boldsymbol{a}) > 2^{\ell + 1}$, and
  \begin{equation}
    \trinorm{ \mathcal{A}^{\boldsymbol{a}}_{\boldsymbol{b}} (\dot{G}_{\sigma},
    \mathcal{F}^{\boldsymbol{b}}_{\sigma}) }_{\sigma} \lesssim
    \llbracket \sigma \rrbracket^{[\boldsymbol{a}] - 1}, \qquad \trinorm{
    \mathcal{B}^{\boldsymbol{a}}_{\boldsymbol{b}, \boldsymbol{c}} (\dot{G}_{\sigma},
    \mathcal{F}^{\boldsymbol{b}}_{\sigma}, \mathcal{F}^{\boldsymbol{c}}_{\sigma})
    }_{\sigma} \lesssim \llbracket \sigma \rrbracket^{[\boldsymbol{a}]
    - 1} . \label{eq:bound-ind}
  \end{equation}
  In the rest of the proof of the induction step we treat separately two
  cases: $\mathcal{F}^{\boldsymbol{a}}_{\sigma}$ is irrelevant, \textit{i.e.}
  $[\boldsymbol{a}] > 0$, or $\mathcal{F}^{\boldsymbol{a}}_{\sigma}$ is relevant,
  {\textit{i.e.}}, $[\boldsymbol{a}] < 0$. Note that, by
  \eqref{eq:no-marginal-cumulants}, there are no marginal cumulants at level
  $\ell > 0$.
  
  If $\mathcal{F}^{\boldsymbol{a}}$ is irrelevant, {\textit{i.e.}}
  $[\boldsymbol{a}] > 0$, then we can bound it by integrating the flow
  equation~\eqref{eq:cumulant-flow} \textit{backward} from the final
  condition at $\sigma = 1$. We stress that for any irrelevant cumulant this
  final condition is vanishing, namely $\mathcal{F}^{\boldsymbol{a}}_1 = 0$.
  Indeed, for $1 \leqslant \ell \leqslant \hat{\ell}$ the only nonzero
  positive-level final kernel is
  \[
    F_1^{[\ell], (1)}(z, z_1) = r^{[\ell]}_{\varepsilon, M} \delta (z, z_1),
  \]
  and its one-point cumulant has homogeneity
  \[
    \delta \ell - 2 \beta < 0.
  \]
  For $\ell > \hat{\ell}$, the final kernel vanishes by
  \eqref{eq:initial-cond-induction-lemma}. Moreover, any cumulant of order at
  least two containing a deterministic positive-level final kernel vanishes.
  Therefore, for every irrelevant cumulant with $L (\boldsymbol{a}) > 0$, we
  have $\mathcal{F}^{\boldsymbol{a}}_1 = 0$.
  Using~\eqref{eq:cumulant-flow}, Remark~\ref{rem:norm-scale-monotonicity}
  and~\eqref{eq:bound-ind} we show easily that
  \[ \trinorm{ \mathcal{F}^{\boldsymbol{a}}_{\mu} }_{\mu} \lesssim
     \int_{\mu}^1 \llbracket \sigma \rrbracket^{[\boldsymbol{a}] - 1} \mathd
     \sigma \lesssim \llbracket \mu \rrbracket^{[\boldsymbol{a}]} . \]
  On the other hand, if $\mathcal{F}^{\boldsymbol{a}}$ is relevant,
  {\textit{i.e.}} $[\boldsymbol{a}] < 0$, then on account of the discussion of
  Section~\ref{sec:class-rel-cumulants} $\mathcal{F}^{\boldsymbol{a}}$ is of the
  form $\bar{F}_{\mu}^{[\ell], (1)} =\mathfrak{K}_1 (F_{\mu}^{[\ell],(1)})$
  with $\ell \in \{ 1, \ldots, \hat{\ell} \}$. The treatment of the relevant
  cumulant proceeds via a localisation procedure, which shows that only the
  local part of the cumulant requires renormalisation. However, in our
  fractional setting, where the kernels exhibit only limited polynomial decay,
  this must be done with care. Specifically, we first introduce a preliminary
  truncation before performing a Taylor expansion to localise the fields.
  Consequently, to bound $\bar{F}_{\mu}^{[\ell], (1)}$, we shall use the
  following decomposition
  \begin{equation}
    \begin{aligned}
      \bar{F}_{\mu}^{[\ell], (1)} & = \bar{F}_1^{[\ell], (1)} - \int_{\mu}^1
      [(1 - K_{\sigma}^{1, 1}) \dot{\bar{F}}_{\sigma}^{[\ell], (1)}]\,\mathd
      \sigma \\ &- \int_{\mu}^1 [K^{1, 1}_{\sigma} \dot{\bar{F}}_{\sigma}^{[\ell],
      (1)}] \cdot (1 - h_{\mu})\,\mathd \sigma\\
      & + \boldsymbol{L} \int_{\mu}^1 [K^{1, 1}_{\sigma}
      \dot{\bar{F}}_{\sigma}^{[\ell], (1)}] \cdot (1 - h_{\mu})\,\mathd \sigma
      \\&- \boldsymbol{L} \int_{\mu}^1 [K^{1, 1}_{\sigma}
      \dot{\bar{F}}_{\sigma}^{[\ell], (1)}]\,\mathd \sigma\\
      & - \boldsymbol{R} \int_{\mu}^1 [K^{1, 1}_{\sigma}
      \dot{\bar{F}}_{\sigma}^{[\ell], (1)}] \cdot h_{\mu} \mathd \sigma .
    \end{aligned} \label{eq:flow-decomposition}
  \end{equation}
  This decomposition implements the renormalisation procedure in stages. It
  separates the part not captured by the regularising kernel, the contribution
  away from the diagonal, the local Taylor component, the renormalised local
  coefficient, and the Taylor remainder near the diagonal. The nonlocal and
  remainder terms gain either a power of the ratio of scales or a power of the
  distance from the diagonal, while the purely local coefficient is fixed by
  the renormalisation condition.
  Here
  \[ \dot{\bar{F}}_{\sigma}^{[\ell], (1)} \assign \partial_{\sigma}
     \bar{F}_\sigma^{[\ell], (1)}, \]
  the operator $\boldsymbol{L}$ is a localisation operator and
  $\boldsymbol{R}$ is a remainder operator, both defined in
  Appendix~\ref{app:localization}. The
  polynomial weight $h_{\mu}$ was introduced in
  Definition~\ref{def:weights-kernels}. Recall that, by
  Definition~\ref{def:regularizing-kernels-for-kernels}, in the present case
  $K_{\sigma}^{\boldsymbol{a}} = K^{1, 1}_{\sigma} = 1 \otimes K_{\sigma}$. We
  impose the following renormalisation condition
  \begin{equation}\label{eq:ren-cond} \bar{F}_1^{[\ell], (1)} =  \boldsymbol{L} \int_{1 / 2}^1 [K^{1,
     1}_{\sigma} \dot{\bar{F}}_{\sigma}^{[\ell], (1)}]\,\mathd \sigma . \end{equation}
   Since
  \[
  \bar F^{[\ell],(1)}_1(z,z_1)
  =r^{[\ell]}_{\varepsilon,M}\delta(z-z_1),
  \] the above condition fixes the
  counterterm $r^{[\ell]}_{\varepsilon, M}$.  
  The renormalisation condition is non-circular. Indeed, by the triangular
  structure of the cumulant flow,
  \[
  \partial_\sigma \bar F_\sigma^{[\ell],(1)}
  \]
  depends only on cumulants of total level $\ell-1$. Hence, at the $\ell$-th
  step of the induction, it depends only on the counterterms
  $r_{\varepsilon,M}^{[1]},\ldots,r_{\varepsilon,M}^{[\ell-1]}$, which have
  already been fixed, and is independent of $r_{\varepsilon,M}^{[\ell]}$.
  Consequently, the right-hand side of the renormalisation condition is already
  known and uniquely determines the counterterm
  $r_{\varepsilon,M}^{[\ell]}$.
  
  Note that the integral in~\eqref{eq:ren-cond} is finite for every $\varepsilon\in(0,1)$.
  Indeed, let $\boldsymbol a$ be the one-element tuple whose unique entry is
  $(\ell,1)$. Applying $\boldsymbol L K^{1,1}_\sigma$ to the cumulant flow
  equation~\eqref{eq:cumulant-flow} and using
  Lemma~\ref{lem:local-B-cancellation}, we obtain
  \[
    \boldsymbol L
    (K^{1,1}_\sigma\dot{\bar F}^{[\ell],(1)}_\sigma)
    =-
    \boldsymbol L\left(
      K^{1,1}_\sigma
      \sum_{\boldsymbol b}
      \mathcal A^{\boldsymbol a}_{\boldsymbol b}
      (\dot G_\sigma,\mathcal F^{\boldsymbol b}_\sigma)
    \right).
  \]
  The kernel on the right-hand side is invariant under spatial reflections.
  Consequently, by~\eqref{eq:L-operator}, as in the proof of
  Lemma~\ref{lem:estimates-local}, its local projection contains only the
  zeroth moment. In particular, writing $V_{\boldsymbol b}$ for the summands,
  \[
    \left\|\boldsymbol L\sum_{\boldsymbol b}V_{\boldsymbol b}\right\|
    \leqslant
    \left\|\sum_{\boldsymbol b}V_{\boldsymbol b}\right\|
    \leqslant
    \sum_{\boldsymbol b}\|V_{\boldsymbol b}\|.
  \]
  Since $n(\boldsymbol a)=1$, the basic cumulant norm
  \eqref{eq:basic-norm-cumulant} coincides with the kernel
  norm~\eqref{eq:basic-norm}.
  Combining the ordinary estimate~\eqref{eq:bound-ind}, the
  unweighted estimate~\eqref{eq:A-epsilon-bound-corrected}, and the induction
  hypothesis~\eqref{eq:induction}, we conclude that
  \begin{equation}
    \left\|
      \boldsymbol L
      (K^{1,1}_\sigma\dot{\bar F}^{[\ell],(1)}_\sigma)
    \right\|
    \lesssim
    \llbracket\sigma\rrbracket^{[\boldsymbol a]-1}
    \wedge
    \varepsilon^{-\flat}
    \llbracket\sigma\rrbracket^{[\boldsymbol a]+\flat-1}.
    \label{eq:expected-epsilon}
  \end{equation}
  Notice that~\eqref{eq:expected-epsilon} is asserted only for the local
  projection. Such an estimate may generally be false for the full kernel
  $K^{1,1}_\sigma\dot{\bar F}^{[\ell],(1)}_\sigma$, because its
  $\mathcal B$-contribution need not have a sublattice total variation gain.
  Consequently,
  the renormalisation condition~\eqref{eq:ren-cond} and~\eqref{eq:expected-epsilon} yield
  \begin{align*}
    |r^{[\ell]}_{\varepsilon,M}|
    &\leqslant
    \int_{1/2}^1
    \left\|
      \boldsymbol L
      (K^{1,1}_\sigma\dot{\bar F}^{[\ell],(1)}_\sigma)
    \right\|\,\mathd\sigma
    \\
    &\lesssim
    \int_0^\varepsilon
    \varepsilon^{-\flat}
    \llbracket\sigma\rrbracket^{[\boldsymbol a]+\flat-1}
    \,\mathd\llbracket\sigma\rrbracket
    +
    \int_\varepsilon^{1/2}
    \llbracket\sigma\rrbracket^{[\boldsymbol a]-1}
    \,\mathd\llbracket\sigma\rrbracket
    \\
    &\lesssim \varepsilon^{[\boldsymbol a]},
  \end{align*}
  where we used
  \begin{equation}
    [\boldsymbol a]+\flat
    \geqslant
    \flat-\varrho+\theta+\delta+\beta
    =\flat-2\beta+\delta
    =\flat-2s+2\delta+\kappa>0
    \label{eq:proof-a-bound}
  \end{equation}
  by~\eqref{eq:homogeneity-cum},~\eqref{eq:values-param}
  and~\eqref{eq:choice-of-flat}. 
  
  Taking into account the renormalisation
  condition~\eqref{eq:ren-cond}, we obtain
  \[ \begin{aligned}
       \trinorm{ \bar{F}_{\mu}^{[\ell], (1)} }_{\mu} & \leqslant
       \int_{\mu}^1 \trinorm{ (1 - K^{1, 1}_{\sigma})
       \dot{\bar{F}}_{\sigma}^{[\ell], (1)} }_{\mu} \mathd \sigma +
       \int_{\mu}^1 \trinorm{ [K^{1, 1}_{\sigma}
       \dot{\bar{F}}_{\sigma}^{[\ell], (1)}] \cdot (1 - h_{\mu})
       }_{\mu} \mathd \sigma\\
       &   + \int_{\mu}^1 \trinorm{ \boldsymbol{L} ([K^{1, 1}_{\sigma}
       \dot{\bar{F}}_{\sigma}^{[\ell], (1)}] \cdot (1 - h_{\mu}))
       }_{\mu} \mathd \sigma
          + \int_{1 / 2}^{\mu} \trinorm{
       \boldsymbol{L} (K^{1, 1}_{\sigma} \dot{\bar{F}}_{\sigma}^{[\ell], (1)})
       }_{\mu} \mathd \sigma\\
       &   + \int_{\mu}^1 \trinorm{ \boldsymbol{R} ([K^{1, 1}_{\sigma}
       \dot{\bar{F}}_{\sigma}^{[\ell], (1)}] \cdot h_{\mu}) }_{\mu}
       \mathd \sigma .
     \end{aligned} \]
  Finally, applying
  Lemmas~\ref{lem:estimates-1-h},~\ref{lem:estimates-1-K},~\ref{lem:estimates-local}~and~\ref{lem:estimates-remainder},
  we arrive at
  \[ \begin{aligned}
       \trinorm{ \bar{F}_{\mu}^{[\ell], (1)} }_{\mu} & \lesssim
       \llbracket \mu \rrbracket^{- 2 s} \int_{\mu}^1 \llbracket \sigma
       \rrbracket^{[\boldsymbol{a}] + 2 s - 1} \mathd \sigma + \llbracket \mu
       \rrbracket^{- \flat} \int_{\mu}^1 \llbracket \sigma
       \rrbracket^{[\boldsymbol{a}] + \flat - 1} \mathd \sigma\\
       &   + \llbracket \mu \rrbracket^{- \flat} \int_{\mu}^1 \llbracket
       \sigma \rrbracket^{[\boldsymbol{a}] + \flat - 1} \mathd \sigma + \int_{1
       / 2}^{\mu} \llbracket \sigma \rrbracket^{[\boldsymbol{a}] - 1} \mathd
       \sigma\\
       &   + \llbracket \mu \rrbracket^{- \flat} \int_{\mu}^1 \llbracket \sigma
       \rrbracket^{[\boldsymbol{a}] + \flat-1} \mathd \sigma\\
       & \lesssim  \llbracket \mu \rrbracket^{[\boldsymbol{a}]},
     \end{aligned} \]
  where we used again~\eqref{eq:proof-a-bound}. This concludes the proof of
  the induction step.
\end{proof}

\subsection{Local estimates for the flow of
kernels}\label{sec:localization-of-random-kernels}

In this section, we prove an auxiliary result showing that the norms of the
effective force kernels $F^{\mathfrak{a}}_{\mu}$ can be controlled in terms of
the norms of the localised kernels $F^{\mathfrak{a}}_{\mu} \cdot
v^{\mathfrak{a}}_{\mu}$ restricted to a neighbourhood of the diagonal,
where $v^{\mathfrak{a}}_{\mu}$ is the weight introduced in
Definition~\ref{def:weights-kernels}. We use the notation
\[ \left\| X \right\|_{L_{\mathbb{P}}^N} \assign (\mathbb{E} | X |^N)^{1
   / N} . \]
\begin{lemma}
  \label{lem:nonrelevant-final-data}
  For every $\mathfrak{a} \in \mathfrak{A}$,
  \[
    [\mathfrak{a}] \geqslant 0
    \quad \Longrightarrow \quad
    F_1^{\mathfrak{a}} = 0.
  \]
\end{lemma}

\begin{proof}
  By~\eqref{eq:initial-cond-induction-lemma}, the only nonzero final kernels at
  level $\ell = 0$ are those with $k = 0, 1, 3$. Their homogeneities are
  \[
    [0, 0] = - \alpha < 0,
    \qquad
    [0, 1] = - \alpha + \beta = - 2 \beta - \kappa < 0,
    \qquad
    [0, 3] = - \alpha + 3 \beta = - \kappa < 0.
  \]
  At positive level, the only nonzero final kernels are
  \[
    F_1^{[\ell], (1)}(z,z_1) =
    r_{\varepsilon, M}^{[\ell]} \delta(z-z_1),
    \qquad 1 \leqslant \ell \leqslant \hat{\ell}.
  \]
  Their kernel homogeneity is
  \[
    [\ell, 1] = - \alpha + \delta \ell + \beta
    = \delta \ell - 2 \beta - \kappa
    < - \kappa < 0.
  \]
  For $\ell > \hat{\ell}$, the final datum is zero. Therefore every kernel
  with nonnegative homogeneity has zero final datum.
\end{proof}

\begin{lemma}
  \label{lem:loc-kernels}For every \(\ell\in\{0,\ldots,\bar\ell\}\) and \(N\in\mathbb N_+\), we have
  \begin{equation}
    \begin{aligned}
      &\sup_{\mathfrak{a}\mid \ell(\mathfrak{a}) = \ell}
      \sup_{\mu \in (1 / 2, 1)}
      \left\| \llbracket \mu \rrbracket^{- [\mathfrak{a}] - \kappa}
      \| F^{\mathfrak{a}}_{\mu} \|_{\mu} \right\|_{L_{\mathbb{P}}^N}\\
      &\quad \lesssim_{\ell} 1
      + \sup_{\substack{\mathfrak{a}\mid
      \ell(\mathfrak{a}) \leqslant \ell\\ [\mathfrak{a}] < 0}}
      \sup_{\mu \in (1 / 2, 1)}
      \left\| \llbracket \mu \rrbracket^{- [\mathfrak{a}] - \kappa}
      \| F^{\mathfrak{a}}_{\mu} \cdot v^{\mathfrak{a}}_{\mu}
      \|_{\mu} \right\|_{L_{\mathbb{P}}^{2^{\ell} N}}^{2^{\ell}},
    \end{aligned}
    \label{eq:lemma-F-vF}
  \end{equation}
  and
  \begin{equation}
    \begin{aligned}
      &\sup_{\mathfrak{a}\mid \ell(\mathfrak{a}) = \ell}
      \sup_{\mu \in (1 / 2, 1)}
      \left\| \llbracket \mu \rrbracket^{- [\mathfrak{a}] - \kappa + 1}
      \| \partial_{\mu} F^{\mathfrak{a}}_{\mu} \|_{\mu}
      \right\|_{L_{\mathbb{P}}^N}\\
      &\quad \lesssim_{\ell} 1
      + \sup_{\substack{\mathfrak{a}\mid
      \ell(\mathfrak{a}) \leqslant \ell\\ [\mathfrak{a}] < 0}}
      \sup_{\mu \in (1 / 2, 1)}
      \left\| \llbracket \mu \rrbracket^{- [\mathfrak{a}] - \kappa}
      \| F^{\mathfrak{a}}_{\mu} \cdot v^{\mathfrak{a}}_{\mu}
      \|_{\mu} \right\|_{L_{\mathbb{P}}^{2^{\ell} N}}^{2^{\ell}} .
    \end{aligned}
    \label{eq:lemma-dot-F-vF}
  \end{equation}
\end{lemma}

  \begin{proof}
  We shall prove the lemma by induction on $\ell \in \mathbb{N}_0$. First, we
  discuss the case $\ell = 0$ for which we have $\partial_{\mu}
  F^{\mathfrak{a}}_{\mu} = 0$ and thus $F^{\mathfrak{a}}_{\mu} =
  F^{\mathfrak{a}}_1$ for all $\mu \in (1 / 2, 1)$. Observe that the kernels
  $F^{\mathfrak{a}}_1$ are local, that is, they are supported on the diagonal,
  as can be seen from~\eqref{eq:force}. Hence,
  \[ F^{\mathfrak{a}}_{\mu} \cdot v^{\mathfrak{a}}_{\mu} = F^{\mathfrak{a}}_1
     \cdot v^{\mathfrak{a}}_{\mu} = F^{\mathfrak{a}}_1 =
     F^{\mathfrak{a}}_{\mu}, \qquad \ell (\mathfrak{a}) = 0, \]
  for all $\mu \in (1 / 2, 1)$ since $v^{\mathfrak{a}}_{\mu} = 1$ on the
  diagonal. Moreover, for kernels with $[\mathfrak{a}] \geqslant 0$,
  Lemma~\ref{lem:nonrelevant-final-data} gives $F^{\mathfrak{a}}_1 = 0$, hence
  $F^{\mathfrak{a}}_{\mu} = 0$ for all $\mu \in (1 / 2, 1)$. This proves the
  bounds~\eqref{eq:lemma-F-vF} and~\eqref{eq:lemma-dot-F-vF} for $\ell = 0$.
  
  Let $\ell_0 \in \mathbb{N}_0$ and suppose that the
  bounds~\eqref{eq:lemma-F-vF} and~\eqref{eq:lemma-dot-F-vF} are true for all
  $\ell \leqslant \ell_0$. We shall prove the bounds for $\ell = \ell_0 + 1$.
  It follows from~\eqref{eq:flow-kernels} and
  Lemma~\ref{lem:estimate-B-kernels} that
  \begin{align*}
    \llbracket \mu \rrbracket^{- [\mathfrak{a}] - \kappa + 1}
    \| \partial_{\mu} F^{\mathfrak{a}}_{\mu} \|_{\mu}
    &\lesssim \sum_{\mathfrak{b}, \mathfrak{c}}
    \llbracket \mu \rrbracket^{- [\mathfrak{a}] - 2 \kappa + 1}
    \| B^{\mathfrak{a}}_{\mathfrak{b}, \mathfrak{c}}
    (\dot{G}_{\mu}, F^{\mathfrak{b}}_{\mu},
    F^{\mathfrak{c}}_{\mu}) \|_{\mu}\\
    &\lesssim \sum_{\mathfrak{b}, \mathfrak{c}}
    \llbracket \mu \rrbracket^{- [\mathfrak{b}] - \kappa}
    \| F^{\mathfrak{b}}_{\mu} \|_{\mu}
    \llbracket \mu \rrbracket^{- [\mathfrak{c}] - \kappa}
    \left\| F^{\text{} \mathfrak{c}}_{\mu} \right\|_{\mu} .
  \end{align*}
  Consequently, by H{\"o}lder's inequality
  \[ \left\| \llbracket \mu \rrbracket^{- [\mathfrak{a}] - \kappa + 1} \|
     \partial_{\mu} F^{\mathfrak{a}}_{\mu} \|_{\mu} \right\|_{L_{\mathbb{P}}^N}
     \lesssim \sum_{\mathfrak{b}, \mathfrak{c}} \| \llbracket \mu
     \rrbracket^{- [\mathfrak{b}] - \kappa} \left\| F^{\mathfrak{b}}_{\mu} \|_{\mu}
     \right\|_{L_{\mathbb{P}}^{2 N}}  \left\| \llbracket \mu \rrbracket^{-
     [\mathfrak{c}] - \kappa} \left\| F^{\text{} \mathfrak{c}}_{\mu}
     \right\|_{\mu} \right\|_{L_{\mathbb{P}}^{2 N}} . \]
  Since $\ell (\mathfrak{b}) \vee \ell (\mathfrak{c}) \leqslant \ell
  (\mathfrak{a}) - 1 \leqslant \ell_0$ and $\ell = \ell (\mathfrak{a}) =
  \ell_0 + 1$ the above bound together with the induction hypothesis imply the
  bound~\eqref{eq:lemma-dot-F-vF}. The dependence on $\ell$ in the
  bound~\eqref{eq:lemma-dot-F-vF} comes from the estimate on the (finite)
  number of terms in the sum over $\mathfrak{b}, \mathfrak{c}$.
  
  Let $F^{\mathfrak{a}}_{\mu}$ be non-relevant, so that
  $[\mathfrak{a}] \geqslant 0$. Then
  Lemma~\ref{lem:nonrelevant-final-data} yields
  $F^{\mathfrak{a}}_1 = 0$ and therefore
  \[
    F^{\mathfrak{a}}_{\mu} = - \int_{\mu}^1
    \partial_{\eta} F^{\mathfrak{a}}_{\eta}\,\mathd \eta .
  \]
  Using the monotonicity of the kernel norms~\eqref{eq:norm-mon}, we obtain
  \[ \| F^{\mathfrak{a}}_{\mu} \|_{\mu} \lesssim \int_{\mu}^1 \| \partial_{\eta} F^{\mathfrak{a}}_{\eta} \|_{\eta}\,\mathd \eta . \] By Minkowski's inequality, this gives
  \[ 
  \begin{aligned} 
  	\left\| \| F^{\mathfrak{a}}_{\mu} \|_{\mu} \right\|_{L_{\mathbb P}^N} &\lesssim \int_{\mu}^1 \left\| \| \partial_{\eta} F^{\mathfrak{a}}_{\eta} \|_{\eta} \right\|_{L_{\mathbb P}^N} \,\mathd \eta \\ &\lesssim \int_{\mu}^1 \llbracket \eta \rrbracket^{[\mathfrak a]+\kappa-1} \,\mathd \eta\, \sup_{\eta\in(\mu,1)} \llbracket \eta \rrbracket^{-[\mathfrak a]-\kappa+1} \left\| \| \partial_{\eta} F^{\mathfrak{a}}_{\eta} \|_{\eta} \right\|_{L_{\mathbb P}^N}.
  	\end{aligned} 
  \]
  Since $[\mathfrak a]+\kappa\geqslant \kappa>0$, we have
  \[
    \int_{\mu}^1 \llbracket \eta \rrbracket^{[\mathfrak a]+\kappa-1}
    \,\mathd \eta \lesssim \llbracket \mu \rrbracket^{[\mathfrak a]+\kappa}.
  \]
  Consequently,
  \[
    \begin{aligned}
      \llbracket \mu \rrbracket^{-[\mathfrak a]-\kappa}
      \left\| \| F^{\mathfrak{a}}_{\mu} \|_{\mu} \right\|_{L_{\mathbb P}^N}
      &\lesssim \sup_{\eta\in(\mu,1)} \llbracket \eta
      \rrbracket^{-[\mathfrak a]-\kappa+1} \left\| \|
      \partial_{\eta} F^{\mathfrak{a}}_{\eta} \|_{\eta}
      \right\|_{L_{\mathbb P}^N}.
    \end{aligned}
  \]
  Taking the supremum over $\mu \in (1 / 2, 1)$, the required bound for
  non-relevant kernels in~\eqref{eq:lemma-F-vF} follows from
  \eqref{eq:lemma-dot-F-vF}. This argument applies equally to marginal and
  irrelevant kernels in the non-relevant branch.
  
  Let us now consider relevant kernels $F^{\mathfrak{a}}$, for which
  $[\mathfrak{a}] < 0$. If $k (\mathfrak{a}) = 0$, then
  $v^{\mathfrak{a}}_{\mu} F^{\mathfrak{a}}_{\mu} = F^{\mathfrak{a}}_{\mu}$ and
  the bound~\eqref{eq:lemma-F-vF} is trivial. For $k (\mathfrak{a}) > 0$ we
  shall use the following decomposition
  \[
    \begin{aligned}
      F^{\mathfrak{a}}_{\mu} \cdot (1 - v^{\mathfrak{a}}_{\mu}) &=
      F^{\mathfrak{a}}_1 \cdot (1 - v^{\mathfrak{a}}_{\mu}) - \int_{\mu}^1
      \partial_{\eta} F^{\mathfrak{a}}_{\eta} \cdot (1 -
      v^{\mathfrak{a}}_{\mu})\,\mathd \eta\\
      & = - \int_{\mu}^1 [(1 - \tilde{K}^{\mathfrak{a}}_{\eta})
      \partial_{\eta} F^{\mathfrak{a}}_{\eta}] \cdot (1 -
      v^{\mathfrak{a}}_{\mu})\,\mathd \eta - \int_{\mu}^1
      [\tilde{K}^{\mathfrak{a}}_{\eta} \partial_{\eta}
      F^{\mathfrak{a}}_{\eta}] \cdot (1 - v^{\mathfrak{a}}_{\mu})\,\mathd \eta
      .
    \end{aligned}
  \]
  Note that $F^{\mathfrak{a}}_1 \cdot (1 -
  v^{\mathfrak{a}}_{\mu}) = 0$, which follows from the fact that the kernels
  $F_1^{\mathfrak{a}}$ are local, that is, they are supported on the diagonal, and
  $v^{\mathfrak{a}}_{\mu} = 1$ on the diagonal. Then
  \begin{align*}
     &\| F^{\mathfrak{a}}_{\mu} - F^{\mathfrak{a}}_{\mu} \cdot
     v^{\mathfrak{a}}_{\mu} \|_{\mu}\\
     &\qquad \leqslant \int_{\mu}^1 \| [(1 -
     \tilde{K}^{\mathfrak{a}}_{\eta}) \partial_{\eta} F^{\mathfrak{a}}_{\eta}]
     \cdot (1 - v^{\mathfrak{a}}_{\mu}) \|_{\mu} \mathd \eta\\
     &\qquad\quad + \int_{\mu}^1 \|
     [\tilde{K}^{\mathfrak{a}}_{\eta} \partial_{\eta} F^{\mathfrak{a}}_{\eta}]
     \cdot (1 - v^{\mathfrak{a}}_{\mu}) \|_{\mu} \mathd \eta .
  \end{align*}
  We now observe that, on account of Lemma~\ref{lem:estimates-1-v} and
  Lemma~\ref{lem:estimates-1-K-coeff}, we have
  \begin{align*}
     &\| F^{\mathfrak{a}}_{\mu} - F^{\mathfrak{a}}_{\mu} \cdot
     v^{\mathfrak{a}}_{\mu} \|_{\mu}\\
     &\qquad \leqslant \int_{\mu}^1 \left(
     \frac{\llbracket \eta \rrbracket^{2 s}}{\llbracket \mu \rrbracket^{2 s}}
     + \frac{\llbracket \eta \rrbracket^{\flat - \ell
     (\mathfrak{a}) \kappa_{\mathfrak{o}}}}{\llbracket \mu \rrbracket^{\flat -
     \ell (\mathfrak{a}) \kappa_{\mathfrak{o}}}} \right)  \| \partial_{\eta}
     F^{\mathfrak{a}}_{\eta} \|_{\eta} \mathd \eta\\
     &\qquad \lesssim \int_{\mu}^1
     \frac{\llbracket \eta \rrbracket^{\flat - \ell
     (\mathfrak{a}) \kappa_{\mathfrak{o}}}}{\llbracket \mu \rrbracket^{\flat -
     \ell (\mathfrak{a}) \kappa_{\mathfrak{o}}}}  \| \partial_{\eta}
     F^{\mathfrak{a}}_{\eta} \|_{\eta} \mathd \eta .
  \end{align*}
  Consequently, by Minkowski's inequality
  \[ \left\| \| F^{\mathfrak{a}}_{\mu} - F^{\mathfrak{a}}_{\mu} \cdot
     v^{\mathfrak{a}}_{\mu} \|_{\mu} \right\|_{L_{\mathbb{P}}^N} \lesssim
     \int_{\mu}^1 \frac{\llbracket \eta \rrbracket^{\flat -
     \ell (\mathfrak{a}) \kappa_{\mathfrak{o}}}}{\llbracket \mu
     \rrbracket^{\flat - \ell (\mathfrak{a}) \kappa_{\mathfrak{o}}}}
     \left\| \| \partial_{\eta} F^{\mathfrak{a}}_{\eta} \|_{\eta}
     \right\|_{L_{\mathbb{P}}^N}
     \mathd \eta . \]
  Note that by \eqref{eq:homogeneity-of-kernels}, \eqref{eq:values-param},
  \eqref{eq:choice-of-flat} and~\eqref{eq:cond-kappa-o},
  \begin{align*}
    [\mathfrak{a}] + \flat - \ell(\mathfrak{a})\kappa_{\mathfrak{o}}
    &\geqslant \flat - \alpha
    + (\delta - \kappa_{\mathfrak{o}})\ell(\mathfrak{a}) + \beta\\
    &= \flat + (\delta - \kappa_{\mathfrak{o}})\ell(\mathfrak{a})
    - 2 \beta - \kappa\\
    &= (\delta - \kappa_{\mathfrak{o}})\ell(\mathfrak{a})
    + \delta + \flat - 2 s > 0,
  \end{align*}
  for $k (\mathfrak{a}) > 0$ and $\ell (\mathfrak{a}) > 0$. Thus, the previous
  estimate and the bound~\eqref{eq:lemma-dot-F-vF} imply the bound
  \eqref{eq:lemma-F-vF} for the relevant kernels. This concludes the proof.
\end{proof}

\subsection{From cumulants to random
kernels}\label{sec:concl-of-cumul-analysis}

The bound~\eqref{eq:stochastic-estimates-kernels} stated in
Theorem~\ref{thm:effective-force} is proved by combining the cumulant
estimates from Lemma~\ref{lem:cumul-est-1} with the auxiliary bounds
established in Lemma~\ref{lem:loc-kernels}, via the following lemma, which
shows how estimates on the cumulants translate into estimates on the effective
force kernels. The remaining part of Theorem~\ref{thm:effective-force},
namely, the proof that the estimates~\eqref{eq:est-F} hold with $C_F = 1 + \|
F^{\mathfrak{A}} \|^2$, is presented in Section \ref{sec:post-processing}.

Our analysis relies on the decomposition of the effective force kernels into
mean and fluctuation parts:
\[ F^{\mathfrak{a}}_{\mu} = {\bar{F}_{\mu}^{\mathfrak{a}}}  +
   \tilde{F}_{\mu}^{\mathfrak{a}}, \qquad \qquad
   {\bar{F}_{\mu}^{\mathfrak{a}}}  \assign \mathbb{E}F^{\mathfrak{a}}_{\mu},
   \qquad \tilde{F}_{\mu}^{\mathfrak{a}} \assign F^{\mathfrak{a}}_{\mu}
   -\mathbb{E}F^{\mathfrak{a}}_{\mu} . \]
\begin{lemma}
  \label{lem:kolmogorov} Let $(F_{\mu}^{(\varepsilon,M)})_{\mu \in
  (1 / 2, 1)}$ be the solution of the approximate flow equation with final
  condition~\eqref{eq:initial-cond-induction-lemma} with $r_{\varepsilon, M}$
  as in Lemma \ref{lem:cumul-est-1}. For every $N \in \mathbb{N}_+$, we have
  \begin{equation}
    \mathbb{E} \| F^{\mathfrak{A}} \|^N = \left\| \| F^{\mathfrak{A}} \|
    \right\|^N_{L^N_{\mathbb{P}}} \lesssim 1, \label{eq:bound-FA}
  \end{equation}
  where $\| F^{\mathfrak{A}} \|$ was defined in~\eqref{eq:norm-family-F}.
\end{lemma}

\begin{proof}
  If $[\mathfrak{a}] \geqslant 0$, then
  $F_1^{\mathfrak{a}} = 0$ by Lemma~\ref{lem:nonrelevant-final-data}. Hence
  \[
    F^{\mathfrak{a}}_{\mu}
    = - \int_{\mu}^1 \partial_{\eta} F^{\mathfrak{a}}_{\eta}\,\mathd \eta
  \]
  and
  \begin{align*}
     \llbracket \mu \rrbracket^{- [\mathfrak{a}]}
     \| F^{\mathfrak{a}}_{\mu} \|_{\mu}
     &\lesssim \llbracket \mu \rrbracket^{- [\mathfrak{a}]}
     \int_{\mu}^1 \| \partial_{\eta} F^{\mathfrak{a}}_{\eta} \|_{\mu}
     \mathd \eta\\
     &\lesssim \int_{\mu}^1
     \llbracket \eta \rrbracket^{- [\mathfrak{a}]}
     \| \partial_{\eta} F^{\mathfrak{a}}_{\eta} \|_{\eta} \mathd \eta\\
     &\lesssim \int_{1 / 2}^1
     \llbracket \eta \rrbracket^{- [\mathfrak{a}]}
     \| \partial_{\eta} F^{\mathfrak{a}}_{\eta} \|_{\eta} \mathd \eta .
  \end{align*}
  On the other hand, if the kernel $F^{\mathfrak{a}}_{\mu}$ is relevant, that
  is $[\mathfrak{a}] < 0$, then $\llbracket \mu \rrbracket^{- [\mathfrak{a}]} = 0$ at $\mu =
  1$. Consequently, we have $\llbracket \mu \rrbracket^{- [\mathfrak{a}]}
  F^{\mathfrak{a}}_{\mu} = 0$ at $\mu = 1$. From
  \begin{equation}
    \partial_{\mu} (\llbracket \mu \rrbracket^{- [\mathfrak{a}]}
    F^{\mathfrak{a}}_{\mu}) = \llbracket \mu \rrbracket^{- [\mathfrak{a}]}
    \partial_{\mu} F^{\mathfrak{a}}_{\mu} + [\mathfrak{a}] 
    \llbracket \mu \rrbracket^{- [\mathfrak{a}] - 1} F^{\mathfrak{a}}_{\mu},
    \label{eq:F-rho-derivative}
  \end{equation}
  we deduce that
  \begin{align*}
    \llbracket \mu \rrbracket^{- [\mathfrak{a}]}
    \| F^{\mathfrak{a}}_{\mu} \|_{\mu}
    &\lesssim \int_{\mu}^1
    \bigl(\llbracket \eta \rrbracket^{- [\mathfrak{a}]}
    \| \partial_{\eta} F^{\mathfrak{a}}_{\eta} \|_{\mu}
    + \llbracket \eta \rrbracket^{- [\mathfrak{a}] - 1}
    \| F^{\mathfrak{a}}_{\eta} \|_{\mu}\bigr)\,\mathd \eta\\
    &\lesssim \int_{1 / 2}^1
    \bigl(\llbracket \eta \rrbracket^{- [\mathfrak{a}]}
    \| \partial_{\eta} F^{\mathfrak{a}}_{\eta} \|_{\eta}
    + \llbracket \eta \rrbracket^{- [\mathfrak{a}] - 1}
    \| F^{\mathfrak{a}}_{\eta} \|_{\eta}\bigr)\,\mathd \eta .
  \end{align*}
  Next, Minkowski's inequality gives
  \begin{align*}
    &\left\| \sup_{\mu \in (1 / 2, 1)}
    \llbracket \mu \rrbracket^{- [\mathfrak{a}]}
    \| F^{\mathfrak{a}}_{\mu} \|_{\mu} \right\|_{L^N_{\mathbb{P}}}\\
    &\quad \lesssim \int_{1 / 2}^1
    \left\| \llbracket \eta \rrbracket^{- [\mathfrak{a}]}
    \| \partial_{\eta} F^{\mathfrak{a}}_{\eta} \|_{\eta}
    \right\|_{L^N_{\mathbb{P}}}\,\mathd \eta
    + \int_{1 / 2}^1 \left\| \llbracket \eta \rrbracket^{- [\mathfrak{a}] - 1}
    \| F^{\mathfrak{a}}_{\eta} \|_{\eta}
    \right\|_{L^N_{\mathbb{P}}}\,\mathd \eta,
  \end{align*}
  for every kernel. Since $\int_{1 / 2}^1 \llbracket
  \eta \rrbracket^{\kappa - 1} \mathd \eta < \infty$, we arrive at
  \begin{align*}
    \left\| \sup_{\mu \in (1 / 2, 1)} \llbracket \mu \rrbracket^{- [\mathfrak{a}]}
    \| F^{\mathfrak{a}}_{\mu} \|_{\mu} \right\|_{L^N_{\mathbb{P}}}
    &\lesssim \sup_{\mu \in (1 / 2, 1)}
    \left\| \llbracket \mu \rrbracket^{- [\mathfrak{a}] + 1 - \kappa}
    \| \partial_{\mu} F^{\mathfrak{a}}_{\mu} \|_{\mu}
    \right\|_{L^N_{\mathbb{P}}}\\
    &\quad + \sup_{\mu \in (1 / 2, 1)}
    \left\| \llbracket \mu \rrbracket^{- [\mathfrak{a}] - \kappa}
    \| F^{\mathfrak{a}}_{\mu} \|_{\mu}
    \right\|_{L^N_{\mathbb{P}}} .
  \end{align*}
  Using Lemma~\ref{lem:loc-kernels} and the decomposition
  $F_{\mu}^{\mathfrak{a}} = \bar{F}_{\mu}^{\mathfrak{a}} +
  \tilde{F}_{\mu}^{\mathfrak{a}}$ we obtain
  \begin{align*}
    &\left\| \sup_{\mu \in (1 / 2, 1)} \llbracket \mu \rrbracket^{- [\mathfrak{a}]}
    \| F^{\mathfrak{a}}_{\mu} \|_{\mu} \right\|_{L^N_{\mathbb{P}}}
    \lesssim 
    1 + \sup_{\substack{\mathfrak{b}\mid
    \ell (\mathfrak{b}) \leqslant \ell(\mathfrak{a})\\ [\mathfrak{b}] < 0}}
    \sup_{\mu \in (1 / 2, 1)}
    \left\| \llbracket \mu \rrbracket^{- [\mathfrak{b}] - \kappa}
    \| v^{\mathfrak{b}}_{\mu} F^{\mathfrak{b}}_{\mu}
    \|_{\mu} \right\|_{L_{\mathbb{P}}^{2^{\ell(\mathfrak{a})} N}}^{2^{\ell(\mathfrak{a})}}\\
    &\lesssim 1 + \sup_{\substack{\mathfrak{b}\mid
    \ell (\mathfrak{b}) \leqslant \ell(\mathfrak{a})\\ [\mathfrak{b}] < 0}}
    \sup_{\mu \in (1 / 2, 1)}
    \left( \left\| \llbracket \mu \rrbracket^{- [\mathfrak{b}] - \kappa}
    \| v^{\mathfrak{b}}_{\mu} \tilde{F}^{\mathfrak{b}}_{\mu} \|_{\mu}
    \right\|^{2^{\ell(\mathfrak{a})}}_{L_{\mathbb{P}}^{2^{\ell(\mathfrak{a})} N}}
    +
    \bigl(\llbracket \mu \rrbracket^{- [\mathfrak{b}] - \kappa}
     \| v^{\mathfrak{b}}_{\mu} \bar{F}^{\mathfrak{b}}_{\mu}
     \|_{\mu}\bigr)^{2^{\ell(\mathfrak{a})}} \right) .
  \end{align*}
  Applying Lemmas~\ref{lem:v-and-K} and~\ref{lem:control-of-avgs} and
  recalling that $\trinorm{ \mathcal{F}^{\boldsymbol{A}} } \lesssim 1$, on account of Lemma~\ref{lem:cumul-est-1}, we conclude that
  \[ \left\| \sup_{\sigma \in (1 / 2, 1)}  \llbracket \sigma \rrbracket^{-
     [\mathfrak{a}]}  \| F^{\mathfrak{a}}_{\sigma} \|_{\sigma}
     \right\|_{L^N_{\mathbb{P}}} \lesssim 1 . \]
  Finally, to establish the slightly stronger bound
  \[ \left\| \sup_{\sigma \in (1 / 2, 1)}  \llbracket \sigma \rrbracket^{d / 2  + s
     + 2 \kappa}  \| F^{[0], (0)}_{\sigma} \|_{\sigma} \right\|_{L^N_{\mathbb{P}}} \lesssim 1
  \]
  for the noise $F^{[0], (0)}_{\sigma} = \xi^{(\varepsilon, M)}$, we proceed as before, and invoke the
  second bound stated in Lemma~\ref{lem:v-and-K}.
\end{proof}

The following factorisation is chosen to work uniformly above and below the
lattice scale. When $\llbracket\mu\rrbracket\geqslant\varepsilon$, both the
temporal and spatial variables can be smoothed as in the continuum. When
$\llbracket\mu\rrbracket<\varepsilon$, the spatial structure is atomic and is
therefore retained explicitly in one factor.

\begin{definition}\label{def:K_hat_check}
	Let $\hat\kappa\in[0,1]$. For \(\llbracket\mu\rrbracket\geq\varepsilon\), we set
	\begin{equation}
		\label{eq:partial-kernel-above-lattice}
		\begin{aligned}
			\check K_\mu
			&=
			(1+\llbracket\mu\rrbracket^{2s}\partial_t)^{-\hat\kappa}
			(1-\llbracket\mu\rrbracket^2\Delta_\varepsilon)^{-\hat\kappa},
			\\
			\hat K_\mu
			&=
			(1+\llbracket\mu\rrbracket^{2s}\partial_t)^{\hat\kappa-1}
			(1-\llbracket\mu\rrbracket^2\Delta_\varepsilon)^{\hat\kappa-2},
		\end{aligned}
	\end{equation}
	For \(\llbracket\mu\rrbracket<\varepsilon\), we set
	\begin{equation}
		\label{eq:partial-kernel-below-lattice}
		\begin{aligned}
			\check K_\mu
			&=
			(1+\llbracket\mu\rrbracket^{2s}\partial_t)^{-\hat\kappa}
			(1-\llbracket\mu\rrbracket^2\Delta_\varepsilon)^{-2},
			\\
			\hat K_\mu
			&=
			(1+\llbracket\mu\rrbracket^{2s}\partial_t)^{\hat\kappa-1}
			\otimes 1.
		\end{aligned}
	\end{equation}
	Observe that in both regimes the product is exactly
	\begin{equation}
		\label{eq:partial-kernel}
		\check K_\mu\hat K_\mu
		=
		(1+\llbracket\mu\rrbracket^{2s}\partial_t)^{-1}
		(1-\llbracket\mu\rrbracket^2\Delta_\varepsilon)^{-2}
		=
		K_\mu .
	\end{equation}
\end{definition}

\begin{lemma}
  \label{lem:v-and-K}For all $\mathfrak{a} \in \mathfrak{A}$ and $N \in
  \mathbb{N}_+$, we have
  \begin{align*}
    & \sup_{\mu \in (1 / 2, 1)} \left\| \llbracket \mu \rrbracket^{-
    [\mathfrak{a}] - \kappa} \| \tilde{F}^{\mathfrak{a}}_{\mu} \cdot
    v^{\mathfrak{a}}_{\mu} \|_{\mu} \right\|_{L^{2 N}_{\mathbb{P}}} \lesssim_N
    \trinorm{ \mathcal{F}^{\boldsymbol{A}} }, \\
    & \sup_{\mu \in (1 / 2, 1)} \left\| \llbracket \mu
    \rrbracket^{d / 2 + s + \kappa}  \| F^{[0], (0)}_{\mu} \|_{\mu}
    \right\|_{L^{2 N}_{\mathbb{P}}} \lesssim_N 1. 
  \end{align*}
\end{lemma}

\begin{proof}
  Let $k = k (\mathfrak{a})$. By Lemma~\ref{lem:aux-weight-improved}, there exists a
  weight $\tilde{V}^{\mathfrak{a}}_{\mu}$ with support properties analogous to
  those of $v^{\mathfrak{a}}_{\mu}$ such that
  \[ \| \tilde{F}^{\mathfrak{a}}_{ \mu} \cdot v^{\mathfrak{a}}_{\mu} \|_{\mu}
     = \| \mathfrak{o}^{\mathfrak{a}} \cdot [\tilde{K}^{\mathfrak{a}}_{\mu}
     (\tilde{F}^{\mathfrak{a}}_{ \mu} \cdot v^{\mathfrak{a}}_{\mu})] \cdot
     \tilde{w}^{\mathfrak{a}}_{\mu} \| \lesssim  \|
     \mathfrak{o}^{\mathfrak{a}} \cdot [\tilde{K}^{\mathfrak{a}}_{\mu}
     \tilde{F}_{\mu}^{\mathfrak{a}}] \cdot \tilde{V}^{\mathfrak{a}}_{\mu}
     \tilde{w}^{\mathfrak{a}}_{\mu} \| . \]
  Moreover, we have
  \[ \sup_{z \in \Lambda} \int_{\Lambda^k}  \tilde{w}^{\mathfrak{a}}_{\mu} (z,
     Z^{\mathfrak{a}}) \tilde{V}^{\mathfrak{a}}_{\mu} (z, Z^{\mathfrak{a}}) 
     \textrm{} \mathd Z^{\mathfrak{a}} \lesssim \llbracket \mu
     \rrbracket^{2 s k}  (\varepsilon \vee \llbracket \mu \rrbracket)^{d k},
     \qquad Z^{\mathfrak{a}} = (z_1, \ldots, z_k) \in \Lambda^k . \]
  Thus, we have
  \[ \| \tilde{F}^{\mathfrak{a}}_{ \mu} \cdot v^{\mathfrak{a}}_{\mu} \|_{\mu}
     \lesssim \| \mathfrak{o}^{\mathfrak{a}} \cdot
     [\tilde{K}^{\mathfrak{a}}_{\mu} \tilde{F}_{\mu}^{\mathfrak{a}}] \cdot
     \tilde{V}^{\mathfrak{a}}_{\mu} \tilde{w}^{\mathfrak{a}}_{\mu} \| \lesssim
     \llbracket \mu \rrbracket^{2 s k}  (\varepsilon \vee
     \llbracket \mu \rrbracket)^{d k}  \| \mathfrak{o}^{\mathfrak{a}} \cdot
     [\tilde{K}^{\mathfrak{a}}_{\mu} \tilde{F}_{\mu}^{\mathfrak{a}}]
     \|_{L^{\infty}}, \]
  where, here and below, $\| \Cdot \|_{L^{\infty}}$ denotes the supremum
  norm in all the variables. Hence,
  \begin{equation}
    \left\| \| \tilde{F}^{\mathfrak{a}}_{ \mu} \cdot v^{\mathfrak{a}}_{\mu}
    \|_{\mu}  \right\|_{L^{2 N}_{\mathbb{P}}} \lesssim \llbracket
    \mu \rrbracket^{2 s k}  (\varepsilon \vee \llbracket \mu \rrbracket)^{d
    k}  \left\| \| \mathfrak{o}^{\mathfrak{a}} \cdot
    [\tilde{K}^{\mathfrak{a}}_{\mu} \tilde{F}_{\mu}^{\mathfrak{a}}]
    \|_{L^{\infty}}  \right\|_{L^{2 N}_{\mathbb{P}}} . \label{eq:nuF-F}
  \end{equation}
  To analyse the last term, we decompose the operator $K_{\mu}$ as
  \[ K_{\mu} = \check{K}_{\mu}  \hat{K}_{\mu}, \]
  where the kernels $\check K_\mu$ and $\hat K_\mu$ are introduced in Definition~\ref{def:K_hat_check} and $\hat{\kappa} > 0$ is a small constant to be fixed such that
  \[
  0<\hat\kappa<\frac1{2+2\bar k}.
  \]
  By
  Lemma~\ref{lem:check-K-factor},
  \begin{equation}
    \| \check{K}_{\mu} \|_{L^{2 N / (2 N - 1)} (1
    /\mathfrak{o}^{\mathfrak{a}}) } \lesssim \llbracket \mu \rrbracket^{- (d
    + 2 s) / 2 N}, \label{eq:bound-K-check}
  \end{equation}
  provided \[4 N \hat{\kappa} > 3.\]
  Note that
  \[ \tilde{K}^{\mathfrak{a}}_{\mu} = \check{K}^{\otimes (1 +
     k)}_{\mu} \hat{K}^{\otimes (1 + k)}_{\mu} K_{\mu}^{1, k} \]
  by Definitions~\ref{def:regularizing-kernels}
  and~\ref{def:regularizing-kernels-for-kernels}. Hence, by weighted Young's
  inequality, Fubini's theorem and~\eqref{eq:bound-K-check}, we have
  \[ \begin{aligned}
       &\mathbb{E} [\| \mathfrak{o}^{\mathfrak{a}} \cdot
       [\tilde{K}^{\mathfrak{a}}_{\mu} \tilde{F}_{\mu}^{\mathfrak{a}}]
       \|_{L^{\infty}}^{2 N}]\\[3pt]
        &\lesssim \llbracket \mu \rrbracket^{- (d + 2 s) (1 + k)}
       \int_{\Lambda^{1 + k}} \mathbb{E} [(\mathfrak{o}^{\mathfrak{a}} (z) 
       (\hat{K}^{\otimes (1 + k)}_{\mu} K_{\mu}^{1, k}
       \tilde{F}_{\mu}^{\mathfrak{a}}) (z ; Z^{\mathfrak{a}}))^{2 N}]\,\mathd z
       \mathd Z^{\mathfrak{a}}\\[3pt]
       &\lesssim \llbracket \mu \rrbracket^{- (d + 2 s) (1 + k)}
       \int_{\Lambda^{1 + k}} \mathbb{E} [((\hat{K}^{\otimes (1 + k)}_{\mu}
       K_{\mu}^{1, k} \tilde{F}_{\mu}^{\mathfrak{a}}) (z ; Z^{\mathfrak{a}})
       w^{\mathfrak{a}}_{\mu} (z, Z^{\mathfrak{a}}))^{2 N}]{\left(
       \frac{\mathfrak{o}^{\mathfrak{a}} (z)}{w^{\mathfrak{a}}_{\mu} (z,
       Z^{\mathfrak{a}})} \right) }^{2 N} \mathd z \mathd Z^{\mathfrak{a}}\\[3pt]
       &\lesssim \llbracket \mu \rrbracket^{- (d + 2 s)
       (1 + k)}
       \left\|
       \mathbb{E} [([\hat{K}^{\otimes (1 + k)}_{\mu} K_{\mu}^{1, k}
       \tilde{F}_{\mu}^{\mathfrak{a}}]
       \cdot w^{\mathfrak{a}}_{\mu})^{2 N}]
       \right\|_{L^{\infty}},
     \end{aligned} \]
  where we used that,
  \[ \int_{\Lambda^{1 + k}} {\left( \frac{\mathfrak{o}^{\mathfrak{a}}
     (z)}{w^{\mathfrak{a}}_{\mu} (z, Z^{\mathfrak{a}})} \right) }^{2 N} \mathd
     z \mathd Z^{\mathfrak{a}} \lesssim 1, \]
  for $N \in \mathbb{N}_+$ large enough. Thus,
  \begin{equation}
    \left\| \| \mathfrak{o}^{\mathfrak{a}} \cdot [\tilde{K}^{\mathfrak{a}}_{\mu}
    \tilde{F}_{\mu}^{\mathfrak{a}}] \|_{L^{\infty}}  \right\|_{L^{2 N}_{\mathbb{P}}}
    \lesssim \llbracket \mu \rrbracket^{- \frac{(d + 2 s) (1 +
    k)}{2 N} }  \| \mathbb{E} [([\hat{K}^{\otimes (1 +
    k)}_{\mu} K_{\mu}^{1, k} \tilde{F}_{\mu}^{\mathfrak{a}}] \cdot
    w^{\mathfrak{a}}_{\mu})^{2 N}] \|_{L^{\infty}}^{1 / 2 N} .
    \label{eq:steptwo-Fbound}
  \end{equation}
  We need to control the $L^{2 N}_{\mathbb{P}}$ norm of the centred random
  variable
  \[ (\hat{K}^{\otimes (1 + k)}_{\mu} K_{\mu}^{1, k}
     \tilde{F}_{\mu}^{\mathfrak{a}}) (z ; Z^{\mathfrak{a}})
     w^{\mathfrak{a}}_{\mu} (z, Z^{\mathfrak{a}}) . \]
  To this end, we observe that, by induction, it is easy to show that
  $\tilde{F}_{\mu}^{\mathfrak{a}}$ is a polynomial of the Gaussian random
  field $\xi^{(\varepsilon, M)}$. Consequently, by Nelson's hypercontractivity
  estimate, all higher moments can be controlled by the second one. In
  particular, the $2 N$-th moment can be bounded in terms of the second-order
  cumulant
  \[ \mathbb{E} (\tilde{F}_{\mu}^{\mathfrak{a}} \otimes
     \tilde{F}_{\mu}^{\mathfrak{a}}) =\mathfrak{K}_2 (F_{\mu}^{\mathfrak{a}},
     F_{\mu}^{\mathfrak{a}}) =\mathcal{F}^{(\mathfrak{a}, \mathfrak{a})}_{\mu} . \]
  Using the notation
  \[ w_{\mu}^{(\mathfrak{a}, \mathfrak{a})} \assign w_{\mu}^{\mathfrak{a}}
     \otimes w_{\mu}^{\mathfrak{a}}, \qquad
     K^{(\mathfrak{a}, \mathfrak{a})}_{\mu} \assign K_{\mu}^{1, k} \otimes
     K_{\mu}^{1, k}, \]
  we have
  \begin{equation}
    \begin{aligned}
      &\| \mathbb{E} [([\hat{K}^{\otimes (1 + k)}_{\mu} K_{\mu}^{1, k}
      \tilde{F}_{\mu}^{\mathfrak{a}}] \cdot w^{\mathfrak{a}}_{\mu})^{2 N}]
      \|_{L^{\infty}}\\
      &\quad \lesssim_N \| \mathbb{E} [([\hat{K}^{\otimes (1 + k)}_{\mu}
      K_{\mu}^{1, k} \tilde{F}_{\mu}^{\mathfrak{a}}] \cdot
      w^{\mathfrak{a}}_{\mu})^2] \|_{L^{\infty}}^N\\
      &\quad \lesssim \| \mathfrak{K}_2 ([\hat{K}^{\otimes (1 + k)}_{\mu}
      K_{\mu}^{1, k} F_{\mu}^{\mathfrak{a}}] \cdot w^{\mathfrak{a}}_{\mu},
      [\hat{K}^{\otimes (1 + k)}_{\mu} K_{\mu}^{1, k} F_{\mu}^{\mathfrak{a}}]
      \cdot w^{\mathfrak{a}}_{\mu}) \|^N_{L^{\infty}} .
    \end{aligned} \label{eq:hypercontractivity}
  \end{equation}
  Since
  \[ \trinorm{ \mathcal{F}^{(\mathfrak{a}, \mathfrak{a})}_{\mu}
     }_{\mu}  = \trinorm{ [K^{(\mathfrak{a}, \mathfrak{a})}_{\mu}
     \mathcal{F}^{(\mathfrak{a}, \mathfrak{a})}_{\mu}] \cdot
     w_{\mu}^{(\mathfrak{a}, \mathfrak{a})} } , \]
  by Lemma~\ref{lem:kolmogorov-aux} below, we obtain
  \begin{align*}
     &\| \mathfrak{K}_2 ([\hat{K}^{\otimes (1 + k)}_{\mu} K_{\mu}^{1, k}
     F_{\mu}^{\mathfrak{a}}] \cdot w^{\mathfrak{a}}_{\mu},
     [\hat{K}^{\otimes (1 + k)}_{\mu} K_{\mu}^{1, k}
     F_{\mu}^{\mathfrak{a}}] \cdot w^{\mathfrak{a}}_{\mu}) \| _{L^{\infty}}\\
     &\quad \lesssim
     \llbracket \mu \rrbracket^{- 2 s (1 + 2 k)}
     (\varepsilon \vee \llbracket \mu \rrbracket)^{- d (1 + 2 k)}
     \trinorm{
     \mathcal{F}^{(\mathfrak{a}, \mathfrak{a})}_{\mu} }_{\mu}  .
  \end{align*}
  Hence, by~\eqref{eq:hypercontractivity},
  \[ \| \mathbb{E} ([\hat{K}^{\otimes (1 + k)}_{\mu} K_{\mu}^{1, k}
     \tilde{F}_{\mu}^{\mathfrak{a}}] \cdot w^{\mathfrak{a}}_{\mu})^{2 N}
     \|_{L^{\infty}}^{\frac{1}{2 N}} \lesssim \llbracket \mu
     \rrbracket^{- s (1 + 2 k)}  (\varepsilon \vee \llbracket \mu
     \rrbracket)^{- d (1 + 2 k) / 2}  \trinorm{
     \mathcal{F}^{(\mathfrak{a}, \mathfrak{a})}_{\mu} }_{\mu}^{1 / 2}
     . \]
  Overall, using~\eqref{eq:nuF-F},~\eqref{eq:steptwo-Fbound} together with the
  Definition~\eqref{eq:norm-family-cumulants} of $\trinorm{
  \mathcal{F}^{\boldsymbol{A}} }$ and the fact that
  \[ [(\mathfrak{a}, \mathfrak{a})] = - \varrho + 2 (\theta + \alpha) + 2
     [\mathfrak{a}], \]
  we arrive at
  \begin{equation}
    \sup_{\mu} \| \llbracket \mu \rrbracket^{- [\mathfrak{a}] - \kappa} \|
    v^{\mathfrak{a}}_{\mu} \tilde{F}^{\mathfrak{a}}_{\mu} \|_{\mu} \|_{L^{2
    N}_{\mathbb{P}}} \lesssim_N \trinorm{ \mathcal{F}^{\boldsymbol{A}}
    }\, \sup_{\mu}\, \llbracket \mu \rrbracket^{-
    \frac{\varrho}{2} + \theta + \alpha - \frac{d + 2 s}{2} - \frac{(d + 2 s)
    (1 + k)}{2 N} - \kappa} . \label{eq:exp-bound-2}
  \end{equation}
  In order to conclude, we need the right-hand side above to be finite, that
  is, we require
  \begin{equation}
    \theta + \alpha - \frac{\varrho}{2} - \frac{d + 2 s}{2} - \frac{(d + 2 s)
    (1 + \bar{k})}{2 N} - \kappa \geqslant 0.
    \label{eq:wanted-bound}
  \end{equation}
  Given our choices of parameters~\eqref{eq:values-param} this bound has the
  form
  \begin{equation}
    \frac{3}{4} \delta_{\star} - \frac{(d + 2 s)
    (1 + \bar{k})}{2 N} - \frac{3}{2} \kappa \geqslant 0.
    \label{eq:smallness-conditions}
  \end{equation}
  Since $\kappa \in (0, \delta_{\star} / 2)$ by~\eqref{eq:cond-kappa}, it
  suffices to take $N \in \mathbb{N}_+$ large enough to satisfy this
  inequality and conclude our bound. The bound for smaller values of $N \in
  \mathbb{N}_+$ then follows immediately by Jensen's inequality.
  
  The second of the stated bounds, namely the estimate for the noise $F^{[0],
  (0)} = \xi^{(\varepsilon, M)} $, is proved by the same argument as
  for~\eqref{eq:exp-bound-2}, except that we now use the precise estimate
  \[ \trinorm{ \mathcal{F}^{(\mathfrak{a}, \mathfrak{a})}_{\mu}}_{\mu} = 1,\qquad \mathfrak{a}=
  (0, 0), \]
  for the covariance of the noise
  \[ \mathcal{F}^{(\mathfrak{a}, \mathfrak{a})}_{\mu} =\mathfrak{K}_2(F_\mu^{[0], (0)},F_\mu^{[0], (0)})
  =
     \mathfrak{K}_2(\xi^{(\varepsilon, M)} , \xi^{(\varepsilon, M)} ),
      \]
  which follows from~\eqref{eq:est-cov-noise}, instead of the weaker bound
  \[ \trinorm{ \mathcal{F}^{(\mathfrak{a}, \mathfrak{a})}_{\mu}}_{\mu}
     \leqslant 
     \llbracket
     \mu \rrbracket^{[(\mathfrak{a}, \mathfrak{a})]}\,
     \trinorm{ \mathcal{F}^{\boldsymbol{A}} }^2 . \]
  Specifically, in place of~\eqref{eq:exp-bound-2} we obtain
  \begin{align*}
     &\sup_{\mu}\, \left\| \llbracket \mu \rrbracket^{d / 2 + s + \kappa}
     \| F_\mu^{[0], (0)} \|_{\mu} \right\|_{L^{2 N}_{\mathbb{P}}}\\
     &\quad =
     \sup_{\mu}\, \left\| \llbracket \mu \rrbracket^{d / 2 + s + \kappa}
     \| \mathfrak{o}\cdot K_{\mu} F_\mu^{[0], (0)} \| \right\|_{L^{2 N}_{\mathbb{P}}}\\
     &\quad \lesssim \sup_{\mu}\,
     \llbracket \mu \rrbracket^{- \frac{(d + 2 s)}{2 N} + \kappa}
     \lesssim 1,
  \end{align*}
  for $N \in \mathbb{N}_+$ large enough, depending on $\kappa > 0$.
\end{proof}

\begin{lemma}
  \label{lem:kolmogorov-aux}Let $F$ be a random distributional field over $\Lambda^{1+k}$ admitting second moment.
  Assume that the law of $F$ is invariant under spacetime translations and such that
  \begin{equation}
    z \mapsto F (z, z + z_1, \ldots, z + z_k)
    \label{eq:F-periodicity-lem}
  \end{equation}
  is a.s. periodic in space with period $M$ for all $z_1, \ldots, z_k \in
  \Lambda$. For $z, z_1, \ldots, z_{2 k} \in \Lambda$, define
  \begin{align*}
    \mathcal{F} (z, z_1, \ldots, z_{2 k})
    &=\mathfrak{K}_2 (F (0, z_1, \ldots, z_k),
    F (z, z_{k + 1}, \ldots, z_{2 k}))
    .
  \end{align*} 
  For all $\omega \geqslant 0$ and $\hat{\kappa} \in [0, 1 / (2 k + 2))$, we
  have
  \begin{align*}
    &\| \mathfrak{K}_2 ((\hat{K}^{\otimes (1 + k)}_{\mu} F),
    (\hat{K}^{\otimes (1 + k)}_{\mu} F))
    \cdot (w^{(1 + k), \omega}_{\mu}
    \otimes w^{(1 + k), \omega}_{\mu}) \|_{L^{\infty}}\\
    &\quad \lesssim \llbracket \mu \rrbracket^{- 2 s (1 + 2 k)}
    (\varepsilon \vee \llbracket \mu \rrbracket)^{- d (1 + 2 k)}\\
    &\qquad \times \int_{\Lambda_M \times \Lambda^{2 k}}
    w^{(1 + k), \omega}_{\mu} (0, z_1, \ldots, z_k)
    w^{(1 + k), \omega}_{\mu} (z, z_{k + 1}, \ldots, z_{2 k})\\
    &\qquad \times
    | \mathcal{F} (z, z_1, \ldots, z_{2 k}) |
    \,\mathd z \mathd z_1 \ldots \mathd z_{2 k},
  \end{align*}
  uniformly over $F$ and $\mu \in [1 / 2, 1)$, where \(\hat K_\mu\) is the kernel defined in
  Definition~\ref{def:K_hat_check}.
\end{lemma}

\begin{proof}
  We claim that
  \[ \mathfrak{K}_2 ((\hat{K}^{\otimes (1 + k)}_{\mu} F) (0 , z_1, \ldots,
     z_k), (\hat{K}^{\otimes (1 + k)}_{\mu} F) (z , z_{k + 1}, \ldots, z_{2
     k})) = (\mathcal{K} _{\mu} \ast  \mathcal{F}) (z, z_1, \ldots, z_{2 k})
  \]
  for all $z, z_1, \ldots, z_{2 k} \in \Lambda$, where $\ast $ stands for the
  convolution on $\Lambda^{1 + 2 k}$ and
  \[ \mathcal{K} _{\mu} (z, z_1, \ldots, z_{2 k}) \assign \int_{\Lambda}
     \hat{K}_{\mu} (z'')  \hat{K} _{\mu} (z + z'')  \hat{K}_{\mu} (z_1 + z'')
     \ldots \hat{K}_{\mu} (z_{2 k} + z'') \: \mathd z'' . \]
  Indeed, we have
  \begin{align*}
 &\quad (\mathcal{K}_{\mu} \ast  \mathcal{F}) (z, z_1, \ldots, z_{2 k})\\
 &\quad = \int_{\Lambda^{2 + 2 k}} \hat{K}_{\mu} (z'')  \hat{K}_{\mu} (z' +
    z'')  \hat{K}_{\mu} (z'_1 + z'') \ldots \hat{K}_{\mu} (z'_{2 k} + z'')\\
 &\quad \qquad \times \mathcal{F} (z - z', z_1 - z'_1, \ldots, z_{2 k} -
    z'_{2 k})\,\mathd z' \mathd z'' \mathd z'_1 \ldots \mathd z'_{2 k}\\
 &\quad = \int_{\Lambda^{2 + 2 k}} \hat{K}_{\mu} (z'')  \hat{K}_{\mu} (z') 
    \hat{K}_{\mu} (z'_1) \ldots \hat{K}_{\mu} (z'_{2 k})\\
 &\quad \qquad \times \mathfrak{K}_2 \bigl(
    F (0 - z'' , z_1 - z'_1, \ldots, z_k - z'_k),\\
 &\quad \qquad \qquad
    F (z - z', z_{k + 1} - z'_{k + 1}, \ldots, z_{2 k} - z'_{2 k})\bigr)\\
 &\quad \qquad \times
    \mathd z' \mathd z'' \mathd z'_1 \ldots \mathd z'_{2 k}\\
 &\quad =\mathfrak{K}_2 \bigl(
    (\hat{K}^{\otimes (1 + k)}_{\mu} F) (0 , z_1, \ldots, z_k),
    (\hat{K}^{\otimes (1 + k)}_{\mu} F)
    (z , z_{k + 1}, \ldots, z_{2 k})\bigr),
  \end{align*}
  where the second equality above follows from the translation invariance of
  the law of $F$.
  
  It remains to show that
  \begin{align*}
    &| (\mathcal{K}_{\mu} \ast \mathcal{F}) (z, z_1, \ldots, z_{2 k}) |
    w^{(1 + k), \omega}_{\mu} (0, z_1, \ldots, z_k)
    w^{(1 + k), \omega}_{\mu} (z, z_{k + 1}, \ldots, z_{2 k})\\
    &\quad \lesssim \llbracket \mu \rrbracket^{- 2 s (1 + 2 k)}
    (\varepsilon \vee \llbracket \mu \rrbracket)^{- d (1 + 2 k)}\\
    &\qquad \times
    \int_{\Lambda_M \times \Lambda^{2 k}}
    w^{(1 + k), \omega}_{\mu} (0, z_1, \ldots, z_k)
    w^{(1 + k), \omega}_{\mu} (z, z_{k + 1}, \ldots, z_{2 k})\\
    &\qquad \times | \mathcal{F} (z, z_1, \ldots, z_{2 k}) |
    \mathd z \mathd z_1 \ldots \mathd z_{2 k} .
  \end{align*}
  We have
  \begin{align*}
    &| (\mathcal{K}_{\mu} \ast  \mathcal{F}) (z, z_1, \ldots, z_{2 k}) |\\
    &\quad \leqslant \sum_{n \in \mathbb{Z}^d}
    \int_{\Lambda_M \times \Lambda^{2 k}}
    | \mathcal{K}_{\mu} (z' + M n, z'_1, \ldots, z'_{2 k}) |\\
    &\qquad \times
	    | \mathcal{F} (z - z' - M n, z_1 - z'_1, \ldots, z_{2 k} - z'_{2 k}) |
    \mathd z' \mathd z'_1 \ldots \mathd z'_{2 k} .
  \end{align*}
  Thus, by the same argument as in~\eqref{eq:weight-w-2}, we obtain
  \begin{align*}
    &| (\mathcal{K}_{\mu} \ast  \mathcal{F}) (z, z_1, \ldots, z_{2 k}) |
    w^{(1 + k), \omega}_{\mu} (0, z_1, \ldots, z_k)
    w^{(1 + k), \omega}_{\mu} (z , z_{k + 1}, \ldots, z_{2 k})\\
    &\quad \lesssim \sum_{n \in \mathbb{Z}^d}
    \int_{\Lambda_M \times \Lambda^{2 k}}
    w^{\omega}_{\mu} (z' + M n) w^{\omega}_{\mu} (z_1') \cdots
    w^{\omega}_{\mu} (z_{2 k}')
    | \mathcal{K}_{\mu} (z' + M n, z'_1, \ldots, z'_{2 k}) |\\
    &\qquad \times
    w^{(1 + k), \omega}_{\mu} (0, z_1 - z'_1, \ldots, z_k - z'_k)\\
    &\qquad \times
    w^{(1 + k), \omega}_{\mu}
	    (z - z' - M n, z_{k + 1} - z'_{k + 1}, \ldots, z_{2 k} - z'_{2 k})\\
    &\qquad \times
	    | \mathcal{F} (z - z' - M n, z_1 - z'_1, \ldots, z_{2 k}
	    - z'_{2 k}) |\,\mathd z' \mathd z'_1 \ldots \mathd z'_{2 k},
  \end{align*}
  where for $z = (t, x) \in \Lambda$ and $n \in \mathbb{Z}^d$ we write $z + M
  n = (t, x + M n) \in \Lambda$. Using Definition~\ref{def:weights-kernels} we
  arrive at
  \begin{align*}
       &\int_{\Lambda_M \times \Lambda^{2 k}}
       w^{(1 + k), \omega}_{\mu} (0, z_1 - z'_1, \ldots, z_k - z'_k)\\
       &\quad \times
	       w^{(1 + k), \omega}_{\mu} (z - z' - M n,
	       z_{k + 1} - z'_{k + 1}, \ldots, z_{2 k} - z'_{2 k})\\
       &\quad \times
       \left|
	       \mathcal{F} (z - z' - M n, z_1 - z'_1, \ldots, z_{2 k} - z'_{2 k})
       \right|
       \mathd z' \mathd z'_1 \ldots \mathd z'_{2 k}\\
       &= \int_{\Lambda_M \times \Lambda^{2 k}}
       w^{(1 + k), \omega}_{\mu} (0, z'_1, \ldots, z'_k)\,
       w^{(1 + k), \omega}_{\mu} (0, z'_{k + 1}, \ldots, z'_{2 k})\\
       &\quad \times
       \left|
	       \mathcal{F} (z - z' - M n, z'_1, \ldots, z'_k, z - z' - M n + z'_{k + 1}, \ldots,
	       z - z' - M n + z'_{2 k})
       \right|\\
       &\quad \times \mathd z' \mathd z'_1 \ldots \mathd z'_{2 k} .
  \end{align*}
  Note that by~\eqref{eq:F-periodicity-lem}, the above expression can be
  rewritten as
  \begin{align*}
       &\int_{\Lambda_M \times \Lambda^{2 k}}
       w^{(1 + k), \omega}_{\mu} (0, z'_1, \ldots, z'_k)
       w^{(1 + k), \omega}_{\mu} (0, z'_{k + 1}, \ldots, z'_{2 k})\\
       &\quad \times
       | \mathcal{F} (z', z'_1, \ldots, z'_k,
       z' + z'_{k + 1}, \ldots, z' + z'_{2 k}) |
       \mathd z' \mathd z'_1 \ldots \mathd z'_{2 k}\\
       &= \int_{\Lambda_M \times \Lambda^{2 k}}
       w^{(1 + k), \omega}_{\mu} (0, z'_1, \ldots, z'_k)
       w^{(1 + k), \omega}_{\mu} (z', z'_{k + 1}, \ldots, z'_{2 k})\\
       &\quad \times
       | \mathcal{F} (z', z'_1, \ldots, z'_k,
       z'_{k + 1}, \ldots, z'_{2 k}) |
       \mathd z' \mathd z'_1 \ldots \mathd z'_{2 k} .
  \end{align*}
  Consequently,
  \begin{align*}
       &| (\mathcal{K}_{\mu} \ast  \mathcal{F}) (z, z_1, \ldots, z_{2 k}) |
       w^{(1 + k), \omega}_{\mu} (0, z_1, \ldots, z_k)
       w^{(1 + k), \omega}_{\mu} (z , z_{k + 1}, \ldots, z_{2 k})\\
       &\quad \lesssim \sup_{z', z'_1, \ldots, z'_{2 k} \in \Lambda}
       \sum_{n \in \mathbb{Z}^d} w^{\omega}_{\mu} (z' + M n)
       w^{\omega}_{\mu} (z_1') \cdots w^{\omega}_{\mu} (z_{2 k}')\\
       &\qquad \times
       | \mathcal{K}_{\mu} (z' + M n, z'_1, \ldots, z'_{2 k}) |\\
       &\qquad \times \int_{\Lambda_M \times \Lambda^{2 k}}
       w^{(1 + k), \omega}_{\mu}
       (0, z''_1, \ldots, z''_k) w^{(1 + k), \omega}_{\mu}
       (z'', z''_{k + 1}, \ldots, z''_{2 k})\\
       &\qquad \times
       | \mathcal{F} (z'', z''_1, \ldots, z''_{2 k}) |\,
       \mathd z'' \mathd z''_1 \ldots \mathd z''_{2 k} .
  \end{align*}
  To complete the proof it remains to show that
  \begin{align*}
    &\sup_{z', z'_1, \ldots, z'_{2 k} \in \Lambda}
    \sum_{n \in \mathbb{Z}^d}
    w^{\omega}_{\mu} (z' + M n) w^{\omega}_{\mu} (z_1') \cdots
    w^{\omega}_{\mu} (z_{2 k}')\\
    &\quad \times
    | \mathcal{K}_{\mu} (z' + M n, z'_1, \ldots, z'_{2 k}) |
    \lesssim \llbracket \mu \rrbracket^{- 2 s (1 + 2 k)}
    (\varepsilon \vee \llbracket \mu \rrbracket)^{- d (1 + 2 k)}.
  \end{align*}
  
  To this end, first note that by Lemma~\ref{lem:hat-K-factor},
  \[
  \prod_{i=0}^{2k}w_\mu^{\omega+d+1}(z_i)
  \left|\mathcal K_\mu(z_0,\ldots,z_{2k})\right|
  \lesssim
  \llbracket \mu \rrbracket^{-2s(1+2k)}
  (\varepsilon\vee\llbracket \mu \rrbracket)^{-d(1+2k)}.
  \]
  Consequently,
  \begin{align*}
  	&\sum_{n\in\mathbb Z^d}
  	w_\mu^\omega(z'+Mn)
  	w_\mu^\omega(z'_1)\cdots w_\mu^\omega(z'_{2k})
  	\left|\mathcal K_\mu(z'+Mn,z'_1,\ldots,z'_{2k})\right|
  	\\
  	&\quad\lesssim
  	\llbracket \mu \rrbracket^{-2s(1+2k)}
  	(\varepsilon\vee\llbracket \mu \rrbracket)^{-d(1+2k)}
  	\sum_{n\in\mathbb Z^d}w_\mu^{-d-1}(z'+Mn)
  	\\
  	&\quad\lesssim
  	\llbracket \mu \rrbracket^{-2s(1+2k)}
  	(\varepsilon\vee\llbracket \mu \rrbracket)^{-d(1+2k)},
  \end{align*}
  where we discarded the factors
  \(\prod_{i=1}^{2k}w_\mu^{-d-1}(z'_i)\leq1\).
  This proves the required bound, and finishes the proof.  
\end{proof}

\begin{lemma}
  \label{lem:control-of-avgs} For every \(\mathfrak a\in\mathfrak A\), we have
  \[ \llbracket \mu \rrbracket^{- [\mathfrak{a}] - \kappa}  \| 
     \bar{F}^{\mathfrak{a}}_{\mu} \cdot v^{\mathfrak{a}}_{\mu} \|_{\mu}
     \lesssim \trinorm{ {\mathcal{F}^{\boldsymbol{A}}} 
     }, \]
  uniformly in $\mu \in [1 / 2, 1)$.
\end{lemma}

\begin{proof}
  By Lemma~\ref{lem:aux-weight-improved}, we have
  \[
  \|\bar F_\mu^{\mathfrak a}\cdot v_\mu^{\mathfrak a}\|_\mu
  \lesssim
  \|\bar F_\mu^{\mathfrak a}\|_\mu
  =
  \|
  \mathfrak o^{\mathfrak a}
  \cdot
  [\tilde K_\mu^{\mathfrak a}\bar F_\mu^{\mathfrak a}]
  \cdot
  \tilde w_\mu^{\mathfrak a}
  \|.
  \]
  Moreover, by Definition~\ref{def:weights-kernels},
  \[ \mathfrak{o}^{\mathfrak{a}} \tilde{w}_{\mu}^{\mathfrak{a}} \leqslant
     w_{\mu}^{\mathfrak{a}} . \]
  Thus, by Lemma~\ref{lem:bounds-K-J},~\eqref{eq:homogeneity-of-kernels}
  and~\eqref{eq:homogeneity-cum} we obtain
  \[ \begin{aligned}
       \llbracket \mu \rrbracket^{- [\mathfrak{a}] - \kappa}  \|
       \bar{F}_{\mu}^{\mathfrak{a}} \cdot v^{\mathfrak{a}}_{\mu} \|_{\mu} 
       &\lesssim  \llbracket \mu \rrbracket^{- [\mathfrak{a}] - \kappa}  \|
       [\tilde{K}^{\mathfrak{a}}_{\mu} \bar{F}_{\mu}^{\mathfrak{a}}] \cdot
       w_{\mu}^{\mathfrak{a}} \| \\
       & =  \llbracket \mu \rrbracket^{- [\mathfrak{a}] - \kappa}  \|
       [K_{\mu}^{\otimes (1 + k (\mathfrak{a}))} K^{1, k (\mathfrak{a})}_{\mu}
       \bar{F}_{\mu}^{\mathfrak{a}}] \cdot w_{\mu}^{\mathfrak{a}} \| \\
       & \lesssim  \llbracket \mu \rrbracket^{- [\mathfrak{a}] - \kappa}  \|
       K_{\mu} \|_{\operatorname{TV} (w^{\flat}_{\mu})}^{1 + k (\mathfrak{a})}  \|
       [K^{1, k (\mathfrak{a})}_{\mu} \bar{F}_{\mu}^{\mathfrak{a}}] \cdot
       w_{\mu}^{\mathfrak{a}} \| \\
       & \lesssim  \trinorm{ \mathcal{F}^{\boldsymbol{A}} }
       \llbracket \mu \rrbracket^{[\boldsymbol{a}] - [\mathfrak{a}] - \kappa}\\
       & \lesssim  \trinorm{ \mathcal{F}^{\boldsymbol{A}} }
       \llbracket \mu \rrbracket^{\alpha - \varrho + \theta - \kappa} ,
     \end{aligned} \] 
     where $\boldsymbol{a}=(\mathfrak{a})$.
  To conclude we use the inequality
  \begin{equation}
    \alpha - \varrho + \theta - \kappa \geqslant 0,
    \label{eq:constraint-f-bar}
  \end{equation}
  which coincides with the constraint~$[\bar{F}]$ stated
  in~\eqref{eq:summary-ineq}.
\end{proof}

\subsection{Post-processing}\label{sec:post-processing}

It remains to establish the bounds~\eqref{eq:est-F} with
$C_F = 1 + \| F^{\mathfrak{A}} \|^2$ used above in the proof of
Theorem~\ref{thm:effective-force}, based on the analysis of the flow equation
carried out so far. This is the content of
Lemmas~\ref{lem:overall-difference-bound}, \ref{lem:estimates-DF}
and~\ref{lem:packaging} below.

Before continuing with the specific computations leading to~\eqref{eq:est-F}
it will be useful to discuss how to fix the values of the parameters
\[ \gamma, \vartheta, \bar{\kappa} , \bar{\ell}, \bar{k}, a,
   \nu, \kappa_{\mathfrak{o}}, \kappa . \]
The validity of the
Lemmas~\ref{lem:overall-difference-bound}--\ref{lem:difference-bound} below,
which together yield~\eqref{eq:est-F}, depends on a series of conditions on
these parameters,
namely~\eqref{eq:cond-kappabar-one},~\eqref{eq:choice-zeta},~\eqref{eq:xi-on-zeta},~$\eqref{eq:DF-bound},$~\eqref{eq:choice-of-K-and-l}
and~\eqref{eq:one-J-bound} as well as the
constraints~\eqref{eq:thm-constraints} in Theorem~\ref{thm:apriori-system}.
The parameters ${\kappa, \kappa_{\mathfrak{o}}} $ are further constrained by
conditions~\eqref{eq:cond-kappa} and~\eqref{eq:cond-kappa-o}.

For the reader's convenience, these conditions are summarised in the following
table:
\[ \begin{array}{lc}
     {}[\mathbb{A}] & \vartheta \leqslant 3 (\gamma - \beta) + \delta - a
     \kappa_{\mathfrak{o}} (1 + \bar{\ell}) - \kappa\\
     {}[\mathbb{B}] & \vartheta \leqslant 2 \gamma \wedge \left( 3 \gamma -
     \frac{d + 2 s}{2} - 2 \kappa - a \kappa_{\mathfrak{o}} \right)\\
     {}[\mathbb{C}] & \vartheta \leqslant \delta \bar{\ell} - \alpha - a
     \kappa_{\mathfrak{o}} (2 \bar{\ell} + 2)\\
     {}[\mathbb{D}] & \bar{\kappa} \geqslant \bar{k} \nu + (1 + \bar{\ell})
     \kappa_{\mathfrak{o}}\\
     {}[\mathbb{E}] & \vartheta \leqslant \gamma \wedge (2 s - 2 \gamma)\\
     {}[\mathbb{F}] & \vartheta \geqslant \frac{4 (s + \gamma) 
     \bar{\kappa}}{1 - \bar{\kappa}}\\
     {}[\mathbb{G}] & \kappa \in (0, \delta_{\star} / 2), \qquad
     \frac{2 \beta}{\delta}\notin\mathbb N_+\\
     {}[\mathbb{H}] & \kappa_{\mathfrak{o}} \leqslant \delta
     \wedge (2-\flat) / (2 \bar{\ell} + 2)
   \end{array} \]
In order to satisfy all these constraints, together with the basic bounds
$\vartheta > 0$ and $\delta > 0$, we proceed as follows. First, recall that
\[
\beta=s-\frac{\delta_{\star}}{4}-\frac{\kappa}{2},
\qquad
\delta = \frac{\delta_\star}2,
\qquad
\delta_{\star} = \frac{4 s - d}{3} > 0,
\]
and set
\begin{equation}
  \gamma \assign \beta - \vartheta,
  \qquad
  \vartheta \assign \frac{1}{16} \delta_{\star}.
  \label{eq:fix-gamma-zeta}
\end{equation}
Next, we {\underline{fix $\bar{\kappa} \in(0,1/3)$ small enough}} so that
\[ \frac{(8s + \delta_\star)  \bar{\kappa}}{1 - \bar{\kappa}} \leqslant \frac{1}{16}
   \delta_{\star}. \]
Since 
$4(s+\gamma)<8s+\delta_\star$, this allows us to satisfy~$[\mathbb{F}]$. The constraint~$[\mathbb{A}]$ is
satisfied provided
\begin{equation}
  a \kappa_{\mathfrak{o}} (1 + \bar{\ell}) + \kappa \leqslant \frac{1}{4}
  \delta_{\star} . \label{eq:const-kappa-two}
\end{equation}
Since 
\begin{align*}
   &3 \gamma - \frac{d + 2 s}{2} - 2 \kappa - a \kappa_{\mathfrak{o}}\\
   &\qquad = 3(\gamma - \beta) + 3 \beta
   - \frac{d + 2 s}{2} - 2 \kappa - a \kappa_{\mathfrak{o}}\\
   &\qquad = 3(\gamma - \beta) + \frac{3}{4} \delta_{\star}
   - \frac{7}{2} \kappa - a \kappa_{\mathfrak{o}},
\end{align*}
the condition~$[\mathbb{B}]$ is satisfied provided
\[ 4 \vartheta \leqslant \frac{3}{4} \delta_{\star} - \frac{7}{2} \kappa - a
   \kappa_{\mathfrak{o}}, \]
and
\[ 3 \vartheta \leqslant 2 \beta = 2 s - \delta - \kappa = 2 s -
   \frac{\delta_{\star}}{2} - \kappa, \]
that is, whenever
\begin{equation}
  \frac{7}{2} \kappa + a \kappa_{\mathfrak{o}} \leqslant
  \frac{\delta_{\star}}{2}, \qquad \kappa \leqslant 2 s - \frac{11}{16}
  \delta_{\star} . \label{eq:const-kappa-three}
\end{equation}
Note that
\[ 2 s - \frac{11}{16} \delta_{\star} \geqslant 2 s - \delta_{\star} \geqslant
   2 s - \frac{4 s - d}{3} = \frac{2 s + d}{3} > 0, \]
so the second condition in~\eqref{eq:const-kappa-three} can indeed be
satisfied. To meet condition~$[\mathbb{C}]$, we shall later impose
\begin{equation}
	2a\kappa_{\mathfrak{o}}\leqslant \frac{\delta_\star}{4}.
	\label{eq:const-kappa-four}
\end{equation}
Under this constraint,
\[
\left(
\frac{\delta_\star}{2}-2a\kappa_{\mathfrak{o}}
\right)\bar\ell
-2a\kappa_{\mathfrak{o}}
\geqslant
\frac{\delta_\star}{4}\bar\ell
-\frac{\delta_\star}{4}.
\]
We therefore choose $\bar\ell$ large enough so that
\[
\alpha+\frac{1}{16}\delta_\star
\leqslant
\frac{\delta_\star}{4}\bar\ell
-\frac{\delta_\star}{4}.
\]
Then~$[\mathbb{C}]$ follows once~\eqref{eq:const-kappa-four} is imposed.
Given $\bar{\ell}$ as above, we \underline{fix $\bar{k}$ large enough} as discussed
in Section~\ref{sec:random-flow}. Then we \underline{fix $a$ sufficiently large} so that
\[ a \geqslant \bar{k} (2 s - \delta) / \bar{\kappa}
   \geqslant 2 \bar{k} \beta / \bar{\kappa} \geqslant 2 \bar{k} \gamma /
   \bar{\kappa}. \]
Consequently, relation~\eqref{eq:a-nu-gamma} compels us to \underline{set
$\nu = \gamma / a$}. To ensure condition~$[\mathbb{D}]$ we
introduce a~constraint
\begin{equation}
  (1 + \bar{\ell}) \kappa_{\mathfrak{o}} \leqslant \bar{\kappa} / 2.
  \label{eq:const-kappa-one}
\end{equation}
Condition~$[\mathbb{E}]$ is met if
\begin{equation}
	\kappa\leqslant 2s-\frac34\delta_\star
	=s+\frac d4 .
	\label{eq:const-kappa-five}
\end{equation}
Indeed, the inequality $\vartheta\leqslant 2s-2\gamma$ is automatic, since
\[
2s-2\gamma
=
2s-2\beta+2\vartheta
=
\frac{5}{8}\delta_\star+\kappa
>
\vartheta .
\]
It remains to ensure $\vartheta\leqslant\gamma$. Since
$\gamma=\beta-\vartheta$, this is equivalent to
$2\vartheta\leqslant\beta$, which follows from
\[
2\vartheta
=
\frac{\delta_\star}{8}
\leqslant
s-\frac{\delta_\star}{4}-\frac{\kappa}{2}
=
\beta
\]
provided~\eqref{eq:const-kappa-five} holds.
To conclude, we {\underline{fix $\kappa_{\mathfrak{o}}$ and $\kappa$ small
enough}} as follows. We fix $\kappa_{\mathfrak{o}} > 0$ and choose
\[
  \kappa \in (0, \delta_{\star} / 2)
  \setminus
  \{ 2 s - \delta (\ell + 1) : \ell \in \mathbb{N}_+ \}
\]
sufficiently small so that the remaining constraints, namely
$[\mathbb{G}]$ and $[\mathbb{H}]$, together
with~\eqref{eq:const-kappa-two},~\eqref{eq:const-kappa-three},
~\eqref{eq:const-kappa-four},~\eqref{eq:const-kappa-one}
and~\eqref{eq:const-kappa-five}, as well as~\eqref{eq:cond-kappa-o}, are
satisfied.

Let us now detail the analysis. As a first step, we extract the coercive
contribution from the effective force by defining
\begin{equation}
  \begin{aligned}
    Q_{\sigma} (\psi_\sigma) & \assign J_{\sigma} F_{\sigma}
    (\psi_{\sigma}) - (- \lambda \psi^3_{\sigma})\\
    & = J_{\sigma} F_{\sigma} (\psi_{\sigma}) - J_{\sigma} (- \lambda
    \psi^3_{\sigma}) + (1 - J_{\sigma}) (\lambda \psi^3_{\sigma})\\
    & = J_{\sigma} F_{\sigma}^{[> 0]} (\psi_{\sigma}) + J_{\sigma} F^{[0],
    (1)} \psi_{\sigma} + J_{\sigma} F^{[0], (0)} + (1 - J_{\sigma}) (\lambda
    \psi^3_{\sigma}),
  \end{aligned} , \label{eq:difference-decomposition}
\end{equation}
where $\psi$ is a generic field,
\[ \psi_{\sigma} \assign J_{\sigma} \psi, \]
and
\[ F^{[0], (0)} = \xi^{(\varepsilon, M)}, \qquad F ^{[0], (1)} (\psi) = \bar{r}
   \psi, \qquad F_{\sigma}^{[> 0]} (\psi ) \assign \sum_{\mathfrak{a}
   \mid \ell (\mathfrak{a}) > 0} F_{\sigma}^{\mathfrak{a}}
   (\psi^{\otimes k (\mathfrak{a})}) . \]
Recall that the norms $\trinorm{ \Cdot } =
\trinorm{ \Cdot }_{\bar{\mu}}$ and $\trinorm{ \Cdot
}_{\#} = \trinorm{ \Cdot }_{\#, \bar{\mu}}$ depending
on the terminal scale $\bar{\mu} \in [1 / 2, 1)$ were introduced in
Definition~\ref{def:triple-norms}.

\begin{lemma}
  \label{lem:overall-difference-bound}Assume~\eqref{eq:choice-zeta},
  \eqref{eq:xi-on-zeta}, \eqref{eq:one-J-bound} and
  \begin{equation}
    \bar{k} \geqslant 5, \qquad \bar{\kappa} \geqslant \bar{k} \nu + (1 +
    \bar{\ell}) \kappa_{\mathfrak{o}} \label{eq:cond-kappabar-one} .
  \end{equation}
  Then the bounds
  \begin{equation}
      \| \zeta_{\mu}^{\bar{\kappa}} Q_{\sigma} (\psi_{\sigma}) \|  \lesssim 
      \llbracket \sigma \rrbracket^{- 3 \gamma + \vartheta}  [\|
      F^{\mathfrak{A}} \|  (1 + \trinorm{ \psi })^{\bar{k}} + (1 +
      \trinorm{ \psi })^2 \trinorm{(\mathcal L\psi)_\Cdot
      }_{\#}],
\label{eq:bound-Q}
  \end{equation}
  \begin{equation}
  		\left\|
  		\zeta_{\mu}^{\bar{\kappa}} K_{\sigma}
  		\bigl(F_{\sigma}(\psi_{\sigma})
  		+ \lambda \psi_{\sigma}^{3}\bigr)
  		\right\|
  		\lesssim
  		\llbracket \sigma \rrbracket^{-3\gamma+\vartheta}
  		\|F^{\mathfrak A}\|
  		(1+\trinorm{\psi}_{\bar\mu})^{\bar k},
  	\label{eq:K-force-coercive-remainder}
  \end{equation}
  \begin{equation}
      \| \zeta_{\mu}^{\bar{\kappa}} K_{\sigma} F_{\sigma} (\psi_{\sigma}) \| 
      \lesssim  \llbracket \sigma \rrbracket^{- 3 \gamma}  \|
      F^{\mathfrak{A}} \|  (1 + \trinorm{ \psi })^{\bar{k}},
       \label{eq:bound-JtildeF}
  \end{equation}
  hold uniformly in $\bar{\mu}, \mu \in [1 / 2, 1),$ $\sigma \in [\mu \vee
  \bar{\mu}, 1) $ and $\psi \in \mathcal{S}' (\Lambda)$.
\end{lemma}

\begin{proof}
  First, observe that, by the triangle inequality applied to the
  decomposition~\eqref{eq:difference-decomposition}, we have
  \begin{align*}
    \| \zeta_{\mu}^{\bar{\kappa}} Q_{\sigma} (\psi_{\sigma}) \|
    &\leqslant
    \| \zeta_{\mu}^{\bar{\kappa}} J_{\sigma}
    F_{\sigma}^{[> 0]}(\psi_{\sigma}) \|
    + \| \zeta_{\mu}^{\bar{\kappa}} J_{\sigma}
    (F^{[0], (1)} \psi_{\sigma} + F^{[0], (0)}) \|\\
    &\quad + \| \zeta_{\mu}^{\bar{\kappa}}
    (1 - J_{\sigma})(\lambda \psi^3_{\sigma}) \|\\
    &\lesssim
    \| \zeta_{\mu}^{\bar{\kappa}} K_{\sigma}
    F_{\sigma}^{[>0]}(\psi_{\sigma}) \|
    + \| \zeta_{\mu}^{\bar{\kappa}} K_{\sigma}
    (F^{[0], (1)} \psi_{\sigma} + F^{[0], (0)}) \|\\
    &\quad + \| \zeta_{\mu}^{5 \nu}
    (1 - J_{\sigma})(\lambda \psi^3_{\sigma}) \| .
  \end{align*}
  To obtain the second estimate, we used that $\bar{\kappa} \geqslant 5 \nu$,
  which is a consequence of~\eqref{eq:cond-kappabar-one}, together with
  \[
    J_{\sigma} = L_{\sigma} J_{\sigma} K_{\sigma}, \qquad
    \| L_{\sigma} J_{\sigma} \|_{\operatorname{TV} (\zeta^{- 1})}
    \lesssim 1 .
  \]
  The latter bound is a consequence of Lemma~\ref{lem:bounds-K-J}. The
  bound~\eqref{eq:bound-Q} now follows from
  Lemmas~\ref{lem:control-of-forcing},~\ref{lem:l-0-kernels}
  and~\ref{lem:difference-bound}. 
  
  Next, by the decomposition of the force,
  \[
  F_{\sigma}(\psi_{\sigma})+\lambda\psi_{\sigma}^{3}
  =
  F_{\sigma}^{[>0]}(\psi_{\sigma})
  +F^{[0],(1)}\psi_{\sigma}
  +F^{[0],(0)}.
  \]
  Hence, by Lemmas~\ref{lem:control-of-forcing}
  and~\ref{lem:l-0-kernels},
  \begin{align*}
  	&\left\|
  	\zeta_{\mu}^{\bar\kappa}K_{\sigma}
  	\bigl(F_{\sigma}(\psi_{\sigma})
  	+\lambda\psi_{\sigma}^{3}\bigr)
  	\right\|\\
  	&\quad\leqslant
  	\left\|
  	\zeta_{\mu}^{\bar\kappa}K_{\sigma}
  	F_{\sigma}^{[>0]}(\psi_{\sigma})
  	\right\|
  	+
  	\left\|
  	\zeta_{\mu}^{\bar\kappa}K_{\sigma}
  	\bigl(F^{[0],(1)}\psi_{\sigma}
  	+F^{[0],(0)}\bigr)
  	\right\|\\
  	&\quad\lesssim
  	\llbracket\sigma\rrbracket^{-3\gamma+\vartheta}
  	\|F^{\mathfrak A}\|
  	(1+\trinorm{\psi}_{\bar\mu})^{\bar k}.
  \end{align*}
  This proves~\eqref{eq:K-force-coercive-remainder}.

  Since $\bar{\kappa} \geqslant 3 \nu$ and $\|
  K_{\sigma} \|_{\operatorname{TV} (\zeta^{- 1})} \lesssim 1$ by
  Lemmas~\ref{lem:bounds-K-J} and~\ref{lem:triple-phi}, we have
  \[ \| \zeta_{\mu}^{\bar{\kappa}} K_{\sigma} \psi^3_{\sigma} \| \lesssim \|
     \zeta_{\mu}^{3 \nu} K_{\sigma} \psi^3_{\sigma} \| \lesssim \| K_{\sigma}
     \|_{\operatorname{TV} (\zeta^{- 1})}  \| \zeta_{\mu}^{3 \nu} \psi^3_{\sigma} \|
     \lesssim \| \zeta_{\mu}^{3 \nu} \psi^3_{\sigma} \| \lesssim \llbracket
     \sigma \rrbracket^{- 3 \gamma}  \trinorm{ \psi }^3 . \]
  The bound~\eqref{eq:bound-JtildeF} now follows from the estimate
  \[ \| \zeta_{\mu}^{\bar{\kappa}} K_{\sigma} F_{\sigma} (\psi_{\sigma}) \|
     \lesssim 
     \left\|
     \zeta_{\mu}^{\bar\kappa}K_{\sigma}
     \bigl(F_{\sigma}(\psi_{\sigma})
     +\lambda\psi_{\sigma}^{3}\bigr)
     \right\|
      + \| \zeta_{\mu}^{\bar{\kappa}}
     K_{\sigma} \psi^3_{\sigma} \| \]
  and ~\eqref{eq:K-force-coercive-remainder}.
\end{proof}

\begin{lemma}
  \label{lem:control-of-forcing}Assume~\eqref{eq:cond-kappabar-one}.
  Then for every $\mathfrak{a} \in \mathfrak{A}$ such that $\ell
  (\mathfrak{a}) > 0$, the bound
  \begin{equation}
    \| \zeta_{\mu}^{\bar{\kappa}} K_{\sigma} F^{\mathfrak{a}}_{\sigma}
    (\psi^{\otimes k (\mathfrak{a})}_{\sigma}) \| \lesssim \llbracket \sigma
    \rrbracket^{- \alpha + (\beta - \gamma) k (\mathfrak{a}) + \delta \ell
    (\mathfrak{a}) - a \kappa_{\mathfrak{o}} (1 + \bar{\ell})}  \|
    F^{\mathfrak{A}} \|  {\trinorm{ \psi }^{k (\mathfrak{a})}} 
    \label{eq:est-F-a}
  \end{equation}
  holds uniformly in $\bar{\mu}, \mu \in [1 / 2, 1),$ $\sigma \in [\mu \vee
  \bar{\mu}, 1) $ and $\psi \in \mathcal{S}' (\Lambda)$. Moreover, provided
  \begin{equation}
    \vartheta \leqslant 3 (\gamma - \beta) + \delta - a \kappa_{\mathfrak{o}}
    (1 + \bar{\ell}) - \kappa, \label{eq:choice-zeta}
  \end{equation}
  the following bound
  \[ \| \zeta_{\mu}^{\bar{\kappa}} K_{\sigma} F^{[> 0]}_{\sigma}
     (\psi_{\sigma}) \| \lesssim \llbracket \sigma \rrbracket^{- 3 \gamma +
     \vartheta}  \| F^{\mathfrak{A}} \|  (1 + \trinorm{ \psi }
     )^{\bar{k}} . \]
  holds uniformly in $\bar{\mu}, \mu \in [1 / 2, 1),$ $\sigma \in [\mu \vee
  \bar{\mu}, 1) $ and $\psi \in \mathcal{S}' (\Lambda)$.
\end{lemma}

\begin{proof}
  We first observe that, since $\tilde{J}_{\sigma} J_{\sigma}  = J_{\sigma} $ and
  $K_{\sigma} L_{\sigma} = 1$, we have
  \[ \psi_{\sigma} = \tilde{J}_{\sigma} \psi_{\sigma} = K_{\sigma}^2
     L_{\sigma}^2 \tilde{J}_{\sigma} \psi_{\sigma} . \]
  Consequently, it follows from $\tilde{K}^{\mathfrak{a}}_{\mu} = K_{\mu} 
  \otimes (K^2_{\mu})^{\otimes k (\mathfrak{a})}$ and
  Definition~\ref{def:norm-kernels} of the kernel norm that
  \[ \| \zeta_{\mu}^{\bar{\kappa}} K_{\sigma} F^{\mathfrak{a}}_{\sigma}
     (\psi^{\otimes k (\mathfrak{a})}_{\sigma}) \| \lesssim \|
     \mathfrak{w}_{\mu, \sigma} \|  \| \mathfrak{o}^{\mathfrak{a}} \cdot
     [\tilde{K}^{\mathfrak{a}}_{\sigma} F^{\mathfrak{a}}_{\sigma}] \cdot
     \tilde{w}_{\sigma}^{\mathfrak{a}} \|  \| \rho_{\mu} L_{\sigma}^2
     \tilde{J}_{\sigma} \psi_{\sigma} \|^{k (\mathfrak{a})}, \]
  where
  \begin{align*}
    \mathfrak{w}_{\mu, \sigma}
    \left( z, z_1 {, \ldots, z_{k(\mathfrak{a})}} \right)
    &\assign
    \frac{\zeta_{\mu}^{\bar{\kappa}}(z)}
    {\mathfrak{o}^{\mathfrak{a}}(z)
    \tilde{w}_{\sigma}^{\mathfrak{a}}
    \left( z, z_1 {, \ldots, z_{k(\mathfrak{a})}} \right)
    \rho _{\mu}(z_1) \ldots \rho _{\mu}(z_{k(\mathfrak{a})})},\\
    \| \mathfrak{w}_{\mu, \sigma} \|
    &= \| \mathfrak{w}_{\mu, \sigma} \|_{L^{\infty}
    (\Lambda^{1 + k(\mathfrak{a})})} .
  \end{align*}
  Using that the weight $\mathfrak{o}$ is at scale $1$ (cf.
  Definition~\ref{def:weights-kernels}) and that the weight $\zeta_{\mu}$ is at
  scale $\llbracket \mu \rrbracket^{- a}$ (cf. Definition~\ref{def:weights}), we
  deduce
  \begin{equation}
    \mathfrak{o} (z)^{- 1} = \langle z \rangle_s^{\kappa_{\mathfrak{o}}}
    \leqslant \llbracket \mu \rrbracket^{- a \kappa_{\mathfrak{o}}} \langle
    \llbracket \mu \rrbracket^a .z \rangle_s^{\kappa_{\mathfrak{o}}} =
    \llbracket \mu \rrbracket^{- a \kappa_{\mathfrak{o}}} \zeta_{\mu}^{-
    \kappa_{\mathfrak{o}}} (z) . \label{eq:bound-o-zeta-weights}
  \end{equation}
  Noting that the Steiner diameter of a collection of points is always at
  least as large as the distance between any two points, we exploit the weight
  $(\tilde{w}_{\sigma}^a)^{- 1}$ to propagate the weights $\rho_{\mu}^{- 1}$
  to the output variable. Using furthermore that $\rho_{\mu}  =
  \zeta_{\mu}^{\nu}$ and $\mathfrak{o}^{\mathfrak{a}} =\mathfrak{o}^{1 + \ell
  (\mathfrak{a})}$ as well as $\ell (\mathfrak{a}) \leqslant \bar{\ell}$, $k
  (\mathfrak{a}) \leqslant \bar{k}$ and~$\eqref{eq:cond-kappabar-one}$, we
  obtain
  \begin{equation}
    \begin{aligned}
      \| \mathfrak{w}_{\mu, \sigma} \| & \lesssim \left\|
      \zeta_{\mu}^{\bar{\kappa}}  \left( {\mathfrak{o}^{\mathfrak{a}}} 
      {\rho^{k (\mathfrak{a})}_{\mu}}  \right)^{- 1} \right\|_{L^{\infty}
      (\Lambda)}\\
      & \lesssim \llbracket \mu \rrbracket^{- a \kappa_{\mathfrak{o}} (1 +
      \ell (\mathfrak{a}))}  \| \zeta_{\mu}^{\bar{\kappa} -
      \kappa_{\mathfrak{o}} (1 + \ell (\mathfrak{a})) - \nu k (\mathfrak{a})}
      \|_{L^{\infty} (\Lambda)}\\
      & \lesssim \llbracket \mu \rrbracket^{- a \kappa_{\mathfrak{o}} (1 +
      \bar{\ell})}
    \end{aligned} \label{eq:weight-bound-1}
  \end{equation}
  We also observe that, by Lemmas~\ref{lem:bounds-K-J}
  and~\ref{lem:triple-phi},
  \[ \| \rho_{\mu} L_{\sigma}^2 \tilde{J}_{\sigma} \psi_{\sigma} \| \lesssim
     \| L_{\sigma}^2 \tilde{J}_{\sigma} \|_{\operatorname{TV} (\zeta^{- 1}) } \|
     \rho_{\mu} \psi_{\sigma} \|_{L^{\infty}} \lesssim \llbracket \sigma
     \rrbracket^{- \gamma}  \trinorm{ \psi }  . \]
  Moreover, by Definition~\ref{def:norm-kernels},
  \[ \| \mathfrak{o}^{\mathfrak{a}} \cdot [\tilde{K}^{\mathfrak{a}}_{\sigma}
     F^{\mathfrak{a}}_{\sigma}] \cdot \tilde{w}_{\sigma}^{\mathfrak{a}} \| =
     \| F^{\mathfrak{a}}_{\sigma} \|_{\sigma} \lesssim \llbracket \sigma
     \rrbracket^{[\mathfrak{a}]} \| F^{\mathfrak{A}} \| . \]
  Combining the above estimates and using~\eqref{eq:homogeneity-of-kernels},
  we conclude that
  \[ \| \zeta_{\mu}^{\bar{\kappa}} K_{\sigma} F^{\mathfrak{a}}_{\sigma}
     (\psi^{\otimes k (\mathfrak{a})}_{\sigma}) \| \lesssim \llbracket \sigma
     \rrbracket^{- \alpha + (\beta - \gamma) k (\mathfrak{a}) + \delta \ell
     (\mathfrak{a}) - a \kappa_{\mathfrak{o}} (1 + \bar{\ell})}  \|
     F^{\mathfrak{A}} \|  {\trinorm{ \psi }^{k (\mathfrak{a})}} ,
  \]
  which proves~\eqref{eq:est-F-a}.
  
  To prove the second part, we observe that
  \begin{align*}
    \| \zeta_{\mu}^{\bar{\kappa}} K_{\sigma} F^{[> 0]}_{\sigma}
    (\psi_{\sigma}) \| &\leqslant  \sum_{\mathfrak{a} \mid \ell
    (\mathfrak{a}) > 0, k (\mathfrak{a}) \leqslant \bar{k}} \|
    \zeta_{\mu}^{\bar{\kappa}} K_{\sigma} F^{\mathfrak{a}}_{\sigma}
    (\psi_{\sigma}) \|\\
 &\lesssim  \sum_{\mathfrak{a} \mid \ell (\mathfrak{a}) > 0, k
    (\mathfrak{a}) \leqslant \bar{k}} \llbracket \sigma \rrbracket^{- \alpha +
    (\beta - \gamma) k (\mathfrak{a}) + \delta \ell (\mathfrak{a}) - a
    \kappa_{\mathfrak{o}} (1 + \bar{\ell})}  \| F^{\mathfrak{A}} \| 
    {\trinorm{ \psi }^{k (\mathfrak{a})}} \\
 &\lesssim  \llbracket \sigma \rrbracket^{- \alpha + \delta - a
    \kappa_{\mathfrak{o}} (1 + \bar{\ell})}  \| F^{\mathfrak{A}} \|  (1 +
    \trinorm{ \psi } )^{\bar{k}},
  \end{align*}
  where we used that $\beta \geqslant \gamma$, $1 \leqslant \ell
  (\mathfrak{a}) \leqslant \bar{\ell}$ and $k (\mathfrak{a}) \leqslant
  \bar{k}$. Since, by $\alpha = 3 \beta + \kappa$ and~\eqref{eq:choice-zeta},
  we have
  \[ 3 \gamma - \alpha + \delta - a \kappa_{\mathfrak{o}} (1 + \bar{\ell}) = 3
     (\gamma - \beta) + \delta - a \kappa_{\mathfrak{o}} (1 + \bar{\ell}) -
     \kappa \geqslant \vartheta, \]
  the proof is complete.
\end{proof}

\begin{lemma}
  \label{lem:l-0-kernels}Assume~\eqref{eq:cond-kappabar-one} and
  \begin{equation}
    \vartheta \leqslant 2 \gamma \wedge \left( 3 \gamma - \frac{d + 2 s}{2} -
    2 \kappa - a \kappa_{\mathfrak{o}} \right) . \label{eq:xi-on-zeta}
  \end{equation}
  Then the bound
  \[ \| \zeta_{\mu}^{\bar{\kappa}} K_{\sigma} F^{[0], (0)}_{\sigma} \| \vee \|
     \zeta_{\mu}^{\bar{\kappa}} K_{\sigma} F^{[0], (1)}_{\sigma}
     (\psi_{\sigma}) \| \lesssim \llbracket \sigma \rrbracket^{- 3 \gamma +
     \vartheta}  \| F^{\mathfrak{A}} \|  (1 + \trinorm{ \psi }),
  \]
  holds uniformly in $\bar{\mu}, \mu \in [1 / 2, 1), \sigma \in [\mu \vee
  \bar{\mu}, 1) $ and $\psi \in \mathcal{S}' (\Lambda)$.
\end{lemma}

\begin{proof}
  On account of~\eqref{eq:flow-kernels} and~\eqref{eq:B-form},
  $\partial_{\sigma} F^{[0]}_{\sigma} = 0$. Thus,
  \[ F^{[0], (0)}_{\sigma} (z) = F^{[0], (0)}_1 (z) = \xi^{(\varepsilon, M)}
     (z), \qquad F^{[0], (1)}_{\sigma} (z, z_1) = F^{[0], (1)}_1 (z, z_1) =
     \bar{r} \delta (z, z_1) . \]
  Recall the definition \eqref{eq:norm-family-F} of $\| F^{\mathfrak{A}} \|$.
  By Lemmas~\ref{lem:bounds-K-J} and~\ref{lem:triple-phi}, $\bar{\kappa}
  \geqslant \nu$ and $\vartheta \leqslant 2 \gamma$, we obtain
  \begin{align*}
    \| \zeta_{\mu}^{\bar{\kappa}} K_{\sigma} F^{[0], (1)}_{\sigma}
    (\psi_{\sigma}) \| &=  | \bar{r} |  \|
    \zeta_{\mu}^{\bar{\kappa}} K_{\sigma} \psi_{\sigma} \| \lesssim \|
    F^{\mathfrak{A}} \|  \| K_{\sigma} \|_{\operatorname{TV} (\zeta^{- 1})}  \|
    \zeta_{\mu}^{\bar{\kappa}} \psi_{\sigma} \|\\
 &\lesssim  \llbracket \sigma \rrbracket^{- \gamma}  \| F^{\mathfrak{A}}
    \|  \trinorm{ \psi } \lesssim \llbracket \sigma \rrbracket^{-
    3 \gamma + \vartheta}  \| F^{\mathfrak{A}} \|  \trinorm{ \psi
    } .
  \end{align*}
  Moreover, by~\eqref{eq:xi-on-zeta}, we have
  \begin{align*}
     \| \zeta_{\mu}^{\bar{\kappa}} K_{\sigma} F^{[0], (0)}_{\sigma} \|
     &\lesssim \| \zeta_{\mu}^{\bar{\kappa}} \mathfrak{o}^{- 1} \|
     \| \mathfrak{o}K_{\sigma} F^{[0], (0)} \|\\
     &\lesssim \llbracket \mu \rrbracket^{- a \kappa_{\mathfrak{o}}}
     \llbracket \sigma \rrbracket^{- (d + 2 s) / 2 - 2 \kappa}
     \| F^{\mathfrak{A}} \|\\
     &\lesssim \llbracket \sigma \rrbracket^{- 3 \gamma + \vartheta}
     \| F^{\mathfrak{A}} \|,
  \end{align*}
  where the bound $\| \zeta_{\mu}^{\bar{\kappa}} \mathfrak{o}^{- 1} \|
  \lesssim \llbracket \mu \rrbracket^{- a
  \kappa_{\mathfrak{o}}}$ follows
  from~\eqref{eq:bound-o-zeta-weights}, since $\bar{\kappa} \geqslant
  \kappa_{\mathfrak{o}}$ by~\eqref{eq:cond-kappabar-one}.
\end{proof}

\begin{lemma}
  \label{lem:estimates-DF}Assume~\eqref{eq:cond-kappabar-one} and
  \begin{equation}
    \vartheta \leqslant \delta - a \kappa_{\mathfrak{o}} (1 + \bar{\ell}) .
    \label{eq:DF-bound}
  \end{equation}
  Then the bound
  \[ \| \zeta_{\mu} K_{\sigma} (\mathD F_{\sigma} (\psi_{\sigma})
     (\dot{G}_{\sigma} \hat{\psi})) \| \lesssim \llbracket \sigma
     \rrbracket^{\vartheta - 1}  \| F^{\mathfrak{A}} \|  (1 + \trinorm{ \psi
     })^{\bar{k}}  \| \zeta_{\mu}^{1 - \bar{\kappa}} \tilde{J
     }_{\sigma}^2 \hat{\psi} \|, \]
  holds uniformly in $\bar{\mu} \in [1 / 2, 1),$ $\mu \in [1 / 2, 1)$, $\sigma
  \in [\mu \vee \bar{\mu}, 1)$ and $\psi, \hat{\psi} \in \mathcal{S}'
  (\Lambda)$.
\end{lemma}

\begin{proof}
  We have
  \[ \left\| \zeta_{\mu} K_{\sigma} \left( {\mathD F_{\sigma}} 
     (\psi_{\sigma}) (\dot{G}_{\sigma} \hat{\psi}) \right) \right\| \lesssim
     \sum_{\mathfrak{a} \mid k (\mathfrak{a}) > 0} \| \zeta_{\mu}
     K_{\sigma} F_{\sigma}^{\mathfrak{a}} (\psi_{\sigma}^{\otimes (k
     (\mathfrak{a}) - 1)} \otimes \dot{G}_{\sigma} \hat{\psi}) \| . \]
  Observe that since $\dot{G}_{\sigma} \assign \mathcal{L}^{- 1}_{\varepsilon}
  \dot{J}_{\sigma}$, $\tilde{J}_{\sigma}^2 \dot{J}_{\sigma}  =
  \dot{J}_{\sigma} $, $\tilde{J}_{\sigma} J_{\sigma}  = J_{\sigma}$ and
  $K_{\sigma} L_{\sigma} = 1$,
  \[ \psi_{\sigma} = K_{\sigma}^2 L_{\sigma}^2 \tilde{J}_{\sigma}
     \psi_{\sigma}, \qquad \dot{G}_{\sigma} \hat{\psi} = \tilde{J}_{\sigma}^2
     \dot{G}_{\sigma} \hat{\psi} = K_{\sigma}^2 L_{\sigma}^2 \dot{G}_{\sigma}
     \tilde{\psi}, \qquad \tilde{\psi} = \tilde{J }_{\sigma}^2 \hat{\psi}
     . \]
  Using Definition~\ref{def:norm-kernels} of the kernel norm, we obtain
  \begin{align*}
    &\| \zeta_{\mu} K_{\sigma} F_{\sigma}^{\mathfrak{a}}
    (\psi_{\sigma}^{\otimes (k(\mathfrak{a}) - 1)}
    \otimes \dot{G}_{\sigma}\tilde{\psi}) \|\\
    &\quad \lesssim
    \| \mathfrak{w}_{\mu, \sigma} \|
    \| \mathfrak{o}^{\mathfrak{a}}
    \cdot [\tilde{K}_{\sigma}^{\mathfrak{a}}F_{\sigma}^{\mathfrak{a}}]
    \cdot \tilde{w}_{\sigma}^{\mathfrak{a}} \|\\
    &\qquad \times
    \| \rho_{\mu} L^2_{\sigma} \tilde{J}_{\sigma} \psi_{\sigma}
    \|^{k(\mathfrak{a}) - 1}
    \| \zeta_{\mu}^{1 - \bar{\kappa}} L_{\sigma}^2
    \dot{G}_{\sigma} \tilde{\psi} \|,
  \end{align*}
  where
  \begin{align*}
    \mathfrak{w}_{\mu, \sigma} \left( z, z_1 {, \ldots, z_{k
     (\mathfrak{a})}}  \right)
     &\assign \frac{\zeta_{\mu} (z)}
     {\mathfrak{o}^{\mathfrak{a}} (z)  \tilde{w}_{\sigma}^{\mathfrak{a}}
     \left( z, z_1 {, \ldots, z_{k (\mathfrak{a})}}  \right) \zeta_{\mu}^{1 -
     \bar{\kappa}} (z_1) \rho _{\mu} (z_2) \ldots \rho _{\mu} (z_{k
     (\mathfrak{a})})},\\
    \| \mathfrak{w}_{\mu, \sigma} \|
     &\assign \| \mathfrak{w}_{\mu, \sigma} \|_{L^{\infty} (\Lambda^{1 + k
     (\mathfrak{a})})} .
  \end{align*}
  Using furthermore that $\rho_{\mu} = \zeta_{\mu}^{\nu}$ and
  $\mathfrak{o}^{\mathfrak{a}} =\mathfrak{o}^{1 + \ell (\mathfrak{a})}$ as
  well as $\ell (\mathfrak{a}) \leqslant \bar{\ell}$, $k (\mathfrak{a})
  \leqslant \bar{k}$ and~$\eqref{eq:cond-kappabar-one}$, we obtain
  \begin{equation}
    \begin{aligned}
      \| \mathfrak{w}_{\mu, \sigma} \|_{L^{\infty}} & \lesssim \|
      \zeta_{\mu} (\mathfrak{o}^{\mathfrak{a}} \zeta^{1 - \bar{\kappa}}_{\mu}
      \rho_{\mu}^{(k (\mathfrak{a}) - 1)})^{- 1} \|_{L^{\infty} (\Lambda)}\\
      & \lesssim \llbracket \mu \rrbracket^{- a
      \kappa_{\mathfrak{o}} (1 + \ell (\mathfrak{a}))}  \| \zeta^{\bar{\kappa}
      - \kappa_{\mathfrak{o}} (1 + \ell (\mathfrak{a})) - \nu (k
      (\mathfrak{a}) - 1)}_{\mu} \|_{L^{\infty} (\Lambda)}\\
      & \lesssim \llbracket \mu \rrbracket^{- a
      \kappa_{\mathfrak{o}} (1 + \bar{\ell})} .
    \end{aligned} \label{eq:weight-bound-2}
  \end{equation}
  We also observe that, by Lemmas~\ref{lem:bounds-K-J}
  and~\ref{lem:triple-phi},
  \[ \| \rho_{\mu} L_{\sigma}^2 \tilde{J}_{\sigma} \psi_{\sigma} \| \lesssim
     \| L_{\sigma}^2 \tilde{J}_{\sigma} \|_{\operatorname{TV} (\zeta^{- 1}) } \|
     \rho_{\mu} \psi_{\sigma} \|_{L^{\infty}} \lesssim \llbracket \sigma
     \rrbracket^{- \gamma}  \trinorm{ \psi }  \]
  and
  \[ \| \zeta_{\mu}^{1 - \bar{\kappa}} L_{\sigma}^2 \dot{G}_{\sigma}
     \tilde{\psi} \| \lesssim \| L_{\sigma}^2 \dot{G}_{\sigma} \|_{\operatorname{TV}
     (\zeta^{- 1}) }  \| \zeta_{\mu}^{1 - \bar{\kappa}} \tilde{\psi} \|
     \lesssim \llbracket \sigma \rrbracket^{2 s - 1}  \| \zeta_{\mu}^{1 -
     \bar{\kappa}} \tilde{\psi} \| . \]
  Moreover, by Definition~\ref{def:norm-kernels},
  \[ \| \mathfrak{o}^{\mathfrak{a}} \cdot [\tilde{K}^{\mathfrak{a}}_{\sigma}
     F^{\mathfrak{a}}_{\sigma}] \cdot \tilde{w}_{\sigma}^{\mathfrak{a}} \| =
     \| F^{\mathfrak{a}}_{\sigma} \|_{\sigma} \lesssim \llbracket \sigma
     \rrbracket^{[\mathfrak{a}]} \| F^{\mathfrak{A}} \| . \]
  Combining the above estimates and using~\eqref{eq:homogeneity-of-kernels},
  we conclude that
  \begin{align*}
    &\| \zeta_{\mu} K_{\sigma} F_{\sigma}^{\mathfrak{a}}
    (\psi_{\sigma}^{\otimes (k(\mathfrak{a}) - 1)}
    \otimes \dot{G}_{\sigma} \tilde{\psi}) \|\\
    &\quad \lesssim
    \llbracket \sigma \rrbracket^{
    2 s - \alpha + \beta + (\beta - \gamma)(k(\mathfrak{a}) - 1)
    - a \kappa_{\mathfrak{o}}(1 + \bar{\ell}) - 1}
    \| F^{\mathfrak{A}} \|
    (1 + \trinorm{ \psi })^{\bar{k}}
    \| \zeta^{1 - \bar{\kappa}}_{\mu} \tilde{\psi} \|\\
    &\quad \lesssim
    \llbracket \sigma \rrbracket^{\delta
    - a \kappa_{\mathfrak{o}}(1 + \bar{\ell}) - 1}
    \| F^{\mathfrak{A}} \|
    (1 + \trinorm{ \psi })^{\bar{k}}
    \| \zeta^{1 - \bar{\kappa}}_{\mu} \tilde{\psi} \|,
  \end{align*}
  where in the last step we used $2 s - \alpha + \beta = \delta$ and that
  $\beta \geqslant \gamma$. This finishes the proof.
\end{proof}

\begin{lemma}
  \label{lem:packaging}Assume that
  \begin{equation}
    \vartheta \leqslant \delta \bar{\ell} - \alpha - a \kappa_{\mathfrak{o}}
    (2 \bar{\ell} + 2) . \label{eq:choice-of-K-and-l}
  \end{equation}
  Then
  \[ H_{\sigma} (\psi) \assign \partial_{\sigma} F_{\sigma} (\psi) + \mathD
     F_{\sigma} (\psi)  (\dot{G}_{\sigma} F_{\sigma} (\psi)), \qquad \psi \in
     \mathcal{S}' (\Lambda), \]
  satisfies the bound
  \[ \| \zeta_{\mu} K_{\sigma} H_{\sigma} (\psi_{\sigma}) \| \lesssim
     \llbracket \sigma \rrbracket^{\vartheta - 1}  \| F^{\mathfrak{A}} \|^2 
     (1 + \trinorm{ \psi })^{2 \bar{k}} . \]
  uniformly in $\bar{\mu}, \mu \in [1 / 2, 1),$ $\sigma \in [\mu \vee
  \bar{\mu}, 1)$ and $\psi \in \mathcal{S}' (\Lambda)$.
\end{lemma}

\begin{proof}
  We observe that, on account of the perturbative flow equation, we have
  \[ H_{\sigma}  (\psi_{\sigma}) = \sum_{\mathfrak{a} \in
     \mathfrak{A}_{\mathrm{ext}} \setminus \mathfrak{A}}
     \sum_{\mathfrak{b}, \mathfrak{c} \in \mathfrak{A}}
     B^{\mathfrak{a}}_{\mathfrak{b}, \mathfrak{c}}
     (\dot{G}_{\sigma}, F^{\mathfrak{b}}_{\sigma}, F^{\mathfrak{c}}_{\sigma})
     (\psi^{\otimes k (\mathfrak{a})}_{\sigma}) . \]
  Then, working
  as in the proof of Lemma~\ref{lem:control-of-forcing} and using
  Lemma~\ref{lem:estimate-B-kernels} as well as $k (\mathfrak{a}) = k
  (\mathfrak{b}) + k (\mathfrak{c}) - 1 \leqslant 2 \bar{k}$, we obtain
  \begin{align*}
    \| \zeta_{\mu} K_{\sigma} H_{\sigma} (\psi_{\sigma}) \| &\leqslant 
    \sum_{\mathfrak{a} \in \mathfrak{A}_{\mathrm{ext}} \setminus \mathfrak{A}}
    \sum_{\mathfrak{b}, \mathfrak{c} \in \mathfrak{A}} \|
    \zeta_{\mu} K_{\sigma} B^{\mathfrak{a}}_{\mathfrak{b}, \mathfrak{c}}
    (\dot{G}_{\sigma}, F^{\mathfrak{b}}_{\sigma}, F^{\mathfrak{c}}_{\sigma})
    (\psi^{\otimes k (\mathfrak{a})}_{\sigma}) \|\\
 &\lesssim  \sum_{\mathfrak{a} \in \mathfrak{A}_{\mathrm{ext}} \setminus \mathfrak{A}}
    \sum_{\mathfrak{b}, \mathfrak{c} \in \mathfrak{A}}
    \llbracket \mu \rrbracket^{- a \kappa_{\mathfrak{o}}
    (1 + \ell (\mathfrak{a}))}\\
 &\quad \times
    \llbracket \sigma \rrbracket^{- \gamma k (\mathfrak{a})}
    \| B^{\mathfrak{a}}_{\mathfrak{b}, \mathfrak{c}} (\dot{G}_{\sigma},
    F^{\mathfrak{b}}_{\sigma}, F^{\mathfrak{c}}_{\sigma}) \|_{\sigma}
    \trinorm{ \psi }^{k (\mathfrak{a})}\\
 &\lesssim  \sum_{\mathfrak{a} \in \mathfrak{A}_{\mathrm{ext}} \setminus \mathfrak{A}}
    \sum_{\mathfrak{b}, \mathfrak{c} \in \mathfrak{A}}
    \llbracket \mu \rrbracket^{- a \kappa_{\mathfrak{o}}
    (1 + \ell (\mathfrak{a}))}\\
 &\quad \times
    \llbracket \sigma \rrbracket^{[\mathfrak{a}] - [\mathfrak{b}]
    - [\mathfrak{c}] - \gamma k (\mathfrak{a})-1}
    \| F^{\mathfrak{b}}_{\sigma} \|_{\sigma}
    \| F^{\mathfrak{c}}_{\sigma} \|_{\sigma}
    \trinorm{ \psi }^{k (\mathfrak{a})}\\
 &\lesssim  \sum_{\mathfrak{a} \in \mathfrak{A}_{\mathrm{ext}} \setminus \mathfrak{A}}
    \sum_{\mathfrak{b}, \mathfrak{c} \in \mathfrak{A}}
    \llbracket \sigma \rrbracket^{[\mathfrak{a}] - \gamma k (\mathfrak{a})
    - a \kappa_{\mathfrak{o}} (1 + \ell (\mathfrak{a})) - 1}\\
 &\quad \times
    \| F^{\mathfrak{A}} \|^2
    (1 + \trinorm{ \psi })^{2 \bar{k}} .
  \end{align*}
  The factor $\llbracket \mu \rrbracket^{- a
  \kappa_{\mathfrak{o}} (1 + \ell (\mathfrak{a}))}$ appearing in the second
  line arises due to the estimate~\eqref{eq:bound-o-zeta-weights} by a bound
  similar to~\eqref{eq:weight-bound-1}. Now since $\beta \geqslant \gamma$
  and $\ell (\mathfrak{a}) \in [\bar{\ell}+1, 2 \bar{\ell}+1]$, we have
  \begin{align*}
    &[\mathfrak{a}] - \gamma k (\mathfrak{a})
    - a \kappa_{\mathfrak{o}} (1 + \ell(\mathfrak{a}))\\
    &\quad = - \alpha + (\beta - \gamma) k(\mathfrak{a})
    + \delta \ell(\mathfrak{a})
    - a \kappa_{\mathfrak{o}} (1 + \ell(\mathfrak{a}))\\
    &\quad \geqslant - \alpha + \delta \bar{\ell}
    - a \kappa_{\mathfrak{o}} (2 \bar{\ell} + 2) \geqslant \vartheta,
  \end{align*}
  where we exploited~\eqref{eq:homogeneity-of-kernels}
  and~\eqref{eq:choice-of-K-and-l}. This proves the claim.
\end{proof}

\begin{lemma}
  \label{lem:difference-bound}For all
  \begin{equation}
    \vartheta \leqslant 2 s - 2 \gamma \label{eq:one-J-bound}
  \end{equation}
  the following bound
  \[ \| \zeta_{\mu}^{5 \nu}  (1 - J_{\sigma}) (\psi_{\sigma}^3) \| \lesssim
     \llbracket \sigma \rrbracket^{- 3 \gamma + \vartheta}  \trinorm{ \psi
     }^2  \trinorm{ \mathcal{L} \psi_{\Cdot} }_{\#} .
  \]
  holds uniformly in $\bar{\mu}, \mu \in [1 / 2, 1),$ $\sigma \in [\mu \vee
  \bar{\mu}, 1) $ and $\psi \in \mathcal{S}' (\Lambda)$.
\end{lemma}

\begin{proof}
  Let $J_{> \eta} \assign 1 - J_{\eta}$. Then
  \begin{align*}
     (1 - J_{\sigma}) (\psi_{\sigma}^3)
     &= (1 - J_{\sigma}) [(J_{\eta} \psi_{\sigma})^3
     + (J_{> \eta} \psi_{\sigma})^3\\
     &\quad + 3 (J_{\eta} \psi_{\sigma})^2
     (J_{> \eta} \psi_{\sigma})
     + 3 (J_{\eta} \psi_{\sigma})
     (J_{> \eta} \psi_{\sigma})^2] .
  \end{align*}
  Choosing $\eta  = \frac{\sigma}{6 - 5 \sigma}$, which implies $6 \eta (1 -
  \eta )^{- 1} = \sigma (1 - \sigma)^{- 1}$, the first contribution $(1 -
  J_{\sigma}) (J_{\eta } \psi_{\sigma})^3$ vanishes by the Fourier space
  support property of the product. As for the other contributions, we have,
  for example,
  \begin{align*}
     &\| \zeta_{\mu}^{5 \nu}  (1 - J_{\sigma})
     (J_{> \eta} \psi_{\sigma})^3 \|
     = \| \rho_{\mu}^5  (1 - J_{\sigma})
     (J_{> \eta} \psi_{\sigma})^3 \|\\
     &\quad \lesssim \| \rho_{\mu}^5
     (J_{> \eta} \psi_{\sigma})^3 \|
     \lesssim \| \rho_{\mu} \psi_{\sigma} \|^2
     \| \rho_{\mu}^3  (J_{> \eta} \psi_{\sigma}) \| .
  \end{align*}
  Since $\sigma \geqslant \mu$ and $\llbracket \eta \rrbracket \approx
  \llbracket \sigma \rrbracket$, by Lemmas~\ref{lem:triple-phi}
  and~\ref{lem:schauder-two}, we have
  \[ \| \zeta^{5 \nu}_{\mu}  (1 - J_{\sigma}) (J_{> \eta} \psi_{\sigma})^3 \|
     \lesssim \llbracket \eta \rrbracket^{2 s}  \llbracket \sigma
     \rrbracket^{- 5 \gamma}  \trinorm{ \psi }^2  \trinorm{
     \mathcal{L} \psi_{\Cdot} }_{\#} \lesssim \llbracket \sigma
     \rrbracket^{2 s - 5 \gamma}  \trinorm{ \psi }^2  \trinorm{
     \mathcal{L} \psi_{\Cdot} }_{\#} . \]
  The contributions coming from the terms $(J_{\eta} \psi_{\sigma})^2  (J_{>
  \eta} \psi_{\sigma})$ and $(J_{\eta} \psi_{\sigma})  (J_{> \eta }
  \psi_{\sigma})^2$ satisfy the same bound. Since $\vartheta \leqslant 2 s - 2
  \gamma$, this finishes the proof.
\end{proof}

Lemmas~\ref{lem:overall-difference-bound} and~\ref{lem:difference-bound}
provide the first estimate in~\eqref{eq:est-F}; Lemma~\ref{lem:packaging}
provides the second; the second bound in
Lemma~\ref{lem:overall-difference-bound} provides the third; and
Lemma~\ref{lem:estimates-DF} provides the fourth. Together with the random
kernel estimates from Section~\ref{sec:concl-of-cumul-analysis}, this completes
the proof of Theorem~\ref{thm:effective-force}.

\appendix\section{Auxiliary estimates}\label{app:technical}

We collect in this appendix various technical estimates of general character. Recall that $d=3$ and $s\in(3/4,1)$. As in the main text, all estimates are understood to hold uniformly in $\varepsilon \in 2^{-\mathbb{N}_0}$ and $M \in \mathbb{N}_+$.

\subsection{Kernel estimates}

\begin{definition}
	\label{def:derivatives}
	Recall from Section~\ref{sec:notations} that for a spacetime multi-index
	\[
	A \in \mathbb N_0^{\{ 0, 1\pm, \ldots, d\pm\}},
	\]
	we denote its spatial part by
	\[
	\bar A\assign (A_{i,+},A_{i,-})_{i=1}^d
	\in \mathbb N_0^{\{1\pm, \ldots, d\pm\}}
	\]
	and write
	\[
	|\bar A|
	\assign
	\sum_{i=1}^d(A_{i,+}+A_{i,-}),
	\qquad
	|A|_s
	\assign
	2sA_0+|\bar A|.
	\]
	For \((k_0,k)\in\mathbb R\times\mathbb R^d\), let
	\[
	d_\varepsilon^{i,\pm}(k)
	\assign
	\pm\frac{e^{\pm i\varepsilon k_i}-1}{\varepsilon},
	\]
	and
	\[
	d_\varepsilon^{\bar A}(k)
	\assign
	\prod_{i=1}^d
	\bigl(d_\varepsilon^{i,+}(k)\bigr)^{A_{i,+}}
	\bigl(d_\varepsilon^{i,-}(k)\bigr)^{A_{i,-}},
	\qquad
	d_\varepsilon^A(k_0,k)
	\assign
	(ik_0)^{A_0}d_\varepsilon^{\bar A}(k).
	\]
	Thus
	\[
	\widehat{\partial_\varepsilon^Af}(k_0,k)
	=
	d_\varepsilon^A(k_0,k)\hat f(k_0,k).
	\]
\end{definition}

\begin{lemma}
	\label{lem:q-epsilon-bounds}
	Recall that
	\[
	q_\varepsilon(k)^2
	=
	\frac4{\varepsilon^2}
	\sum_{j=1}^d
	\sin^2\left(\frac{\varepsilon k_j}{2}\right)
	=
	\frac2{\varepsilon^2}
	\sum_{j=1}^d
	\bigl(1-\cos(\varepsilon k_j)\bigr).
	\]
	The following bounds hold uniformly in \(\varepsilon\).
	
	\begin{enumerate}[(a)]
		\item For every \(k\in(\mathbb R_\varepsilon^d)^\ast\),
		\[
		q_\varepsilon(k)\asymp |k|.
		\]
		
		\item For every \(\beta\geqslant0\) and
		\(\tilde a\in\mathbb N_0^d\),
		\begin{equation}
			\label{eq:q-epsilon-derivatives}
			\left|
			\partial_k^{\tilde a}q_\varepsilon^\beta(k)
			\right|
			\lesssim_{\beta,\tilde a}
			q_\varepsilon^{\beta-|\tilde a|}(k),
			\qquad q_\varepsilon(k)>0.
		\end{equation}
		In the special case \(\beta=2\), one has the stronger estimate
		\begin{equation}
			\label{eq:q-epsilon-square-derivatives}
			\left|
			\partial_k^{\tilde a}q_\varepsilon^2(k)
			\right|
			\lesssim_{\tilde a}
			\varepsilon^{(|\tilde a|-2)_+}
			q_\varepsilon^{(2-|\tilde a|)_+}(k)
		\end{equation}
		for every \(k\in(\mathbb R_\varepsilon^d)^\ast\), including the
		points at which \(q_\varepsilon(k)=0\). Equivalently,
		\[
		\left|
		\partial_k^{\tilde a}q_\varepsilon^2(k)
		\right|
		\lesssim_{\tilde a}
		\begin{cases}
			q_\varepsilon^2(k),
			& |\tilde a|=0,\\
			q_\varepsilon(k),
			& |\tilde a|=1,\\
			1,
			& |\tilde a|=2,\\
			\varepsilon^{|\tilde a|-2},
			& |\tilde a|\geqslant3.
		\end{cases}
		\]
		
		\item For every
		\(\bar A\in\mathbb N_0^{\{1\pm,\ldots,d\pm\}}\) and
		\(\tilde a\in\mathbb N_0^d\),
		\begin{equation}
			\label{eq:d-epsilon-derivatives}
			\left|
			\partial_k^{\tilde a}d_\varepsilon^{\bar A}(k)
			\right|
			\lesssim_{\bar A,\tilde a}
			q_\varepsilon^{|\bar A|-|\tilde a|}(k),
			\qquad q_\varepsilon(k)>0.
		\end{equation}
		Consequently, if
		\(A\in\mathbb N_0^{\{0,1\pm,\ldots,d\pm\}}\), then, for
		\(a_0\in\mathbb N_0\) and \(\tilde a\in\mathbb N_0^d\),
		\begin{equation}
			\label{eq:d-epsilon-full-derivatives}
			\left|
			\partial_{k_0}^{a_0}\partial_k^{\tilde a}
			d_\varepsilon^A(k_0,k)
			\right|
			\lesssim_{A,\tilde a,a_0}
			|k_0|^{(A_0-a_0)_+}
			q_\varepsilon^{|\bar A|-|\tilde a|}(k),
			\qquad q_\varepsilon(k)>0,
		\end{equation}
		with the left-hand side equal to zero if \(a_0>A_0\).
		
		\item Let \(\chi\in C_c^\infty((0,\infty))\). For every
		\(\beta\geqslant0\), \(\tilde a\in\mathbb N_0^d\), and
		\(\eta\geqslant1/2\),
		\begin{equation}
			\label{eq:q-epsilon-cutoff-derivatives}
			\left|
			\partial_k^{\tilde a}
			\left(
			q_\varepsilon^\beta(k)
			\chi(\eta q_\varepsilon(k))
			\right)
			\right|
			\lesssim_{\beta,\tilde a,\chi}
			\eta^{|\tilde a|-\beta}
			1_{\operatorname{supp}\chi}
			(\eta q_\varepsilon(k)).
		\end{equation}
		More generally, for every
		\(\bar A\in\mathbb N_0^{\{1\pm,\ldots,d\pm\}}\),
		\begin{equation}
			\label{eq:q-epsilon-cutoff-d-derivatives}
			\left|
			\partial_k^{\tilde a}
			\left(
			q_\varepsilon^\beta(k)
			\chi(\eta q_\varepsilon(k))
			d_\varepsilon^{\bar A}(k)
			\right)
			\right|
			\lesssim_{\beta,\tilde a,\bar A,\chi}
			\eta^{|\tilde a|-\beta-|\bar A|}
			1_{\operatorname{supp}\chi}
			(\eta q_\varepsilon(k)).
		\end{equation}
		
		\item For every \(R>0\) and \(\gamma<d\),
		\begin{equation}
			\label{eq:q-epsilon-volume-singular}
			\int_{\{q_\varepsilon\leqslant R\}}
			q_\varepsilon^{-\gamma}(k)\,\mathd k
			\lesssim_\gamma R^{d-\gamma}.
		\end{equation}
		In particular, for every \(r<d+4s\),
		\begin{equation}
			\label{eq:q-epsilon-check-integrability}
			\int_{\{q_\varepsilon\leqslant2\}}
			\left(
			1\vee q_\varepsilon^{4s-r}(k)
			\right)\,\mathd k
			\lesssim_r1.
		\end{equation}
	\end{enumerate}
\end{lemma}

\begin{proof}
	The comparison \(q_\varepsilon(k)\asymp|k|\) follows from
	\[
	|\sin x|\asymp|x|,
	\qquad |x|\leqslant\pi/2.
	\]
	
	We first prove the strengthened estimate
	\eqref{eq:q-epsilon-square-derivatives}. Since
	\[
	q_\varepsilon^2(k)
	=
	\frac2{\varepsilon^2}
	\sum_{j=1}^d
	\bigl(1-\cos(\varepsilon k_j)\bigr),
	\]
	one has
	\[
	\partial_k^{\tilde a}q_\varepsilon^2=0
	\]
	whenever \(\tilde a\) has nonzero components in at least two
	distinct coordinate directions. It therefore suffices to consider
	\(\tilde a=ne_j\) for some \(j\in\{1,\ldots,d\}\) and \(n\geqslant1\).
	Direct differentiation gives
	\[
	\left|
	\partial_{k_j}^n q_\varepsilon^2(k)
	\right|
	\leqslant
	2\varepsilon^{n-2}
	\begin{cases}
		|\sin(\varepsilon k_j)|,
		& n \text{ odd},\\
		|\cos(\varepsilon k_j)|,
		& n \text{ even}.
	\end{cases}
	\]
	For \(n=1\), using
	\[
	|\sin(\varepsilon k_j)|
	\leqslant
	2\left|
	\sin\left(\frac{\varepsilon k_j}{2}\right)
	\right|
	\leqslant
	\varepsilon q_\varepsilon(k),
	\]
	we obtain
	\[
	\left|
	\partial_{k_j}q_\varepsilon^2(k)
	\right|
	\lesssim q_\varepsilon(k).
	\]
	For \(n=2\), the derivative is uniformly bounded, while for
	\(n\geqslant3\),
	\[
	\left|
	\partial_{k_j}^n q_\varepsilon^2(k)
	\right|
	\lesssim_n\varepsilon^{n-2}.
	\]
	Together with the case \(n=0\), this proves
	\eqref{eq:q-epsilon-square-derivatives}.
	
	We next derive the weaker estimate needed for arbitrary powers.
	Since
	\[
	q_\varepsilon(k)\lesssim\varepsilon^{-1},
	\]
	estimate \eqref{eq:q-epsilon-square-derivatives} implies, whenever
	\(q_\varepsilon(k)>0\),
	\[
	\left|
	\partial_k^{\tilde a}q_\varepsilon^2(k)
	\right|
	\lesssim_{\tilde a}
	q_\varepsilon^{2-|\tilde a|}(k).
	\]
	Indeed, this is immediate for \(|\tilde a|\leqslant2\), whereas for
	\(|\tilde a|\geqslant3\) one uses
	\[
	\varepsilon^{|\tilde a|-2}
	\lesssim
	q_\varepsilon^{2-|\tilde a|}(k).
	\]
	
	Now write
	\[
	q_\varepsilon^\beta
	=
	(q_\varepsilon^2)^{\beta/2}.
	\]
	By Faà di Bruno's formula, every term in
	\(\partial_k^{\tilde a}q_\varepsilon^\beta\), with
	\(\tilde a\neq0\), is of the form
	\[
	C
	(q_\varepsilon^2)^{\beta/2-m}
	\prod_{i=1}^m
	\partial_k^{\tilde a_i}q_\varepsilon^2,
	\qquad
	\tilde a_1+\cdots+\tilde a_m=\tilde a,
	\qquad
	|\tilde a_i|\geqslant1.
	\]
	Consequently,
	\[
	\left|
	(q_\varepsilon^2)^{\beta/2-m}
	\prod_{i=1}^m
	\partial_k^{\tilde a_i}q_\varepsilon^2
	\right|
	\lesssim
	(q_\varepsilon^2)^{\beta/2-m}
	\prod_{i=1}^m
	q_\varepsilon^{2-|\tilde a_i|}
	=
	q_\varepsilon^{\beta-|\tilde a|}.
	\]
	This proves \eqref{eq:q-epsilon-derivatives}.
	
	We next prove \eqref{eq:d-epsilon-derivatives}. For each
	\(i\in\{1,\ldots,d\}\), each sign, and each
	\(\tilde r\in\mathbb N_0^d\),
	\[
	\left|
	\partial_k^{\tilde r}d_\varepsilon^{i,\pm}(k)
	\right|
	\lesssim_{\tilde r}
	q_\varepsilon^{1-|\tilde r|}(k),
	\qquad q_\varepsilon(k)>0.
	\]
	For \(\tilde r=0\), this follows from
	\[
	|d_\varepsilon^{i,\pm}(k)|
	\lesssim q_\varepsilon(k).
	\]
	For \(|\tilde r|=1\), it is immediate from the explicit formula for
	\(d_\varepsilon^{i,\pm}\). For \(|\tilde r|\geqslant2\), the
	derivative vanishes unless \(\tilde r=ne_i\), in which case it is
	bounded by \(\varepsilon^{|\tilde r|-1}\). Since
	\(q_\varepsilon(k)\lesssim\varepsilon^{-1}\), this gives
	\[
	\varepsilon^{|\tilde r|-1}
	\lesssim
	q_\varepsilon^{1-|\tilde r|}(k).
	\]
	The product rule now yields
	\eqref{eq:d-epsilon-derivatives}. Estimate
	\eqref{eq:d-epsilon-full-derivatives} follows by differentiating
	the factor \((ik_0)^{A_0}\).
	
	Estimate \eqref{eq:q-epsilon-cutoff-derivatives} follows from
	\eqref{eq:q-epsilon-derivatives} and the chain rule, since on
	\(\operatorname{supp}\chi(\eta q_\varepsilon)\) one has
	\[
	q_\varepsilon\asymp\eta^{-1}.
	\]
	Combining the same argument with
	\eqref{eq:d-epsilon-derivatives} gives
	\eqref{eq:q-epsilon-cutoff-d-derivatives}.
	
	Finally, \eqref{eq:q-epsilon-volume-singular} follows from
	\(q_\varepsilon\asymp|k|\). Taking
	\(\gamma=(r-4s)_+\) gives
	\eqref{eq:q-epsilon-check-integrability}.
\end{proof}

\begin{lemma}[Ordinary lattice heat-kernel estimates]
	\label{lem:sublattice-heat-kernel}
	Let \(p_t^\varepsilon\) be the kernel of
	\(e^{t\Delta_\varepsilon}\), with respect to
	\(\varepsilon^d\)-counting measure.  Then \(p_t^\varepsilon\geq0\),
	\[
	\int_{\mathbb R_\varepsilon^d}p_t^\varepsilon(x)\,\mathd x=1,
	\qquad
	\int_{\mathbb R_\varepsilon^d}|x|^2
	p_t^\varepsilon(x)\,\mathd x=2dt,
	\]
	and
	\begin{equation}
		\label{eq:sublattice-heat-Linfty}
		\|p_t^\varepsilon\|_{L^\infty}
		\lesssim(\sqrt t+\varepsilon)^{-d}.
	\end{equation}
	For every \(q\geq2\),
	\begin{equation}
		\label{eq:sublattice-heat-moments}
		\int_{\mathbb R_\varepsilon^d}|x|^q
		p_t^\varepsilon(x)\,\mathd x
		\lesssim_q t^{q/2}+t\varepsilon^{q-2}.
	\end{equation}
	Moreover, for every \(\tau\in(0,1]\) and
	\(i\in\{1,\ldots,d\}\),
	\begin{equation}
		\label{eq:sublattice-heat-gradient}
		\begin{aligned}
			&\|\partial^{i\pm}p_t^\varepsilon\|
			_{\operatorname{TV}(1+|\Cdot|^2/\tau^2)}
			\\
			&\quad\lesssim
			(\sqrt t+\varepsilon)^{-1}
			\left(1+\frac{t+\varepsilon^2}{\tau^2}\right)
			+\frac{\sqrt t}{\tau^2}.
		\end{aligned}
	\end{equation}
	Finally, for \(p\in[1,\infty)\) and \(\omega\geq0\), there exists
	\(C=C(p,\omega,d)\) such that
	\begin{equation}
		\label{eq:sublattice-heat-weighted-Lp}
		\|(1+|\Cdot|)^\omega p_t^\varepsilon\|_{L^p}
		\lesssim
		(\sqrt t+\varepsilon)^{-d(1-1/p)}(1+t)^C.
	\end{equation}
\end{lemma}

\begin{proof}
	Positivity and conservation of mass follow from the Markov property.
	If \(X_t^\varepsilon\) denotes the continuous-time nearest-neighbour
	random walk with generator \(\Delta_\varepsilon\), then
	\[
	\Delta_\varepsilon|x|^2=2d.
	\]
	Dynkin's formula therefore gives
	\(\mathbb E|X_t^\varepsilon|^2=2dt\), proving the exact second-moment
	identity.
	
	The Fourier representation is
	\[
	p_t^\varepsilon(x)
	=
	\int_{(-\pi/\varepsilon,\pi/\varepsilon]^d}
	e^{-tq_\varepsilon(k)^2}e^{ik\cdot x}
	\frac{\mathd k}{(2\pi)^d}.
	\]
	On the Brillouin zone,
	\(q_\varepsilon(k)^2\gtrsim|k|^2\).  Integrating the Gaussian gives
	\(Ct^{-d/2}\), while bounding the integrand by one gives
	\(C\varepsilon^{-d}\).  Taking the minimum proves
	\eqref{eq:sublattice-heat-Linfty}.
	
	Let \(N_t\) be the number of jumps up to time \(t\).  It is Poisson with
	parameter \(2dt/\varepsilon^2\), and the quadratic variation of the
	martingale \(X^\varepsilon\) satisfies
	\([X^\varepsilon]_t=\varepsilon^2N_t\).  By the
	Burkholder--Davis--Gundy inequality and the standard Poisson moment
	estimate,
	\[
	\mathbb E|X_t^\varepsilon|^q
	\lesssim_q
	\varepsilon^q\mathbb E N_t^{q/2}
	\lesssim_q
	t^{q/2}+t\varepsilon^{q-2}.
	\]
	This proves \eqref{eq:sublattice-heat-moments}.
	
	We prove \eqref{eq:sublattice-heat-gradient} for the forward difference.
	In one dimension let \(P_t(n)\) be the probability that
	\(X_t^\varepsilon=n\varepsilon\).  The sequence \(P_t(n)\) is symmetric
	and non-increasing for \(n\geq0\).  Discrete summation by parts gives
	\[
	\sum_{n\in\mathbb Z}
	\left(1+\frac{n^2\varepsilon^2}{\tau^2}\right)
	|P_t(n+1)-P_t(n)|
	\lesssim
	P_t(0)\left(1+\frac{\varepsilon^2}{\tau^2}\right)
	+\frac{\varepsilon^2}{\tau^2}
	\sum_n|n|P_t(n).
	\]
	The on-diagonal estimate yields
	\(P_t(0)\lesssim\varepsilon/(\sqrt t+\varepsilon)\), and
	Cauchy--Schwarz yields
	\[\sum_n|n|P_t(n)\lesssim\sqrt t/\varepsilon.\]
	Dividing by \(\varepsilon\), as required by the discrete difference,
	proves the one-dimensional estimate.  In dimension \(d\), the walk is a
	product of one-dimensional walks.  The coordinates different from \(i\)
	contribute their second moments, of total size \(O(t/\tau^2)\), multiplied
	by the unweighted variation in the \(i\)-th coordinate.  This gives
	\eqref{eq:sublattice-heat-gradient}.  The backward difference is identical.
	
	For \eqref{eq:sublattice-heat-weighted-Lp}, choose an even integer
	\(q\geq\omega p\).  Positivity gives, for \(1\leqslant p<\infty\),
	\[
	\begin{aligned}
		\|(1+|\Cdot|)^\omega p_t^\varepsilon\|_{L^p}^p
		&\leq
		\|p_t^\varepsilon\|_{L^\infty}^{p-1}
		\int(1+|x|)^{\omega p}p_t^\varepsilon(x)\,\mathd x
		\\
		&\lesssim
		(\sqrt t+\varepsilon)^{-d(p-1)}(1+t)^C.
	\end{aligned}
	\]
	Taking the \(p\)-th root concludes the proof.
\end{proof}

\begin{lemma}
	\label{lem:K-kernel}
	Let \(k\in\mathbb N_0\) and
	\(A\in\mathbb N_0^{\{0,1\pm,\ldots,d\pm\}}\) satisfy
	\(A_0\in\{0,1\}\) and \(|\bar A|\leq2\).  Then
	\begin{equation}
		\label{eq:K-kernel-corrected-derivative}
		\|\partial^AK_\mu\|_{\operatorname{TV}(w_\mu^2)}
		\lesssim
		\llbracket\mu\rrbracket^{-|A|_s}
		\left(
		1+\frac{\varepsilon}{\llbracket\mu\rrbracket}
		\right)^{\mathbbm{1}_{|\bar A|=1}}.
	\end{equation}
	Moreover, uniformly in
	\(\mu\in[0,1)\) and \(\eta\in[\mu,1)\),
	\begin{align}
		\|K_{\mu,\eta}\|_{\operatorname{TV}(w_\mu^2)}
		&\lesssim1,
		\label{eq:K-kernel-cross-corrected}
		\\
		\|(L_\eta-1)K_\mu\|_{\operatorname{TV}(w_\mu^2)}
		&\lesssim
		\llbracket\mu\rrbracket^{-2s}
		\llbracket\eta\rrbracket^{2s},
		\label{eq:K-kernel-Lminus1-corrected}
		\\
		\|(\tilde L_\eta^{1,k}-1)\tilde K_\mu^{1,k}\|
		_{\operatorname{TV}((w_\mu^2)^{\otimes(1+k)})}
		&\lesssim
		\llbracket\mu\rrbracket^{-2s}
		\llbracket\eta\rrbracket^{2s}.
		\label{eq:K-kernel-tensor-corrected}
	\end{align}
\end{lemma}

\begin{proof}
	Set \(\tau=\llbracket\mu\rrbracket\). Using the notation introduced in Lemma~\ref{lem:K-kernel-aux} below, we write
	\[
		K_\mu=\mathring K_\tau\otimes(\bar K_\tau*\bar K_\tau),
	\]
	and observe that
	\[
	w_\mu^2(t,x)
	\lesssim
		\mathring w_\tau(t)\bar w_\tau(x).
	\]
	The estimates with no spatial difference and with exactly one spatial
	difference follow from weighted Young's inequality and
	Lemma~\ref{lem:K-kernel-aux}.  A time derivative costs
	\(\tau^{-2s}\).
	
	Consider now two spatial differences.  Since difference operators commute
	with convolution, every such term can be written as a finite linear
	combination of kernels
	\[
	(\partial^{i,\square}\bar K_\tau)
	*
	(\partial^{j,\triangle}\bar K_\tau),
	\qquad
	\square,\triangle\in\{+,-\}.
	\]
	The weight satisfies
	\[
	\bar w_\tau(x+y)
	\lesssim
	\bar w_\tau(x)+\bar w_\tau(y).
	\]
	Weighted Young's inequality therefore yields
	\[
	\begin{aligned}
		&\tau^2
		\|
		(\partial^{i,\square}\bar K_\tau)
		*
		(\partial^{j,\triangle}\bar K_\tau)
		\|_{\operatorname{TV}(\bar w_\tau)}
		\\
		&\quad\lesssim
		\bigl[
		\tau\|\partial^{i,\square}\bar K_\tau\|
		_{\operatorname{TV}(\bar w_\tau)}
		\bigr]
		\bigl[
		\tau\|\partial^{j,\triangle}\bar K_\tau\|
		_{\operatorname{TV}}
		\bigr]
		\\
		&\qquad+
		\bigl[
		\tau\|\partial^{i,\square}\bar K_\tau\|
		_{\operatorname{TV}}
		\bigr]
		\bigl[
		\tau\|\partial^{j,\triangle}\bar K_\tau\|
		_{\operatorname{TV}(\bar w_\tau)}
		\bigr].
	\end{aligned}
	\]
	By \eqref{eq:K-aux-gradient-unweighted} and
	\eqref{eq:K-aux-gradient-weighted}, each product on the right is bounded by
	\[
	\left(
	1+\frac{\varepsilon}{\tau}
	\right)
	\min\left\{1,\frac{\tau}{\varepsilon}\right\}
	\lesssim1.
	\]
	This proves \eqref{eq:K-kernel-corrected-derivative} for
	\(|\bar A|=2\), including mixed spacetime derivatives.
	
	We turn to the estimates involving two scales.  Since
	\(\eta\geq\mu\), one has
	\(\llbracket\eta\rrbracket\leq\llbracket\mu\rrbracket\).
	Expand \(K_{\mu,\eta}=L_\eta K_\mu\) into temporal and spatial factors.
	Every nonidentity temporal factor is controlled by
	\eqref{eq:K-aux-time-cross}, whereas every nonidentity spatial factor is a complete
	scaled Laplacian block and is controlled by
	\eqref{eq:K-aux-space-cross}.  This proves
	\eqref{eq:K-kernel-cross-corrected}.
	
	Similarly, every term in \((L_\eta-1)K_\mu\) contains either
	\(\llbracket\eta\rrbracket^{2s}\partial_t\) or
	\(\llbracket\eta\rrbracket^2\Delta_\varepsilon\).  The temporal contribution
	is bounded by
	\[
	\llbracket\eta\rrbracket^{2s}
	\llbracket\mu\rrbracket^{-2s}.
	\]
	The spatial contribution is bounded by
	\[
	\llbracket\eta\rrbracket^2
	\llbracket\mu\rrbracket^{-2}
	\leq
	\llbracket\eta\rrbracket^{2s}
	\llbracket\mu\rrbracket^{-2s},
	\]
	because \(s\leq1\) and
	\(\llbracket\eta\rrbracket/\llbracket\mu\rrbracket\leq1\).
	Mixed terms are smaller.  This proves
	\eqref{eq:K-kernel-Lminus1-corrected}.
	
	Finally, expand
	\(\tilde L_\eta^{1,k}-1\) as a finite sum of tensor products in which at
	least one factor is \(L_\eta-1\).  Apply
	\eqref{eq:K-kernel-Lminus1-corrected} to that factor and
	\eqref{eq:K-kernel-cross-corrected} to all remaining factors.  This proves
	\eqref{eq:K-kernel-tensor-corrected}.
\end{proof}

\begin{lemma}
	\label{lem:K-kernel-aux}
	For \(\tau\in(0,1]\) and \(\eta\in[\tau,1]\), define
	\begin{align*}
			\mathring L_\tau&\assign1+\tau^{2s}\partial_t,
		&
		\bar L_\tau&\assign1-\tau^2\Delta,
		\\
			\mathring K_\tau&\assign\mathring L_\tau^{-1},
		&
		\bar K_\tau&\assign\bar L_\tau^{-1},
		\\
			\mathring K_{\eta,\tau}&\assign\mathring L_\tau\mathring K_\eta,
		&
		\bar K_{\eta,\tau}&\assign\bar L_\tau\bar K_\eta,
	\end{align*}
	and
	\[
		\mathring w_\tau(t)\assign1+\frac{|t|^2}{\tau^{4s}},
	\qquad
	\bar w_\tau(x)\assign1+\frac{|x|^2}{\tau^2}.
	\]
	Then the following estimates hold uniformly in
	\(\varepsilon\in2^{-\mathbb N_0}\).
	
	\begin{enumerate}[(a)]
		\item Uniformly in \(\tau\in(0,1]\),
		\begin{align}
				\|\mathring K_\tau\|_{\operatorname{TV}(\mathring w_\tau)}
				+
				\tau^{2s}\|\partial_t\mathring K_\tau\|
				_{\operatorname{TV}(\mathring w_\tau)}
			&\lesssim1,
			\label{eq:K-aux-time}
			\\
			\|\bar K_\tau\|_{\operatorname{TV}(\bar w_\tau)}
			+
			\tau^2\|\Delta\bar K_\tau\|
			_{\operatorname{TV}(\bar w_\tau)}
			&\lesssim1.
			\label{eq:K-aux-space-zero-two}
		\end{align}
		The first spatial differences satisfy
		\begin{align}
			\tau\|\partial^{i\pm}\bar K_\tau\|_{\operatorname{TV}}
			&\lesssim
			\min\left\{1,\frac{\tau}{\varepsilon}\right\},
			\label{eq:K-aux-gradient-unweighted}
			\\
			\tau\|\partial^{i\pm}\bar K_\tau\|
			_{\operatorname{TV}(\bar w_\tau)}
			&\lesssim
			1+\frac{\varepsilon}{\tau}
			\lesssim
			\frac{\varepsilon\vee\tau}{\tau}.
			\label{eq:K-aux-gradient-weighted}
		\end{align}
		
		\item Uniformly in \(\tau\in(0,1]\) and \(\eta\in[\tau,1]\),
		\begin{align}
				\|(1-\mathring L_\tau)\mathring K_\eta\|
				_{\operatorname{TV}(\mathring w_\eta)}
			&\lesssim\tau^{2s}\eta^{-2s},
			\label{eq:K-aux-time-cross}
			\\
			\|(1-\bar L_\tau)\bar K_\eta\|
			_{\operatorname{TV}(\bar w_\eta)}
			&\lesssim\tau^2\eta^{-2}
			\leq\tau^{2s}\eta^{-2s}.
			\label{eq:K-aux-space-cross}
		\end{align}
		
		\item Uniformly in the same range,
		\begin{equation}
				\|\mathring K_{\eta,\tau}\|_{\operatorname{TV}(\mathring w_\eta)}
			+
			\|\bar K_{\eta,\tau}\|_{\operatorname{TV}(\bar w_\eta)}
			\lesssim1.
			\label{eq:K-aux-cross-kernels}
		\end{equation}
	\end{enumerate}
\end{lemma}

\begin{proof}
	The temporal kernel is explicit:
	\[
		\mathring K_\tau(t)
	=
	\tau^{-2s}\mathbbm{1}_{\{t\geq0\}}
	e^{-t/\tau^{2s}}.
	\]
	The first term in \eqref{eq:K-aux-time} follows by the change of variables
	\(t=\tau^{2s}u\).  The distributional derivative is
	\[
		\partial_t\mathring K_\tau
	=
	\tau^{-2s}\delta_0
	-
	\tau^{-4s}\mathbbm{1}_{\{t\geq0\}}
	e^{-t/\tau^{2s}},
	\]
	and the same change of variables proves the second term.
	
	For the spatial kernel, the resolvent formula gives
	\begin{equation}
		\label{eq:K-aux-resolvent}
		\bar K_\tau
		=
		\int_0^\infty e^{-\lambda}
		e^{\lambda\tau^2\Delta_\varepsilon}\,\mathd\lambda,
		\qquad
		\bar K_\tau(x)
		=
		\int_0^\infty e^{-\lambda}
		p_{\lambda\tau^2}^\varepsilon(x)\,\mathd\lambda.
	\end{equation}
	The kernel is nonnegative and has mass one.  Using the exact
	second-moment identity from Lemma~\ref{lem:sublattice-heat-kernel},
	\[
	\begin{aligned}
		\|\bar K_\tau\|_{\operatorname{TV}(\bar w_\tau)}
		&=
		\int_0^\infty e^{-\lambda}
		\left(
		1+\tau^{-2}
		\int|x|^2p_{\lambda\tau^2}^\varepsilon(x)\,\mathd x
		\right)\mathd\lambda
		\\
		&=
		\int_0^\infty e^{-\lambda}(1+2d\lambda)\,\mathd\lambda
		\lesssim1.
	\end{aligned}
	\]
	Since
	\[
	(1-\tau^2\Delta_\varepsilon)\bar K_\tau=\delta_0,
	\qquad
	\tau^2\Delta_\varepsilon\bar K_\tau
	=
	\bar K_\tau-\delta_0,
	\]
	the scaled-Laplacian estimate in
	\eqref{eq:K-aux-space-zero-two} follows.
	
	We prove the weighted first-difference estimate.  If
	\(\tau\geq\varepsilon\), Minkowski's inequality,
	\eqref{eq:K-aux-resolvent}, and
	\eqref{eq:sublattice-heat-gradient} give
	\[
	\begin{aligned}
		&\tau\|\partial^{i\pm}\bar K_\tau\|
		_{\operatorname{TV}(\bar w_\tau)}
		\\
		&\quad\lesssim
		\int_0^\infty e^{-\lambda}
		\left[
		\frac{\tau}{\tau\sqrt\lambda+\varepsilon}
		\left(
		1+\lambda+\frac{\varepsilon^2}{\tau^2}
		\right)
		+\sqrt\lambda
		\right]\mathd\lambda.
	\end{aligned}
	\]
	On \(0<\lambda<(\varepsilon/\tau)^2\), use
	\[
	\frac{\tau}{\tau\sqrt\lambda+\varepsilon}
	\leq\frac{\tau}{\varepsilon};
	\]
	the length of this interval is
	\((\varepsilon/\tau)^2\).  On its complement, use
	\[
	\frac{\tau}{\tau\sqrt\lambda+\varepsilon}
	\leq\lambda^{-1/2}.
	\]
	Both integrals are bounded uniformly, so
	\[
	\tau\|\partial^{i\pm}\bar K_\tau\|
	_{\operatorname{TV}(\bar w_\tau)}
	\lesssim1
	\qquad\text{when }\tau\geq\varepsilon.
	\]
	
	Suppose now that \(\tau<\varepsilon\).  Positivity implies
	\[
	|\partial^{i+}\bar K_\tau(x)|
	\leq
	\varepsilon^{-1}
	\bigl(
	\bar K_\tau(x+\varepsilon e_i)+\bar K_\tau(x)
	\bigr).
	\]
	Furthermore,
	\[
	\bar w_\tau(y-\varepsilon e_i)
	\lesssim
	\left(1+\frac{\varepsilon^2}{\tau^2}\right)
	\bar w_\tau(y).
	\]
	Consequently,
	\[
	\begin{aligned}
		\tau\|\partial^{i+}\bar K_\tau\|
		_{\operatorname{TV}(\bar w_\tau)}
		&\lesssim
		\frac{\tau}{\varepsilon}
		\left(1+\frac{\varepsilon^2}{\tau^2}\right)
		\|\bar K_\tau\|_{\operatorname{TV}(\bar w_\tau)}
		\\
		&\lesssim
		1+\frac{\varepsilon}{\tau}.
	\end{aligned}
	\]
	The backward difference is identical.  This proves
	\eqref{eq:K-aux-gradient-weighted}.
	
	For the unweighted estimate, the preceding diffusive argument gives
	\(O(1)\) when \(\tau\geq\varepsilon\).  When
	\(\tau<\varepsilon\), positivity and mass one give
	\[
	\|\partial^{i\pm}\bar K_\tau\|_{\operatorname{TV}}
	\leq2\varepsilon^{-1}.
	\]
	Multiplication by \(\tau\) proves
	\eqref{eq:K-aux-gradient-unweighted}.
	
	Finally,
	\[
		(1-\mathring L_\tau)\mathring K_\eta
	=
		-\tau^{2s}\partial_t\mathring K_\eta,
	\]
	and the explicit temporal kernel at scale \(\eta\) gives
	\eqref{eq:K-aux-time-cross}.  Similarly,
	\[
	(1-\bar L_\tau)\bar K_\eta
	=
	\tau^2\Delta_\varepsilon\bar K_\eta
	=
	\frac{\tau^2}{\eta^2}
	(\bar K_\eta-\delta_0),
	\]
	which proves \eqref{eq:K-aux-space-cross}.  Since
	\[
		\mathring K_{\eta,\tau}
		=
		\mathring K_\eta-(1-\mathring L_\tau)\mathring K_\eta,
	\qquad
	\bar K_{\eta,\tau}
	=
	\bar K_\eta-(1-\bar L_\tau)\bar K_\eta,
	\]
	\eqref{eq:K-aux-cross-kernels} follows from the preceding estimates.
\end{proof}

\begin{remark}
	The factor \(1+\varepsilon/\tau\) in
	\eqref{eq:K-aux-gradient-weighted} is necessary.  Indeed, when
	\(\tau/\varepsilon\) is small,
	\[
	\bar K_\tau(0)
	=
	\varepsilon^{-d}
	\left(
	1+O(\tau^2/\varepsilon^2)
	\right),
	\qquad
	\bar K_\tau(-\varepsilon e_i)
	=
	O\left(
	\varepsilon^{-d}\tau^2/\varepsilon^2
	\right).
	\]
	Therefore
	\[
	|\partial^{i+}\bar K_\tau(-\varepsilon e_i)|
	\gtrsim\varepsilon^{-d-1}.
	\]
	Retaining this single nonnegative summand in the weighted
	total-variation norm gives
	\[
	\tau\|\partial^{i+}\bar K_\tau\|
	_{\operatorname{TV}(\bar w_\tau)}
	\gtrsim\frac{\varepsilon}{\tau}.
	\]
\end{remark}

\begin{lemma}
	\label{lem:check-K-factor}
	Let \(N\in\mathbb N_+\), set
	\[
	p\assign\frac{2N}{2N-1},
	\]
	and assume
	\[
	\hat\kappa>
	\max\left\{\frac1{2N},\frac d{4N}\right\}.
	\]
	Recall Definition~\ref{def:K_hat_check} of \(\check K_\mu\). Then
	\begin{equation}
		\label{eq:check-K-piecewise-bound}
		\|\check K_\mu\|_{L^p(\zeta^{-1})}
		\lesssim
		\llbracket\mu\rrbracket^{-(d+2s)/(2N)}
	\end{equation}
	uniformly in \(\mu\in(0,1)\) and
	\(\varepsilon\in2^{-\mathbb N_0}\).
\end{lemma}

\begin{proof}
	Put \(\tau=\llbracket\mu\rrbracket\).  The temporal kernel of
	\((1+\tau^{2s}\partial_t)^{-\hat\kappa}\) is
	\[
	\breve K_\tau(t)
	=
	\frac{\mathbbm{1}_{\{t>0\}}}{\Gamma(\hat\kappa)}
	\tau^{-2s\hat\kappa}
	t^{\hat\kappa-1}e^{-t/\tau^{2s}}
	=
	\tau^{-2s}\breve K_1(t/\tau^{2s}).
	\]
	For every \(\omega\geq0\),
	\[
	\|(1+|\Cdot|)^\omega\breve K_\tau\|_{L^p(\mathbb R)}
	\lesssim_\omega
	\tau^{-2s(1-1/p)},
	\]
	provided
	\(p(\hat\kappa-1)>-1\), equivalently
	\(\hat\kappa>1-1/p=1/(2N)\).
	
	Suppose first that \(\tau\geq\varepsilon\).  The spatial factor is
	\[
	(1-\tau^2\Delta_\varepsilon)^{-\hat\kappa}
	=
	\frac1{\Gamma(\hat\kappa)}
	\int_0^\infty
	e^{-\lambda}
	e^{\lambda\tau^2\Delta_\varepsilon}
	\lambda^{\hat\kappa-1}\,\mathd\lambda.
	\]
	By Minkowski's inequality and
	\eqref{eq:sublattice-heat-weighted-Lp},
	\[
	\begin{aligned}
		&\|(1+|\Cdot|)^\omega
		(1-\tau^2\Delta_\varepsilon)^{-\hat\kappa}
		\|_{L^p}
		\\
		&\quad\lesssim
		\int_0^\infty
		e^{-\lambda}\lambda^{\hat\kappa-1}
		(\tau\sqrt\lambda+\varepsilon)^{-d(1-1/p)}
		(1+\tau^2\lambda)^C\,\mathd\lambda.
	\end{aligned}
	\]
	Since \(\tau\leq1\), the last polynomial factor is bounded by
	\((1+\lambda)^C\).  Factoring out \(\tau^{-d(1-1/p)}\) and setting
	\(\delta=\varepsilon/\tau\in(0,1]\), it remains to control
	\[
	\int_0^\infty
	e^{-\lambda}(1+\lambda)^C
	\lambda^{\hat\kappa-1}
	(\sqrt\lambda+\delta)^{-d(1-1/p)}\,\mathd\lambda.
	\]
	Split the integral into
	\((0,\delta^2)\), \((\delta^2,1)\), and \((1,\infty)\).  On the first interval the
	integral is bounded by \(C\delta^{2\hat\kappa-d(1-1/p)}\); on the second it is
	bounded by
	\[
	C\int_{\delta^2}^1
	\lambda^{\hat\kappa-1-d(1-1/p)/2}\,\mathd\lambda;
	\]
	the last interval is controlled by the exponential.  All three terms are
	uniformly bounded precisely when
	\(2\hat\kappa>d(1-1/p)\), i.e.
	\(\hat\kappa>d/(4N)\).  Hence
	\[
	\|(1+|\Cdot|)^\omega
	(1-\tau^2\Delta_\varepsilon)^{-\hat\kappa}
	\|_{L^p}
	\lesssim_\omega
	\tau^{-d/(2N)}.
	\]
	
	Suppose now that \(\tau<\varepsilon\).  The spatial factor is
	\[
	(1-\tau^2\Delta_\varepsilon)^{-2}.
	\]
	It is a positive probability kernel.  Its \(L^\infty\)-norm is bounded by
	\(C\varepsilon^{-d}\), because its Fourier multiplier is bounded by one.
	Interpolation between \(L^1\) and \(L^\infty\) gives
	\[
	\|(1-\tau^2\Delta_\varepsilon)^{-2}\|_{L^p}
	\lesssim
	\varepsilon^{-d(1-1/p)}
	\leq
	\tau^{-d(1-1/p)}.
	\]
	To include polynomial weights, use
	\[
	(1-\tau^2\Delta_\varepsilon)^{-2}
	=
	\int_0^\infty e^{-\lambda}\lambda
	e^{\lambda\tau^2\Delta_\varepsilon}\,\mathd\lambda.
	\]
	By \eqref{eq:sublattice-heat-moments}, for every fixed \(q\),
	\[
	\int(1+|x|)^q
	(1-\tau^2\Delta_\varepsilon)^{-2}(x)\,\mathd x
	\lesssim_q 1
	\]
	uniformly when \(\tau<\varepsilon\leq1\).  Positivity therefore gives
	\[
	\begin{aligned}
		&\|(1+|\Cdot|)^\omega
		(1-\tau^2\Delta_\varepsilon)^{-2}
		\|_{L^p}^p
		\\
		&\quad\leq
		\|(1-\tau^2\Delta_\varepsilon)^{-2}\|_{L^\infty}^{p-1}
		\int(1+|x|)^{\omega p}
		(1-\tau^2\Delta_\varepsilon)^{-2}(x)\,\mathd x
		\\
		&\quad\lesssim
		\varepsilon^{-d(p-1)}
		\leq
		\tau^{-d(p-1)}.
	\end{aligned}
	\]
	Thus the same spatial estimate holds below the lattice scale.
	
	Choose \(\omega\) large enough that
	\[
	\zeta(t,x)^{-1}
	\lesssim
	(1+|t|)^\omega(1+|x|)^\omega.
	\]
	The temporal and spatial kernels factorise, so Fubini gives
	\[
	\|\check K_\mu\|_{L^p(\zeta^{-1})}
	\lesssim
	\tau^{-2s/(2N)}\tau^{-d/(2N)}
	=
	\tau^{-(d+2s)/(2N)}.
	\]
	This finishes the proof.
\end{proof}

\begin{lemma}
	\label{lem:hat-K-factor}
	Recall that \(d=3\).  Let \(m\in\mathbb N_+\) and
	\begin{equation}
		\label{eq:cal-K-factor}
		\mathcal K_\mu(z_1,\ldots,z_m)
		\assign
		\int_{\Lambda_\varepsilon}
		\hat K_\mu(z)\hat K_\mu(z_1+z)\cdots
		\hat K_\mu(z_m+z)\,\mathd z,
	\end{equation}
	where $\hat K_\mu$ is the kernel introduced in Definition~\ref{def:K_hat_check}.
	If
	\[
	0\leq\hat\kappa<\frac1{m+1},
	\]
	then \(\mathcal K_\mu\) is continuous in the time variables and, for every
	\(\omega\geq0\),
	\begin{equation}
		\label{eq:hat-K-piecewise-bound}
		\|
		(w_\mu^\omega)^{\otimes m}\mathcal K_\mu
		\|_{L^\infty}
		\lesssim
		\llbracket\mu\rrbracket^{-2sm}
		(\varepsilon\vee\llbracket\mu\rrbracket)^{-dm},
	\end{equation}
	uniformly in \(\mu\in(0,1)\) and
	\(\varepsilon\in2^{-\mathbb N_0}\).
\end{lemma}

\begin{proof}
	Put
	$
	\tau\assign\llbracket\mu\rrbracket
	$.
	The kernel \(\hat K_\mu\) factorises into its temporal and spatial
	parts:
	\[
	\hat K_\mu(t,x)=H_\tau(t)\bar H_{\tau,\varepsilon}(x).
	\]
	Accordingly, the kernel in
	\eqref{eq:cal-K-factor} factorises as
	\[
	\mathcal K_\mu(t_1,x_1,\ldots,t_m,x_m)
	=
	\mathcal {\mathring K}_\mu(t_1,\ldots,t_m)
	\bar{\mathcal K}_\mu(x_1,\ldots,x_m).
	\]
	
	The temporal kernel is the kernel of
	\((1+\tau^{2s}\partial_t)^{-(1-\hat\kappa)}\), namely
	\[
	H_\tau(t)
	=
	\frac{\mathbbm{1}_{\{t>0\}}}{\Gamma(1-\hat\kappa)}
	\tau^{-2s(1-\hat\kappa)}
	t^{-\hat\kappa}e^{-t/\tau^{2s}}
	=
	\tau^{-2s}H_1(t/\tau^{2s}).
	\]
	Since \((m+1)\hat\kappa<1\), one has
	\(H_1\in L^{m+1}(\mathbb R)\). Hence Hölder's inequality gives
	\[
	\sup_{t_1,\ldots,t_m\in\mathbb R}
	\int_{\mathbb R}
	|H_1(t)H_1(t+t_1)\cdots H_1(t+t_m)|
	\,\mathd t
	<\infty.
	\]
	The exponential decay of \(H_1\), together with the
	submultiplicativity of polynomial weights, gives, for every
	\(\omega\geq0\),
	\[
	\sup_{t_1,\ldots,t_m\in\mathbb R}
	\prod_{i=1}^m
	\left(
	1+\frac{|t_i|^{1/(2s)}}{\tau}
	\right)^\omega
	\left|
	\mathcal {\mathring K}_\mu(t_1,\ldots,t_m)
	\right|
	\lesssim_\omega
	\tau^{-2sm}.
	\]
	Continuity in \(t_1,\ldots,t_m\) follows from the continuity of
	translations in \(L^{m+1}(\mathbb R)\).
	
	We next estimate the spatial correlation kernel. Suppose first that
	\(\tau\geq\varepsilon\), and put
	\[
	\delta\assign\frac{\varepsilon}{\tau}\in(0,1].
	\]
	The spatial part of \(\hat K_\mu\) is
	\[
	(1-\tau^2\Delta_\varepsilon)^{-(2-\hat\kappa)}.
	\]
	Its kernel has the Fourier representation
	\begin{align*}
		\bar H_{\tau,\varepsilon}(x)
		&=
		\int_{\mathbb T_\varepsilon^d}
		e^{i\tilde k\cdot x}
		(1+\tau^2q_\varepsilon(\tilde k)^2)^{-(2-\hat\kappa)}
		\,\mathd\tilde k
		\\
		&=
		\tau^{-d}
		\int_{\mathbb T_\delta^d}
		e^{i\tilde q\cdot(x/\tau)}
		a_\delta(\tilde q)
		\,\mathd\tilde q,
	\end{align*}
	where
	\[
	a_\delta(\tilde q)
	\assign
	(1+q_\delta(\tilde q)^2)^{-(2-\hat\kappa)}.
	\]
	By the strengthened \(\beta=2\) estimate in
	Lemma~\ref{lem:q-epsilon-bounds}, for every multi-index \(A\in\mathbb N_0^d\),
	\begin{equation}
		\label{eq:hat-K-symbol-L1-factor}
		\sup_{\delta\in(0,1]}
		\|
		\partial_{\tilde q}^A a_\delta
		\|_{L^1(\mathbb T_\delta^d)}
		<\infty.
	\end{equation}
	Thus, for every \(\omega\in\mathbb N_0\),
	\[
	\left(
	1+\frac{|x|^2}{\tau^2}
	\right)^\omega
	\bar H_{\tau,\varepsilon}(x)
	=
	\tau^{-d}
	\int_{\mathbb T_\delta^d}
	e^{i\tilde q\cdot(x/\tau)}
	(1-\Delta_{\tilde q})^\omega a_\delta(\tilde q)
	\,\mathd\tilde q.
	\]
	Using \eqref{eq:hat-K-symbol-L1-factor}, we conclude that, for every
	\(\omega\geq0\),
	\begin{equation}
		\label{eq:hat-K-spatial-decay-factor}
		|
		\bar H_{\tau,\varepsilon}(x)
		|
		\lesssim_\omega
		\tau^{-d}
		\left(
		1+\frac{|x|}{\tau}
		\right)^{-\omega},
		\qquad
		\tau\geq\varepsilon.
	\end{equation}
	
	Let
	\[
	\bar{\mathcal K}_\mu(x_1,\ldots,x_m)
	\assign
	\int_{\mathbb R_\varepsilon^d}
	\bar H_{\tau,\varepsilon}(x)
	\bar H_{\tau,\varepsilon}(x+x_1)\cdots
	\bar H_{\tau,\varepsilon}(x+x_m)
	\,\mathd x.
	\]
	For every \(\bar\omega>d\), the discrete Riemann-sum estimate
	\begin{equation}
		\label{eq:hat-K-lattice-sum-factor}
		\int_{\mathbb R_\varepsilon^d}
		\left(
		1+\frac{|x|}{\tau}
		\right)^{-\bar\omega}
		\,\mathd x
		\lesssim_{\bar\omega}
		\tau^d,
		\qquad
		\tau\geq\varepsilon,
	\end{equation}
	holds uniformly in \(\varepsilon\) and \(\tau\). Moreover,
	\[
	1+\frac{|x_i|}{\tau}
	\lesssim
	\left(
	1+\frac{|x|}{\tau}
	\right)
	\left(
	1+\frac{|x+x_i|}{\tau}
	\right).
	\]
	Choose \(\bar\omega>m\omega+d\). Inserting
	\eqref{eq:hat-K-spatial-decay-factor}, transferring the weights by
	the preceding inequality, and then using
	\eqref{eq:hat-K-lattice-sum-factor}, we obtain
	\begin{align*}
		&\prod_{i=1}^m
		\left(
		1+\frac{|x_i|}{\tau}
		\right)^\omega
		|
		\bar{\mathcal K}_\mu(x_1,\ldots,x_m)
		|
		\\
		&\quad\lesssim
		\tau^{-d(m+1)}
		\int_{\mathbb R_\varepsilon^d}
		\left(
		1+\frac{|x|}{\tau}
		\right)^{-\bar\omega+m\omega}
		\prod_{i=1}^m
		\left(
		1+\frac{|x+x_i|}{\tau}
		\right)^{-\bar\omega+\omega}
		\,\mathd x
		\\
		&\quad\lesssim
		\tau^{-d(m+1)}
		\int_{\mathbb R_\varepsilon^d}
		\left(
		1+\frac{|x|}{\tau}
		\right)^{-\bar\omega+m\omega}
		\,\mathd x
		\lesssim
		\tau^{-dm}.
	\end{align*}
	Hence
	\begin{equation}
		\label{eq:hat-K-spatial-correlation-above-factor}
		\|
		((1+|\Cdot|/\tau)^\omega)^{\otimes m}
		\bar{\mathcal K}_\mu
		\|_{L^\infty}
		\lesssim
		\tau^{-dm},
		\qquad
		\tau\geq\varepsilon.
	\end{equation}
	
	Suppose now that \(\tau<\varepsilon\). The spatial part of
	\(\hat K_\mu\) is the identity, whose kernel is
	\[
	\delta_0^\varepsilon(x)
	=
	\varepsilon^{-d}\mathbbm{1}_{\{x=0\}}.
	\]
	Its spatial correlation kernel is therefore
	\[
	\bar{\mathcal K}_\mu(x_1,\ldots,x_m)
	=
	\varepsilon^{-dm}
	\prod_{j=1}^m
	\mathbbm{1}_{\{x_j=0\}}.
	\]
	The spatial weights equal one on its support, and consequently
	\begin{equation}
		\label{eq:hat-K-spatial-correlation-below-factor}
		\|
		((1+|\Cdot|/\tau)^\omega)^{\otimes m}
		\bar{\mathcal K}_\mu
		\|_{L^\infty}
		=
		\varepsilon^{-dm},
		\qquad
		\tau<\varepsilon.
	\end{equation}
	
	Finally, by Definition~\ref{def:weight-w},
	\[
	w_\mu^\omega(t_i,x_i)
	\lesssim
	\left(
	1+\frac{|t_i|^{1/(2s)}}{\tau}
	\right)^\omega
	\left(
	1+\frac{|x_i|}{\tau}
	\right)^\omega.
	\]
	Combining the temporal estimate with
	\eqref{eq:hat-K-spatial-correlation-above-factor} and
	\eqref{eq:hat-K-spatial-correlation-below-factor} gives
	\[
	\|
	(w_\mu^\omega)^{\otimes m}\mathcal K_\mu
	\|_{L^\infty}
	\lesssim
	\tau^{-2sm}
	\begin{cases}
		\tau^{-dm},&\tau\geq\varepsilon,\\
		\varepsilon^{-dm},&\tau<\varepsilon.
	\end{cases}
	\]
	This is precisely \eqref{eq:hat-K-piecewise-bound}.
\end{proof}

\begin{lemma}
	\label{lem:J-decay}
	Let $\ell\in\mathbb N_0$. Put
	\[
	\tau
	\assign \llbracket\sigma\rrbracket,
	\qquad
	\rho
	\assign \tau\vee\varepsilon.
	\]
	The following estimates hold uniformly in
	\(\varepsilon\in2^{-\mathbb N_0}\),
	\(\sigma\in[1/2,1)\), and \((t,x)\in\Lambda_\varepsilon\).
	
	\begin{enumerate}[(a)]
		\item For every \(N\in\mathbb N_0\) and
		\(A\in\mathbb N_0^{\{0,1\pm,\ldots,d\pm\}}\),
		\begin{align}
			|\partial^A\tilde J_{\sigma,\ell}(t,x)|
			&\lesssim_{A,N,\ell}
			\tau^{-2s(1+A_0)}
			\rho^{-d-|\bar A|}
			\left(1+\frac{|t|}{\tau^{2s}}+\frac{|x|}{\rho}\right)^{-N},
			\label{eq:J-decay-corrected-pointwise}
			\\
			|\partial^A\partial_\sigma
			\tilde J_{\sigma,\ell}(t,x)|
			&\lesssim_{A,N,\ell}
			\llbracket\sigma\rrbracket^{-1}
			\tau^{-2s(1+A_0)}
			\rho^{-d-|\bar A|}
			\left(1+\frac{|t|}{\tau^{2s}}+\frac{|x|}{\rho}\right)^{-N}.
			\label{eq:J-decay-corrected-sigma}
		\end{align}
		
		\item Let \(a\in\mathbb N_0\), let \(\bar A\) be a spatial
		finite-difference multi-index, and let
		\(\omega\geq0\).  Define
		\[
		\mathfrak l_{\bar A,\omega}(u)
		\assign
		\begin{cases}
			1, & \bar A=0,\\
			(1+u)^{(\omega-|\bar A|)_+}, & \bar A\neq0.
		\end{cases}
		\]
		Then
		\begin{align}
			\|\partial_t^a\partial^{\bar A}
			\tilde J_{\sigma,\ell}\|_{\operatorname{TV}(w_\sigma^\omega)}
			&\lesssim_{a,\bar A,\omega,\ell}
			\llbracket\sigma\rrbracket^{-2sa-|\bar A|}
			\mathfrak l_{\bar A,\omega}
			\left(\frac{\varepsilon}{\llbracket\sigma\rrbracket}\right),
			\label{eq:J-decay-corrected-TV}
			\\
			\|\partial_t^a\partial^{\bar A}\partial_\sigma
			\tilde J_{\sigma,\ell}\|_{\operatorname{TV}(w_\sigma^\omega)}
			&\lesssim_{a,\bar A,\omega,\ell}
			\llbracket\sigma\rrbracket^{-2sa-|\bar A|-1}
			\mathfrak l_{\bar A,\omega}
			\left(\frac{\varepsilon}{\llbracket\sigma\rrbracket}\right).
			\label{eq:J-decay-corrected-sigma-TV}
		\end{align}
		In particular, there is no sublattice loss when \(\bar A=0\), nor when
		\(|\bar A|\geq1\) and \(\omega\leq |\bar A|\).
	\end{enumerate}
	Consequently, for every fixed \(m\in\mathbb N_0\),
	\begin{align}
		\|L_\sigma^m\tilde J_{\sigma,\ell}\|
		_{\operatorname{TV}(w_\sigma^2)}
		&\lesssim_{m,\ell}1,
		\label{eq:J-decay-corrected-Lm}
		\\
		\|\partial_t\tilde J_{\sigma,\ell}\|
		_{\operatorname{TV}(w_\sigma^2)}
		&\lesssim_\ell\llbracket\sigma\rrbracket^{-2s},
		\label{eq:J-decay-corrected-time}
		\\
		\|\nabla_\varepsilon\tilde J_{\sigma,\ell}\|
		_{\operatorname{TV}(w_\sigma^1)}
		&\lesssim_\ell\llbracket\sigma\rrbracket^{-1}.
		\label{eq:J-decay-corrected-gradient}
	\end{align}
	Since \(\zeta^{-1}\lesssim w_\sigma^1\), the last estimate implies
	\[
	\|\nabla_\varepsilon\tilde J_{\sigma,\ell}\|
	_{\operatorname{TV}(\zeta^{-1})}
	\lesssim_\ell\llbracket\sigma\rrbracket^{-1}.
	\]
\end{lemma}

\begin{proof}
	Fix \(\ell\in\mathbb N_0\), and set
	\[
	c_{\sigma,\ell}\assign 2^{-\ell}\sigma^{-1}.
	\]
	Thus the symbol of \(\tilde J_{\sigma,\ell}\) is
	\[
	j\bigl(c_{\sigma,\ell}\tau |k_0|^{1/(2s)}\bigr)
	j\bigl(c_{\sigma,\ell}\tau q_\varepsilon(k)\bigr).
	\]
	Since \(\sigma\in[1/2,1)\), one has
	\begin{equation}
		c_{\sigma,\ell}\asymp_\ell 1,
		\qquad
		c_{\sigma,\ell}\tau\asymp_\ell\tau,
		\qquad
		\left|
		\frac{\partial_\sigma(c_{\sigma,\ell}\tau)}
		{c_{\sigma,\ell}\tau}
		\right|
		=\frac1{\sigma\tau}
		\lesssim\tau^{-1}.
		\label{eq:J-decay-corrected-tau}
	\end{equation}
	
	Denote by \(\mathring J_\tau\) and
	\(\bar J_{\tau,\varepsilon}\) the temporal and spatial kernels whose
	Fourier multipliers are, respectively,
	\[
	k_0\longmapsto
	j\bigl(c_{\sigma,\ell}\tau |k_0|^{1/(2s)}\bigr),
	\qquad
	k\longmapsto
	j\bigl(c_{\sigma,\ell}\tau q_\varepsilon(k)\bigr).
	\]
	Here and below, the fixed index \(\ell\) is suppressed from the notation.
	With the usual convention of denoting a convolution operator and its
	kernel by the same symbol, one then has
	\begin{equation}
		\tilde J_{\sigma,\ell}(t,x)
		=\mathring J_\tau(t)\bar J_{\tau,\varepsilon}(x).
		\label{eq:J-decay-corrected-factorisation}
	\end{equation}
	
	The function
	\[
	k_0\longmapsto j\bigl(|k_0|^{1/(2s)}\bigr)
	\]
	is smooth and compactly supported, since it is constant in a
	neighbourhood of the origin.  Its inverse Fourier transform is therefore
	Schwartz.  Using \eqref{eq:J-decay-corrected-tau}, we obtain, for every
	\(a,N\in\mathbb N_0\),
	\begin{align}
		|\partial_t^a\mathring J_\tau(t)|
		&\lesssim_{a,N,\ell}
		\tau^{-2s(1+a)}
		\left(1+\frac{|t|}{\tau^{2s}}\right)^{-N},
		\label{eq:J-decay-corrected-time-proof}
		\\
		|\partial_t^a\partial_\sigma\mathring J_\tau(t)|
		&\lesssim_{a,N,\ell}
		\tau^{-2s(1+a)-1}
		\left(1+\frac{|t|}{\tau^{2s}}\right)^{-N}.
		\label{eq:J-decay-corrected-time-sigma-proof}
	\end{align}
	Indeed, differentiation with respect to \(\sigma\) replaces \(j(u)\)
	by a constant multiple of \(u j'(u)\), together with the logarithmic
	derivative appearing in \eqref{eq:J-decay-corrected-tau}.
	
	We next consider the spatial factor.  Suppose first
	that \(\tau\geq\varepsilon\), set
	\[
	\delta\assign\frac{\varepsilon}{\tau}\leq1,
	\]
	and make the change of variables \(\tilde q=\tau\tilde k\) in Fourier space.
	Then
	\[
	\partial^{\bar A}\bar J_{\tau,\varepsilon}(x)
	=\tau^{-d-|\bar A|}
	\int_{(\mathbb R_\delta^d)^*}
	d_\delta^{\bar A}(\tilde q)
	j\bigl(c_{\sigma,\ell}q_\delta(\tilde q)\bigr)
	e^{i\tilde q\cdot x/\tau}
	\frac{\mathd\tilde q}{(2\pi)^d}.
	\]
	On the support of the integrand,
	\[
	q_\delta(\tilde q)\lesssim_\ell1,
	\]
	and hence \(|\tilde q|\lesssim_\ell1\), uniformly in \(\delta\leq1\).
	The estimates of Lemma~\ref{lem:q-epsilon-bounds}, together with the fact
	that \(j\) is constant near the origin, show that the Fourier derivatives
	of the integrand have uniformly bounded \(L^1\)-norms.  Repeated
	integration by parts therefore gives
	\begin{equation}
		|\partial^{\bar A}\bar J_{\tau,\varepsilon}(x)|
		\lesssim_{\bar A,N,\ell}
		\tau^{-d-|\bar A|}
		\left(1+\frac{|x|}{\tau}\right)^{-N}.
		\label{eq:J-decay-corrected-space-above}
	\end{equation}
	
	Suppose next that \(\tau<\varepsilon\), set
	\[
	\lambda\assign\frac{\tau}{\varepsilon}\in(0,1),
	\]
	and rescale \(\tilde q=\varepsilon\tilde k\).  This gives
	\[
	\partial^{\bar A}\bar J_{\tau,\varepsilon}(x)
	=\varepsilon^{-d-|\bar A|}
	\int_{(-\pi,\pi]^d}
	d_1^{\bar A}(\tilde q)
	j\bigl(c_{\sigma,\ell}\lambda q_1(\tilde q)\bigr)
	e^{i\tilde q\cdot x/\varepsilon}
	\frac{\mathd\tilde q}{(2\pi)^d}.
	\]
	As \(\sigma\in[1/2,1)\) and \(\lambda\in[0,1]\), the multipliers in the
	last display form a bounded subset of
	\(C^\infty((\mathbb R/2\pi\mathbb Z)^d)\).  Integration by parts on the
	torus yields
	\begin{equation}
		|\partial^{\bar A}\bar J_{\tau,\varepsilon}(x)|
		\lesssim_{\bar A,N,\ell}
		\varepsilon^{-d-|\bar A|}
		\left(1+\frac{|x|}{\varepsilon}\right)^{-N}.
		\label{eq:J-decay-corrected-space-below}
	\end{equation}
	The same arguments after differentiation with respect to \(\sigma\)
	give an additional factor \(\tau^{-1}\), by
	\eqref{eq:J-decay-corrected-tau}.
	Combining \eqref{eq:J-decay-corrected-factorisation}--
	\eqref{eq:J-decay-corrected-space-below}, and recalling that
	\(\rho=\tau\vee\varepsilon\), proves
	\eqref{eq:J-decay-corrected-pointwise} and
	\eqref{eq:J-decay-corrected-sigma}.
	
	It remains to establish the weighted total-variation estimates.  The
	temporal bounds imply, for every \(\omega\geq0\),
	\begin{align}
		\int_{\mathbb R}
		\left(1+\frac{|t|^{1/(2s)}}{\tau}\right)^\omega
		|\partial_t^a\mathring J_\tau(t)|\,\mathd t
		&\lesssim_{a,\omega,\ell}\tau^{-2sa},
		\label{eq:J-decay-corrected-time-TV-proof}
		\\
		\int_{\mathbb R}
		\left(1+\frac{|t|^{1/(2s)}}{\tau}\right)^\omega
		|\partial_t^a\partial_\sigma
		\mathring J_\tau(t)|\,\mathd t
		&\lesssim_{a,\omega,\ell}\tau^{-2sa-1}.
		\label{eq:J-decay-corrected-time-sigma-TV-proof}
	\end{align}
	
	We now distinguish the deep sublattice regime.  Since
	\[
	q_\varepsilon(\tilde k)\leq2\sqrt d\,\varepsilon^{-1},
	\qquad
	c_{\sigma,\ell}\leq2^{1-\ell},
	\]
	there exists \(c_{d,\ell}>0\) such that
	\[
	j\bigl(c_{\sigma,\ell}\tau q_\varepsilon(\tilde k)\bigr)=1
	\quad\text{for every }\tilde k\in(\mathbb R_\varepsilon^d)^*
	\]
	whenever \(\tau\leq c_{d,\ell}\varepsilon\).  Consequently,
	\begin{equation}
		\bar J_{\tau,\varepsilon}=\delta_0^\varepsilon,
		\qquad
		\partial_\sigma\bar J_{\tau,\varepsilon}=0
		\qquad
		\text{if }\tau\leq c_{d,\ell}\varepsilon.
		\label{eq:J-decay-corrected-delta}
	\end{equation}
	
	If \(|\bar A|\geq1\), the kernel
	\(\partial^{\bar A}\delta_0^\varepsilon\) is supported in a ball of
	radius \(C_{\bar A}\varepsilon\), and
	\[
	\|\partial^{\bar A}\delta_0^\varepsilon\|_{\operatorname{TV}}
	\lesssim_{\bar A}\varepsilon^{-|\bar A|}.
	\]
	Hence, in the regime \(\tau\leq c_{d,\ell}\varepsilon\),
	\begin{align*}
		\|\partial^{\bar A}\delta_0^\varepsilon\|
		_{\operatorname{TV}((1+|\Cdot|/\tau)^\omega)}
		&\lesssim_{\bar A,\omega,\ell}
		\varepsilon^{-|\bar A|}
		\left(1+\frac{\varepsilon}{\tau}\right)^\omega
		\\
		&\lesssim_{\bar A,\omega,\ell}
		\tau^{-|\bar A|}
		\left(1+\frac{\varepsilon}{\tau}\right)^{(\omega-|\bar A|)_+}.
	\end{align*}
	For \(|\bar A|=0\), the corresponding weighted norm of
	\(\delta_0^\varepsilon\) is exactly one.
	
	In the complementary regime \(\tau>c_{d,\ell}\varepsilon\), one has
	\[
	\rho=\tau\vee\varepsilon\lesssim_{d,\ell}\tau.
	\]
	The pointwise bounds, with \(N\) sufficiently large, therefore imply
	\begin{align*}
		\|\partial^{\bar A}\bar J_{\tau,\varepsilon}\|
		_{\operatorname{TV}((1+|\Cdot|/\tau)^\omega)}
		&\lesssim_{\bar A,\omega,\ell}\tau^{-|\bar A|},
		\\
		\|\partial^{\bar A}\partial_\sigma
		\bar J_{\tau,\varepsilon}\|
		_{\operatorname{TV}((1+|\Cdot|/\tau)^\omega)}
		&\lesssim_{\bar A,\omega,\ell}\tau^{-|\bar A|-1}.
	\end{align*}
	
	Using the factorisation
	\eqref{eq:J-decay-corrected-factorisation}, the product rule for
	\(\partial_\sigma\), and
	\[
	w_\sigma^\omega(t,x)
	\lesssim_\omega
	\left(1+\frac{|t|^{1/(2s)}}{\tau}\right)^\omega
	\left(1+\frac{|x|}{\tau}\right)^\omega,
	\]
	the preceding estimates prove
	\eqref{eq:J-decay-corrected-TV} and
	\eqref{eq:J-decay-corrected-sigma-TV}.
	
	Finally, since \(\tau=\llbracket\sigma\rrbracket\),
	\[
	L_\sigma^m
	=(1+\tau^{2s}\partial_t)^m
	(1-\tau^2\Delta_\varepsilon)^{2m}.
	\]
	Every term in this expansion is a constant multiple of
	\[
	\tau^{2sa+2b}
	\partial_t^a(-\Delta_\varepsilon)^b.
	\]
	If \(b=0\), apply \eqref{eq:J-decay-corrected-TV} with \(|\bar A|=0\) and
	\(\omega=2\).  If \(b\geq1\), write
	\((-\Delta_\varepsilon)^b\) as a finite linear combination of spatial
	finite differences of order \(2b\), and apply the same estimate with
	\(|\bar A|=2b\) and \(\omega=2\).  Since \(2\leq2b\), no sublattice factor is
	lost.  This proves \eqref{eq:J-decay-corrected-Lm}.
	
	The estimates \eqref{eq:J-decay-corrected-time} and
	\eqref{eq:J-decay-corrected-gradient} follow from
	\eqref{eq:J-decay-corrected-TV} with
	\[
	(a,|\bar A|,\omega)=(1,0,2)
	\qquad\text{and}\qquad
	(a,|\bar A|,\omega)=(0,1,1),
	\]
	respectively.
\end{proof}

\begin{lemma}
	\label{lem:Gdot-decay}
	Assume \(d\geq1\) and \(s\in(3/4,1]\).  Let
	\(A\in\mathbb N_0^{\{0,1\pm,\ldots,d\pm\}}\),
	\(\flat\in(0,2s)\), and let \(p\in\mathbb N_0\) be fixed.  Then,
	uniformly in \(\varepsilon\in2^{-\mathbb N_0}\),
	\(\sigma\in[1/2,1)\), and every admissible spatial period \(M\geq1\),
	\begin{align}
		\|\partial^A\dot G_\sigma\|_{\operatorname{TV}(\zeta^{-1})}
		&\lesssim_A
		\llbracket\sigma\rrbracket^{2s-1-|A|_s},
		\label{eq:Gdot-corrected-fixed-weight}
		\\
		\|L_\sigma^p\dot G_\sigma\|_{\operatorname{TV}(w_\sigma^\flat)}
		&\lesssim_{p,\flat}
		\llbracket\sigma\rrbracket^{2s-1},
		\label{eq:Gdot-corrected-weighted-TV}
		\\
		\left\|
		\left|(L_\sigma^p\dot G_\sigma)w_\sigma^\flat\right|_M
		\right\|_{L^\infty}
		&\lesssim_{p,\flat}
		\llbracket\sigma\rrbracket^{-d-1},
		\label{eq:Gdot-corrected-periodised-Linfty}
		\\
		\left\|\left|L_\sigma^p\dot G_\sigma\right|_M\right\|_{L^\infty}
		&\lesssim_p
		\llbracket\sigma\rrbracket^{-d-1}
		\wedge
		\varepsilon^{-d}\llbracket\sigma\rrbracket^{-1}.
		\label{eq:Gdot-corrected-periodised-unweighted}
	\end{align}
	Moreover,
	\begin{equation}
		\|\partial^A G_{1/2}\|_{\operatorname{TV}(\zeta^{-1})}
		\lesssim_A1,
		\qquad
		\|G\|_{\operatorname{TV}(w_0^\flat)}
		\lesssim_\flat1.
		\label{eq:Gdot-corrected-one-half}
	\end{equation}
\end{lemma}

\begin{proof}
	Put \(r\assign\llbracket\sigma\rrbracket\).  Lemma~\ref{lem:Gdot-decay-old}
	implies, for every \(\omega\in[0,2s)\),
	\begin{align}
		\|\dot G_\sigma\|_{\operatorname{TV}(w_\sigma^\omega)}
		&\lesssim_\omega r^{2s-1},
		\label{eq:Gdot-corollary-basic-TV}
		\\
		\left\|\left|\dot G_\sigma w_\sigma^\omega\right|_M\right\|_{L^\infty}
		&\lesssim_\omega r^{-d-1},
		\label{eq:Gdot-corollary-basic-Linfty}
		\\
		\left\|\left|\dot G_\sigma\right|_M\right\|_{L^\infty}
		&\lesssim
		r^{-d-1}\wedge r^{-1}\varepsilon^{-d}.
		\label{eq:Gdot-corollary-basic-sublattice}
	\end{align}
	Indeed, for the first estimate choose
	\(\bar\alpha\in(\omega,2s)\) and \(\mathring\alpha\) sufficiently large
	in~\eqref{eq:Gdot-bound}, and integrate separately in time and space.
	The second estimate follows in the same way, since
	\(\bar\alpha>\omega\) makes the spatial translates summable.  For the
	last estimate, fix \(\bar\alpha\in(0,2s)\) and use
	\eqref{eq:Gdot-bound} separately in the regimes
	\(r\geq\varepsilon\) and \(r<\varepsilon\); in the latter regime the value
	at the origin is bounded by \(r^{-1}\varepsilon^{-d}\), and the same bound
	controls all the spatial translates.
	
	The support properties of the cutoff imply
	\begin{equation}
		\tilde J_\sigma\dot G_\sigma=\dot G_\sigma,
		\qquad
		\tilde J_{1/2}G_{1/2}=G_{1/2}.
		\label{eq:Gdot-corollary-cutoff}
	\end{equation}
	Consequently,
	\[
	\partial^A\dot G_\sigma
	=(\partial^A\tilde J_\sigma)*\dot G_\sigma,
	\qquad
	L_\sigma^p\dot G_\sigma
	=(L_\sigma^p\tilde J_\sigma)*\dot G_\sigma.
	\]
	By Lemma~\ref{lem:J-decay},
	\[
	\|\partial^A\tilde J_\sigma\|_{\operatorname{TV}(\zeta^{-1})}
	\lesssim_A r^{-|A|_s},
	\qquad
	\|L_\sigma^p\tilde J_\sigma\|_{\operatorname{TV}(w_\sigma^\flat)}
	\lesssim_{p,\flat}1.
	\]
	Here the first estimate follows from
	\eqref{eq:J-decay-corrected-TV} with \(\omega=1\): there is no sublattice
	loss when \(|\bar A|=0\), nor when \(|\bar A|\geq1\), since then
	\(1\leq|\bar A|\).  Weighted Young inequalities, together with
	\eqref{eq:Gdot-corollary-basic-TV}, prove
	\eqref{eq:Gdot-corrected-fixed-weight} and
	\eqref{eq:Gdot-corrected-weighted-TV}.  The corresponding periodised
	Young inequality, combined with
	\eqref{eq:Gdot-corollary-basic-Linfty} and
	\eqref{eq:Gdot-corollary-basic-sublattice}, proves
	\eqref{eq:Gdot-corrected-periodised-Linfty} and
	\eqref{eq:Gdot-corrected-periodised-unweighted}.
	
	The first estimate in~\eqref{eq:Gdot-corrected-one-half} follows in the
	same way from \eqref{eq:G-one-half}, the second identity in
	\eqref{eq:Gdot-corollary-cutoff}, and Lemma~\ref{lem:J-decay}.  Moreover,
	\eqref{eq:G-one-half} gives
	\(\|G_{1/2}\|_{\operatorname{TV}(w_0^\flat)}\lesssim_\flat1\).
	Finally, using
	\[
	G=G_{1/2}+\int_{1/2}^1\dot G_\sigma\,\mathd\sigma,
	\qquad
	w_0^\flat\leq w_\sigma^\flat,
	\]
	and integrating~\eqref{eq:Gdot-corollary-basic-TV} in \(\sigma\) proves
	the second estimate in~\eqref{eq:Gdot-corrected-one-half}.
\end{proof}

\begin{lemma}
  \label{lem:Gdot-decay-old}
  Fix \(d\ge1\) and \(s\in(3/4,1]\). For every $\mathring\alpha\geq 0$ and $\bar\alpha \in [0, 2 s)$, we have
  \begin{equation}
  	\left|G_{1/2}(t,x)\right|
  	\lesssim_{\mathring\alpha,\bar\alpha}
  	(1+|t|)^{-\mathring\alpha}
  	(1+|x|)^{-d-\bar\alpha}
  	\label{eq:G-one-half}
  \end{equation}
  and
  \begin{equation}
  	\left|\dot G_\sigma(t,x)\right|
  	\lesssim_{\mathring\alpha,\bar\alpha}
  	\llbracket\sigma\rrbracket^{-d-1}
  	\left(
  	1+\frac{|t|}{\llbracket\sigma\rrbracket^{2s}}
  	\right)^{-\mathring\alpha}
  	\left(
  	1+\frac{\varepsilon\vee|x|}
  	{\llbracket\sigma\rrbracket}
  	\right)^{-d}
  	\left(
  	1+\frac{|x|}{\llbracket\sigma\rrbracket}
  	\right)^{-\bar\alpha}
  	\label{eq:Gdot-bound}
  \end{equation}
  uniformly in $(t, x) \in \Lambda$ and $\sigma \in [1 / 2, 1) .$
\end{lemma}

\begin{proof}
  We only prove~\eqref{eq:Gdot-bound} since the proof of~\eqref{eq:G-one-half}
  follows the same lines. Observe that
  \[ \dot{G}_{\sigma} (t, x) = \int_{\Lambda^{\ast}_{\varepsilon}}
     \frac{\partial_{\sigma} (j_{\sigma} (| k'_0 |^{1 / 2 s}) j_{\sigma}
     (q_{\varepsilon} (\tilde{k}')))}{i k_0' +
     m^2 + q_{\varepsilon}^{2 s} (\tilde{k}')} e^{i (k_0' t + \tilde{k}' x)}
     \frac{\mathd k_0' \mathd \tilde{k}'}{(2 \pi)^{d + 1}} . \]
  To prove a bound for the $L^{\infty}$-norm of
  \[ (t, x) \mapsto t^{a_0} x^{\tilde{a}} \dot{G}_{\sigma} (t, x),
  \]
  it suffices to control the $L^1$-norm of
  \begin{equation}
    (k'_0, \tilde{k}') \mapsto \partial_{k_0'}^{a_0} \partial_{\tilde{k}'}^{\tilde{a}}
    \left( \frac{\partial_{\sigma} (j_{\sigma} (| k'_0 |^{1 / 2 s}) j_{\sigma}
    (q_{\varepsilon} (\tilde{k}')))}{i k_0' +
    m^2 + q_{\varepsilon}^{2 s} (\tilde{k}')} \right) .
    \label{eq:G-dot-proof-Fourier-1}
  \end{equation}
  Since for noninteger parameters $s$ the bound
  \[ | \partial_{\tilde{k}'}^{\tilde{a}} q_{\varepsilon}^{2 s} (\tilde{k}') |
     \lesssim | \tilde{k}' |^{2 s - | \tilde{a} |}, \]
  uniform in $\varepsilon \in (0, 1)$, is optimal, the $L^1$-norm of the
  function~\eqref{eq:G-dot-proof-Fourier-1} is bounded uniformly in
  $\varepsilon \in (0, 1)$ only if $| \tilde{a} | < d + 2 s$. As a result, the
  above simple strategy can only be used to prove the lemma for $\bar\alpha \in
  [0, \lfloor 2 s \rfloor]$. To establish the claim for every $\bar\alpha \in [0,
  2 s)$ a more refined argument is needed. To this end, let $h : \mathbb{R}
  \rightarrow \mathbb{R}$ be defined by
  \[ h (\omega) \assign - \omega \partial_{\omega} j (\omega) \]
  for $\omega \in \mathbb{R}$. Moreover, let the families of kernels
  $(\tilde{G}_{\tau})_{\tau \in (0, \infty)}$, $(\check{G}_{\tau})_{\tau \in
  (0, \infty)}$, $(\hat{G}_{\tau, \eta})_{\tau, \eta \in (0, \infty)}$ be
  defined by
  \begin{align*}
       \tilde{G}_{\tau} (t, x)
       &\assign \tau^{- 1}
       \int_{\Lambda^{\ast}_{\varepsilon}}
       \frac{j \left( \tau | k'_0 |^{1 / 2 s}  \right)
       h (\tau q_{\varepsilon} (\tilde{k}') )
       + h \left( \tau | k'_0 |^{1 / 2 s}  \right)
       j (\tau q_{\varepsilon} (\tilde{k}') )}
       {i k_0' + m^2 + q^{2 s}_{\varepsilon}(\tilde{k}')}\\
       &\qquad \times
       e^{i (k'_0 t + \tilde{k}' x)}
       \frac{\mathd k_0' \mathd \tilde{k}'}{(2 \pi)^{d + 1}},
  \end{align*}
  \begin{align*}
       \hat{G}_{\tau, \eta} (t, x)
       &\assign \tau^{- 1} \eta^{- 1}
       \int_{\Lambda^{\ast}_{\varepsilon}}
       \frac{q^{2 s}_{\varepsilon}(\tilde{k}')
       h \left( \tau | k'_0 |^{1 / 2 s}  \right)
       j(\tau q_{\varepsilon}(\tilde{k}') )
       h(\tau \eta q_{\varepsilon}(\tilde{k}'))}
       {(k_0' - i m^2)^2 + q^{4 s}_{\varepsilon}(\tilde{k}')}\\
       &\qquad \times
       e^{i (k'_0 t + \tilde{k}' x)}
       \frac{\mathd k_0' \mathd \tilde{k}'}{(2 \pi)^{d + 1}},
  \end{align*}
  \begin{align*}
       \check{G}_{\tau} (t, x)
       &\assign \tau^{- 1}
       \int_{\Lambda^{\ast}_{\varepsilon}}
       \left( \frac{j \left( \tau | k'_0 |^{1 / 2 s}  \right)
       h(\tau q_{\varepsilon}(\tilde{k}'))}
       {i k_0' + m^2 + q^{2 s}_{\varepsilon}(\tilde{k}')} \right.\\
       &\qquad \left.
       - \frac{(i k_0' + m^2)
       h \left( \tau | k'_0 |^{1 / 2 s}  \right)
       j(\tau q_{\varepsilon}(\tilde{k}') )}
       {(k_0' - i m^2)^2 + q^{4 s}_{\varepsilon}(\tilde{k}')}
       \right)\,
       e^{i (k'_0 t + \tilde{k}' x)}
       \frac{\mathd k_0' \mathd \tilde{k}'}{(2 \pi)^{d + 1}} .
  \end{align*}
  We claim that
  \begin{equation}
    \dot{G}_{\sigma} (t, x) = \sigma^{- 2} {\tilde{G}_{\llbracket \sigma
    \rrbracket  / \sigma }}  (t, x), \qquad \tilde{G}_{\tau} =
    \check{G}_{\tau} + \int_{1 / 2}^{\infty} \hat{G}_{\tau, \eta} \mathd \eta
    . \label{eq:G-dot-proof-identities}
  \end{equation}
  The first identity follows from
  \[ \partial_{\sigma} j_{\sigma} (\omega) = \llbracket \sigma \rrbracket^{-
     1} \sigma^{- 1} h (\llbracket \sigma \rrbracket  \sigma^{- 1} \omega) .
  \]
  To verify the second, we use the identities
  \begin{align*}
     &\frac{h \left( \tau  | k'_0 |^{1 / 2 s}  \right)
     j (\tau q_{\varepsilon} (\tilde{k}') )}
     {i k_0' + m^2 + q^{2 s}_{\varepsilon}(\tilde{k}')}\\
     &\quad =
     \frac{q^{2 s}_{\varepsilon} (\tilde{k}') h \left( \tau
     | k'_0 |^{1 / 2 s}  \right) j (\tau q_{\varepsilon}(\tilde{k}') ) }
     {(k_0' - i m^2)^2 + q^{4 s}_{\varepsilon}(\tilde{k}')}
     - \frac{(i k_0' + m^2) h \left( \tau | k'_0
     |^{1 / 2 s}  \right) j (\tau q_{\varepsilon}(\tilde{k}') )}
     {(k_0' - i m^2)^2 + q^{4 s}_{\varepsilon}(\tilde{k}')}
  \end{align*}
  as well as
  \[ \int^{\infty}_{1 / 2} \eta^{- 1} h (\tau \eta q _{\varepsilon}
     (\tilde{k}'))\,\mathd \eta = - \int^{\infty}_{1 / 2} \partial_{\eta} j
     (\tau \eta q _{\varepsilon} (\tilde{k}'))\,\mathd \eta = j (\tau q
     _{\varepsilon} (\tilde{k}') / 2) \]
  and
  \[ j (\tau q _{\varepsilon} (\tilde{k}')) j (\tau q _{\varepsilon}
     (\tilde{k}') / 2) = j (\tau q _{\varepsilon} (\tilde{k}')) . \]
  Let us motivate the usefulness of the representation of $\dot{G}_{\sigma}$
  given by~\eqref{eq:G-dot-proof-identities}. First, note that if $| \tilde{a}
  | < d + 4 s$, then the $L^1$-norm of the function
  \begin{align*}
    (k'_0, \tilde{k}') \mapsto
    \partial_{k_0'}^{a_0}
    \partial_{\tilde{k}'}^{\tilde{a}}
    \left( \left(
    \frac{j \left( \tau | k'_0 |^{1 / 2 s} \right)
    h(\tau q _{\varepsilon}(\tilde{k}'))}
    {i k_0' + m^2 + q^{2 s}_{\varepsilon}(\tilde{k}')}\right.\right.\\
    \left.\left.
    - \frac{(i k_0' + m^2)
    h \left( \tau | k'_0 |^{1 / 2 s} \right)
    j(\tau q _{\varepsilon}(\tilde{k}'))}
    {(k_0' - i m^2)^2 + q^{4 s}_{\varepsilon}(\tilde{k}')}
    \right)  \right),
  \end{align*}
  is bounded uniformly in $\varepsilon \in (0, 1)$. This stands in contrast to
  the function in~\eqref{eq:G-dot-proof-Fourier-1}, whose $L^1$-norm remains
  uniformly bounded in $\varepsilon \in (0, 1)$ only if $| \tilde{a} | < d + 2
  s$. As a result, it is possible to control the $L^{\infty}$-norm of the
  function
  \[ (t, x) \mapsto t^{a_0} x^{\tilde{a}} \check{G}_{\tau} (t, x)
  \]
  for $| \tilde{a} | < d + 4 s$. Since $\lfloor 4 s \rfloor \geqslant 3 > 2 s
  > \bar\alpha$ this indicates that $\check{G}_{\tau}$ indeed exhibits
  the desired decay at infinity. On the other hand, the $L^1$-norm of the
  function
  \[ (k'_0, \tilde{k}') \mapsto \partial_{k_0'}^{a_0}
     \partial_{\tilde{k}'}^{\tilde{a}} \left( \frac{q^{2 s}_{\varepsilon}
     (\tilde{k}') h \left( \tau  | k'_0 |^{1 / 2 s}  \right) j
     (\tau^{} q _{\varepsilon} (\tilde{k}') ) h (\tau  \eta  q _{\varepsilon}
     (\tilde{k}'))}{(k_0' - i m^2)^2 + q^{4 s}_{\varepsilon} (\tilde{k}')}
      \right) \]
  is bounded uniformly in $\varepsilon \in (0, 1)$ for all $a_0, \tilde{a}$
  and consequently $\hat{G}_{\tau, \eta}$ has good decay
  properties. However, the bound for the $L^{\infty}$-norm of the function
  \[ (t, x) \mapsto t^{a_0} x^{\tilde{a}} \hat{G}_{\tau, \eta} (t,
     x) \]
  depends on $\eta$. To complete the argument, we must therefore control the
  above norms uniformly in $\varepsilon \in (0, 1)$, $\tau \in (0, 1)$ and
  $\eta \in (1 / 2, \infty)$.
  
  Fix an integer $N\geqslant\mathring\alpha$. We claim that the following
  bounds imply~\eqref{eq:Gdot-bound},
  \begin{equation}
    \begin{aligned}
      | \tilde{G}_{\tau} (t, 0) | & \lesssim \tau^{- 1}
      \varepsilon^{- d} (1 + | t |  / \tau^{2 s})^{- N},\\
      | \check{G}_{\tau} (t, x) | & \lesssim \tau^{- d - 1}
      (1 + | x |  / \tau)^{- d - 3} (1 + | t |  / \tau^{2 s})^{- N},\\
      | \hat{G}_{\tau, \eta} (t, x) | & \lesssim \eta^{- 2 s - 1}
      \tau^{- d - 1}  (1 + | x |  / \tau)^{- d}  (1 + | t |  / \tau^{2
      s})^{- N},\\
      | \hat{G}_{\tau, \eta} (t, x) | & \lesssim \eta^{1 - 2 s}
      \tau^{- d - 1}  (1 + | x |  / \tau)^{- d - 2}  (1 + | t |  /
      \tau^{2 s})^{- N} .
    \end{aligned} \label{eq:G-dot-proof-bounds}
  \end{equation}
  To prove this claim, first use the last two bounds to conclude that for all
  $\bar\alpha \in [0, 2]$ we have
  \[ | \hat{G}_{\tau, \eta} (t, x) | \lesssim \eta^{\bar\alpha - 1 - 2
     s} \tau^{- d - 1}  (1 + | x |  / \tau)^{- d - \bar\alpha}  {(1 + | t
     |  / \tau^{2 s})^{- N}}  . \]
  Consequently, the second identity in~\eqref{eq:G-dot-proof-identities} and
  the above bound yield
  \[ | \tilde{G}_{\tau} (t, x) | \lesssim_{\bar\alpha} \tau^{- d - 1}
     (1 + | x | / \tau)^{- d - \bar\alpha}
     (1 + | t | / \tau^{2 s})^{- N} \]
  for $\bar\alpha \in [0, 2 s)$. If $x\neq0$, then $|x|\geqslant\varepsilon$,
  while for $x=0$ we combine the above estimate with the first estimate in
  \eqref{eq:G-dot-proof-bounds}. We obtain
  \[ | \tilde{G}_{\tau} (t, x) | \lesssim_{\bar\alpha} \tau^{- d - 1}
     (1 + (\varepsilon \vee |x|) / \tau)^{- d}
     (1 + |x| / \tau)^{- \bar\alpha}
     (1 + |t| / \tau^{2 s})^{- N} . \]
  Since $N\geqslant\mathring\alpha$ and
  $\llbracket\sigma\rrbracket/\sigma\asymp\llbracket\sigma\rrbracket$ for
  $\sigma\in[1/2,1)$, the bound~\eqref{eq:Gdot-bound} follows from the first
  identity in~\eqref{eq:G-dot-proof-identities}.
  
  It remains to prove~\eqref{eq:G-dot-proof-bounds}. Starting from the bound
  for $\hat{G}_{\tau, \eta} (t, x)$, we have
  \begin{align*}
       (\hat{G}_{\tau, \eta}) (t, x)
       &= \tau^{- 1} \eta^{- 1} \int_{\Lambda^{\ast}_{\varepsilon}}
       \frac{q^{2 s}_{\varepsilon} (\tilde{k}')}
       {(k_0' - i m^2)^2 + q^{4 s}_{\varepsilon}(\tilde{k}')}
       h \left( \tau | k'_0 |^{1 / 2 s}  \right)
       j (\tau q_{\varepsilon} (\tilde{k}') )\\
       &\qquad \times
       h (\tau \eta q_{\varepsilon}(\tilde{k}'))
       e^{i (k'_0 t + \tilde{k}' x)}
       \frac{\mathd k_0' \mathd \tilde{k}'}{(2 \pi)^{d + 1}}\\
       &= \tau^{- d - 1} \eta^{- 1}
       \int_{\Lambda^{\ast}_{\varepsilon / \tau}}
       \frac{q^{2 s}_{\varepsilon / \tau} (\tilde{k})}
       {(k_0 - i \tau^{2 s} m^2)^2
       + q^{4 s}_{\varepsilon / \tau}(\tilde{k})}
       \\
       &\qquad \times h \big( | k_0 |^{1 / 2 s}  \big)
       j (q_{\varepsilon / \tau} (\tilde{k}))\\
       &\qquad \times
       h (\eta q_{\varepsilon / \tau}(\tilde{k}))
       e^{i(k_0 t / \tau^{2 s} + \tilde{k} x / \tau)}
       \frac{\mathd k_0 \mathd \tilde{k}}{(2 \pi)^{d + 1}},
  \end{align*}
  where we set $k_0 = \tau^{2 s} {k'_0} $ and $\tilde{k} = \tau \tilde{k}'$.
  It follows that
  \begin{equation}
    \begin{aligned}
      &(\tau^{- 2 s} | t |)^{a_0} (\tau^{- 1} | x |)^{\tilde{a}}
      | (\hat{G}_{\tau, \eta}) (t, x) |\\
      &\quad \lesssim \tau^{- d - 1} \eta^{- 1}
      \int_{\Lambda^{\ast}_{\varepsilon / \tau}}
      \Bigg| \partial_{k_0}^{a_0}
      \partial_{\tilde{k}}^{\tilde{a}} \Bigg[
      \frac{q^{2 s}_{\varepsilon / \tau}(\tilde{k})
      h \left( | k_0 |^{1 / 2 s}  \right)
      j(q_{\varepsilon / \tau}(\tilde{k}))
      h(\eta q_{\varepsilon / \tau}(\tilde{k}))}
      {(k_0 - i \tau^{2 s} m^2)^2
      + q^{4 s}_{\varepsilon / \tau}(\tilde{k})}
       \Bigg] \Bigg|
      \frac{\mathd k_0 \mathd \tilde{k}}{(2 \pi)^{d + 1}} .
    \end{aligned} \label{eq:integral-to-control}
  \end{equation}
  Observe the following:
  \begin{itemize}
    \item In~\eqref{eq:integral-to-control}, the factor $h \left( | k_0
    |^{1 / 2 s}  \right)$ restricts the integration domain to $| k_0 |^{1 / 2 s}  \in [1, 2]$. In this region, the denominator is never
    vanishing. More precisely
    \[ | (k_0 - i \tau^{2 s} m^2)^2 + q^{4 s}_{\varepsilon / \tau}
       (\tilde{k}) |^2 \geqslant c_m>0. \]
    \item Similarly, the factor $h (\eta q _{\varepsilon / \tau} (\tilde{k}))$
    restricts the integration domain to those $\tilde{k}$ such that $q
    _{\varepsilon / \tau} (\tilde{k}) \in [\eta^{- 1}, 2 \eta^{- 1}]$. Thus,
    by~\eqref{eq:q-epsilon-derivatives}, on the integration domain
    \[ | \partial_{\tilde{k}}^{\tilde{a}} q^{2 s}_{\varepsilon / \tau}
       (\tilde{k}) | \lesssim q^{2 s - | \tilde{a} |}_{\varepsilon / \tau}
       (\tilde{k}) \approx \eta^{| \tilde{a} | - 2 s} . \]
  \end{itemize}
  On account of these comments and~\eqref{eq:q-epsilon-cutoff-derivatives},
  we have the following bound for the integrand
  \begin{align*}
    &\left| \partial_{k_0}^{a_0} \partial_{\tilde{k}}^{\tilde{a}} \left(
    \frac{q^{2 s}_{\varepsilon / \tau} (\tilde{k})
    h \left( | k_0 |^{1 / 2 s} \right)
    j(q _{\varepsilon / \tau} (\tilde{k}))
    h(\eta q _{\varepsilon / \tau} (\tilde{k}))
    }
    {(k_0 - i \tau^{2 s} m^2)^2
    + q^{4 s}_{\varepsilon / \tau}(\tilde{k})}
    \right) \right|\\
    &\quad \lesssim \eta^{| \tilde{a} | - 2 s}
    1_{[0, 2]} (| k_0 |^{1 / 2 s})
    1_{[1, 2]} (\eta q _{\varepsilon / \tau }(\tilde{k})),
  \end{align*}
  for all $(a_0, \tilde{a}) \in \mathbb{N}_0^{1 + d}$. Taking $a_0=p$
  and summing over all $\tilde a$ with $|\tilde a|=q$, this implies that
  \begin{align*}
    &(\tau^{- 2 s} | t |)^p (\tau^{- 1} | x |)^q
    | (\hat{G}_{\tau, \eta}) (t, x) |\\
    &\quad \lesssim \tau^{- d - 1} \eta^{q - 2 s - 1}
    \int_{\Lambda^{\ast}_{\varepsilon / \tau}}
    1_{[0, 2]} (| k_0 |^{1 / 2 s})
    1_{[1, 2]} (\eta q _{\varepsilon / \tau } (\tilde{k}))
    \frac{\mathd k_0 \mathd \tilde{k}}{(2 \pi)^{d + 1}}\\
    &\quad \lesssim \tau^{- d - 1} \eta^{q - 2 s - d - 1}
    \int_{\Lambda^{\ast}_{\varepsilon / (\tau \eta)}}
    1_{[0, 2]} (| k_0 |^{1 / 2 s})
    1_{[1, 2]} (q _{\varepsilon / (\tau \eta)} (\tilde{k}))
    \frac{\mathd k_0 \mathd \tilde{k}}{(2 \pi)^{d + 1}}\\
    &\quad \lesssim \tau^{- d - 1} \eta^{q - 2 s - d - 1},
  \end{align*}
  where the last estimate follows from~\eqref{eq:q-epsilon-volume-singular}.
  For any $p, q \in \mathbb{N}_0$, taking $p\in\{0,N\}$ and
  $q\in\{0,d,d+2\}$ proves the last two bounds in
  \eqref{eq:G-dot-proof-bounds}. The proof for $\check{G}_{\tau}
  (t, x)$ is similar. Working as in the previous case, we have
  \begin{align*}
       &(\tau^{- 2 s} | t |)^p (\tau^{- 1} | x |)^q
       | \check{G}_{\tau} (t, x) |\\
       &\quad \lesssim \tau^{- d - 1}
       \int_{\Lambda^{\ast}_{\varepsilon / \tau}} \biggl|
       \partial_{k_0}^{a_0} \partial_{\tilde{k}}^{\tilde{a}}
       \biggl[ \biggl(
       \frac{j \left( | k_0 |^{1 / 2 s}  \right)
       h (q_{\varepsilon / \tau} (\tilde{k}))}
       {i k_0 + \tau^{2 s} m^2 + q^{2 s}_{\varepsilon / \tau} (\tilde{k})}
       \\
       &\qquad
       - \frac{(i k_0 + \tau^{2 s} m^2)
       h \left( | k_0 |^{1 / 2 s}  \right)
       j (q_{\varepsilon / \tau} (\tilde{k}))}
       {(k_0 - i \tau^{2 s} m^2)^2
       + q^{4 s}_{\varepsilon / \tau} (\tilde{k})} \biggr)
        \biggr] \biggr|\,
       \frac{\mathd k_0 \mathd \tilde{k}}{(2 \pi)^{d + 1}} .
  \end{align*}
  As above, also here the denominators are non-vanishing on the integration
  domain, more precisely,
  \[ | i k_0 + \tau^{2 s} m^2 + q^{2 s}_{\varepsilon / \tau} (\tilde{k}) |
     \geqslant 1 \qquad \operatorname{and} \qquad | (k_0 - i \tau^{2 s} m^2)^2 +
     q^{4 s}_{\varepsilon / \tau} (\tilde{k}) | \geqslant 1. \]
  In the first term, this is due to the factor $h (q_{\varepsilon /
  \tau} (\tilde{k}))$, which restricts the integration domain to $q_{
  \varepsilon / \tau} (\tilde{k}) \in [1, 2]$, while in the second one
  this is due to the factor $h \left( | k_0 |^{1 / 2 s} 
  \right)$, which restricts the integration domain to $| k_0 |^{1 / 2 s} 
  \in [1, 2]$. By~\eqref{eq:q-epsilon-derivatives}, since
  $| \partial_{\tilde{k}}^{\tilde{a}} q^{4 s}_{\varepsilon
  / \tau} (\tilde{k}) | \lesssim q^{4 s - | \tilde{a} |}_{\varepsilon / \tau}
  (\tilde{k})$, we get
  \begin{align*}
    &\biggl| \partial_{k_0}^{a_0} \partial_{\tilde{k}}^{\tilde{a}} \biggl[
    \biggl( \frac{j \left( | k_0 |^{1 / 2 s} \right)
    h(q _{\varepsilon / \tau} (\tilde{k}))}
    {i k_0 + \tau^{2 s} m^2
    + q^{2 s}_{\varepsilon / \tau}(\tilde{k})}\\
    &\qquad
    - \frac{(i k_0 + \tau^{2 s} m^2)
    h \left( | k_0 |^{1 / 2 s} \right)
    j(q _{\varepsilon / \tau}(\tilde{k}))}
    {(k_0 - i \tau^{2 s} m^2)^2
    + q^{4 s}_{\varepsilon / \tau}(\tilde{k})} \biggr)
     \biggr] \biggr|\\
    &\quad \lesssim
    1_{[0, 2]} (| k_0 |^{1 / 2 s})
    1_{[0, 2]}(q _{\varepsilon / \tau }(\tilde{k}))\,
    (1 \vee q^{4 s - | \tilde{a} |}_{\varepsilon / \tau}(\tilde{k})) .
  \end{align*}
  Taking $a_0=p$ and summing over all $\tilde a$ with $|\tilde a|=q$,
  this implies that, for $q < d + 4 s$,
  \begin{align*}
    &(\tau^{- 2 s} | t |)^p (\tau^{- 1} | x |)^q
    | \check{G}_{\tau} (t, x) |\\
    &\quad \lesssim \tau^{- d - 1}
    \int_{\Lambda^{\ast}_{\varepsilon / \tau}}
    1_{[0, 2]} (| k_0 |^{1 / 2 s})
    1_{[0, 2]} (q _{\varepsilon / \tau } (\tilde{k}))
    (1 \vee q^{4 s - q}_{\varepsilon / \tau}(\tilde{k}))
    \mathd k_0 \mathd \tilde{k}\\
    &\quad \lesssim \tau^{- d - 1} .
  \end{align*}
  Here the last estimate is~\eqref{eq:q-epsilon-check-integrability}.
  Since $d + 4 s > d + 3$, taking $p\in\{0,N\}$ and
  $q\in\{0,d+3\}$ proves the bound for $\check{G}_{\tau}$.
  Finally, we discuss the first bound in \eqref{eq:G-dot-proof-bounds}. By
  changing the variables as above, we arrive at
  \begin{align*}
    &(\tau^{- 2 s} | t |)^p
    | \tilde{G}_{\tau} (t, 0) |
    \leqslant \tau^{- 1 - d}
    \\
    &\quad\times 
    \int_{\Lambda^{\ast}_{\varepsilon / \tau}}
    \left|
    \partial_{k_0}^p \left[
    \frac{j \left( | k_0 |^{1 / 2 s} \right)
    h(q_{\varepsilon / \tau}(\tilde{k}))
    + h \left( | k_0 |^{1 / 2 s} \right)
    j(q_{\varepsilon / \tau}(\tilde{k}))}
    {i k_0 + \tau^{2 s} m^2
    + q^{2 s}_{\varepsilon / \tau}(\tilde{k})}
     \right] \right|
    \mathd k_0 \mathd \tilde{k} .
  \end{align*}
  We use the same argument as for $\check{G}_{\tau}$: indeed both
  terms in the integrand have non-vanishing denominators due to the factors $h
  (q _{\varepsilon / \tau} (\tilde{k}) )$ and $h \left( | k_0 |^{1
  / 2 s}  \right)$, respectively. Analogously, we obtain
  \begin{align*}
    &\left| \partial_{k_0}^p \left[
    \frac{j \left( | k_0 |^{1 / 2 s} \right)
    h(q _{\varepsilon / \tau}(\tilde{k}))
    + h \left( | k_0 |^{1 / 2 s} \right)
    j(q _{\varepsilon / \tau}(\tilde{k}))}
    {i k_0 + \tau^{2 s} m^2
    + q^{2 s}_{\varepsilon / \tau}(\tilde{k})}
     \right] \right|\\
    &\quad \lesssim 1_{[0, 2]} (| k_0 |^{1 / 2 s}),
  \end{align*}
  and thus
  \begin{align*}
    &(\tau^{- 2 s} | t |)^p
    | \tilde{G}_{\tau} (t, 0) |\\
    &\quad \lesssim \tau^{- 1 - d}
    \int_{\Lambda^{\ast}_{\varepsilon / \tau}}
    1_{[0, 2]} (| k_0 |^{1 / 2 s})
    \frac{\mathd k_0 \mathd \tilde{k}}{(2 \pi)^{d + 1}}\\
    &\quad \lesssim
    \tau^{- 1 - d}(\varepsilon \tau^{- 1})^{- d}
    = \tau^{- 1} \varepsilon^{- d} .
  \end{align*}
  Taking $p\in\{0,N\}$ proves the first bound in
  \eqref{eq:G-dot-proof-bounds}. This concludes the proof.
\end{proof}

\subsection{Properties of the weights}\label{sec:weight-properties}

We collect here some important properties of the various weights introduced
throughout the paper.

\begin{lemma}
  \label{lem:fractional-bound}For $\alpha \geqslant 0$, we have
  \[ \| (- \Delta_{\varepsilon})^s \zeta_{\sigma}^{\alpha} \| \lesssim
     \llbracket \sigma \rrbracket^{2 s a} \]
  uniformly in $\sigma \in (0, 1)$.
\end{lemma}

\begin{proof}
  Let $f \in C^2 (\mathbb{R}^d)$. Recall that
  by~\eqref{eq:integral-rep-frac-lap-discrete}, we have
  \begin{align*}
    (- \Delta_{\varepsilon})^s f (x)
    &= \int_{\mathbb{R}_{\varepsilon}^d}
    H^{(\varepsilon)}_s (y) (f (x) - f (x - y))\,\mathd y,
    \qquad x \in \mathbb{R}_{\varepsilon}^d,
  \end{align*}
  where the kernel $H^{(\varepsilon)}_s : \mathbb{R}_{\varepsilon}^d
  \rightarrow \mathbb{R}$ is positive and such that $H^{(\varepsilon)}_s (0) =
  0$, $H^{(\varepsilon)}_s (x) = H^{(\varepsilon)}_s (- x)$ and
  \[ | H^{(\varepsilon)}_s (x) | \lesssim | x |^{- d - 2 s}, \qquad x \in
     \mathbb{R}_{\varepsilon}^d . \]
  Let $B (\delta)$ be the ball in $\mathbb{R}_{\varepsilon}^d$ of radius $\delta
  > 0$ centred at the origin and $B^c (\delta) \assign
  \mathbb{R}_{\varepsilon}^d \setminus B (\delta)$. For any $\delta > 0$, we
  have
\[
\begin{aligned}
	|(-\Delta_\varepsilon)^sf(x)|
	&\le \frac12\int_{B(\delta)}H_s^{(\varepsilon)}(y)\,
	|2f(x)-f(x+y)-f(x-y)|\,\mathrm dy\\
	&\quad+\int_{B^c(\delta)}H_s^{(\varepsilon)}(y)\,
	|f(x)-f(x-y)|\,\mathrm dy,
\end{aligned}
\]
  To bound the first term observe that
  \[ | 2 f (x) - f (x + y) - f (x - y) | \lesssim \| \nabla \nabla f \|  | y
     |^2, \]
  by Taylor's theorem, where $\| \nabla \nabla f \|$ denotes the supremum norm
  of the Hessian of $f$ in the continuum. As a result, we obtain

  \begin{align*}
    &\int_{B(\delta)} H^{(\varepsilon)}_s (y)
    \,|2 f (x) - f (x + y) - f (x - y)|\,\mathd y \\
    &\quad \lesssim \| \nabla \nabla f \|
    \int_{B(\delta)} | y |^{2 - d - 2 s} \mathd y
    \lesssim \| \nabla \nabla f \| \delta^{2 - 2 s} .
  \end{align*}
  Moreover, we have
  \[  \int_{B^c (\delta)} H^{(\varepsilon)}_s (y)  \,|f (x) - f (x - y)|\,
     \mathd y  \lesssim \| f \|  \int_{B^c (\delta)} | y |^{- d - 2 s}
     \mathd y \lesssim \| f \| \delta^{- 2 s} . \]
  To conclude, we apply the above estimates with $f (x) = \zeta^{\alpha}_{\sigma}
  (t, x)$ and choose $\delta = \llbracket \sigma \rrbracket^{- a}$.
\end{proof}

\begin{lemma}
  \label{lem:D-zeta}We have
  \[ \| \zeta_{\sigma}^{1 / 3} \mathfrak{D}_s
     (\zeta_{\sigma}^{- 1 / 3}) \|^2 \lesssim \llbracket \sigma \rrbracket^{2
     s a}, \]
  uniformly in $\sigma \in (0, 1)$.
\end{lemma}

\begin{proof}
  Recall that $\mu_s$ denotes the kernel defined by~\eqref{eq:int-kernel},
  which is local in time. Let $B (z, \delta)$ be the ball of radius $\delta > 0$
  centred at $z \in \Lambda$. By~\eqref{eq:D}, we have
  \[ \begin{aligned}
       {}[\zeta_{\sigma}^{1 / 3} \mathfrak{D}_s
       (\zeta_{\sigma}^{- 1 / 3}) (z)]^2 & =  \zeta_{\sigma}^{2 / 3} (z) \int
       [\zeta_{\sigma}^{- 1 / 3} (z') - \zeta_{\sigma}^{- 1 / 3} (z)]^2 \mu_s
       (z, \mathd z')\\
       & =  \zeta_{\sigma}^{2 / 3} (z) \int_{B^c (z,
       \delta)} [\zeta_{\sigma}^{- 1 / 3} (z') - \zeta_{\sigma}^{- 1 / 3}
       (z)]^2 \mu_s (z, \mathd z')\\
       &   + \zeta_{\sigma}^{2 / 3} (z) \int_{B (z,
       \delta)} [\zeta_{\sigma}^{- 1 / 3} (z') - \zeta_{\sigma}^{- 1 / 3}
       (z)]^2 \mu_s (z, \mathd z')\\
       & =  \mathbb{I}+\mathbb{I}\mathbb{I}.
     \end{aligned} \]
  On $B^c (z, \delta) \assign \Lambda_{\varepsilon} \setminus B (z, \delta)$
  we have
  \begin{align*}
       \mathbb{I}
       &= \int_{B^c (z, \delta)}
       [\zeta_{\sigma}^{1 / 3} (z)
       \zeta_{\sigma}^{- 1 / 3} (z') - 1]^2 \mu_s (z, \mathd z')\\
       &\lesssim
       \int_{B^c (z, \delta)} [\zeta_{\sigma}^{2 / 3} (z)
       \zeta_{\sigma}^{- 2 / 3} (z') + 1] \mu_s (z, \mathd z')\\
       &\lesssim \int_{B^c (z, \delta)} \zeta_{\sigma}^{2 / 3} (z)
       \zeta_{\sigma}^{- 2 / 3} (z') \mu_s (z, \mathd z') + \delta^{- 2 s} .
  \end{align*}
  Recall that
  \[ \zeta_{\sigma}^{2 / 3} (z) \zeta_{\sigma}^{- 2 / 3} (z') \lesssim
     \zeta^{- 2 / 3}_{\sigma} (z - z') = \langle \llbracket
     \sigma \rrbracket^a . (z - z') \rangle_s^{2 / 3} = (1 + \llbracket \sigma
     \rrbracket^{2 a} | \bar{z} - \bar{z}' |^2)^{1 / 3}, \]
  uniformly in $z = (z_0, \bar{z}) \in \Lambda$ and $z' = (z_0, \bar{z}') \in
  \Lambda$. Since $2 s > 2 / 3$,
  \begin{align*}
    \int_{B^c (z, \delta)} \zeta_{\sigma}^{2 / 3} (z) \zeta_{\sigma}^{- 2 / 3}
    (z') \mu_s (z, \mathd z') &\lesssim  \int_{B^c (z, \delta)} (1 +
    \llbracket \sigma \rrbracket^{2 a} | \bar{z} - \bar{z}' |^2)^{1 / 3} \mu_s
    (z, \mathd z')\\
 &\lesssim  \delta^{- 2 s} + \llbracket \sigma \rrbracket^{2 a / 3}
    \delta^{- 2 s + 2 / 3} .
  \end{align*}
  Let us bound the second term. Thanks to the smoothness of $\zeta_{\sigma}^{-
  1 / 3}$ in space, by Taylor's theorem
  \[ \mathbb{I}\mathbb{I} \leqslant \zeta_{\sigma}^{2 / 3} (z) \int_{B (z,
     \delta)} \| \nabla  \zeta_{\sigma}^{- 1 / 3} \|^2  | \bar{z}' - \bar z |^2
     \mu_s (z, \mathd z') . \]
  Since $\| \nabla  \zeta_{\sigma}^{- 1 / 3} \|^2 \lesssim \llbracket \sigma
  \rrbracket^{2 a}$ and $\| \zeta_{\sigma}^{2 / 3} \| \lesssim 1$, we conclude
  \[ \mathbb{I}\mathbb{I} \lesssim \llbracket \sigma \rrbracket^{2 a} \int_{B
     (z, \delta)} | \bar{z}' - \bar{z} |^2 \mu_s (z, \mathd z') \lesssim
     \delta^{2 - 2 s} \llbracket \sigma \rrbracket^{2 a} = \llbracket \sigma
     \rrbracket^{2 sa} . \]
  Choosing now $\delta = \llbracket \sigma \rrbracket^{- a}$, we get
  \[ \mathbb{I} \lesssim \llbracket \sigma \rrbracket^{2 as} + \llbracket
     \sigma \rrbracket^{2 a / 3} \llbracket \sigma \rrbracket^{2 s a - 2 a /
     3} \lesssim \llbracket \sigma \rrbracket^{2 as}, \qquad
     \mathbb{I}\mathbb{I} \lesssim \delta^{2 - 2 s} \llbracket \sigma
     \rrbracket^{2 a} = \llbracket \sigma \rrbracket^{2 sa} . \]
  This finishes the proof.
\end{proof}

\begin{lemma}
  \label{lem:weight-interaction} Recall that $1<\flat<2s<2$. The bound
  \[ \frac{(1 - h_{\mu}) w^{(2), \flat}_{\mu}}{w^{(2), \flat}_{\sigma}}
     \lesssim \llbracket \sigma \rrbracket^{\flat} \llbracket \mu
     \rrbracket^{- \flat} \]
  holds uniformly in $\mu, \sigma \in (0, 1)$.
\end{lemma}

\begin{proof}
  We have
  \[ \begin{aligned}
       &\frac{(1 - h_{\mu} (z, z_1)) w^{(2), \flat}_{\mu} (z, z_1)}{w^{(2),
       \flat}_{\sigma} (z, z_1)}  =  \frac{w^{(2), \flat}_{\mu} (z,
       z_1)}{w^{(2), \flat}_{\sigma} (z, z_1)} \frac{\llbracket \mu
       \rrbracket^{- 2} | z - z_1 |_s^2}{1 + \llbracket \mu \rrbracket^{- 2} |
       z - z_1 |_s^2}\\
       &\quad =  \frac{w^{(2), \flat}_{\mu} (z, z_1)}{w^{(2), \flat}_{\sigma} (z,
       z_1)} \frac{(\llbracket \mu \rrbracket^{- 1} | z - z_1 |_s)^{\flat}
       (\llbracket \mu \rrbracket^{- 1} | z - z_1 |_s)^{2 - \flat}}{1 +
       \llbracket \mu \rrbracket^{- 2} | z - z_1 |_s^2}\\
       &\quad =  \llbracket \sigma \rrbracket^{\flat} \llbracket \mu
       \rrbracket^{- \flat} \frac{(1 + \llbracket \mu \rrbracket^{- 1} | z -
       z_1 |_s)^{\flat}}{(1 + \llbracket \sigma \rrbracket^{- 1} | z - z_1
       |_s)^{\flat}} \frac{(\llbracket \sigma \rrbracket^{- 1} | z - z_1
       |_s)^{\flat} (\llbracket \mu \rrbracket^{- 1} | z - z_1 |_s)^{2 -
       \flat}}{1 + \llbracket \mu \rrbracket^{- 2} | z - z_1 |_s^2}\\
       &\quad \leqslant  \llbracket \sigma \rrbracket^{\flat} \llbracket \mu
       \rrbracket^{- \flat} \frac{(1 + \llbracket \mu \rrbracket^{- 1} | z -
       z_1 |_s)^{\flat} (\llbracket \mu \rrbracket^{- 1} | z - z_1 |_s)^{2 -
       \flat}}{1 + \llbracket \mu \rrbracket^{- 2} | z - z_1 |_s^2}\\
       &\quad \leqslant  \llbracket \sigma \rrbracket^{\flat} \llbracket \mu
       \rrbracket^{- \flat} \frac{(1 + \llbracket \mu \rrbracket^{- 1} | z -
       z_1 |_s)^2}{1 + \llbracket \mu \rrbracket^{- 2} | z - z_1 |_s^2}
       \lesssim \llbracket \sigma \rrbracket^{\flat} \llbracket \mu
       \rrbracket^{- \flat},
     \end{aligned} \]
  which proves the claim.
\end{proof}

\begin{lemma}
  \label{lem:weight-interaction-v}For all $\omega > 0$ and $m\in\mathbb N_0$, the following bound
  \[ \frac{(1 - v^{(1 + m)}_{\mu}) w^{(1 + m), \omega}_{\mu}}{w^{(1 + m),
     \omega}_{\sigma}} \lesssim \llbracket \sigma \rrbracket^{\omega}
     \llbracket \mu \rrbracket^{- \omega} \]
  holds uniformly in $\mu, \sigma \in (0, 1)$.
\end{lemma}

\begin{proof}
  By Definition~\ref{def:weights-kernels}, $1 - v^{(m + 1)}_{\mu} = 0$ if
  $\llbracket \mu \rrbracket^{- 1} \operatorname{St} (z, y_1, \ldots, y_m) \leqslant
  1$. As a consequence,
  \[ \begin{aligned}
       \frac{(1 - v^{(1 + m)}_{\mu}) w^{(1 + m), \omega}_{\mu}}{w^{(1 + m),
       \omega}_{\sigma}} & =  \frac{\left( 1 + \llbracket \mu \rrbracket^{-
       1} \operatorname{St}(z, y_1, \ldots, y_m) 
       \right)^{\omega}}{\left( 1 + \llbracket \sigma \rrbracket^{- 1}
       \operatorname{St}(z, y_1, \ldots, y_m)  \right)^{\omega}}
       (1 - v^{(1 + m)}_{\mu})\\
       & =  \llbracket \mu \rrbracket^{- \omega} \llbracket \sigma
       \rrbracket^{\omega} \frac{\left( \llbracket \mu \rrbracket  +
       \operatorname{St}(z, y_1, \ldots, y_m) \right)^{\omega}}{\left(
       \llbracket \sigma \rrbracket  + \operatorname{St}(z, y_1, \ldots, y_m)
       \right)^{\omega}} (1 - v^{(1 + m)}_{\mu}),
     \end{aligned} \]
  is non-vanishing only if $\operatorname{St} (z, y_1, \ldots, y_m) > \llbracket \mu
  \rrbracket$. Hence
  \[ \frac{(1 - v^{(1 + m)}_{\mu}) w^{(1 + m), \omega}_{\mu}}{w^{(1 + m),
     \omega}_{\sigma}} \leqslant \llbracket \mu \rrbracket^{- \omega}
     \llbracket \sigma \rrbracket^{\omega} \frac{\left( 2
     \operatorname{St}(z, y_1, \ldots, y_m) \right)^{\omega}}{\left(
     \operatorname{St}(z, y_1, \ldots, y_m) \right)^{\omega}} (1 - v^{(m +
     1)}_{\mu}) \lesssim \llbracket \mu \rrbracket^{- \omega} \llbracket
     \sigma \rrbracket^{\omega}, \]
  which proves the claim.
\end{proof}

\subsection{Norm estimates}

\begin{lemma}
  \label{lem:estimates-1-h}The cumulant norms introduced in
  Definition~\ref{def:norm-cumulants} satisfy the bound
  \[ \trinorm{ [K^{1, 1}_{\sigma} \dot{\bar{F}}_{\sigma}^{[\ell], (1)}]
     \cdot (1 - h_{\mu}) }_{\mu} \lesssim \llbracket \mu
     \rrbracket^{- \flat} \llbracket \sigma \rrbracket^{\flat}  \trinorm{
     \dot{\bar{F}}_{\sigma}^{[\ell], (1)} }_{\sigma}, \]
  uniformly over $\sigma, \mu \in [1 / 2, 1)$ and
  $\dot{\bar{F}}_{\sigma}^{[\ell], (1)}$.
\end{lemma}

\begin{proof}
  By~\eqref{eq:weight-w-2}, Young's inequality and
  Lemma~\ref{lem:weight-interaction} we obtain
  \[ \begin{aligned}
       \trinorm{ [K^{1, 1}_{\sigma} \dot{\bar{F}}_{\sigma}^{[\ell], (1)}]
       \cdot (1 - h_{\mu}) }_{\mu} & =  \trinorm{ [K^{1, 1}_{\mu}
       ([K^{1, 1}_{\sigma} \dot{\bar{F}}_{\sigma}^{[\ell], (1)}] \cdot (1 -
       h_{\mu}))] \cdot w_{\mu}^{(2), \flat} }\\
       & \lesssim  \| K^{1, 1}_{\mu} \|_{\operatorname{TV}
       ((w_{\mu}^{\flat})^{\otimes 2})}  \trinorm{ [K^{1, 1}_{\sigma}
       \dot{\bar{F}}_{\sigma}^{[\ell], (1)}] \cdot (1 - h_{\mu}) w_{\mu}^{(2),
       \flat} }\\
       & \lesssim  \| w_{\mu}^{(2), \flat} (1 - h_{\mu})  (w_{\sigma}^{(2),
       \flat})^{- 1} \|_{L^{\infty}}  \trinorm{
       \dot{\bar{F}}_{\sigma}^{[\ell], (1)} }_{\sigma}\\
       & \lesssim  \llbracket \mu \rrbracket^{- \flat} \llbracket \sigma
       \rrbracket^{\flat} \trinorm{ \dot{\bar{F}}_{\sigma}^{[\ell], (1)}
       }_{\sigma} .
     \end{aligned} \]
  This proves the statement.
\end{proof}

\begin{lemma}
  \label{lem:estimates-1-v} For every \(\mathfrak a\in\mathfrak A\), the kernel norms introduced in
  Definition~\ref{def:norm-kernels} satisfy the bound
  \[ \| [\tilde{K}^{\mathfrak{a}}_{\sigma} \dot{F}^{\mathfrak{a}}_{\sigma}]
     \cdot (1 - v^{\mathfrak{a}}_{\mu}) \|_{\mu} \lesssim \llbracket \mu
     \rrbracket^{\ell (\mathfrak{a}) \kappa_{\mathfrak{o}} - \flat} \llbracket
     \sigma \rrbracket^{\flat - \ell (\mathfrak{a}) \kappa_{\mathfrak{o}}}  \|
     \dot{F}^{\mathfrak{a}}_{\sigma} \|_{\sigma}, \]
  uniformly over $\sigma, \mu \in [1 / 2, 1)$ and
  $\dot{F}^{\mathfrak{a}}_{\sigma}$.
\end{lemma}

\begin{proof}
  By~\eqref{eq:weight-w-o-2}, Young's inequality,
  Lemma~\ref{lem:weight-interaction-v} and $\tilde{w}^{\mathfrak{a}}_{\mu}
  \assign w^{(1 + k (\mathfrak{a})), \flat - \ell (\mathfrak{a})
  \kappa_{\mathfrak{o}}}_{\mu}$ we obtain
  \[ \begin{aligned}
       \| [\tilde{K}^{\mathfrak{a}}_{\sigma} \dot{F}^{\mathfrak{a}}_{\sigma}]
       \cdot (1 - v^{\mathfrak{a}}_{\mu}) \|_{\mu} & =  \|
       \mathfrak{o}^{\mathfrak{a}} \cdot [\tilde{K}^{\mathfrak{a}}_{\mu}
       ([\tilde{K}^{\mathfrak{a}}_{\sigma} \dot{F}^{\mathfrak{a}}_{\sigma}]
       \cdot (1 - v^{\mathfrak{a}}_{\mu}))] \cdot
       \tilde{w}_{\mu}^{\mathfrak{a}} \|\\
       & \lesssim  \| \tilde{K}^{\mathfrak{a}}_{\mu} \|_{\operatorname{TV}
       (\mathfrak{w}_{\mu})}  \| \mathfrak{o}^{\mathfrak{a}} \cdot
       [\tilde{K}^{\mathfrak{a}}_{\sigma} \dot{F}^{\mathfrak{a}}_{\sigma}]
       \cdot (1 - v^{\mathfrak{a}}_{\mu}) \tilde{w}_{\mu}^{\mathfrak{a}} \|\\
       & \lesssim  \| (1 - v^{\mathfrak{a}}_{\mu})
       \tilde{w}_{\mu}^{\mathfrak{a}}  (\tilde{w}_{\sigma}^{\mathfrak{a}})^{-
       1} \|_{L^{\infty}}  \| \dot{F}^{\mathfrak{a}}_{\sigma} \|_{\sigma}\\
       & \lesssim  \llbracket \mu \rrbracket^{\ell (\mathfrak{a})
       \kappa_{\mathfrak{o}} - \flat} \llbracket \sigma \rrbracket^{\flat -
       \ell (\mathfrak{a}) \kappa_{\mathfrak{o}}} \|
       \dot{F}^{\mathfrak{a}}_{\sigma} \|_{\sigma},
     \end{aligned} \]
  where
  \[ \mathfrak{w}_{\mu} (z, z_1, \ldots, z_{k
     (\mathfrak{a})}) {= w_{\mu}^2}  (z) {w_{\mu}^{\flat}}  (z_1) {\cdots
     w_{\mu}^{\flat}}  (z_{k (\mathfrak{a})}) . \]
  This proves the statement.
\end{proof}

\begin{lemma}[Auxiliary localisation]
	\label{lem:aux-weight-improved}
	Let $\mathfrak a\in\mathfrak A$ and $\mu\in[1/2,1)$, and set
	\[
	k\assign k(\mathfrak a),
	\qquad
	\tau\assign\llbracket\mu\rrbracket,
	\qquad
		\rho\assign\tau\vee\varepsilon,
		\qquad
		\delta\assign\frac{\tau}{\rho}\in(0,1].
	\]
	There exist nonnegative smooth functions
	\[
	\tilde v_{\mu,0}^{\mathfrak a},
	\qquad
	\tilde v_{\mu,1}^{\mathfrak a}
	\colon\Lambda_0^{1+k}\longrightarrow[0,\infty)
	\]
	such that, uniformly in $\mu$, $\varepsilon$, and $M$,
	\begin{equation}
		\label{eq:aux-weight-combined-Linfty}
		\|\tilde v_{\mu,0}^{\mathfrak a}\|_{L^\infty}
		+\|\tilde v_{\mu,1}^{\mathfrak a}\|_{L^\infty}
		\lesssim1.
	\end{equation}
	Set
	\(\rho_0\assign\tau\) and \(\rho_1\assign\rho\). Then, for \(j\in\{0,1\}\),
	\begin{equation}
		\label{eq:aux-weight-combined-support}
		\operatorname{supp} \tilde v_{\mu,j}^{\mathfrak a}
		\subset \left\{
		(z,z_1,\ldots,z_k):
		\begin{array}{l}
			|z_0-z_{i,0}|\leqslant c\tau^{2s},\\[2pt]
				|\bar z-\bar z_i|\leqslant c\rho_j,
			\quad i=1,\ldots,k,
		\end{array}
		\right\},
	\end{equation}
	where $c$ is independent of $\mu$, $\varepsilon$, and $M$.
	
	Define
	\begin{equation}
		\label{eq:aux-weight-combined-total}
		\tilde V_\mu^{\mathfrak a}
		\assign
		\tilde v_{\mu,0}^{\mathfrak a}
			+\delta^2\tilde v_{\mu,1}^{\mathfrak a}.
	\end{equation}
	Then, for every kernel $F^{\mathfrak a}$,
	\begin{equation}
		\label{eq:aux-weight-combined-bound}
		\|F^{\mathfrak a}v_\mu^{\mathfrak a}\|_\mu
		\lesssim
		\left\|
		\mathfrak o^{\mathfrak a}\cdot
		[\tilde K_\mu^{\mathfrak a}F^{\mathfrak a}]\cdot
		\tilde V_\mu^{\mathfrak a}
		\tilde w_\mu^{\mathfrak a}
		\right\|
		\lesssim
		\|F^{\mathfrak a}\|_\mu.
	\end{equation}
	Moreover, for every $\omega\in[0,2]$,
	\begin{equation}
		\label{eq:aux-weight-combined-volume}
		\sup_{z\in\Lambda}
		\int_{\Lambda^k}
		w_\mu^{(1+k),\omega}(z,Z^{\mathfrak a})
		\tilde V_\mu^{\mathfrak a}(z,Z^{\mathfrak a})
		\,\mathd Z^{\mathfrak a}
		\lesssim
			\tau^{2sk}\rho^{dk}.
	\end{equation}
	In particular, taking
	$\omega=\flat-\ell(\mathfrak a)\kappa_{\mathfrak o}$ in
	\eqref{eq:aux-weight-combined-volume} gives the estimate
	with $w_\mu^{(1+k),\omega}=\tilde w_\mu^{\mathfrak a}$.
\end{lemma}

\begin{remark}[Meaning of the auxiliary localisation]
    Lemma~\ref{lem:aux-weight-improved} allows one to transfer the original localisation
	\(v_\mu^{\mathfrak a}\) from the kernel
	\(F^{\mathfrak a}\) to the regularised kernel
	\(\tilde K_\mu^{\mathfrak a}F^{\mathfrak a}\). More precisely,
	\eqref{eq:aux-weight-combined-bound} shows that, up to a uniform
	constant, multiplication by \(v_\mu^{\mathfrak a}\) before applying
	the norm may be replaced by multiplication of
	\(\tilde K_\mu^{\mathfrak a}F^{\mathfrak a}\) by the auxiliary
	localisation
	\[
	\tilde V_\mu^{\mathfrak a}
	=
	\tilde v_{\mu,0}^{\mathfrak a}
	+\delta^2\tilde v_{\mu,1}^{\mathfrak a}.
	\]
	The first term is localised at the natural scale \(\tau\) in space
	and \(\tau^{2s}\) in time. The second term accounts for lattice
	translations generated by the spatial operator when \(\tau\) lies
	below the lattice spacing. Its spatial support may therefore extend
	to the effective lattice scale
	\(\rho=\tau\vee\varepsilon\), but it comes with the factor
	\[
	\delta^2
	=
	\left(\frac{\tau}{\rho}\right)^2.
	\]
	In the regime \(\tau\gtrsim\varepsilon\), one has
	\(\rho\asymp\tau\) and \(\delta\asymp1\), so both terms are
	localised at the same scale. In the sublattice regime
	\(\tau\ll\varepsilon\), the support of
	\(\tilde v_{\mu,0}^{\mathfrak a}\) is contained in the exact spatial
	diagonal, whereas \(\tilde v_{\mu,1}^{\mathfrak a}\) is supported on
	finitely many neighbouring lattice diagonals at distance
	\(O(\varepsilon)\).
	
	The estimate \eqref{eq:aux-weight-combined-volume} shows that this
	enlargement of the support does not change the natural weighted
	volume of the localisation. Indeed, a spatial displacement of size
	\(O(\varepsilon)\) may cost at most
	\[
	\left(\frac{\rho}{\tau}\right)^\omega
	=
	\delta^{-\omega}
	\]
	in the Steiner weight, while the translated part carries the
	compensating factor \(\delta^2\). The restriction
	\(\omega\leqslant2\) therefore yields
	\[
	\delta^2\delta^{-\omega}
	=
	\delta^{2-\omega}
	\leqslant1.
	\]
	Thus \(\tilde V_\mu^{\mathfrak a}\) has, uniformly at and below the
	lattice scale, the same effective weighted volume
	\(\tau^{2sk}\rho^{dk}\) as a cutoff localised in \(k\) input
	variables at temporal scale \(\tau^{2s}\) and spatial scale \(\rho\).
\end{remark}

\begin{proof}
	The case $k=0$ is immediate: take
	\[
	\tilde v_{\mu,0}^{\mathfrak a}=1,
	\qquad
	\tilde v_{\mu,1}^{\mathfrak a}=0.
	\]
	We henceforth assume $k\geqslant1$.
	
	We first record a weighted Young estimate used below. Set
	\[
	\mathfrak w_\mu
	\assign
	w_\mu^2\otimes(w_\mu^\flat)^{\otimes k}.
	\]
	By \eqref{eq:weight-w-o-2}, whenever $Q$ is a tensor-product convolution
	kernel satisfying
	\[
	\|Q\|_{\operatorname{TV}(\mathfrak w_\mu)}\lesssim1,
	\]
	one has
	\begin{equation}
		\label{eq:aux-weight-combined-young}
		\left\|
		\mathfrak o^{\mathfrak a}\cdot(Q*G)\cdot
		\tilde w_\mu^{\mathfrak a}
		\right\|
		\lesssim
		\left\|
		\mathfrak o^{\mathfrak a}\cdot G\cdot
		\tilde w_\mu^{\mathfrak a}
		\right\|.
	\end{equation}
	
	Put
	\[
	h\assign\tilde K_\mu^{\mathfrak a}F^{\mathfrak a}.
	\]
	For the tensor variables (the output variable when $j=0$, and the
	$j$-th input variable when $j=1,\ldots,k$), set
	\[
	q_0=1,
	\qquad q_j=2\quad(j\geqslant1),
	\qquad
	p_0=2,
	\qquad p_j=4\quad(j\geqslant1),
	\]
	and define
	\[
	\mathring{\mathcal L}_\tau
	\assign
	\bigotimes_{j=0}^k
	(1+\tau^{2s}\partial_{t_j})^{q_j},
	\qquad
	\mathring{\mathcal K}_\tau
	\assign\mathring{\mathcal L}_\tau^{-1},
	\]
	\[
	\bar{\mathcal L}_\tau
	\assign
	\bigotimes_{j=0}^k
	(1-\tau^2\Delta_{\varepsilon,j})^{p_j},
	\qquad
	\bar{\mathcal K}_\tau
	\assign\bar{\mathcal L}_\tau^{-1}.
	\]
	By Definitions~\ref{def:regularizing-kernels} and
	\ref{def:regularizing-kernels-for-kernels},
	\begin{equation}
		\label{eq:aux-weight-combined-factorisation}
		\tilde L_\mu^{\mathfrak a}
		=\mathring{\mathcal L}_\tau\bar{\mathcal L}_\tau,
		\qquad
		\tilde K_\mu^{\mathfrak a}
		=\mathring{\mathcal K}_\tau\bar{\mathcal K}_\tau,
		\qquad
		F^{\mathfrak a}
		=\mathring{\mathcal L}_\tau\bar{\mathcal L}_\tau h.
	\end{equation}
	
	\smallskip
	\noindent\emph{Temporal reduction.}
	We first commute the temporal part of the regularising operator through
	the cutoff. Let $b=\tau^{2s}$, $q\in\{1,2\}$, and
	$R_{\tau,q}=(1+b\partial_t)^{-q}$. For smooth compactly supported $g,u$,
	\begin{equation}
		\label{eq:aux-weight-combined-temporal-one}
		R_{\tau,q}[g(1+b\partial_t)^qu]
		=
		\sum_{m=0}^q\binom qm
		\sum_{r=0}^m(-1)^r\binom mr
		[(b\partial_t)^{m-r}R_{\tau,q}]*
		[u(b\partial_t)^rg].
	\end{equation}
	Indeed, expand $(1+b\partial_t)^q$ and integrate by parts $m$ times in
	the term containing $b^m\partial_t^m u$. Moreover,
	\begin{equation}
		\label{eq:aux-weight-combined-temporal-algebra}
		(b\partial_t)^nR_{\tau,q}
		=
		\sum_{\ell=0}^n
		(-1)^{n-\ell}\binom n\ell
		(1+b\partial_t)^{-(q-\ell)},
		\qquad 0\leqslant n\leqslant q.
	\end{equation}
	The term with exponent zero is the identity kernel. The explicit
	exponential formula for $(1+b\partial_t)^{-1}$ shows that every kernel on
	the right-hand side of
	\eqref{eq:aux-weight-combined-temporal-algebra} has uniformly bounded
	weighted total variation at scale $\tau$.
	
	Iterating \eqref{eq:aux-weight-combined-temporal-one} over the tensor
	variables yields a finite set $\Gamma$, scaled time-derivative operators
	\[
	D_{\tau,t}^{\gamma}
	\assign
	\prod_{j=0}^k
	(\tau^{2s}\partial_{t_j})^{\gamma_j},
	\qquad
	0\leqslant\gamma_j\leqslant q_j,
	\]
	and tensor-product temporal kernels $Q_{\tau,\gamma}$ such that
	\begin{equation}
		\label{eq:aux-weight-combined-temporal-tensor}
		\mathring{\mathcal K}_\tau
		[g\mathring{\mathcal L}_\tau u]
		=
		\sum_{\gamma\in\Gamma}
		Q_{\tau,\gamma}*[uD_{\tau,t}^{\gamma}g]
	\end{equation}
	and
	\begin{equation}
		\label{eq:aux-weight-combined-temporal-TV}
		\|Q_{\tau,\gamma}\|_{\operatorname{TV}(\mathfrak w_\mu)}
		\lesssim1.
	\end{equation}
	Since the temporal and spatial operators commute,
	\eqref{eq:aux-weight-combined-factorisation} and
	\eqref{eq:aux-weight-combined-temporal-tensor}, applied with
	$g=v_\mu^{\mathfrak a}$ and $u=\bar{\mathcal L}_\tau h$, give the common
	identity
	\begin{equation}
		\label{eq:aux-weight-combined-common-temporal}
		\tilde K_\mu^{\mathfrak a}
		(F^{\mathfrak a}v_\mu^{\mathfrak a})
		=
		\sum_{\gamma\in\Gamma}
		\bar{\mathcal K}_\tau Q_{\tau,\gamma}*
		[(\bar{\mathcal L}_\tau h)
		D_{\tau,t}^{\gamma}v_\mu^{\mathfrak a}].
	\end{equation}
	Thus the temporal part has been treated completely, and it remains only
	to handle the spatial operator in
	\eqref{eq:aux-weight-combined-common-temporal}.
	
	The finitely many functions
	$D_{\tau,t}^{\gamma}v_\mu^{\mathfrak a}$ are uniformly bounded and are
	supported where
	\begin{equation}
		\label{eq:aux-weight-combined-chi-support}
		|z_0-z_{j,0}|\leqslant C\tau^{2s},
		\qquad
		|\bar z-\bar z_j|\leqslant C\tau,
		\qquad j=1,\ldots,k.
	\end{equation}
	We may therefore choose a nonnegative smooth function $\chi_\tau$,
	uniformly bounded and supported in the same set up to increasing $C$, such
	that
	\begin{equation}
		\label{eq:aux-weight-combined-chi-majorant}
		|D_{\tau,t}^{\gamma}v_\mu^{\mathfrak a}|
		\lesssim\chi_\tau,
		\qquad \gamma\in\Gamma.
	\end{equation}
	By Lemma~\ref{lem:K-kernel-aux} and weighted Young's inequality,
	\begin{equation}
		\label{eq:aux-weight-combined-spatial-TV}
		\|\bar{\mathcal K}_\tau\|_{
		\operatorname{TV}(\mathfrak w_\mu)}
		\lesssim1.
	\end{equation}
	
	Fix a sufficiently small constant $c_0>0$, to be chosen below. We now
	treat the spatial part of
	\eqref{eq:aux-weight-combined-common-temporal} in two regimes.
	
	\smallskip
	\noindent\emph{1. The regime $\tau\geqslant c_0\varepsilon$.}
	In this regime $\rho\asymp\tau$. We use the exact discrete product rule
	\begin{equation}
		\label{eq:aux-weight-combined-discrete-product}
		\Delta_\varepsilon(ug)
		=g\Delta_\varepsilon u+u\Delta_\varepsilon g
		+\sum_{i=1}^d
		\left[
		(\partial^{i+}u)(\partial^{i+}g)
		+(\partial^{i-}u)(\partial^{i-}g)
		\right].
	\end{equation}
	When a mixed first-difference term occurs, discrete summation by parts
	gives, for example,
	\begin{equation}
		\label{eq:aux-weight-combined-paired-differences}
		\begin{aligned}
			&\tau^2\bar K_\tau*
			[(\partial^{i+}u)(\partial^{i+}g)]
			\\
			&\quad=
			[\tau\partial^{i+}\bar K_\tau]*
			[u\,T_i^-(\tau\partial^{i+}g)]
			-\bar K_\tau*
			[u\,\tau^2\partial^{i-}\partial^{i+}g],
		\end{aligned}
	\end{equation}
	where $\bar K_\tau$ is the kernel introduced in
	Lemma~\ref{lem:K-kernel-aux} and $T_i^\pm$ denotes translation of the
	spatial variable by $\pm\varepsilon e_i$. The backward-difference term is
	treated analogously. Thus every isolated first spatial difference
	transferred to a regularising kernel is paired with a first spatial
	difference of the cutoff.
	
	Iterating these identities in all $1+k$ spatial tensor variables in
	\eqref{eq:aux-weight-combined-common-temporal} yields a finite set
	$\mathcal I$, constants $c_{A,B}$, and scaled spatial
	difference-translation operators $D_{\tau,x}^A,D_{\tau,x}^B$ such that
	\begin{equation}
		\label{eq:aux-weight-combined-decomposition-large}
		\tilde K_\mu^{\mathfrak a}
		(F^{\mathfrak a}v_\mu^{\mathfrak a})
		=
		\sum_{\gamma\in\Gamma}
		\sum_{(A,B)\in\mathcal I}
		c_{A,B}
		[Q_{\tau,\gamma}*(D_{\tau,x}^A
		\bar{\mathcal K}_\tau)]*
		[h\,D_{\tau,x}^B D_{\tau,t}^{\gamma}
		v_\mu^{\mathfrak a}].
	\end{equation}
	Every spatial difference in $D_{\tau,x}^A,D_{\tau,x}^B$ is multiplied by
	the corresponding power of $\tau$.
	
	Since $\tau\geqslant c_0\varepsilon$, the kernel estimates in
	Lemmas~\ref{lem:K-kernel} and \ref{lem:K-kernel-aux}, together with
	\eqref{eq:aux-weight-combined-temporal-TV}, imply, for every term in the
	finite expansion,
	\begin{equation}
		\label{eq:aux-weight-combined-kernel-large}
		\|Q_{\tau,\gamma}*(D_{\tau,x}^A
		\bar{\mathcal K}_\tau)\|_{
		\operatorname{TV}(\mathfrak w_\mu)}
		\lesssim_{c_0}1.
	\end{equation}
	Indeed, complete second-order spatial blocks are uniformly bounded, while
	an isolated first difference costs at most
	$1+\varepsilon/\tau\lesssim_{c_0}1$. Since
	$v_\mu^{\mathfrak a}$ is obtained from a fixed smooth function by
	parabolic rescaling at scale $\tau$, all the scaled cutoff derivatives in
	\eqref{eq:aux-weight-combined-decomposition-large} satisfy
	\begin{equation}
		\label{eq:aux-weight-combined-cutoff-large}
		\|D_{\tau,x}^B D_{\tau,t}^{\gamma}
		v_\mu^{\mathfrak a}\|_{L^\infty}
		\lesssim_{c_0}1.
	\end{equation}
	Their supports are contained in a fixed enlargement of the support of
	$v_\mu^{\mathfrak a}$. We may therefore choose a nonnegative smooth
	function $\tilde v_{\mu,0}^{\mathfrak a}$ satisfying
	\[
	\|\tilde v_{\mu,0}^{\mathfrak a}\|_{L^\infty}\lesssim1,
	\qquad
	|D_{\tau,x}^B D_{\tau,t}^{\gamma}v_\mu^{\mathfrak a}|
	\lesssim\tilde v_{\mu,0}^{\mathfrak a},
	\]
	for all terms occurring in
	\eqref{eq:aux-weight-combined-decomposition-large}, and supported in
	\eqref{eq:aux-weight-combined-support}. Set
	\[
	\tilde v_{\mu,1}^{\mathfrak a}=0.
	\]
	Combining \eqref{eq:aux-weight-combined-young},
	\eqref{eq:aux-weight-combined-decomposition-large},
	\eqref{eq:aux-weight-combined-kernel-large}, and
	\eqref{eq:aux-weight-combined-cutoff-large} proves the first inequality in
	\eqref{eq:aux-weight-combined-bound}. The second follows from
	\eqref{eq:aux-weight-combined-Linfty}.
	
	On the support of $\tilde v_{\mu,0}^{\mathfrak a}$ one has
	\[
	w_\mu^{(1+k),\omega}\lesssim1
	\qquad(\omega\in[0,2]),
	\]
	because $\rho\asymp\tau$. At fixed output point, the temporal support has
	volume $O(\tau^{2sk})$, while its spatial support has normalised lattice
	volume $O(\rho^{dk})$. This proves
	\eqref{eq:aux-weight-combined-volume} in the present regime.
	
	\smallskip
	\noindent\emph{2. The regime $\tau<c_0\varepsilon$.}
	We now rewrite the discrete spatial derivatives as finite linear combination of lattice translations and do not commute them through the cutoff. Since $\rho=\varepsilon$ in this regime,
	\[
	\frac{\tau^2}{\varepsilon^2}=\delta^2<1.
	\]
	Put
	\[
	\Delta_j^{(1)}\assign\varepsilon^2\Delta_{\varepsilon,j}.
	\]
		For a tensor point
		\[
		Z=(z,z_1,\ldots,z_k),
		\qquad
		z=(z_0,\bar z),
		\qquad
		z_j=(z_{j,0},\bar z_j),\quad j=1,\ldots,k,
		\]
		and
		\(\boldsymbol\nu=(\nu_0,\ldots,\nu_k)\in(\mathbb Z^d)^{1+k}\),
		define the spatial translation of \(Z\) by
		\[
		\Theta_{\boldsymbol\nu}Z
		\assign
		\bigl(
		(z_0,\bar z+\varepsilon\nu_0),
	(z_{1,0},\bar z_1+\varepsilon\nu_1),
	\ldots,
	(z_{k,0},\bar z_k+\varepsilon\nu_k)
	\bigr).
	\]
	Thus \(\Theta_{\boldsymbol\nu}\) leaves all time coordinates unchanged.  We
	denote by \(T_{\boldsymbol\nu}\) the induced translation operator on
	functions,
		\begin{equation}
			\label{eq:aux-weight-combined-translation-convention}
			(T_{\boldsymbol\nu}f)(Z)
			\assign
			f(\Theta_{\boldsymbol\nu}Z).
		\end{equation}
	We use the analogous notation
	\(\Theta_{j,\nu_j}\) and \(T_{j,\nu_j}\) when only the spatial component
	of the \(j\)-th tensor variable is translated.
	
	For each \(j\), the finite-range operator
	\((1-\delta^2\Delta_j^{(1)})^{p_j}\) admits an expansion
	\[
	(1-\delta^2\Delta_j^{(1)})^{p_j}
	=
	\sum_{\nu_j\in\mathcal N_j}
	c_{j,\nu_j}(\delta^2)T_{j,\nu_j},
	\]
	where \(\mathcal N_j\subset\mathbb Z^d\) is finite and
	\[
	|c_{j,0}(\delta^2)|\lesssim1,
	\qquad
	|c_{j,\nu_j}(\delta^2)|\lesssim\delta^2
	\quad(\nu_j\neq0).
	\]
	Tensorising gives
	\begin{equation}
		\label{eq:aux-weight-combined-stencil}
		\bar{\mathcal L}_\tau
		=
		\sum_{\boldsymbol\nu\in\mathcal N}
		c_{\boldsymbol\nu}(\delta^2)T_{\boldsymbol\nu},
	\end{equation}
	where
	\(\mathcal N\subset(\mathbb Z^d)^{1+k}\) is finite.  With
	\[
	r(\boldsymbol\nu)
	\assign
	\#\{j\in\{0,\ldots,k\}:\nu_j\neq0\},
	\]
	one has
	\begin{equation}
		\label{eq:aux-weight-combined-coefficients}
		|c_{\boldsymbol\nu}(\delta^2)|
		\lesssim
		\delta^{2r(\boldsymbol\nu)}.
	\end{equation}
	Inserting \eqref{eq:aux-weight-combined-stencil} into the common temporal
	reduction \eqref{eq:aux-weight-combined-common-temporal} gives
	\begin{equation}
		\label{eq:aux-weight-combined-exact-expansion}
		\tilde K_\mu^{\mathfrak a}
		(F^{\mathfrak a}v_\mu^{\mathfrak a})
		=
		\sum_{\boldsymbol\nu\in\mathcal N}
		\sum_{\gamma\in\Gamma}
		c_{\boldsymbol\nu}(\delta^2)
		\bar{\mathcal K}_\tau Q_{\tau,\gamma}*
		\bigl[
		(T_{\boldsymbol\nu}h)
		D_{\tau,t}^{\gamma}v_\mu^{\mathfrak a}
		\bigr].
	\end{equation}
	
	Enlarge, if necessary, the constant \(C\) in
	\eqref{eq:aux-weight-combined-chi-support} so that \(\chi_\tau\) is
	supported where
	\[
	|z_0-z_{j,0}|\leqslant C\tau^{2s},
	\qquad
	|\bar z-\bar z_j|\leqslant C\tau,
	\qquad j=1,\ldots,k.
	\]
	Choose \(c_0>0\) so small that \(2Cc_0<1\).  If
	\(\tau<c_0\varepsilon\), then, on the support of \(\chi_\tau\),
	\[
	|\bar z_i-\bar z_j|
	\leqslant
	|\bar z_i-\bar z|+|\bar z-\bar z_j|
	\leqslant2C\tau
	<\varepsilon,
	\]
	where we use the convention \(\bar z_0\assign\bar z\).  Since all spatial
	variables belong to \(\varepsilon\mathbb Z^d\), this implies
	\begin{equation}
		\label{eq:aux-weight-combined-diagonal}
		\bar z=\bar z_1=\cdots=\bar z_k.
	\end{equation}
	
	For \(\boldsymbol\nu\in\mathcal N\), define the translated cutoff
	\begin{equation}
		\label{eq:aux-weight-combined-translated-cutoff}
		\chi_{\tau,\boldsymbol\nu}
		\assign
		T_{-\boldsymbol\nu}\chi_\tau.
	\end{equation}
	Equivalently, for a tensor point
	\(Y=(y,y_1,\ldots,y_k)\),
	\begin{equation}
		\label{eq:aux-weight-combined-translated-cutoff-pointwise}
		\chi_{\tau,\boldsymbol\nu}(Y)
		=
		\chi_\tau(\Theta_{-\boldsymbol\nu}Y).
	\end{equation}
	We claim that
	\begin{equation}
		\label{eq:aux-weight-combined-translation}
		\left\|
		\mathfrak o^{\mathfrak a}\cdot
		(T_{\boldsymbol\nu}h)
		D_{\tau,t}^{\gamma}v_\mu^{\mathfrak a}
		\tilde w_\mu^{\mathfrak a}
		\right\|
		\lesssim
		\left\|
		\mathfrak o^{\mathfrak a}\cdot h\cdot
		\chi_{\tau,\boldsymbol\nu}
		\tilde w_\mu^{\mathfrak a}
		\right\|.
	\end{equation}
	Indeed, by \eqref{eq:aux-weight-combined-chi-majorant}, the left-hand side
	is bounded by the same norm with
	\(D_{\tau,t}^{\gamma}v_\mu^{\mathfrak a}\) replaced by \(\chi_\tau\).
	Make the change of variables
	\[
		Y=\Theta_{\boldsymbol\nu}Z.
		\]
		By \eqref{eq:aux-weight-combined-translation-convention},
		\[
		(T_{\boldsymbol\nu}h)(Z)=h(Y),
		\qquad
		\chi_\tau(Z)
		=
		\chi_\tau(\Theta_{-\boldsymbol\nu}Y)
		=
	\chi_{\tau,\boldsymbol\nu}(Y).
	\]
	The normalised lattice measure and the supremum in the output variable
	are translation invariant.  Writing
	\(Y=(y,y_1,\ldots,y_k)\), the output variables \(z\) and \(y\) differ by
	one of finitely many spatial shifts of size \(O(\varepsilon)\).  Since
	\(\varepsilon\leqslant1\),
	\[
	\mathfrak o^{\mathfrak a}(z)
	\lesssim
	\mathfrak o^{\mathfrak a}(y).
		\]
		
		It remains to compare the Steiner weights.  On the support of
		\(\chi_\tau(Z)\), all spatial components of \(Z\) coincide by
		\eqref{eq:aux-weight-combined-diagonal}.  Let
	\(\operatorname{St}_t(z_0,z_{1,0},\ldots,z_{k,0})\) denote the Steiner
	diameter of the time coordinates for the metric
	\(|t-t'|^{1/(2s)}\).  Since the time projection is \(1\)-Lipschitz for
	the parabolic metric and \(\Theta_{\boldsymbol\nu}\) does not change the
	time coordinates,
	\[
		\operatorname{St}(Z)
		=
		\operatorname{St}_t(z_0,z_{1,0},\ldots,z_{k,0})
		\leqslant
	\operatorname{St}(Y).
	\]
	The exponent in \(\tilde w_\mu^{\mathfrak a}\) is nonnegative, and hence
	\[
		\tilde w_\mu^{\mathfrak a}(Z)
		\leqslant
		\tilde w_\mu^{\mathfrak a}(Y).
	\]
	Together with the preceding observations, this proves
	\eqref{eq:aux-weight-combined-translation}.
	
	Combining \eqref{eq:aux-weight-combined-young},
	\eqref{eq:aux-weight-combined-temporal-TV},
	\eqref{eq:aux-weight-combined-spatial-TV},
	\eqref{eq:aux-weight-combined-exact-expansion}, and
	\eqref{eq:aux-weight-combined-translation}, we obtain
	\begin{equation}
		\label{eq:aux-weight-combined-pre-majorants}
		\|F^{\mathfrak a}v_\mu^{\mathfrak a}\|_\mu
		\lesssim
		\sum_{\boldsymbol\nu\in\mathcal N}
		|c_{\boldsymbol\nu}(\delta^2)|
		\left\|
		\mathfrak o^{\mathfrak a}\cdot h\cdot
		\chi_{\tau,\boldsymbol\nu}
		\tilde w_\mu^{\mathfrak a}
		\right\|.
	\end{equation}
	Define
	\begin{align*}
		\tilde v_{\mu,0}^{\mathfrak a}
		&\assign
		\sum_{r(\boldsymbol\nu)=0}
		|c_{\boldsymbol\nu}(\delta^2)|
		\chi_{\tau,\boldsymbol\nu},
		\\
		\tilde v_{\mu,1}^{\mathfrak a}
		&\assign
		\sum_{r(\boldsymbol\nu)\geqslant1}
		\delta^{-2}
		|c_{\boldsymbol\nu}(\delta^2)|
		\chi_{\tau,\boldsymbol\nu}.
	\end{align*}
	The coefficient estimate
	\eqref{eq:aux-weight-combined-coefficients} and the finiteness of
	$\mathcal N$ imply \eqref{eq:aux-weight-combined-Linfty}. Their supports
	satisfy \eqref{eq:aux-weight-combined-support}, because all shifts have
	size $O(\varepsilon)=O(\rho)$ and
	\eqref{eq:aux-weight-combined-chi-support} has spatial width $O(\tau)$.
	
	Since $\mathcal N$ is finite, after increasing the two majorants by a
	fixed multiplicative constant, each norm in
	\eqref{eq:aux-weight-combined-pre-majorants} is bounded by the norm with
	the corresponding sum of cutoffs. Consequently,
	\eqref{eq:aux-weight-combined-pre-majorants} proves the first inequality
	in \eqref{eq:aux-weight-combined-bound}. The second follows from
	\eqref{eq:aux-weight-combined-Linfty} and $\delta\leqslant1$.
	
	It remains to prove \eqref{eq:aux-weight-combined-volume}. Fix the output
	point $z$. For every $\boldsymbol\nu\in\mathcal N$, the support of
	$\chi_{\tau,\boldsymbol\nu}(z,\cdot)$ has temporal volume
	$O(\tau^{2sk})$. By \eqref{eq:aux-weight-combined-diagonal}, each input
	spatial coordinate is fixed relative to the output coordinate at one of
	finitely many lattice offsets. Its normalised spatial volume is therefore
	$O(\varepsilon^{dk})$. Hence
	\begin{equation}
		\label{eq:aux-weight-combined-chi-volume}
		\sup_z\int_{\Lambda^k}
		\chi_{\tau,\boldsymbol\nu}(z,Z)\,\mathd Z
		\lesssim
		\tau^{2sk}\varepsilon^{dk}.
	\end{equation}
	If $r(\boldsymbol\nu)=0$, then $\boldsymbol\nu=0$ and
	$w_\mu^{(1+k),\omega}\lesssim1$ on the support. If
	$r(\boldsymbol\nu)\geqslant1$, the finite range of the shifts and
	\eqref{eq:aux-weight-combined-chi-support} imply
	\[
	\operatorname{St}(z,z_1,\ldots,z_k)\lesssim\varepsilon.
	\]
	Therefore, for every $\omega\in[0,2]$,
	\begin{equation}
		\label{eq:aux-weight-combined-shifted-weight}
		w_\mu^{(1+k),\omega}(z,z_1,\ldots,z_k)
		\lesssim
		\left(\frac{\varepsilon}{\tau}\right)^\omega
		=\delta^{-\omega}.
	\end{equation}
	Since $\rho=\varepsilon$ in this regime,
	\eqref{eq:aux-weight-combined-chi-volume} and
	\eqref{eq:aux-weight-combined-shifted-weight} yield
	\begin{align*}
		\sup_z\int
		w_\mu^{(1+k),\omega}
		\tilde v_{\mu,0}^{\mathfrak a}
		\,\mathd Z
		&\lesssim
		\tau^{2sk}\rho^{dk},
		\\
		\delta^2\sup_z\int
		w_\mu^{(1+k),\omega}
		\tilde v_{\mu,1}^{\mathfrak a}
		\,\mathd Z
		&\lesssim
		\tau^{2sk}\rho^{dk}\,
		\delta^{2-\omega}.
	\end{align*}
	Since $\omega\in[0,2]$ and $\delta\in(0,1]$, both right-hand sides
	are bounded by $C\tau^{2sk}\rho^{dk}$. This proves
	\eqref{eq:aux-weight-combined-volume} and completes the proof.
\end{proof}

\begin{lemma}
  \label{lem:estimates-1-K-coeff} For every \(\mathfrak a\in\mathfrak A\), the kernel norms introduced in
  Definition~\ref{def:norm-kernels} satisfy the bound
  \[ \| [(1 - \tilde{K}^{\mathfrak{a}}_{\sigma})
     \dot{F}^{\mathfrak{a}}_{\sigma}] \cdot (1 - v_{\mu}^{\mathfrak{a}})
     \|_{\mu} \lesssim \llbracket \mu \rrbracket^{- 2 s}  \llbracket \sigma
     \rrbracket^{2 s}  \| \dot{F}^{\mathfrak{a}}_{\sigma} \|_{\sigma}, \]
  uniformly over $\mu \in [1 / 2, 1)$, $\sigma \in [\mu, 1)$ and
  $\dot{F}^{\mathfrak{a}}_{\sigma}$.
\end{lemma}

\begin{proof}
  By Lemma~\ref{lem:aux-weight-improved} we have
  \[ \| [(1 - \tilde{K}^{\mathfrak{a}}_{\sigma})
     \dot{F}^{\mathfrak{a}}_{\sigma}] \cdot (1 - v_{\mu}^{\mathfrak{a}})
     \|_{\mu} \lesssim \| (1 - \tilde{K}^{\mathfrak{a}}_{\sigma})
     \dot{F}^{\mathfrak{a}}_{\sigma} \|_{\mu} = \| \mathfrak{o}^{\mathfrak{a}}
     \cdot [(\tilde{L}^{\mathfrak{a}}_{\sigma} - 1)
     \tilde{K}^{\mathfrak{a}}_{\mu} \tilde{K}^{\mathfrak{a}}_{\sigma}
     \dot{F}^{\mathfrak{a}}_{\sigma}] \cdot \tilde{w}^{\mathfrak{a}}_{\mu}
     \|. 
  \]
  By~\eqref{eq:weight-w-o-2}, Young's inequality, Lemma~\ref{lem:K-kernel} and
  $\tilde{w}_{\mu}^{\mathfrak{a}} \leqslant \tilde{w}_{\sigma}^{\mathfrak{a}}$
  we obtain
  \begin{align*}
       \| \mathfrak{o}^{\mathfrak{a}} \cdot
       [(\tilde{L}^{\mathfrak{a}}_{\sigma} - 1) \tilde{K}^{\mathfrak{a}}_{\mu}
       \tilde{K}^{\mathfrak{a}}_{\sigma} \dot{F}^{\mathfrak{a}}_{\sigma}]
       \cdot \tilde{w}^{\mathfrak{a}}_{\mu} \| & \lesssim  \|
       (\tilde{L}^{\mathfrak{a}}_{\sigma} - 1) \tilde{K}^{\mathfrak{a}}_{\mu}
       \|_{\operatorname{TV} (\mathfrak{w}_{\mu})}  \|
       \mathfrak{o}^{\mathfrak{a}} \cdot [\tilde{K}^{\mathfrak{a}}_{\sigma}
       \dot{F}^{\mathfrak{a}}_{\sigma}] \cdot \tilde{w}^{\mathfrak{a}}_{\mu}
       \|\\
       & \lesssim  \llbracket \mu \rrbracket^{- 2 s} \llbracket \sigma
       \rrbracket^{2 s}  \| \dot{F}^{\mathfrak{a}}_{\sigma} \|_{\sigma},
  \end{align*}
  where $\mathfrak{w}_{\mu} (z, z_1, \ldots, z_{k
  (\mathfrak{a})}) {= w_{\mu}^2}  (z) {w_{\mu}^{\flat}}  (z_1) {\cdots
  w_{\mu}^{\flat}}  (z_{k (\mathfrak{a})})$.
\end{proof}

\begin{lemma}
  \label{lem:estimates-1-K}The cumulant norms introduced in
  Definition~\ref{def:norm-cumulants} satisfy the bound
  \[ \trinorm{ (1 - K^{1, 1}_{\sigma}) \dot{\bar{F}}_{\sigma}^{[\ell], (1)}
     }_{\mu} \lesssim \llbracket \mu \rrbracket^{- 2 s} \llbracket
     \sigma \rrbracket^{2 s}  \trinorm{ \dot{\bar{F}}_{\sigma}^{[\ell], (1)}
     }_{\sigma}, \]
  uniformly over $\mu \in [1 / 2, 1)$ and $\sigma \in [\mu, 1)$.
\end{lemma}

\begin{proof}
  By Young's inequality, Lemma~\ref{lem:K-kernel} and $w_{\mu}^{\mathfrak{a}}
  \leqslant w_{\sigma}^{\mathfrak{a}}$ we obtain
  \[ \begin{aligned}
       \trinorm{ (1 - K^{1, 1}_{\sigma}) \dot{\bar{F}}_{\sigma}^{[\ell],
       (1)} }_{\mu} & =  \trinorm{ [(1 \otimes (L_{\sigma} - 1)
       K_{\mu}) K^{1, 1}_{\sigma} \dot{\bar{F}}_{\sigma}^{[\ell], (1)}] \cdot
       w_{\mu}^{(2), \flat} }\\
       & \lesssim  \| (L_{\sigma} - 1) K_{\mu} \|_{\operatorname{TV}
       (w_{\mu}^{\flat})}  \trinorm{ [K^{1, 1}_{\sigma}
       \dot{\bar{F}}_{\sigma}^{[\ell], (1)}] \cdot w_{\mu}^{(2), \flat}
       }\\
       & \lesssim  \llbracket \mu \rrbracket^{- 2 s} \llbracket \sigma
       \rrbracket^{2 s}  \trinorm{ \dot{\bar{F}}_{\sigma}^{[\ell], (1)}
       }_{\sigma},
     \end{aligned} \]
  as claimed.
\end{proof}

\subsection{Schauder estimates}

In this section, we establish estimates that characterise the spacetime
regularity of distributions $\phi$ in terms of the size of $\mathcal{L} \phi$.

\begin{lemma}
  \label{lem:LP-block}For all $\bar A \in \mathbb{N}_0^{\{1 \pm, \ldots, d \pm
  \}}$ the following bounds
  \begin{align*}
    \| \partial^{\bar A} \bar{\Delta}_i \|_{\operatorname{TV}(\zeta^{- 1})}
    &\lesssim \llbracket \mu_i \rrbracket^{- | \bar A |},\\
    \| (1 - \Delta_{\varepsilon})^{- 1}
    \bar{\Delta}_i \|_{\operatorname{TV}(\zeta^{- 1})}
    &\lesssim \llbracket \mu_i \rrbracket^2,\\
    \| (-\Delta_{\varepsilon})^s
    \bar{\Delta}_i \|_{\operatorname{TV}(\zeta^{- 1})}
    &\lesssim \llbracket \mu_i \rrbracket^{- 2 s}
  \end{align*}
  hold uniformly in $i \geqslant - 1$.
\end{lemma}

\begin{proof}
  The estimates follow from the identity $\llbracket \mu_i \rrbracket = 2^{-
  i - 2}$ and Definition~\ref{def:space-LP} of the spatial Littlewood-Paley blocks.
\end{proof}

\begin{lemma}
  \label{lem:estimates-for-G}For all $A \in \mathbb{N}_0^{\{ 0, 1 \pm, \ldots,
  d \pm \}}$ the following bounds
  \begin{align*}
     \| \bar{\Delta}_i  G \|_{\operatorname{TV} (\zeta^{- 1})}
     &\lesssim \llbracket \mu_i \rrbracket_{}^{2 s},\\
     \| \partial^A \dot{G}_{\eta} \|_{\operatorname{TV} (\zeta^{- 1})}
     &\lesssim \llbracket \eta \rrbracket^{2 s - 1 - | A |_s},\\
     \| (G - G_{\eta}) \|_{\operatorname{TV} (\zeta^{- 1})}
     &\lesssim \llbracket \eta \rrbracket^{2 s},
  \end{align*}
  hold uniformly in $i \geqslant - 1$ and $\eta \in (1 / 2, 1)$.
\end{lemma}

\begin{proof}
  Recall that $\bar{\Delta}_i $ is the spatial Littlewood-Paley block, $\zeta
  (z) \assign (1 + | z |^2_s)^{- 1 / 2}$ and note that
  \[ G = \int_{\mu_i}^1 \dot{G}_{\mu} \mathd \mu + \int_{1 / 2}^{\mu_i}
     \dot{G}_{\mu} \mathd \mu + G_{1 / 2} . \]
  Taking into account Remark~\ref{rem:weight-properties}~a), we obtain
  \[ \begin{aligned}
       \| \bar{\Delta}_i  G \|_{\operatorname{TV} (\zeta^{- 1})} & \lesssim  \|
       \bar{\Delta}_i  \|_{\operatorname{TV} (\zeta^{- 1})} \int_{\mu_i}^1 \|
       \dot{G}_{\mu} \|_{\operatorname{TV} (\zeta^{- 1})} \mathd \mu\\
       &   + \| (1 - \Delta_{\varepsilon})^{- 1} \bar{\Delta}_i
       \|_{\operatorname{TV} (\zeta^{- 1})} \int_{1 / 2}^{\mu_i} \| (1 -
       \Delta_{\varepsilon}) \dot{G}_{\mu} \|_{\operatorname{TV} (\zeta^{- 1})} \mathd
       \mu\\
       &   + \| (1 - \Delta_{\varepsilon})^{- 1} \bar{\Delta}_i
       \|_{\operatorname{TV} (\zeta^{- 1})}  \| (1 - \Delta_{\varepsilon}) G_{1 / 2}
       \|_{\operatorname{TV} (\zeta^{- 1})} .
     \end{aligned} \]
  Consequently, using Lemma~\ref{lem:Gdot-decay} we arrive at
  \[ \begin{aligned}
       \| \bar{\Delta}_i  G \|_{\operatorname{TV} (\zeta^{- 1})} & \lesssim
       \int_{\mu_i}^1 \llbracket \mu  \rrbracket^{2 s - 1} \mathd \mu +
       \llbracket \mu_i \rrbracket^2 \int_{1 / 2}^{\mu_i} {\llbracket \mu
       \rrbracket }^{2 s - 3} \mathd \mu + \llbracket \mu_i \rrbracket^2
       \times 1\\
       & \lesssim  \llbracket \mu_i \rrbracket_{}^{2 s} + \llbracket \mu_i
       \rrbracket_{}^2  \llbracket \mu_i \rrbracket ^{2 s - 2} + \llbracket
       \mu_i \rrbracket^2 \lesssim \llbracket \mu_i \rrbracket_{}^{2 s} .
     \end{aligned} \]
  This proves the first bound. 
  The second bound is precisely
  \eqref{eq:Gdot-corrected-fixed-weight}.
  To show the third bound notice that
  \[ \| (G - G_{\eta}) \|_{\operatorname{TV} (\zeta^{- 1})} \leqslant \int_{\eta}^1 \|
     \dot{G}_{\sigma} \|_{\operatorname{TV} (\zeta^{- 1})} \mathd \sigma \lesssim
     \int_{\eta}^1 \llbracket \sigma \rrbracket^{2 s - 1} \mathd \sigma
     \lesssim \llbracket \eta \rrbracket^{2 s} . \]
  This completes the proof.
\end{proof}

\begin{lemma}
  \label{lem:schauder-full}
  Recall that $\gamma\in(0,s)$ and $\gamma=a\nu$. Then
  \begin{equation}
    \| \rho_{\eta} \bar{\Delta}_i  \phi \| \lesssim \llbracket \eta
    \rrbracket^{- \gamma} \llbracket \bar{\mu} \rrbracket^{\gamma} \llbracket
    \mu_i \vee \bar{\mu} \rrbracket^{- \gamma} [\trinorm{ \phi } +
    \llbracket \bar{\mu} \rrbracket^{2 s - 2 \gamma} \trinorm{ K_{\Cdot}
    \mathcal{L} \phi }_{\#}], \label{eq:schauder}
  \end{equation}
  uniformly in $\phi \in \mathcal{S}' (\Lambda)$, $i \geqslant - 1$ and $\eta
  \in [\bar{\mu}, 1)$, where the norms $\trinorm{ \Cdot } =
  \trinorm{ \Cdot }_{\bar{\mu}}$ and $\trinorm{ \Cdot
  }_{\#} = \trinorm{ \Cdot }_{\#, \bar{\mu}}$
  are defined by~\eqref{eq:sol-norm} and~\eqref{eq:sharp-norm}. In particular,
  for $\hat{\gamma} > \gamma$, we have
  \[ \sum_i 2^{- \hat{\gamma} i}  \| \rho_{\eta} \bar{\Delta}_i  \phi \|
     \lesssim \llbracket \eta \rrbracket^{- \gamma} [\trinorm{ \phi
     } + \llbracket \bar{\mu} \rrbracket^{2 s - 2 \gamma}
     \trinorm{ K_{\Cdot} \mathcal{L} \phi }_{\#}], \]
  uniformly in $\bar{\mu} \in [1 / 2, 1)$, $\eta \in [\bar{\mu}, 1)$ and $\phi
  \in \mathcal{S}' (\Lambda)$.
\end{lemma}

\begin{proof}
  First observe that
  \begin{equation}
    \begin{aligned}
      \| \rho_{\eta}  \bar{\Delta}_i  \phi \|
      &\lesssim \sup_j \| \rho_{\eta} \chi_j  \bar{\Delta}_i  \phi \|\\
      &\lesssim \sup_j \| \rho_{\eta} \chi_j
      \zeta_{\mu_j \vee \bar{\mu}}^{- 1} \|
      \| \zeta_{\mu_j \vee \bar{\mu}} \bar{\Delta}_i  \phi \|\\
      &\lesssim \llbracket \eta \rrbracket^{- \gamma}
      \sup_j \llbracket \mu_j \vee \bar{\mu} \rrbracket^{\gamma}
      \| \zeta_{\mu_j \vee \bar{\mu}}  \bar{\Delta}_i  \phi \|,
    \end{aligned}
    \label{eq:schauder-change-weights}
  \end{equation}
  uniformly in $i \geqslant - 1$ and $\eta \geqslant \bar{\mu}$. To obtain the
  last estimate we used the fact that
  \begin{equation}
     \| \rho_{\eta} \chi_j \zeta_{\mu_j \vee \bar{\mu}}^{- 1} \|
    \lesssim \| \rho_{\eta} \chi_j \zeta_{\mu_j}^{- 1} \| \lesssim
    \llbracket \eta \rrbracket^{- \gamma} \llbracket \mu_j \vee \bar{\mu}
    \rrbracket^{\gamma} . \label{eq:bound-weight-sch}
  \end{equation}
  Indeed note that if $\bar{\mu} \leqslant \eta \leqslant \mu_j$, we have
  \[ \| \rho_{\eta} \chi_j \zeta_{\mu_j}^{- 1} \| \leqslant (1 + \llbracket
     \eta \rrbracket^a \llbracket \mu_j \rrbracket^{- a})^{- \nu} \lesssim
     \llbracket \eta \rrbracket^{- \gamma} \llbracket \mu_j
     \rrbracket^{\gamma} = \llbracket \eta \rrbracket^{- \gamma} \llbracket
     \mu_j \vee \bar{\mu} \rrbracket^{\gamma} . \]
  On the other hand, if $\eta \geqslant \mu_j \vee \bar{\mu}$, then $\|
  \rho_{\eta} \chi_j \zeta_{\mu_j}^{- 1} \| \leqslant 1$ and $1 \leqslant
  \llbracket \eta \rrbracket^{- \gamma} \llbracket \mu_j \vee \bar{\mu}
  \rrbracket^{\gamma}$. Hence, the bound \eqref{eq:bound-weight-sch} follows.
  By the estimate \eqref{eq:schauder-change-weights} it suffices to establish
  the bound \eqref{eq:schauder} with the weight $\rho_{\eta}$ replaced by
  $\zeta_{\eta}$. To prove the latter bound note that
  \begin{align*}
    \sup_{\eta \geqslant \bar{\mu}}  \llbracket \eta \rrbracket^{\gamma}  \|
    \zeta_{\eta} \bar{\Delta}_i  \phi \| &=  \sup_{\bar{\mu} \leqslant \eta
    \leqslant \mu_i}  \llbracket \eta \rrbracket^{\gamma}  \| \zeta_{\eta}
    \bar{\Delta}_i  \phi \| \vee \sup_{\eta \geqslant \mu_i \vee \bar{\mu}} 
    \llbracket \eta \rrbracket^{\gamma}  \| \zeta_{\eta} \bar{\Delta}_i  \phi
    \|\\
 &\leqslant  \llbracket \bar{\mu} \rrbracket^{\gamma} \| \zeta_{\mu_i
    \vee \bar{\mu}} \bar{\Delta}_i  \phi \| \vee \llbracket \bar{\mu}
    \rrbracket^{\gamma} \llbracket \mu_i \vee \bar{\mu} \rrbracket^{- \gamma}
    \sup_{\eta \geqslant \mu_i \vee \bar{\mu}}  \llbracket \eta
    \rrbracket^{\gamma}  \| \zeta_{\eta} \bar{\Delta}_i  \phi \|\\
 &\leqslant  \llbracket \bar{\mu} \rrbracket^{\gamma} \llbracket \mu_i
    \vee \bar{\mu} \rrbracket^{- \gamma} \sup_{\eta \geqslant \mu_i \vee
    \bar{\mu}}  \llbracket \eta \rrbracket^{\gamma}  \| \zeta_{\eta}
    \bar{\Delta}_i  \phi \| .
  \end{align*}
  In order to bound the supremum in the last line above we shall use the
  following decomposition
  \[ \phi = \phi_{\eta} + (G - G_{\eta}) \mathcal{L} \phi, \qquad \eta
     \geqslant \mu_i \vee \bar{\mu} . \]
  Note that
  \[ \sup_{\eta \geqslant \mu_i \vee \bar{\mu}}  \llbracket \eta
     \rrbracket^{\gamma}  \| \zeta_{\eta} \bar{\Delta}_i  \phi \| \lesssim
     \sup_{\eta \geqslant \mu_i \vee \bar{\mu}}  \llbracket \eta
     \rrbracket^{\gamma}  \| \zeta_{\eta} \bar{\Delta}_i  \phi_{\eta} \| +
     \sup_{\eta \geqslant \mu_i \vee \bar{\mu}} 
     \llbracket \eta \rrbracket^{\gamma}  \|
     \zeta_{\eta} \bar{\Delta}_i  (G - G_{\eta}) \mathcal{L} \phi \|, \]
  and
  \[ \sup_{\eta \geqslant \mu_i \vee \bar{\mu}}  \llbracket \eta
     \rrbracket^{\gamma}  \| \zeta_{\eta} \bar{\Delta}_i  \phi_{\eta} \|
     \lesssim \sup_{\eta \geqslant \mu_i \vee \bar{\mu}}  \llbracket \eta
     \rrbracket^{\gamma}  \| \zeta_{\eta} \phi_{\eta} \| \lesssim \sup_{\eta
     \geqslant \mu_i \vee \bar{\mu}}  \llbracket \eta \rrbracket^{\gamma}  \|
     \zeta^{1 / 3}_{\eta} \phi_{\eta} \| \leqslant \trinorm{ \phi
     } . \]
  To complete the proof it is now enough to show that
  \begin{equation}
    \| \zeta_{\eta} \bar{\Delta}_i  (G - G_{\eta}) \mathcal{L} \phi \|
    \lesssim \llbracket \eta \rrbracket^{2 s - 3 \gamma}  \trinorm{
    K_{\Cdot} \mathcal{L} \phi }_{\#},
    \label{eq:schauder-proof-bound}
  \end{equation}
  uniformly in $\eta \geqslant \mu_i \vee \bar{\mu}$. To this end, recall that
  $L_{\mu} = {(1 + \llbracket \mu \rrbracket^{2 s} \partial_t) }
  \bar{L}_{\mu}$, where $\bar{L}_{\mu} {\assign } (1 - \llbracket
  \mu \rrbracket^2 \Delta)^2$. Now
  \begin{align*}
    \bar{\Delta}_i  (G - G_{\eta}) \mathcal{L} \phi &=  \bar{\Delta}_i  (G - G_{\eta})
    L_{\eta} K_{\eta} \mathcal{L} \phi\\
 &=  \bar{L}_{\eta} \bar{\Delta}_i  (G - G_{\eta}) (1 + \llbracket \eta
    \rrbracket^{2 s} \partial_t) (K_{\eta} \mathcal{L} \phi)\\
 &=  \bar{L}_{\eta} \bar{\Delta}_i  (G - G_{\eta}) K_{\eta} \mathcal{L}
    \phi + \llbracket \eta \rrbracket^{2 s} \bar{L}_{\eta} \bar{\Delta}_i 
    \partial_t G (1 - J_{\eta}) K_{\eta} \mathcal{L} \phi\\
 &=  \mathbb{I}+\mathbb{I}\mathbb{I}.
  \end{align*}
  We have the following bound for the first term
  \begin{align*}
    \| \zeta_{\eta} \mathbb{I} \| = \| \zeta_{\eta} \bar{L}_{\eta}
    \bar{\Delta}_i  (G - G_{\eta}) K_{\eta} \mathcal{L} \phi \| &\lesssim 
    \| \bar{L}_{\eta} \bar{\Delta}_i  \|_{\operatorname{TV} (\zeta^{- 1})}  \| (G -
    G_{\eta}) \|_{\operatorname{TV} (\zeta^{- 1})} \| \zeta_{\eta} K_{\eta}
    \mathcal{L} \phi \|\\
 &\lesssim  \llbracket \eta \rrbracket^{2 s} \| \zeta_{\eta} K_{\eta}
    \mathcal{L} \phi \| .
  \end{align*}
  The above estimates follow from
  Remark~\ref{rem:weight-properties}~a), Lemma~\ref{lem:estimates-for-G} and
  the bound
  \[ \| \bar{L}_{\eta} \bar{\Delta}_i  \|_{\operatorname{TV} (\zeta^{- 1})} \lesssim
     1, \]
  which is a consequence of Lemma~\ref{lem:LP-block}. In order to bound $\|
  \zeta_{\eta} \mathbb{I}\mathbb{I} \|_{L^{\infty}}$ we first note that
  \[ \partial_t G = 1 - (m^2 + (- \Delta_{\varepsilon})^s) G \]
  and write
  \begin{align*}
    \mathbb{I}\mathbb{I} &=  \llbracket \eta \rrbracket^{2 s} \bar{L}_{\eta}
    \bar{\Delta}_i  \partial_t G (1 - J_{\eta}) K_{\eta} \mathcal{L} \phi\\
 &=  \llbracket \eta \rrbracket^{2 s} \bar{L}_{\eta} \bar{\Delta}_i  (1
    - J_{\eta}) K_{\eta} \mathcal{L} \phi - \llbracket \eta \rrbracket^{2 s}
    (m^2 + (- \Delta_{\varepsilon})^s) \bar{L}_{\eta} \bar{\Delta}_i  G (1 -
    J_{\eta}) K_{\eta} \mathcal{L} \phi\\
 &=  \mathbb{I}\mathbb{I}_1 +\mathbb{I}\mathbb{I}_2 .
  \end{align*}
  We have the following bound for $\| \zeta_{\eta} \mathbb{I}\mathbb{I}_1
  \|_{L^{\infty}}$:
  \begin{align*}
    \| \zeta_{\eta} \mathbb{I}\mathbb{I}_1 \|
    &= \llbracket \eta \rrbracket^{2 s}
    \| \zeta_{\eta} \bar{L}_{\eta} \bar{\Delta}_i
    (1 - J_{\eta}) K_{\eta} \mathcal{L} \phi \|\\
    &\lesssim  \llbracket \eta \rrbracket^{2 s}
    \| \bar{L}_{\eta} \bar{\Delta}_i \|_{\operatorname{TV} (\zeta^{- 1})}
    \| (1 - J_{\eta}) \| _{\operatorname{TV} (\zeta^{- 1})}\\
    &\qquad \times
    \| \zeta_{\eta} K_{\eta} \mathcal{L} \phi \|\\
    &\lesssim   \llbracket \eta \rrbracket^{2 s}
    \| \zeta_{\eta} K_{\eta} \mathcal{L} \phi \| .
  \end{align*}
  For $\| \zeta_{\eta}  \mathbb{I}\mathbb{I}_2 \|$ we have, exploiting
  Lemma~\ref{lem:estimates-for-G},
  \begin{align*}
    \| \zeta_{\eta}  \mathbb{I}\mathbb{I}_2 \|
    &\leqslant \llbracket \eta \rrbracket^{2 s}
    \| \zeta_{\eta}  (m^2 + (-\Delta_{\varepsilon})^s)
    \bar{L}_{\eta} (\bar{\Delta}_{i - 1} + \bar{\Delta}_i
    + \bar{\Delta}_{i + 1}) \bar\Delta_i G (1 - J_{\eta})
    K_{\eta} \mathcal{L} \phi \|\\
    &\lesssim \llbracket \eta \rrbracket^{2 s}
    \| (m^2 + (-\Delta_{\varepsilon})^s)\bar{L}_{\eta}
    (\bar{\Delta}_{i - 1} + \bar{\Delta}_i + \bar{\Delta}_{i + 1})
    \|_{\operatorname{TV}(\zeta^{- 1})}\\
    &\quad \times
    \| \bar{\Delta}_i G \|_{\operatorname{TV}(\zeta^{- 1})}
    \| (1 - J_{\eta}) \|_{\operatorname{TV}(\zeta^{- 1})}
    \| \zeta_{\eta} K_{\eta} \mathcal{L} \phi \|\\
    &\lesssim \llbracket \eta \rrbracket^{2 s}
    \llbracket \mu_i \rrbracket^{- 2 s}
    \llbracket \mu_i \rrbracket^{2 s}
    \| \zeta_{\eta} K_{\eta} \mathcal{L} \phi \|\\
    &= \llbracket \eta \rrbracket^{2 s}
    \| \zeta_{\eta} K_{\eta} \mathcal{L} \phi \|,
  \end{align*}
  where we used that
  \begin{align*}
    &\|(m^2 + (- \Delta_{\varepsilon})^s) 
     \bar{L}_{\eta} (\bar{\Delta}_{i - 1}  + \bar{\Delta}_i  + \bar{\Delta}_{i
     + 1}) \|_{\operatorname{TV} (\zeta^{- 1})}\lesssim \llbracket \mu_i \rrbracket^{- 2 s},
  \end{align*}
  which follows from Lemma~\ref{lem:LP-block}. Summing up, we have
  \begin{align*}
    \left\| {\zeta_{\eta} } \bar\Delta_i  (G - G_{\eta}) \mathcal{L} \phi \right\|
     &\lesssim \llbracket \eta \rrbracket^{2 s}  \left\| {\zeta_{\eta} } 
     K_{\eta} \mathcal{L} \phi \right\|\\
     &\lesssim \llbracket \eta \rrbracket^{2
     s - 3 \gamma}  \llbracket \eta \rrbracket^{3 \gamma}  \left\|
     {\zeta_{\eta} }  K_{\eta} \mathcal{L} \phi \right\|\\
     &\lesssim \llbracket
     \eta \rrbracket^{2 s - 3 \gamma}  \trinorm{ K_{\Cdot} \mathcal{L}
     \phi }_{\#} .
  \end{align*}
  This proves the bound~\eqref{eq:schauder-proof-bound} and completes the
  proof of the lemma.
\end{proof}

\begin{lemma}
  \label{lem:schauder-two}Let $J_{> \eta} \assign 1 - J_{\eta}$. We have
  \[ \| \rho_{\mu}^3 J_{> \eta} \phi_{\sigma} \| \lesssim \llbracket \eta
     \rrbracket^{2 s}  \llbracket \sigma \rrbracket^{- 3 \gamma}  \trinorm{
     \mathcal{L} \phi_{\Cdot} }_{\#}, \]
  uniformly in $\bar{\mu}, \mu, \eta \in [1 / 2, 1)$, $\sigma \in [\mu \vee
  \bar{\mu}, 1)$ and $\phi \in \mathcal{S}' (\Lambda)$, where the norm
  $\trinorm{ \Cdot }_{\#} = \trinorm{ \Cdot
  }_{\#, \bar{\mu}}$ is defined by~\eqref{eq:sharp-norm}.
\end{lemma}

\begin{proof}
  Using Remark~\ref{rem:weight-properties}~a) and
  Lemma~\ref{lem:estimates-for-G}, we obtain
  \[ \| \rho_{\mu}^3  (G - G_{\eta}) \mathcal{L} \phi_{\sigma} \| \lesssim \|
     G - G_{\eta} \|_{\operatorname{TV} (\zeta^{- 1})}  \| \rho_{\mu}^3 \mathcal{L}
     \phi_{\sigma} \| \lesssim \llbracket \eta \rrbracket^{2 s}  \|
     \rho_{\mu}^3 \mathcal{L} \phi_{\sigma} \| . \]
  Next, we note that
  \begin{align*}
     \| \rho_{\mu}^3 J_{> \eta} \phi_{\sigma} \|
     &= \| \rho_{\mu}^3  (G - G_{\eta})
     \mathcal{L} \phi_{\sigma} \|
     \lesssim \llbracket \eta \rrbracket^{2 s}
     \| \rho_{\mu}^3 \mathcal{L} \phi_{\sigma} \|\\
     &\lesssim
     \llbracket \eta \rrbracket^{2 s} \| \rho_{\sigma}^3
     \mathcal{L} \phi_{\sigma} \|
     \lesssim \llbracket \eta \rrbracket^{2 s}
     \llbracket \sigma \rrbracket^{- 3 \gamma}
     \trinorm{ \mathcal{L} \phi_{\Cdot} }_{\#},
  \end{align*}
  on account of Lemma~\ref{lem:triple-phi} and the identity $\rho_{\sigma} =
  \zeta_{\sigma}^{\nu}$. This finishes the proof.
\end{proof}

\section{Flow equation estimates}\label{app:flow-eq}

\subsection{Estimates for \texorpdfstring{$\mathcal{A}$ and
$\mathcal{B}$}{A and B}}\label{sec:A-B-cumulants}

The goal of this section is to prove the bounds for the operators
{\textsf{$\mathcal{A}^{\boldsymbol{a}}_{\boldsymbol{b}}$}} and
$\mathcal{B}^{\boldsymbol{a}}_{\boldsymbol{b}, \boldsymbol{c}}$ appearing in the
flow equation~\eqref{eq:cumulant-flow} for cumulants, thereby proving
Lemma~\ref{lem:flow-cumulants}. We begin with the definition of the first
operator:
\[ \sum_{\boldsymbol{b}} \mathcal{A}^{\boldsymbol{a}}_{\boldsymbol{b}}
   (\dot{G}_{\sigma}, \mathcal{F}^{\boldsymbol{b}}_{\sigma}) \assign \sum_{i =
   1}^{n (\boldsymbol{a})} \sum_{\boldsymbol{b}} \mathcal{A}^{\boldsymbol{a},
   (i)}_{\boldsymbol{b}} (\dot{G}_{\sigma},
   \mathcal{F}^{\boldsymbol{b}}_{\sigma}), \]
where

\begin{equation}
  \begin{aligned}
    \sum_{\boldsymbol{b}} \mathcal{A}^{\boldsymbol{a}, (i)}_{\boldsymbol{b}}
    (\dot{G}_{\sigma}, \mathcal{F}^{\boldsymbol{b}}_{\sigma})
    & \assign
    \sum_{\ell' = 0}^{\ell (\mathfrak{a}_i) - 1}
    \sum^{k (\mathfrak{a}_i)}_{k' = 0} (k' + 1)
    \mathfrak{C}_i (\dot{G}_{\sigma})\\
    &\qquad \times \mathfrak{K}_{n(\boldsymbol{a}) + 1}
    (F_{\sigma}^{\mathfrak{a}_1}, \ldots,
    F_{\sigma}^{\mathfrak{a}_{i - 1}},
    F_{\sigma}^{[\ell (\mathfrak{a}_i) - 1 - \ell'], (k' + 1)},\\
    &\qquad\qquad
    F_{\sigma}^{[\ell'], (k (\mathfrak{a}_i) - k')},
    F_{\sigma}^{\mathfrak{a}_{i + 1}}, \ldots,
    F_{\sigma}^{\mathfrak{a}_{n(\boldsymbol{a})}}) .
  \end{aligned} \label{eq:specific-operator-A}
\end{equation}
The summation index $\boldsymbol{b}$ is constrained by the allowed values of
$(\ell', k')$ appearing on the right-hand side. 
For fixed $\boldsymbol a = (\mathfrak a_1,\ldots,\mathfrak a_n)
\in \boldsymbol A$ and $i \in \{1,\ldots,n\}$, we say that
$\boldsymbol b$ is \emph{$(\boldsymbol a,i)$-admissible for
	$\mathcal A$} if there exist
\[
0 \leqslant \ell' \leqslant \ell(\mathfrak a_i)-1,
\qquad
0 \leqslant k' \leqslant k(\mathfrak a_i),
\]
such that
\[
\boldsymbol b =
(\mathfrak a_1,\ldots,\mathfrak a_{i-1},
(\ell(\mathfrak a_i)-1-\ell',k'+1),
(\ell',k(\mathfrak a_i)-k'),
\mathfrak a_{i+1},\ldots,\mathfrak a_n).
\]
If $\boldsymbol b$ is not of this form, we set
$\mathcal A^{\boldsymbol a,(i)}_{\boldsymbol b}=0$.
The operator
$\mathfrak{C}_i (\dot{G}_{\sigma})$ acts by applying
$\dot{G}_{\sigma}$ to the output variable of the $(i + 1)$-th kernel, and by
inserting the resulting function into the last input variable of the $i$-th
kernel. This operator generalises $\mathfrak{C}
(\dot{G}_{\sigma})$ defined in~\eqref{eq:frak-B-operator}, which appears in
the flow equation for the effective force kernels~\eqref{eq:flow-kernels}.

We now proceed with the definition of the second operator:
\begin{align*}
  \sum_{\boldsymbol{b}, \boldsymbol{c}}
  \mathcal{B}^{\boldsymbol{a}}_{\boldsymbol{b}, \boldsymbol{c}}
  (\dot{G}_{\sigma}, \mathcal{F}^{\boldsymbol{b}}_{\sigma},
  \mathcal{F}^{\boldsymbol{c}}_{\sigma})\assign \sum_{i = 1}^{n(\boldsymbol{a})}
  \sum_{\boldsymbol{b}, \boldsymbol{c}}
  \mathcal{B}^{\boldsymbol{a}, (i)}_{\boldsymbol{b}, \boldsymbol{c}}
  (\dot{G}_{\sigma}, \mathcal{F}^{\boldsymbol{b}}_{\sigma},
  \mathcal{F}^{\boldsymbol{c}}_{\sigma}),
\end{align*}
where
\begin{equation}
  \begin{aligned}
    \sum_{\boldsymbol{b}, \boldsymbol{c}}
    \mathcal{B}^{\boldsymbol{a}, (i)}_{\boldsymbol{b}, \boldsymbol{c}}
    (\dot{G}_{\sigma}, \mathcal{F}^{\boldsymbol{b}}_{\sigma},
    \mathcal{F}^{\boldsymbol{c}}_{\sigma})
    & \assign \sum_{I_1, I_2}
    \sum_{\ell' = 0}^{\ell(\mathfrak{a}_i) - 1}
    \sum^{k(\mathfrak{a}_i)}_{k' = 0} (k' + 1)
    \mathfrak{C}_{| I_1 | + 1}(\dot{G}_{\sigma})\\
    &\qquad \bigg(
    \mathfrak{K}_{| I_1 | + 1}
    \Big((F_{\sigma}^{\mathfrak{a}_j})_{j \in I_1},
    F^{[\ell(\mathfrak{a}_i) - 1 - \ell'],(k' + 1)}_{\sigma}\Big)\\
    &\qquad\qquad 
    \mathfrak{K}_{| I_2 | + 1}
    \Big(F_{\sigma}^{[\ell'],(k(\mathfrak{a}_i) - k')},
    (F_{\sigma}^{\mathfrak{a}_j})_{j \in I_2}\Big) \bigg).
  \end{aligned}
  \label{eq:specific-operator-B}
\end{equation}
Here, for a fixed $i \in \{ 1, \ldots, n (\boldsymbol{a}) \}$, the sum
$\sum_{I_1, I_2}$ runs over all the partitions of the set
\[ I_1 \sqcup I_2 = \{ 1, \ldots, i - 1, i + 1, \ldots, n (\boldsymbol{a}) \} .
\]
The summation indices $\boldsymbol{b}$ and $\boldsymbol{c}$ are now constrained by the allowed values of
$(I_1, I_2, \ell', k')$ on the right-hand side of
\eqref{eq:specific-operator-B}. 
We say that $(\boldsymbol b,\boldsymbol c)$ is
\emph{$(\boldsymbol a,i)$-admissible for $\mathcal B$} if there exist a
partition
\[
I_1\sqcup I_2=\{1,\ldots,i-1,i+1,\ldots,n\}
\]
and integers $\ell',k'$ in the above ranges such that, with the elements of
$I_1$ and $I_2$ ordered increasingly,
\[
\boldsymbol b =
\bigl((\mathfrak a_j)_{j\in I_1},
(\ell(\mathfrak a_i)-1-\ell',k'+1)\bigr),
\]
and
\[
\boldsymbol c =
\bigl((\ell',k(\mathfrak a_i)-k'),
(\mathfrak a_j)_{j\in I_2}\bigr).
\]
If $(\boldsymbol b,\boldsymbol c)$ is not of this form, we set
$\mathcal B^{\boldsymbol a,(i)}_{\boldsymbol b,\boldsymbol c}=0$.
The operator
$\mathfrak{C}_{| I_1 | + 1} (\dot{G}_{\sigma})$ acts by applying
$\dot{G}_{\sigma}$ to the output variable of the kernel $F_{\sigma}^{[\ell'],(k(\mathfrak a_i)-k')}$, and by
inserting the resulting function into the last input variable of the
kernel $F_{\sigma}^{[\ell(\mathfrak a_i)-1-\ell'],(k'+1)}$. For precise definitions of the operators
{\textsf{$\mathcal{A}^{\boldsymbol{a}}_{\boldsymbol{b}}$}} and
$\mathcal{B}^{\boldsymbol{a}}_{\boldsymbol{b}, \boldsymbol{c}}$,
see~{\cite{duchFlowEquationApproach2021}}.

\begin{remark}[Periodicity and translation invariance]
	\label{rem:periodicity-translation-AB}
	The estimates for the operators $\mathcal A$ and $\mathcal B$ are applied
	to cumulant kernels which are spatially periodic and invariant under
	simultaneous spacetime translations. These properties are preserved by
	the regularising operators, since their kernels are translation invariant,
	and the weights $w_\sigma^{\boldsymbol a}$ are likewise invariant under
	simultaneous translations. 
	Although the norm in Definition~\ref{def:norm-cumulants} distinguishes the
	first output variable, this choice is immaterial for the cumulants appearing
	below. Indeed, let $F$ be simultaneously translation invariant in
	all its output and input variables and spatially periodic with period $M$.
	Then, for every $r\in\{1,\ldots,n\}$,
	\[
	\sup_{z_r\in\Lambda}
	\int_{\Lambda_M^{n-1}}
	\left(\int_{\Lambda^k}
	\bigl|F
	(z_1,\mathd Z_1,\ldots,z_n,\mathd Z_n)\bigr|
	\right)
	\prod_{q\neq r}\mathd z_q
	=
	\trinorm{F}.
	\]
	To see this, translate all variables simultaneously so that the fixed
	variable $z_r$ is moved to the position of the distinguished output.
	Lebesgue measure is translation invariant, while a translated spatial
	fundamental domain may be replaced by $\Lambda_M$ by periodicity.
	
	In the estimate for $\mathcal A$, spatial
	periodicity is used to decompose the integration over the additional
	output variable into translates of the fundamental domain $\Lambda_M$.
	After translating the corresponding output--input block by an element of
	the spatial period lattice, the cumulant kernel and its weight are
	unchanged, while the sum over lattice translates is absorbed into the
	periodisation
	\[
	\bigl|
	(L_\sigma\dot G_\sigma)w_\sigma^\flat
	\bigr|_M.
	\]
	This yields the factor
	\[
	\left\|
	\bigl|
	(L_\sigma\dot G_\sigma)w_\sigma^\flat
	\bigr|_M
	\right\|_{L^\infty}
	\]
	in Lemma~\ref{lem:estimates-A}.
	Translation invariance also shows that the operators $\mathcal A$ and
	$\mathcal B$ preserve the class of cumulants considered here. Moreover,
	in the estimate for $\mathcal B$, it allows one to choose any output
	variable as the distinguished variable in the cumulant norm, while
	periodicity permits the corresponding translated fundamental domain to
	be replaced by $\Lambda_M$.
\end{remark}

\begin{lemma}
	\label{lem:estimates-A}
	Let $\boldsymbol a \in \boldsymbol A$ and
	$i\leqslant n(\boldsymbol a)$. For every
	$(\boldsymbol a,i)$-admissible $\boldsymbol b$ one has
	\[
	\llbracket\sigma\rrbracket^{-[\boldsymbol a]}
	\trinorm{
		\mathcal A^{\boldsymbol a,(i)}_{\boldsymbol b}
		(\dot G_\sigma,\mathcal F^{\boldsymbol b}_\sigma)
	}_\sigma
	\lesssim
	\llbracket\sigma\rrbracket^{-[\boldsymbol b]-1}
	\trinorm{\mathcal F^{\boldsymbol b}_\sigma}_\sigma ,
	\]
	uniformly in $\sigma\in[1/2,1)$ and in
	$\mathcal F^{\boldsymbol b}_\sigma$.
\end{lemma}

\begin{proof}
  Recalling~\eqref{eq:specific-operator-A}, our task is to analyse the
  expression
  \[  \mathfrak{C}_i (\dot{G}_{\sigma}) \mathfrak{K}_{n
     (\boldsymbol{a}) + 1} (F_{\sigma}^{\mathfrak{a}_1}, \ldots,
     F_{\sigma}^{[\ell (\mathfrak{a}_i) - 1 - \ell'], (k' + 1)},
     F_{\sigma}^{[\ell'], (k (\mathfrak{a}_i) - k')}, \ldots,
     F_{\sigma}^{\mathfrak{a}_{n (\boldsymbol{a})}}) . \]
  Using $L_{\sigma} K_{\sigma} = 1$ and
  Definition~\ref{def:regularizing-kernels-for-kernels}, we obtain
  \begin{align*}
 &K_{\sigma}^{n (\boldsymbol{a}), K (\boldsymbol{a})}
    \bigl(\mathfrak{C}_i (\dot{G}_{\sigma})
    \mathfrak{K}_{n(\boldsymbol{a}) + 1} (F_{\sigma}^{\mathfrak{a}_1}, \ldots,
    F_{\sigma}^{[\ell (\mathfrak{a}_i) - 1 - \ell'], (k' + 1)},
    F_{\sigma}^{[\ell'], (k (\mathfrak{a}_i) - k')}, \ldots,
    F_{\sigma}^{\mathfrak{a}_{n (\boldsymbol{a})}})\bigr)\\
  &= \mathfrak{C}_i (L_{\sigma} \dot{G}_{\sigma})
    K_{\sigma}^{n (\boldsymbol{a}) + 1, K (\boldsymbol{a}) + 1} \mathfrak{K}_{n(\boldsymbol{a}) + 1}
    (F_{\sigma}^{\mathfrak{a}_1}, \ldots
    \\
    &\hspace{4cm}\ldots,
    F_{\sigma}^{[\ell (\mathfrak{a}_i) - 1 - \ell'], (k' + 1)},
    F_{\sigma}^{[\ell'], (k (\mathfrak{a}_i) - k')}, \ldots,
    F_{\sigma}^{\mathfrak{a}_{n (\boldsymbol{a})}}) .
  \end{align*}
  As a consequence,
  \[ K_{\sigma}^{\boldsymbol{a}} \mathcal{A}^{\boldsymbol{a}, (i)}_{\boldsymbol{b}}
     (\dot{G}_{\sigma}, \mathcal{F}^{\boldsymbol{b}}_{\sigma})
     =\mathcal{A}^{\boldsymbol{a}, (i)}_{\boldsymbol{b}} (L_{\sigma}
     \dot{G}_{\sigma}, K_{\sigma}^{\boldsymbol{b}}
     \mathcal{F}^{\boldsymbol{b}}_{\sigma}) . \]
  Hence, by Definition~\ref{def:norm-cumulants} of the cumulant
  norm,
  \[ \begin{aligned}
       \trinorm{ \mathcal{A}^{\boldsymbol{a}, (i)}_{\boldsymbol{b}}
       (\dot{G}_{\sigma}, \mathcal{F}^{\boldsymbol{b}}_{\sigma})
       }_{\sigma} & =  \trinorm{ [K_{\sigma}^{\boldsymbol{a}}
       \mathcal{A}^{\boldsymbol{a}, (i)}_{\boldsymbol{b}} (\dot{G}_{\sigma},
       \mathcal{F}^{\boldsymbol{b}}_{\sigma})] \cdot w_{\sigma}^{\boldsymbol{a}}
       }\\
       & =  \trinorm{ \mathcal{A}^{\boldsymbol{a}, (i)}_{\boldsymbol{b}}
       (L_{\sigma} \dot{G}_{\sigma}, K_{\sigma}^{\boldsymbol{b}}
       \mathcal{F}^{\boldsymbol{b}}_{\sigma}) \cdot w_{\sigma}^{\boldsymbol{a}}
       }\\
       & \leqslant  \trinorm{ \mathcal{A}^{\boldsymbol{a},
       (i)}_{\boldsymbol{b}} (| L_{\sigma} \dot{G}_{\sigma} |, |
       K_{\sigma}^{\boldsymbol{b}} \mathcal{F}^{\boldsymbol{b}}_{\sigma} |) \cdot
       w_{\sigma}^{\boldsymbol{a}} }\\
       & \lesssim  \trinorm{ \mathcal{A}^{\boldsymbol{a},
       (i)}_{\boldsymbol{b}} (| L_{\sigma} \dot{G}_{\sigma} | \cdot
       w_{\sigma}^{\flat}, | K_{\sigma}^{\boldsymbol{b}}
       \mathcal{F}^{\boldsymbol{b}}_{\sigma} | \cdot w_{\sigma}^{\boldsymbol{b}})
       },
     \end{aligned} \]
  where the last estimate follows from~\eqref{eq:weight-w} and $| L_{\sigma}
  \dot{G}_{\sigma} | \cdot w_{\sigma}^{\flat}$ denotes the function
  \begin{equation}
    z \mapsto | (L_{\sigma} \dot{G}_{\sigma}) (z) | w_{\sigma}^{\flat} (z) .
    \label{eq:A-proof-dot-G-w}
  \end{equation}
  As a result, using Definition~\ref{def:norm-cumulants} one proves along the lines
  of Lemma~5.29 in~{\cite{duchFlowEquationApproach2021}} that
  \begin{align*}
     &\trinorm{ \mathcal{A}^{\boldsymbol{a}, (i)}_{\boldsymbol{b}}
     (\dot{G}_{\sigma}, \mathcal{F}^{\boldsymbol{b}}_{\sigma})
     }_{\sigma}\\
     &\quad \lesssim \left\| \big| (L_{\sigma} \dot{G}_{\sigma})
     \cdot w_{\sigma}^{\flat} \big|_M  \right\|_{L^{\infty}}
     \trinorm{ [K_{\sigma}^{\boldsymbol{b}}
     \mathcal{F}^{\boldsymbol{b}}_{\sigma}] \cdot w_{\sigma}^{\boldsymbol{b}}
     }\\
     &\quad = \left\| \big| (L_{\sigma} \dot{G}_{\sigma})
     \cdot w_{\sigma}^{\flat} \big|_M  \right\|_{L^{\infty}}
     \trinorm{ \mathcal{F}^{\boldsymbol{b}}_{\sigma}
     }_{\sigma},
  \end{align*}
  where $| (L_{\sigma} \dot{G }_{\sigma} ) \cdot w_{\sigma}^{\flat} | _M$
  denotes the periodisation in space of the
  function~\eqref{eq:A-proof-dot-G-w} with period $M$. Thus, since by
  Lemma~\ref{lem:Gdot-decay}
  \[ \left\| \big| (L_{\sigma} \dot{G }_{\sigma} ) \cdot w_{\sigma}^{\flat}
     \big| _M  \right\|_{L^{\infty}} \lesssim \llbracket \sigma
     \rrbracket^{- 1 - d}, \]
  we obtain
  \[ \trinorm{ \mathcal{A}^{\boldsymbol{a}, (i)}_{\boldsymbol{b}}
     (\dot{G}_{\sigma}, \mathcal{F}^{\boldsymbol{b}}_{\sigma})
     }_{\sigma} \lesssim \llbracket \sigma \rrbracket^{- 1 - d} 
     \trinorm{ \mathcal{F}^{\boldsymbol{b}}_{\sigma} }_{\sigma} .
  \]
  
  It follows that
  \[ \llbracket \sigma \rrbracket^{- [\boldsymbol{a}]}  \trinorm{
     \mathcal{A}^{\boldsymbol{a}, (i)}_{\boldsymbol{b}} (\dot{G}_{\sigma},
     \mathcal{F}^{\boldsymbol{b}}_{\sigma}) }_{\sigma} \lesssim
     \llbracket \sigma \rrbracket^{- [\boldsymbol{a}] - 1 - d}  \trinorm{
     \mathcal{F}^{\boldsymbol{b}}_{\sigma} }_{\sigma} \lesssim
     \llbracket \sigma \rrbracket^{- [\boldsymbol{b}] - 1}  \trinorm{
     \mathcal{F}^{\boldsymbol{b}}_{\sigma} }_{\sigma} . \]
  In the last step, we used~\eqref{eq:homogeneity-cum}
  and~\eqref{eq:compatibility}, which imply
  \[ - [\boldsymbol{a}] - 1 - d = - [\boldsymbol{b}] - 1 + \theta + \beta - \delta
     - d  \geqslant - [\boldsymbol{b}] - 1. \]
  This finishes the proof.
\end{proof}

\begin{lemma}
	\label{lem:estimates-B}
	Let $\boldsymbol a \in \boldsymbol A$ and
	$i\leqslant n(\boldsymbol a)$. For every
	$(\boldsymbol a,i)$-admissible pair $(\boldsymbol b,\boldsymbol c)$ one has
	\[
	\llbracket\sigma\rrbracket^{-[\boldsymbol a]}
	\trinorm{
		\mathcal B^{\boldsymbol a,(i)}_{\boldsymbol b,\boldsymbol c}
		(\dot G_\sigma,
		\mathcal F^{\boldsymbol b}_\sigma,
		\mathcal F^{\boldsymbol c}_\sigma)
	}_\sigma
	\lesssim
	\llbracket\sigma\rrbracket^{-[\boldsymbol b]-[\boldsymbol c]-1}
	\trinorm{\mathcal F^{\boldsymbol b}_\sigma}_\sigma
	\trinorm{\mathcal F^{\boldsymbol c}_\sigma}_\sigma ,
	\]
	uniformly in $\sigma\in[1/2,1)$ and in
	$\mathcal F^{\boldsymbol b}_\sigma,\mathcal F^{\boldsymbol c}_\sigma$.
\end{lemma}

\begin{proof}
  The proof proceeds along the same lines as the previous one.
  Recalling~\eqref{eq:specific-operator-B}, we now need to analyse the
  expression
  \[ \mathfrak{C}_{| I_1 | + 1} (\dot{G}_{\sigma})  (\mathfrak{K}_{|
     I_1 | + 1} ((F_{\sigma}^{\mathfrak{a}_j})_{j \in I_1}, F_{\sigma}^{[\ell(\mathfrak a_i)-1-\ell'],(k'+1)}) \mathfrak{K}_{| I_2 | +
     1} (F_{\sigma}^{[\ell'],(k(\mathfrak a_i)-k')},
     (F_{\sigma}^{\mathfrak{a}_j})_{j \in I_2})) . \]
  Using $L_{\sigma} K_{\sigma} = 1$ and
  Definition~\ref{def:regularizing-kernels-for-kernels}, we obtain
  \[
  \begin{aligned}
  	&K_{\sigma}^{\boldsymbol{a}}\Bigl(
  	\mathfrak{C}_{|I_1|+1}(\dot{G}_{\sigma})
  	\mathfrak{K}_{|I_1|+1}\bigl((F_{\sigma}^{\mathfrak{a}_j})_{j\in I_1},
  F_{\sigma}^{[\ell(\mathfrak a_i)-1-\ell'],(k'+1)}\bigr)
  	\mathfrak{K}_{|I_2|+1}\bigl(F_{\sigma}^{[\ell'],(k(\mathfrak a_i)-k')},
  	(F_{\sigma}^{\mathfrak{a}_j})_{j\in I_2}\bigr)\Bigr) \\
  	&\qquad =
  	\mathfrak{C}_{|I_1|+1}(L_{\sigma}\dot{G}_{\sigma})
  	\Bigl(
  	K_{\sigma}^{|I_1|+1,k_1}
  	\mathfrak{K}_{|I_1|+1}\bigl((F_{\sigma}^{\mathfrak{a}_j})_{j\in I_1},
  	F_{\sigma}^{[\ell(\mathfrak a_i)-1-\ell'],(k'+1)}\bigr)
  	\\
  	&\hspace{3.5cm}
  	K_{\sigma}^{|I_2|+1,k_2}
  	\mathfrak{K}_{|I_2|+1}\bigl(F_{\sigma}^{[\ell'],(k(\mathfrak a_i)-k')},
  	(F_{\sigma}^{\mathfrak{a}_j})_{j\in I_2}\bigr)
  	\Bigr),
  \end{aligned}
  \]
  where 
  \[
  k_1 \assign \sum_{j\in I_1} k(\mathfrak a_j)+k'+1,
  \qquad
  k_2 \assign \sum_{j\in I_2} k(\mathfrak a_j)
  +k(\mathfrak a_i)-k'.
  \]
  As a consequence,
  \[ K_{\sigma}^{\boldsymbol{a}} \mathcal{B}^{\boldsymbol{a}, (i)}_{\boldsymbol{b},
     \boldsymbol{c}} (\dot{G}_{\sigma}, \mathcal{F}^{\boldsymbol{b}}_{\sigma},
     \mathcal{F}^{\boldsymbol{c}}_{\sigma}) =\mathcal{B}^{\boldsymbol{a},
     (i)}_{\boldsymbol{b}, \boldsymbol{c}} (L_{\sigma} \dot{G}_{\sigma},
     K_{\sigma}^{\boldsymbol{b}} \mathcal{F}^{\boldsymbol{b}}_{\sigma},
     K_{\sigma}^{\boldsymbol{c}} \mathcal{F}^{\boldsymbol{c}}_{\sigma}) . \]
  Hence, by Definition~\ref{def:norm-cumulants} of the cumulant norm,
  \[ \begin{aligned}
       \trinorm{ \mathcal{B}^{\boldsymbol{a}, (i)}_{\boldsymbol{b},
       \boldsymbol{c}} (\dot{G}_{\sigma}, \mathcal{F}^{\boldsymbol{b}}_{\sigma},
       \mathcal{F}^{\boldsymbol{c}}_{\sigma}) }_\sigma & =
       \trinorm{ [K_{\sigma}^{\boldsymbol{a}} \mathcal{B}^{\boldsymbol{a},
       (i)}_{\boldsymbol{b}, \boldsymbol{c}} (\dot{G}_{\sigma},
       \mathcal{F}^{\boldsymbol{b}}_{\sigma},
       \mathcal{F}^{\boldsymbol{c}}_{\sigma})] \cdot w_{\sigma}^{\boldsymbol{a}}
       }\\
       & =  \trinorm{ \mathcal{B}^{\boldsymbol{a}, (i)}_{\boldsymbol{b},
       \boldsymbol{c}} (L_{\sigma} \dot{G}_{\sigma}, K_{\sigma}^{\boldsymbol{b}}
       \mathcal{F}^{\boldsymbol{b}}_{\sigma}, K_{\sigma}^{\boldsymbol{c}}
       \mathcal{F}^{\boldsymbol{c}}_{\sigma}) \cdot w_{\sigma}^{\boldsymbol{a}}
       }\\
       & \leqslant  \trinorm{ \mathcal{B}^{\boldsymbol{a},
       (i)}_{\boldsymbol{b}, \boldsymbol{c}} (| L_{\sigma} \dot{G}_{\sigma} |, |
       K_{\sigma}^{\boldsymbol{b}} \mathcal{F}^{\boldsymbol{b}}_{\sigma} |, |
       K_{\sigma}^{\boldsymbol{c}} \mathcal{F}^{\boldsymbol{c}}_{\sigma} |) \cdot
       w_{\sigma}^{\boldsymbol{a}} }\\
       & \lesssim  \trinorm{ \mathcal{B}^{\boldsymbol{a},
       (i)}_{\boldsymbol{b}, \boldsymbol{c}} (| L_{\sigma} \dot{G}_{\sigma} |
       \cdot w_{\sigma}^{\flat},| K_{\sigma}^{\boldsymbol{b}}
       \mathcal{F}^{\boldsymbol{b}}_{\sigma} | \cdot w_{\sigma}^{\boldsymbol{b}},
       | K_{\sigma}^{\boldsymbol{c}} \mathcal{F}^{\boldsymbol{c}}_{\sigma} | \cdot
       w_{\sigma}^{\boldsymbol{c}}) },
     \end{aligned} \]
  where the last estimate follows from~\eqref{eq:weight-w}. As a result, by
  Definition~\ref{def:norm-cumulants} and Lemma~\ref{lem:Gdot-decay}, we have
  \begin{align*}
    \trinorm{ \mathcal{B}^{\boldsymbol{a}, (i)}_{\boldsymbol{b}, \boldsymbol{c}}
    (\dot{G}_{\sigma}, \mathcal{F}^{\boldsymbol{b}}_{\sigma},
    \mathcal{F}^{\boldsymbol{c}}_{\sigma}) }_{\sigma} &\lesssim  \|
    L_{\sigma} \dot{G}_{\sigma} \|_{\operatorname{TV} (w^{\flat}_{\sigma})} 
    \trinorm{ K_{\sigma}^{\boldsymbol{b}} \mathcal{F}^{\boldsymbol{b}}_{\sigma}
    \cdot w_{\sigma}^{\boldsymbol{b}} }  \trinorm{
    K_{\sigma}^{\boldsymbol{c}} \mathcal{F}^{\boldsymbol{c}}_{\sigma} \cdot
    w_{\sigma}^{\boldsymbol{c}} }\\
 &=  \| L_{\sigma} \dot{G}_{\sigma} \|_{\operatorname{TV} (w^{\flat}_{\sigma})} 
    \trinorm{ \mathcal{F}^{\boldsymbol{b}}_{\sigma} }_{\sigma} 
    \trinorm{ \mathcal{F}^{\boldsymbol{c}}_{\sigma} }_{\sigma}\\
 &\lesssim  \llbracket \sigma \rrbracket^{2 s - 1} \trinorm{
    \mathcal{F}^{\boldsymbol{b}}_{\sigma} }_{\sigma}  \trinorm{
    \mathcal{F}^{\boldsymbol{c}}_{\sigma} }_{\sigma} .
  \end{align*}
  Consequently, we arrive at
  \[ \begin{aligned}
       \llbracket \sigma \rrbracket^{- [\boldsymbol{a}]}  \trinorm{
       \mathcal{B}^{\boldsymbol{a}, (i)}_{\boldsymbol{b}, \boldsymbol{c}}
       (\dot{G}_{\sigma}, \mathcal{F}^{\boldsymbol{b}}_{\sigma},
       \mathcal{F}^{\boldsymbol{c}}_{\sigma}) }_{\sigma} & \lesssim
       \llbracket \sigma \rrbracket^{- [\boldsymbol{a}]}  \llbracket \sigma
       \rrbracket^{2 s - 1}  \trinorm{ \mathcal{F}^{\boldsymbol{b}}_{\sigma}
       }_{\sigma}  \trinorm{ \mathcal{F}^{\boldsymbol{c}}_{\sigma}
       }_{\sigma},\\
       & \leqslant  \llbracket \sigma \rrbracket^{- [\boldsymbol{b}] - [\boldsymbol{c}]
       - 1}  \trinorm{ \mathcal{F}^{\boldsymbol{b}}_{\sigma}
       }_{\sigma}  \trinorm{ \mathcal{F}^{\boldsymbol{c}}_{\sigma}
       }_{\sigma} .
     \end{aligned} \]
  In the last step we used~\eqref{eq:homogeneity-cum}
  and~\eqref{eq:compatibility}, which imply
  \[ - [\boldsymbol{a}] + 2 s - 1 = - [\boldsymbol{b}] - [\boldsymbol{c}] - \varrho
     + \theta + \beta - \delta + 2 s - 1 \geqslant - [\boldsymbol{b}] -
     [\boldsymbol{c}] - 1. \]
  This finishes the proof.
\end{proof}

\begin{remark}
  \label{rem:A-B-epsilon-bounds}
  The unweighted $L^\infty$ estimate used in the proof of
  Lemma~\ref{lem:estimates-A} admits an additional sublattice gain, whereas
  the corresponding total variation estimate used in the proof of
  Lemma~\ref{lem:estimates-B} does not. More precisely, uniformly in
  $\sigma\in[1/2,1)$ and in every admissible spatial period $M\geqslant1$,
  \begin{equation}
    \left\|\left|L_\sigma\dot G_\sigma\right|_M\right\|_{L^\infty}
    \lesssim
    \llbracket\sigma\rrbracket^{-d-1}
    \wedge
    \varepsilon^{-\flat}
    \llbracket\sigma\rrbracket^{\flat-d-1}.
    \label{eq:Gdot-epsilon-Linfty}
  \end{equation}
  Indeed, \eqref{eq:Gdot-corrected-periodised-unweighted} with \(p=1\)
  gives
  \[
    \left\|\left|L_\sigma\dot G_\sigma\right|_M\right\|_{L^\infty}
    \lesssim
    \llbracket\sigma\rrbracket^{-d-1}
    \wedge
    \varepsilon^{-d}\llbracket\sigma\rrbracket^{-1}.
  \]
  Since $d=3$ and $\flat<2s\leqslant2$ by~\eqref{eq:choice-of-flat}, the second
  term is bounded by
  $\varepsilon^{-\flat}\llbracket\sigma\rrbracket^{\flat-d-1}$ when
  $\llbracket\sigma\rrbracket\leqslant\varepsilon$. This proves
  \eqref{eq:Gdot-epsilon-Linfty}.

  The first part of the following estimate follows from
  Lemma~\ref{lem:estimates-A}, since the weights are bounded from below by
  one. Repeating its proof without an output weight and using
  \eqref{eq:Gdot-epsilon-Linfty} gives the second part. Thus, for every
  $(\boldsymbol a,i)$-admissible $\boldsymbol b$,
  \begin{equation}
    \begin{aligned}
      &\llbracket\sigma\rrbracket^{-[\boldsymbol a]}
      \trinorm{
        K^{\boldsymbol a}_\sigma
        \mathcal A^{\boldsymbol a,(i)}_{\boldsymbol b}
        (\dot G_\sigma,\mathcal F^{\boldsymbol b}_\sigma)
      }
      \\
      &\qquad\lesssim
      \left(
        \llbracket\sigma\rrbracket^{-[\boldsymbol b]-1}
        \wedge
        \varepsilon^{-\flat}
        \llbracket\sigma\rrbracket^{\flat-[\boldsymbol b]-1}
      \right)
      \trinorm{\mathcal F^{\boldsymbol b}_\sigma}_\sigma.
    \end{aligned}
    \label{eq:A-epsilon-bound-corrected}
  \end{equation}
  Here we also used~\eqref{eq:characterization-a} and the first compatibility
  condition in~\eqref{eq:compatibility}, as in the last step of the
  proof of Lemma~\ref{lem:estimates-A}.

  There is, in general, no analogue of~\eqref{eq:A-epsilon-bound-corrected}
  for the full $\mathcal B$-operator. Indeed, the atomic sublattice part of
  $L_\sigma\dot G_\sigma$ has total variation of order
  $\llbracket\sigma\rrbracket^{2s-1}$ and does not gain the factor
  $\varepsilon^{-\flat}\llbracket\sigma\rrbracket^\flat$. For the relevant
  one-point cumulants, the required improvement will instead follow from the
  exact cancellation of the local $\mathcal B$-contribution established in
  Lemma~\ref{lem:local-B-cancellation} below.
\end{remark}

\subsection{Localisation}\label{app:localization}

In this section, we introduce the Taylor expansion on the (semi)-discrete
lattice $\Lambda = \Lambda_{\varepsilon}$. At first order, we have
\[ \begin{aligned}
     \psi (z_1) & =  \psi (z) + \sum_{i \in \{ 0, 1 \pm, 2 \pm, \ldots, d \pm
     \}} \int_0^1 [\mathd \rho _{z_1 - z} (t)]^i  (\partial^i \psi) (z +
     \rho_{z_1 - z} (t)).
   \end{aligned} \]
Here $\partial^0$ denotes the time derivative, while
\[ \partial^{k\pm} \psi (z) \assign \pm \varepsilon^{- 1} [\psi (z \pm \varepsilon e_k) -
   \psi (z)], \]
denotes the discrete forward ($k +$) and backward ($k -$) derivatives in the
$k$-th spatial direction. For $h \in \Lambda$, the function $\rho_h : [0, 1]
\rightarrow \Lambda$ is a bounded variation path (left-continuous at its spatial jump) such that $\rho_h (0) = 0$
and $\rho_h (1) = h$. We use the notation 
\[
[\mathd \rho_{z_1 - z} (t)]^0
\assign \mathd \rho_{z_1 - z}^0 (t).
\]
We introduce the decomposition
\[
\mathd\rho_{z_1-z}^k=(\mathd\rho_{z_1-z}^k)_++(\mathd\rho_{z_1-z}^k)_-,
\]
where $(\mathd\rho_{z_1-z}^k)_+$ is nonnegative and
$(\mathd\rho_{z_1-z}^k)_-$ is non-positive. We write
\[
[\mathd \rho_{z_1 - z} (t)]^{k \pm}
\assign (\mathd \rho_{z_1 - z}^k (t))_{\pm}.
\]
Note that the path $\rho_h$ is piecewise
constant in space, so that the signed measure $\mathd \rho_h^i$ is well
defined and given by a sum of delta functions multiplied by the corresponding
increments. We choose it so that its total variation is bounded by a constant times $| h | $ and
\[ \int_0^1 [\mathd \rho_h (t)]^{k \pm} = [h]^{k \pm} \assign h^k
   \mathbbm{1}_{\pm h^k \geqslant 0} . \]
\begin{remark}
  Note that in the continuum we could choose
  \[ \rho_{z_1 - z} (t) = (z_1 - z) t. \]
  Moreover, since in the continuum both the right and left derivatives
  $\partial^{k\pm}$ coincide with $\partial^k$, we have
  \[ \psi (z_1) = \psi (z) + \sum_{k \in \{ 0, \ldots, d \}} \int_0^1
     (\partial^k \psi) (z + \rho_{z_1 - z} (t))  (z_1 - z)^k \mathd t. \]
\end{remark}

Since we intend to use the second-order Taylor expansion, it is convenient to
choose the path such that $\rho^0_{z_1 - z} (u) = 0$ for $u \in [0, 1 / 2]$
and
\[ \rho^i_{z_1 - z} (u) = \rho^i_{z_1 - z} (1 / 2) + \mathbbm{1}_{i =
   0}  (z_1 - z)^0   (2 u - 1) \]
for $u \in [1 / 2, 1]$. This allows us to avoid second-order terms with one
derivative in time and one in space. Note also that
\begin{align*}
  \int_0^{1 / 2} | \mathd \rho _{z_1 - z} (u) |
  \lesssim | \bar{z}_1 - \bar{z} |,
  \qquad
  \int_{1 / 2}^1 | \mathd \rho _{z_1 - z} (u) |
  \leqslant | z_{1,0} - z_0 | .
\end{align*}
We have
\begin{align*}
  \psi (z_1) &=  \psi (z) + \sum_i \int_0^1 [\mathd \rho _{z_1 - z} (t)]^i
  \partial^i \psi (z + \rho_{z_1 - z} (t))\\
 &=  \psi (z) + \sum_{i \neq 0} \int_0^{1 / 2} [\mathd \rho _{z_1 - z}
  (t)]^i \partial^i \psi (z + \rho_{z_1 - z} (t))\\
 &\quad + \int_{1 / 2}^1 \mathd
  \rho_{z_1 - z}^0 (t) \partial^0 \psi (z + \rho_{z_1 - z} (t)) .
\end{align*}
Expanding once more the spatial derivatives, we have
\begin{align*}
  \psi (z_1) &= \psi (z)
  + \sum_{i \neq 0} \partial^i \psi (z) [z_1 - z]^i\\
  &\quad + \sum_{i, j \neq 0} \int_0^{1 / 2}
  [\mathd \rho _{z_1 - z} (t)]^i \int_0^t
  [\mathd \rho _{z_1 - z} (u)]^j\,
  \partial^i \partial^j \psi (z + \rho_{z_1 - z} (u))\\
  &\quad + \int_{1 / 2}^1 \mathd \rho_{z_1 - z}^0 (t)
  \partial^0 \psi (z + \rho_{z_1 - z} (t)) .
\end{align*}
As a consequence, for a generic kernel $V (z, z_1)$, and symmetrising the
factor $\partial^i \psi (z)$ using the relation $\partial^{i+} \psi (z) -
\partial^{i-} \psi (z) = \varepsilon (\partial^{i+} \partial^{i-} \psi)
(z)$ we obtain
  \begin{equation}
  \begin{aligned}
    V (\psi) (z)
    &= \int_{\Lambda} \mathd z_1 V (z, z_1) \psi (z_1)\\
    &= \psi (z) \int_{\Lambda} \mathd z_1 V (z, z_1)
    + \sum_{k \neq 0}
    \frac{\partial^{k+}  \psi (z) + \partial^{k-} \psi (z)}{2}
    \int_{\Lambda} \mathd z_1 (z_1 - z)^k V (z, z_1)\\
    &\quad + \int_{\Lambda} \mathd z_1 \int_{1 / 2}^1
    \mathd \rho_{z_1 - z}^0 (t) V (z, z_1)
    \partial^0 \psi (z + \rho_{z_1 - z} (t))\\
    &\quad + \sum_k \frac{\varepsilon (\partial^{k+} \partial^{k-} \psi)(z)}{2}
    \int_{\Lambda} |(z_1 - z)^k| V (z, z_1)\,\mathd z_1\\
    &\quad + \sum_{i, j \neq 0} \int_{\Lambda} \mathd z_1
    \int_0^{1 / 2}\int_0^t
    [\mathd \rho _{z_1 - z} (u)]^j [\mathd \rho _{z_1 - z} (t)]^i\\
    &\qquad \times V (z, z_1) \partial^i \partial^j
    \psi (z + \rho_{z_1 - z} (u)) .
  \end{aligned} \label{eq:Taylor-tested}
\end{equation}
By duality, this Taylor expansion defines operators $\boldsymbol{L}$ and
$\boldsymbol{R}$ acting on the kernel $V$, where
\begin{equation}
  \boldsymbol{R} \assign \boldsymbol{R}^{(0)} + \sum_{i, j \neq 0}
  \boldsymbol{R}^{(i, j)} + \sum_k \boldsymbol{R}_{\varepsilon}^{(k)},
  \label{eq:decomposition-of-R}
\end{equation}
and these operators have the form
\begin{equation}
  \begin{aligned}
    (\boldsymbol{L} V) (z, z_1)
    &\assign \int_{\Lambda} \mathd \tilde{z} V
    (z, \tilde{z})  \delta (z_1 - z)\\
    &\quad + \sum_{i \neq 0}
    \frac{1}{2} \int_{\Lambda} (\tilde{z} - z)^i V (z, \tilde{z} )
    (\partial^{i+}  + \partial^{i-}) \delta (z_1 - z)
  \end{aligned} \label{eq:L-operator}
\end{equation}
and
\begin{equation}
  \begin{aligned}
    (\boldsymbol{R}^{(0)} V) (z, z_1)
    &\assign \int_{\Lambda} \int_{1 / 2}^1
    \mathd \rho_{\tilde{z} - z}^0 (t) V (z, \tilde{z})\\
    &\quad \times \partial^0 \delta (z_1
    - z - \rho_{\tilde{z} - z} (t))\,\mathd \tilde{z},\\
    (\boldsymbol{R}^{(i, j)} V) (z, z_1)
    &\assign \int_{\Lambda}
    \int_0^{1 / 2} \int_0^t  [\mathd \rho _{\tilde{z} - z} (u)]^j
    [\mathd \rho_{\tilde{z} - z} (t)]^i\\
    &\quad \times V (z, \tilde{z}) \partial^i \partial^j
    \delta (z_1 - z - \rho_{\tilde{z} - z} (u))\,\mathd \tilde{z},\\
    (\boldsymbol{R}_{\varepsilon}^{(k)} V) (z, z_1)
    &\assign 
    \frac{\varepsilon}{2} \int_{\Lambda}  |(\tilde{z} - z)^k| V (z, \tilde{z})
    \mathd \tilde{z}  (\partial^{k+}  \partial^{k-}) \delta (z_1 - z) .
  \end{aligned} \label{eq:R-operators}
\end{equation}
When tested with a smooth function this gives the identity
\[ V (\psi) = (\boldsymbol{L} V) (\psi) + (\boldsymbol{R} V) (\psi), \]
as seen in~\eqref{eq:Taylor-tested}. We use this expansion for the analysis of
relevant cumulants, namely
\[ \mathcal{F}^{\boldsymbol{a}}_{\mu} = \bar{F}_{\mu}^{[\ell], (1)}
   =\mathbb{E}F_{\mu}^{[\ell], (1)}, \qquad \ell \in \{ 1, \ldots, \hat{\ell}
   \} , \qquad \boldsymbol{a} = ((\ell, 1)). \]
Recall that we have the following decomposition~\eqref{eq:flow-decomposition},
\begin{align*}
  \bar{F}_{\mu}^{[\ell], (1)} &=  \bar{F}_1^{[\ell], (1)} - \int_{\mu}^1 [(1
  - K^{1, 1}_{\sigma}) \dot{\bar{F}}_{\sigma}^{[\ell], (1)}]\,\mathd \sigma -
  \int_{\mu}^1 [K^{1, 1}_{\sigma} \dot{\bar{F}}_{\sigma}^{[\ell], (1)}] \cdot
  (1 - h_{\mu})\,\mathd \sigma\\
 &\quad + \boldsymbol{L} \int_{\mu}^1 [K^{1, 1}_{\sigma}
  \dot{\bar{F}}_{\sigma}^{[\ell], (1)}] \cdot (1 - h_{\mu})\,\mathd \sigma -
  \boldsymbol{L} \int_{\mu}^1 [K^{1, 1}_{\sigma} \dot{\bar{F}}_{\sigma}^{[\ell],
  (1)}]\,\mathd \sigma\\
 &\quad - \boldsymbol{R} \int_{\mu}^1 [K^{1, 1}_{\sigma}
  \dot{\bar{F}}_{\sigma}^{[\ell], (1)}] \cdot h_{\mu}\, \mathd \sigma,
\end{align*}
where $\dot{\bar{F}}_{\sigma}^{[\ell], (1)} \assign \partial_{\sigma}
\bar{F}_\sigma^{[\ell], (1)}$.

\begin{lemma}
  \label{lem:estimates-local}The following bounds hold uniformly in $\sigma,
  \mu \in [1 / 2, 1)$.
  \begin{align*}
       \trinorm{ \boldsymbol{L} (K^{1, 1}_{\sigma}
       \dot{\bar{F}}_{\sigma}^{[\ell], (1)}) }_{\mu} & \lesssim 
       \trinorm{ \dot{\bar{F}}_{\sigma}^{[\ell], (1)} }_{\sigma},\\
       \trinorm{ \boldsymbol{L} ([K^{1, 1}_{\sigma}
       \dot{\bar{F}}_{\sigma}^{[\ell], (1)}] \cdot (1 - h_{\mu}))
       }_{\mu} & \lesssim  \llbracket \mu \rrbracket^{- \flat}
       \llbracket \sigma \rrbracket^{\flat}  \trinorm{
       \dot{\bar{F}}_{\sigma}^{[\ell], (1)} }_{\sigma} .
   \end{align*}
\end{lemma}

\begin{proof}
  As we argued in Remark~\ref{rem:symmetries}, the cumulants are invariant under
  translations and spatial reflections. As a result, we conclude that
  \[ (\partial^{i+}  + \partial^{i-}) \delta (z_1 - z) \int_{\Lambda} \mathd
     \tilde{z}  (\tilde{z} - z)^i (K^{1, 1}_{\sigma}
     \dot{\bar{F}}_{\sigma}^{[\ell], (1)}) (z, \tilde{z}) = 0
     . \]
  Hence, the following identity
  \[ (\boldsymbol{L} (K^{1, 1}_{\sigma} \dot{\bar{F}}_{\sigma}^{[\ell], (1)}))
     (z, z_1) = \delta (z_1 - z) \int_{\Lambda} (K^{1, 1}_{\sigma}
     \dot{\bar{F}}_{\sigma}^{[\ell], (1)}) (z, \tilde{z})\,\mathd \tilde{z} \]
  holds. Similarly,
  \[ (\boldsymbol{L} ([K^{1, 1}_{\sigma} \dot{\bar{F}}_{\sigma}^{[\ell], (1)}]
     \cdot (1 - h_{\mu}))) (z, z_1) = \delta (z_1 - z) \int_{\Lambda} (K^{1,
     1}_{\sigma} \dot{\bar{F}}_{\sigma}^{[\ell], (1)}) (z, \tilde{z})  (1 -
     h_{\mu} (z, \tilde{z}))\,\mathd \tilde{z} . \]
  Using the first identity, we obtain
  \begin{align*}
     &w^{(2), \flat}_{\mu} (z, z_1)
     (K^{1, 1}_{\mu} \boldsymbol{L} [K^{1,
     1}_{\sigma} \dot{\bar{F}}_{\sigma}^{[\ell], (1)}]) (z, z_1)\\
     &\qquad = w^{(2), \flat}_{\mu} (z, z_1) K_{\mu} (z_1 - z)
     \int_{\Lambda} (K^{1, 1}_{\sigma}
     \dot{\bar{F}}_{\sigma}^{[\ell], (1)}) (z, \tilde{z}) \,\mathd \tilde{z} .
  \end{align*}
  Consequently,
  \begin{align*}
    \trinorm{ \boldsymbol{L} (K^{1, 1}_{\sigma}
    \dot{\bar{F}}_{\sigma}^{[\ell], (1)}) }_{\mu} &\lesssim 
    \sup_z \int_{\Lambda} | K_{\mu} (z_1 - z) | w_{\mu}^{\flat} (z - z_1)
    \mathd z_1 \int_{\Lambda} | (K^{1, 1}_{\sigma}
    \dot{\bar{F}}_{\sigma}^{[\ell], (1)}) (z, \tilde{z}) |\,\mathd \tilde{z}\\
 &\lesssim  \| K_{\mu} \|_{\operatorname{TV} (w_{\mu}^{\flat})} \sup_z
    \int_{\Lambda} | (K^{1, 1}_{\sigma} \dot{\bar{F}}_{\sigma}^{[\ell], (1)})
    (z, \tilde{z}) |\,\mathd \tilde{z}\\
 &\lesssim  \sup_z \int_{\Lambda} | (K^{1, 1}_{\sigma}
    \dot{\bar{F}}_{\sigma}^{[\ell], (1)}) (z, \tilde{z}) |\,\mathd \tilde{z},
  \end{align*}
  by Lemma~\ref{lem:bounds-K-J}. Since
  \begin{align*} \sup_z \int_{\Lambda} | K^{1, 1}_{\sigma} \dot{\bar{F}}_{\sigma}^{[\ell],
     (1)} (z, \tilde{z}) |\,\mathd \tilde{z} 
     &\lesssim
      \sup_z \int_{\Lambda} |
     K^{1, 1}_{\sigma} \dot{\bar{F}}_{\sigma}^{[\ell], (1)} (z, \tilde{z}) |\,
     w^{(2), \flat}_{\sigma} (z, \tilde{z})\, \mathd \tilde{z} 
     \\
     &= \trinorm{
     \dot{\bar{F}}_{\sigma}^{[\ell], (1)} }_{\sigma},
  \end{align*}
  this proves the first of the bounds stated in the lemma. By a similar
  argument, we get

  \[ \trinorm{ \boldsymbol{L} ([K^{1, 1}_{\sigma}
     \dot{\bar{F}}_{\sigma}^{[\ell], (1)}] \cdot (1 - h_{\mu}))
     }_{\mu} \lesssim
     \sup_{z} \int_{\Lambda} | (K^{1, 1}_{\sigma}
     \dot{\bar{F}}_{\sigma}^{[\ell], (1)}) (z, \tilde{z}) |  (1 - h_{\mu} (z,
     \tilde{z}))\,\mathd \tilde{z} . \]
  Consequently, by Lemma~\ref{lem:weight-interaction}, we obtain
  \begin{align*}
    &\trinorm{ \boldsymbol{L} ([K^{1, 1}_{\sigma}
    \dot{\bar{F}}_{\sigma}^{[\ell], (1)}] \cdot (1 - h_{\mu}))
    }_{\mu}\\
    &\quad \lesssim
    \| w^{(2), \flat}_{\mu} (1 - h_{\mu})
    (w^{(2), \flat}_{\sigma})^{- 1} \|_{L^{\infty}}
    \trinorm{ \dot{\bar{F}}_{\sigma}^{[\ell], (1)}
    }_{\sigma}\\
    &\quad \lesssim
    \llbracket \mu \rrbracket^{- \flat}
    \llbracket \sigma \rrbracket^{\flat}
    \trinorm{ \dot{\bar{F}}_{\sigma}^{[\ell], (1)}
    }_{\sigma} .
  \end{align*}
  This proves the second of the bounds stated in the lemma.
\end{proof}

\begin{lemma}[Cancellation of the local $\mathcal B$-contribution]
  \label{lem:local-B-cancellation}
  Let $\ell\in\{1,\ldots,\hat{\ell}\}$ and let $\boldsymbol a$ be the
  one-element tuple whose unique entry is $(\ell,1)$. Then, uniformly in
  $\sigma\in[1/2,1)$,
  \begin{equation}
    \boldsymbol L\left(
      K^{1,1}_\sigma
      \sum_{\boldsymbol b,\boldsymbol c}
      \mathcal B^{\boldsymbol a}_{\boldsymbol b,\boldsymbol c}
      (\dot G_\sigma,
       \mathcal F^{\boldsymbol b}_\sigma,
       \mathcal F^{\boldsymbol c}_\sigma)
    \right)=0.
    \label{eq:local-B-cancellation}
  \end{equation}
\end{lemma}

\begin{proof}
  Since $n(\boldsymbol a)=1$, in~\eqref{eq:specific-operator-B} we have
  $i=1$ and $I_1=I_2=\varnothing$. Consider first the terms with $k'=1$.
  The two cumulants appearing in such a term are indexed by the one-element
  tuples
  \[
    \boldsymbol b=((\ell-1-\ell',2)),
    \qquad
    \boldsymbol c=((\ell',0)).
  \]
  They vanish by Remark~\ref{rem:symmetries}, since respectively
  $n(\boldsymbol b)+K(\boldsymbol b)=3$ and
  $n(\boldsymbol c)+K(\boldsymbol c)=1$ are odd.

  It remains to consider $k'=0$. In this case
  \[
    \boldsymbol b=((\ell-1-\ell',1)),
    \qquad
    \boldsymbol c=((\ell',1)).
  \]
  As in the proof of Lemma~\ref{lem:estimates-B}, we use
  $L_\sigma K_\sigma=1$ to write
  \begin{align*}
    &K^{1,1}_\sigma
    \mathcal B^{\boldsymbol a,(1)}_{\boldsymbol b,\boldsymbol c}
    (\dot G_\sigma,
     \mathcal F^{\boldsymbol b}_\sigma,
     \mathcal F^{\boldsymbol c}_\sigma)
    \\
    &\qquad=
    \mathcal B^{\boldsymbol a,(1)}_{\boldsymbol b,\boldsymbol c}
    (L_\sigma\dot G_\sigma,
     K^{\boldsymbol b}_\sigma\mathcal F^{\boldsymbol b}_\sigma,
     K^{\boldsymbol c}_\sigma\mathcal F^{\boldsymbol c}_\sigma).
  \end{align*}
  Set
  \[
    W_\sigma\assign
    K^{\boldsymbol b}_\sigma\mathcal F^{\boldsymbol b}_\sigma,
    \qquad
    U_\sigma\assign
    K^{\boldsymbol c}_\sigma\mathcal F^{\boldsymbol c}_\sigma.
  \]
  By~\eqref{eq:one-point-mass-translation}, and since convolution with
  $K_\sigma$ preserves translation invariance, the quantity
  \[
    m_{U_\sigma}\assign
    \int_\Lambda U_\sigma(z'';\mathd z_1)
  \]
  is independent of $z''$. Hence, using the definition
  \eqref{eq:def-B-general} of the contraction operator,
  \begin{align*}
    &\int_\Lambda
    \mathcal B^{\boldsymbol a,(1)}_{\boldsymbol b,\boldsymbol c}
    (L_\sigma\dot G_\sigma,W_\sigma,U_\sigma)
    (z;\mathd z_1)
    \\
    &\qquad=
    \int_{\Lambda^2}
    W_\sigma(z;\mathd z')
    (L_\sigma\dot G_\sigma)(z'-z'')
    \left(\int_\Lambda U_\sigma(z'';\mathd z_1)\right)
    \mathd z''
    \\
    &\qquad=
    m_{U_\sigma}
    \int_\Lambda W_\sigma(z;\mathd z')
    \int_\Lambda (L_\sigma\dot G_\sigma)(z'-z'')\,\mathd z''
    =0.
  \end{align*}
  For the last identity, observe that
  \begin{equation}
    \int_\Lambda L_\sigma\dot G_\sigma(z)\,\mathd z=0.
    \label{eq:L-Gdot-zero-mass}
  \end{equation}
  Indeed, by Definition~\ref{def:function-j} and~\eqref{eq:cutoff}, the Fourier
  multiplier of $\partial_\sigma J_\sigma$ equals zero at the zero spacetime frequency.
  Therefore the zero Fourier mode of
  $L_\sigma\dot G_\sigma=L_\sigma\partial_\sigma(J_\sigma G)$ vanishes.

  We have shown that the zeroth moment in the definition
  \eqref{eq:L-operator} of $\boldsymbol L$ vanishes. Moreover, the kernel on the left-hand side of
  \eqref{eq:local-B-cancellation} is invariant under simultaneous
  translations and spatial reflections, by Remark~\ref{rem:symmetries}, the
  corresponding symmetries of \(K_\sigma\) and \(G_\sigma\), and the explicit
  formula \eqref{eq:specific-operator-B}.  Translation invariance reduces it
  to a kernel of the relative variable, and reflection invariance makes that
  relative kernel even in every spatial coordinate.  Its first spatial
  moments are therefore zero, exactly as in the expanded argument in the
  proof of Lemma~\ref{lem:estimates-local}.
  All coefficients in its local projection are zero, and
  \eqref{eq:local-B-cancellation} follows.
\end{proof}

\begin{lemma}
  \label{lem:estimates-remainder}The following bound
  \begin{align*}
    \trinorm{ \boldsymbol{R} ([K^{1, 1}_{\sigma}
    \dot{\bar{F}}_{\sigma}^{[\ell], (1)}] \cdot h_{\mu})
    }_{\mu}
    &\lesssim \llbracket \mu \rrbracket^{- \flat}
    \llbracket \sigma \rrbracket^{\flat}
    \trinorm{ \dot{\bar{F}}_{\sigma}^{[\ell], (1)}
    }_{\sigma}
  \end{align*}
  holds uniformly in $\mu \in [1 / 2, 1), \sigma \in [\mu, 1)$.
\end{lemma}

\begin{proof}
  By applying the triangle inequality to~\eqref{eq:decomposition-of-R}, we
  arrive at
  \begin{align*}
    \trinorm{ \boldsymbol{R} ([K^{1, 1}_{\sigma} \dot{\bar{F}}_{\sigma}^{[\ell],
    (1)}] \cdot h_{\mu}) }_{\mu} &\lesssim  \trinorm{
    \boldsymbol{R}^{(0)} ([K^{1, 1}_{\sigma} \dot{\bar{F}}_{\sigma}^{[\ell],
    (1)}] \cdot h_{\mu}) }_{\mu}\\
 &\quad + \sum_{i, j \neq 0} \trinorm{ \boldsymbol{R}^{(i, j)} ([K^{1,
    1}_{\sigma} \dot{\bar{F}}_{\sigma}^{[\ell], (1)}] \cdot h_{\mu})
    }_{\mu}\\
 &\quad + \sum_k \trinorm{ \boldsymbol{R}^{(k)}_{\varepsilon} ([K^{1,
    1}_{\sigma} \dot{\bar{F}}_{\sigma}^{[\ell], (1)}] \cdot h_{\mu})
    }_{\mu} .
  \end{align*}
  We have
  \begin{align*}
    &(\boldsymbol{R}^{(i, j)} ([K^{1, 1}_{\sigma}
    \dot{\bar{F}}_{\sigma}^{[\ell], (1)}] \cdot h_{\mu})) (z, z_1)\\
    &\quad = \int_{\Lambda} \int_0^{\frac{1}{2}} \int_0^t
    [\mathd \rho _{\tilde{z} - z} (u)]^j
    [\mathd \rho _{\tilde{z} - z} (t)]^i
    h_{\mu}(z, \tilde{z})\\
    &\qquad \times
    (K^{1, 1}_{\sigma} \dot{\bar{F}}_{\sigma}^{[\ell], (1)})(z, \tilde{z})
    \partial^i \partial^j
    \delta (z_1 - z - \rho_{\tilde{z} - z} (u))
    \,\mathd \tilde{z} .
  \end{align*}
  Thus,
  \begin{align*}
    &| w^{(2), \flat}_{\mu} (z, z_1)
    (\boldsymbol{R}^{(i, j)} ([K^{1, 1}_{\sigma}
    \dot{\bar{F}}_{\sigma}^{[\ell], (1)}] \cdot h_{\mu})) (z, z_1) |\\
    &\quad \leqslant \int_{\Lambda} \int_0^{\frac{1}{2}} \int_0^t
    | \mathd \rho _{\tilde{z} - z} (u) |
    \,| \mathd \rho _{\tilde{z} - z} (t) |\,
    w^{(2), \flat}_{\mu} (z, z_1) h_{\mu} (z, \tilde{z})\\
    &\qquad \times
    | K^{1, 1}_{\sigma} \dot{\bar{F}}_{\sigma}^{[\ell], (1)} (z, \tilde z) |
    \,| \partial^i \partial^j
    K_{\mu} (z_1 - z - \rho_{\tilde{z} - z} (u)) |\,
    \mathd \tilde{z} .
  \end{align*}
  Recall that
  \[
  w_{\mu}^{(2), \flat} (z, z_1) = w_{\mu}^{\flat} (z - z_1)
  = (1 + \llbracket \mu \rrbracket^{- 1} | z - z_1 |_s)^{\flat}.
  \]
  Introduce the point
  \[
  Z_u=Z_u(z,\tilde z) \assign z+\rho_{\tilde z-z}(u).
  \]
  Then
  \begin{align*}
  	w^{\flat}_{\mu} (z - z_1)
  	&\lesssim w^{\flat}_{\mu} (z - Z_u)
  	w^{\flat}_{\mu} (Z_u - z_1)\\
  	&\leqslant w^{\flat}_{\mu} (z - \tilde z)
  	w^{\flat}_{\mu} (Z_u - z_1).
  \end{align*}
  The second inequality above follows from
  \[ w^{\flat}_{\mu} (z - Z_u) \lesssim w^{\flat}_{\mu} (z -
     \tilde{z}), \]
  which holds because, by construction, $Z_u$ lies on the path
  connecting $z$ and $\tilde{z}$. Overall, this yields
  \begin{align*}
    &| w_{\mu}^{\flat} (z - z_1)
    (\boldsymbol{R}^{(i, j)} ([K^{1, 1}_{\sigma}
    \dot{\bar{F}}_{\sigma}^{[\ell], (1)}] \cdot h_{\mu})) (z, z_1) |\\
    &\quad \lesssim \int_{\Lambda} \int_0^{\frac{1}{2}} \int_0^t
    | \mathd \rho _{\tilde{z} - z} (u) |\,
    | \mathd \rho _{\tilde{z} - z} (t) |\,
    w_{\mu}^{\flat} (z - \tilde{z}) h_{\mu} (z, \tilde{z})
    (w_{\sigma}^{\flat} (z - \tilde{z}))^{- 1}\\
    &\qquad \times w_{\sigma}^{\flat} (z - \tilde{z})
    \,| (K^{1, 1}_{\sigma}
    \dot{\bar{F}}_{\sigma}^{[\ell], (1)}) (z, \tilde{z}) |\\
    &\qquad \times
    w^{\flat}_{\mu} (Z_u - z_1)
    \,| \partial^i \partial^j K_{\mu}
    (z_1 - Z_u) |\,\mathd \tilde{z} .
  \end{align*}
  As a consequence, we obtain
  \begin{align*}
    &\trinorm{ \boldsymbol{R}^{(i, j)} ([K^{1, 1}_{\sigma}
    \dot{\bar{F}}_{\sigma}^{[\ell], (1)}] \cdot h_{\mu})
    }_{\mu}\\
    &\quad \lesssim
    \sup_{z, \tilde{z} \in \Lambda}
    \int_0^{\frac{1}{2}} \int_0^t
    | \mathd \rho _{\tilde{z} - z} (u) |\,
    | \mathd \rho _{\tilde{z} - z} (t) |\,
    w_{\mu}^{\flat} (z - \tilde{z}) h_{\mu} (z, \tilde{z})
    (w_{\sigma}^{\flat} (z - \tilde{z}))^{- 1}\\
    &\qquad \times
    \trinorm{ \dot{\bar{F}}_{\sigma}^{[\ell], (1)}
    }_{\sigma}
    \sup_{z, \tilde{z} \in \Lambda} \int_{\Lambda}
    w^{\flat}_{\mu} (Z_u - z_1)\,
    | \partial^i \partial^j K_{\mu}
    (z_1 - Z_u) |\, \mathd z_1 .
  \end{align*}
  Using
  \begin{equation}
    \int_0^{\frac{1}{2}} \int_0^t  | \mathd \rho _{\tilde{z} - z} (u) |  |
    \mathd \rho _{\tilde{z} - z} (t) | \lesssim | \overline{\tilde{z}} -
    \bar{z} |^2, \label{eq:bounded-total-mass}
  \end{equation}
  and
  \begin{equation}
    {\int_{\Lambda}}  w^{\flat}_{\mu} (z ) | \partial^i \partial^j K_{\mu} (z
    ) |\,\mathd z \lesssim \llbracket \mu \rrbracket^{- 2},
    \label{eq:bound-for-deriv-of-K}
  \end{equation}
  as well as
  \begin{align*}
      \frac{| \overline{\tilde{z}} - \bar{z} |^2 w_{\mu}^{\flat} (z -
      \tilde{z}) h_{\mu} (z, \tilde{z})}{w_{\sigma}^{\flat} (z - \tilde{z})} &
      \leqslant  \llbracket \mu \rrbracket^{2 - \flat} \llbracket \sigma
      \rrbracket^{\flat} \frac{(1 + \llbracket \mu \rrbracket^{- 1} | z -
      \tilde{z} |_s )^{\flat}}{(1 + \llbracket \sigma \rrbracket^{- 1} | z -
      \tilde{z} |_s)^{\flat}} \frac{\llbracket \mu \rrbracket^{- (2 - \flat)}
      \llbracket \sigma \rrbracket^{- \flat} | z - \tilde{z} |_s^2}{1 +
      {(\llbracket \mu \rrbracket^{- 1} | z - \tilde{z} |_s)^2} }\\
      & \leqslant  \llbracket \mu \rrbracket^{2 - \flat} \llbracket \sigma
      \rrbracket^{\flat} \frac{\left( 1 + \llbracket \mu \rrbracket^{- 1}
      | z - \tilde{z} |_s  \right)^{\flat} (\llbracket \mu
      \rrbracket^{- 1} | z - \tilde{z} |_s)^{2 - \flat}}{1 + {(\llbracket \mu
      \rrbracket^{- 1} | z - \tilde{z} |_s)^2} }\\
      & \lesssim  \llbracket \mu \rrbracket^{2 - \flat} \llbracket \sigma
      \rrbracket^{\flat},
		\label{eq:absorbing-monomials}
  \end{align*}
  we conclude that
  \[ \trinorm{ \boldsymbol{R}^{(i, j)} ([K^{1, 1}_{\sigma}
     \dot{\bar{F}}_{\sigma}^{[\ell], (1)}] \cdot h_{\mu}) }_{\mu}
     \lesssim \llbracket \sigma \rrbracket^{\flat} \llbracket \mu
     \rrbracket^{- \flat}  \trinorm{ \dot{\bar{F}}_{\sigma}^{[\ell], (1)}
     }_{\sigma} . \]
  Next, note that
  \begin{align*}
     &(\boldsymbol{R}^{(0)} ([K^{1, 1}_{\sigma}
     \dot{\bar{F}}_{\sigma}^{[\ell], (1)}] \cdot h_{\mu})) (z, z_1)\\
     &\quad = \int_{\Lambda} \int_{1 / 2}^1
     \mathd \rho_{\tilde{z} - z}^0 (t)
     (h_{\mu} K^{1, 1}_{\sigma}
     \dot{\bar{F}}_{\sigma}^{[\ell], (1)}) (z, \tilde{z})
     \,\partial^0 \delta (z_1 - z - \rho_{\tilde{z} - z} (t))
     \,\mathd \tilde{z} .
  \end{align*}
  Proceeding as above, we get
  \begin{align*}
    &\trinorm{ \boldsymbol{R}^{(0)} ([K^{1, 1}_{\sigma}
    \dot{\bar{F}}_{\sigma}^{[\ell], (1)}] \cdot h_{\mu}) }_{\mu}\\
    &\quad \lesssim  
    \sup_{z, \tilde{z} \in \Lambda}
    \int_{1 / 2}^1 | \mathd \rho^0_{\tilde{z} - z}  (t)|\,
    w_{\mu}^{\flat} (z - \tilde{z}) h_{\mu} (z, \tilde{z})
    \,(w_{\sigma}^{\flat} (z - \tilde{z}))^{- 1}\\
    &\qquad \times \trinorm{ \dot{\bar{F}}_{\sigma}^{[\ell], (1)}
    }_{\sigma} \sup_{z, \tilde{z} \in \Lambda} \int_{\Lambda}
    w_{\mu}^{\flat} (Z_t - z_1)\,
    | \partial^0 K_{\mu} (z_1 -Z_t) |\,\mathd z_1 .
  \end{align*}
  Using
  \[ \int_{1 / 2}^1 |\mathd \rho_{\tilde{z} - z}^0 (t)| \leqslant | \tilde{z}_0
     - z_0 |, \]
  and
  \[ {\int_{\Lambda}}  w_{\mu}^{\flat} (z  - z') \, | \partial^0 K_{\mu} (z - z'
     ) |\,\mathd z \lesssim \llbracket \mu \rrbracket^{- 2 s}, \]
  as well as
  \begin{align*}
       | \tilde{z}_0 - z_0 | \frac{w_{\mu}^{\flat} (z -
       \tilde{z})}{w_{\sigma}^{\flat} (z - \tilde{z})} h_{\mu} (z, \tilde{z}) &
       \leqslant  \llbracket \mu \rrbracket^{2 s - \flat} \llbracket \sigma
       \rrbracket^{\flat} \frac{\left( 1 + \llbracket \mu \rrbracket^{- 1}
       | z - \tilde{z} |_s  \right)^{\flat}}{(1 + \llbracket
       \sigma \rrbracket^{- 1} | z - \tilde{z} |_s)^{\flat}} \frac{\llbracket
       \mu \rrbracket^{- (2 s - \flat)} \llbracket \sigma \rrbracket^{- \flat}
       | z - \tilde{z} |_s^{2 s}}{1 + {(\llbracket \mu \rrbracket^{- 1} | z -
       \tilde{z} |_s)^2} }\\
       & \leqslant  \llbracket \mu \rrbracket^{2 s - \flat} \llbracket
       \sigma \rrbracket^{\flat} \frac{\left( 1 + \llbracket \mu \rrbracket^{-
       1} | z - \tilde{z} |_s  \right)^{\flat} (\llbracket \mu
       \rrbracket^{- 1} | z - \tilde{z} |_s)^{2 s - \flat}}{1 + {(\llbracket
       \mu \rrbracket^{- 1} | z - \tilde{z} |_s)^2} }\\
       & \lesssim  \llbracket \mu \rrbracket^{2 s - \flat} \llbracket \sigma
       \rrbracket^{\flat} \frac{1}{1 + {(\llbracket \mu \rrbracket^{- 1} | z -
       \tilde{z} |_s)^{2 - 2 s}} }\\
       & \lesssim  \llbracket \mu \rrbracket^{2 s - \flat} \llbracket \sigma
       \rrbracket^{\flat},
     \end{align*}
  we arrive at
  \[ \trinorm{ \boldsymbol{R}^{(0)} ([K^{1, 1}_{\sigma}
     \dot{\bar{F}}_{\sigma}^{[\ell], (1)}] \cdot h_{\mu}) }_{\mu}
     \lesssim \llbracket \sigma \rrbracket^{\flat} \llbracket \mu
     \rrbracket^{- \flat}  \trinorm{ \dot{\bar{F}}_{\sigma}^{[\ell], (1)}
     }_{\sigma} . \]
  Finally, we discuss the term
  \begin{align*}
     &(\boldsymbol{R}_{\varepsilon}^{(k)} ([K^{1, 1}_{\sigma}
     \dot{\bar{F}}_{\sigma}^{[\ell], (1)}] \cdot h_{\mu})) (z, z_1)\\
     &\qquad = \frac{\varepsilon}{2} \int_{\Lambda}
     |(\tilde{z} - z)^k|
     (h_{\mu} K^{1, 1}_{\sigma}
     \dot{\bar{F}}_{\sigma}^{[\ell], (1)}) (z, \tilde{z})\,\mathd \tilde{z}\,     (\partial^{k+}  \partial^{k-}) \delta (z_1 - z) .
  \end{align*}
  In this case, the proof also proceeds along the same lines. We obtain
 \begin{align*}
       \trinorm{ \boldsymbol{R}_{\varepsilon}^{(k)} ([K^{1, 1}_{\sigma}
       \dot{\bar{F}}_{\sigma}^{[\ell], (1)}] \cdot h_{\mu}) }_{\mu}
       \lesssim & \varepsilon \sup_{z, \tilde{z} \in \Lambda} |
       \overline{\tilde{z}} - \bar{z} |\, w_{\mu}^\flat (z - \tilde{z}) h_{\mu} (z,
       \tilde{z})  (w_{\sigma}^{\flat} (z - \tilde{z}))^{- 1}\\
       &   \times \trinorm{ \dot{\bar{F}}_{\sigma}^{[\ell], (1)}
       }_{\sigma}\\
       &   \times \int_{\Lambda} w_{\mu}^{\flat} (z - z_1) \, |
       \partial^{k+}  \partial^{k-} K_{\mu} (z_1 - z) |\,\mathd z_1 .
   \end{align*}
  Note that
  \begin{align*}
       \frac{| \overline{\tilde{z}} - \bar{z} | w_{\mu}^{\flat} (z -
       \tilde{z}) h_{\mu} (z, \tilde{z})}{w_{\sigma}^{\flat} (z - \tilde{z})}
       & \leqslant  \llbracket \sigma \rrbracket \frac{(1 + \llbracket \mu
       \rrbracket^{- 1} | z - \tilde{z} |_s)^{\flat}}{(1 + \llbracket \sigma
       \rrbracket^{- 1} | z - \tilde{z} |_s)^{\flat}} \frac{\llbracket \sigma
       \rrbracket^{- 1} | \tilde{z} - z |_s}{{(1 + (\llbracket \mu
       \rrbracket^{- 1} | z - \tilde{z} |_s)^2) } }\\
       & \lesssim  \llbracket \sigma \rrbracket \frac{1}{(1 + \llbracket
       \sigma \rrbracket^{- 1} | \tilde{z} - z |_s)^{\flat - 1}} \frac{1}{{(1
       + \llbracket \mu \rrbracket^{- 1} | z - \tilde{z} |_s)^{2 - \flat}} },
   \end{align*}
  Since $| \overline{\tilde{z}} - \bar{z} | \neq 0$ implies $| \tilde{z} - z
  |_s \geqslant \varepsilon$ we have
  \begin{align*}
    \frac{| \overline{\tilde{z}} - \bar{z} | w_{\mu}^{\flat} (z - \tilde{z})
    h_{\mu} (z, \tilde{z})}{w_{\sigma}^{\flat} (z - \tilde{z})} &\lesssim 
    \llbracket \sigma \rrbracket \frac{1}{(1 + \llbracket \sigma \rrbracket^{-
    1}  \varepsilon)^{\flat - 1}} \frac{1}{{(1 + \llbracket \mu \rrbracket^{-
    1}  \varepsilon)^{2 - \flat}} }\\
 &\leqslant  \llbracket \sigma \rrbracket  \llbracket \sigma
    \rrbracket^{\flat - 1} \varepsilon^{1 - \flat} \llbracket \mu
    \rrbracket^{2 - \flat} {\varepsilon^{\flat - 2}} \\
 &=  \llbracket \sigma \rrbracket^{\flat} \llbracket \mu \rrbracket^{2 -
    \flat} \varepsilon^{- 1} .
  \end{align*}
  Hence, using again~\eqref{eq:bound-for-deriv-of-K} with $i = k +$ and $j = k
  -$ we get
  \begin{align*}
     &\trinorm{ \boldsymbol{R}_{\varepsilon}^{(k)} ([K^{1, 1}_{\sigma}
     \dot{\bar{F}}_{\sigma}^{[\ell], (1)}] \cdot h_{\mu}) }_{\mu}\\
     &\qquad \lesssim \varepsilon \llbracket \mu \rrbracket^{- 2}
     \llbracket \sigma \rrbracket^{\flat}
     \llbracket \mu \rrbracket^{2 - \flat} \varepsilon^{- 1}
     \trinorm{ \dot{\bar{F}}_{\sigma}^{[\ell], (1)} }_{\sigma}\\
     &\qquad =
     \llbracket \mu \rrbracket^{- \flat} \llbracket \sigma \rrbracket^{\flat}
     \trinorm{ \dot{\bar{F}}_{\sigma}^{[\ell], (1)}
     }_{\sigma} .
  \end{align*}
  This finishes the proof.
\end{proof}

\begin{ack}
We would like to thank F.~de~Vecchi and L.~Fresta for discussion pertaining to the analysis of flow equations.
\end{ack}

\begin{funding}
PD was supported by the National Science Centre, Poland, via the grant `Sonata
Bis' 2019/34/E/ST1/00053.  MG and PR have been partially supported by the German Research Foundation (DFG)
under Germany's Excellence Strategy -- GZ 2047/1, project-id 390685813. MG was partially supported by UK Research and
Innovation via the Frontier grant ``StochFields'' (EP/Z534328/1) and by the Simons
Foundation via the Simons Collaboration ``Probabilistic Paths to Quantum Field
Theory'' (SFI-MPS-PP-00012621-05).  For the purpose of open access, the authors
have applied a CC BY public copyright licence to any author-accepted manuscript
arising from this submission.
\end{funding}

\bibliographystyle{template-jems/emsjems}
\bibliography{fractional-phi43-flow-jems}

\end{document}